\documentclass[preprint,12pt,authoryear]{elsarticle}
\usepackage[a4paper,left=1in,right=1in,top=1in,bottom=1in]{geometry}
\usepackage{graphicx}
\usepackage{subcaption}
\usepackage{booktabs}
\usepackage{multirow}
\usepackage{array}
\usepackage{nicematrix}
\usepackage{amsmath}
\usepackage{bm}
\usepackage{amsfonts}
\usepackage{float}
\usepackage{indentfirst}
\usepackage{tikz}
\usetikzlibrary{shapes,arrows,calc,arrows.meta, positioning}
\usetikzlibrary{shapes.geometric, arrows}
\usepackage{algorithm}
\usepackage{algorithmic}
\usepackage[perpage]{footmisc}
\usepackage{xcolor}
\usepackage[colorlinks=true,citecolor=blue,linkcolor=blue,urlcolor=blue]{hyperref}
\usepackage{doi}
\graphicspath{ {figures/} }
\begin{document}

\begin{frontmatter}

\title{Explicit topology control in isogeometric topology optimization using persistent homology}

\author[zju]{Wenliang Yin}
\author[zju]{Gengchen Li}
\author[zju]{Hongwei Lin\corref{cor}}
\ead{hwlin@zju.edu.cn}

\cortext[cor]{Corresponding author.}
\address[zju]{School of Mathematical Sciences, Zhejiang University, Hangzhou, 310058, China}

\begin{abstract}
Explicitly controlling structural topology remains a major challenge in topology
optimization. While many implicit methods regulate topology only indirectly,
available explicit strategies often focus on two-dimensional hole control or
extend to 3D without a unified treatment of tunnels and cavities.
This study proposes a differentiable, persistent
homology-based framework for explicit topology control in 2D and 3D isogeometric
topology optimization. A NURBS-represented density field induces a single sublevel
filtration whose persistence diagrams quantify the connected components, holes in 2D, and tunnels and cavities in 3D directly.
To enforce prescribed bounds on topological feature counts, weighted topological objectives are applied
in a staged manner that coordinates the different homology dimensions.
Each redundant feature is removed either by cutting or by filling according to the
region of its corresponding persistence pair, avoiding the excessive structural
changes associated with filling alone. The topological inverse mapping yields
analytical sensitivities with respect to the density control coefficients, so that
the resulting objective is readily incorporated into gradient-based
optimization and solved by the method of moving asymptotes. Numerical experiments
on 2D and 3D benchmarks demonstrate explicit and reliable topology
control while maintaining competitive mechanical performance and stable convergence.
\end{abstract}


\end{frontmatter}

\section{Introduction}
\label{sec:introduction}

Topology optimization \citep{bendsoeGeneratingOptimalTopologies1988} is a powerful structural design method for
determining the optimal material distribution within prescribed domains
under a given set of loads and boundary conditions to obtain optimal structural performance. 
It has been widely applied in many industrial fields, such as the architecture industry, the automotive industry and the aerospace industry \citep{ALLAIRE2004363, fd277c3919a2408ab668cf9803ba21ab, sobieszczanski-sobieskiMultidisciplinaryAerospaceDesign1997, ZHANG201275}.\par

Despite these achievements, controlling structural complexity, particularly the
topology of optimized structures, remains a fundamental challenge \citep{zuo2022explicit}. 
Structural topology control is significant in both theoretical analysis and practical application of topology optimization.
From a theoretical perspective, the original continuum topology optimization is known to be ill-posed, as the global optimum may contain an infinite number of holes \citep{sigmundNumericalInstabilitiesTopology1998}. 
Controlling the number of holes mitigates numerical instability and improves the well-posedness of the optimization formulation.
From a manufacturing perspective, excessively complex structures incur high manufacturing energy consumption and time cost, making them
frequently difficult or even infeasible to fabricate using traditional processes.
While advanced additive manufacturing techniques have largely removed the cost imposed by structural complexity, enclosed internal cavities remain prohibited in processes such as selective laser sintering \citep{10.1108/13552540510573365} and fused deposition modeling \citep{ZEIN20021169}.
Therefore, controlling the structural topology is essential for maintaining structural reliability and ensuring manufacturability.\par

Existing approaches for controlling structural topology can generally be categorized into implicit and explicit methods.
Implicit methods do not impose a direct constraint on the number of holes or connected components, but instead influence topology indirectly through auxiliary mechanisms. For example,
the widely adopted sensitivity filter \citep{sigmund1994design} and density filter \citep{bourdin2001filters} can eliminate numerous small holes as a by-product of suppressing checkerboard patterns.
Another indirect approach is to control the minimum length scale, such as through the skeleton-based method \citep{zhang2014length} or the geometric constraint method \citep{zhou2015length}.
The classic level-set method permits only a monotonic decrease in the number of holes, so the resulting topology is strongly dependent on the initial design \citep{wang2003levelset}.
The virtual temperature methods \citep{liu2015voids,li2018virtual,yamada2022cavity} can eliminate enclosed voids by introducing an auxiliary temperature field to fulfill the connectivity requirement, but they lack an explicit relationship between the control variables and the resulting number of holes.\par

In contrast, explicit methods seek to establish a direct relationship between the design variables and the numbers of holes or connected components.
For instance, the moving morphable components method (MMC) \citep{zhang2017mmc} changes the topology directly by adding, removing, and reshaping geometric components, while the skeleton-based level-set method \citep{zhang2017skeleton} assigns an independent level-set function to each hole.
Within the framework of discrete-variable topology optimization, two explicit strategies have been developed for controlling different topological quantities. Based on the Euler–Poincaré formula and a fluid flooding strategy, an explicit formulation enables the direct computation and control of the number of holes in 2D and of cavities in 3D \citep{liang2022complexity}. 
Alternatively, a quotient set-based design variable method (QSDV) \citep{hanInequalityConstraintMaximum2022} imposes an inequality constraint on the maximum allowable number of tunnels in 3D.
In addition, several explicit strategies have been incorporated into BESO-type frameworks for direct topology control. 
Specifically, graph- and set-theory-inspired methods \citep{zhao2020direct, xiongNewApproachEliminating2020a}, the topology-preserving method based on a thinning algorithm \citep{heThinningAlgorithmBased2022}, and approaches employing the digital Gauss–Bonnet formula for counting holes \citep{han2021constraints} constrain the number of holes directly.
A further strategy \citep{heHolefillingBasedApproach2023} in this category additionally constrains the number of holes in 2D and of cavities and tunnels in 3D by filling in cavities and covering tunnels with patches.
Under the SIMP framework, two explicit formulations calculate the number of holes in 2D structures directly based on fire-burning method (FBM) and control this number through tuning the filter radius \citep{wang2022topological} and hole filling \citep{zuo2022explicit}, respectively. The latter has subsequently been extended to 3D, yielding a formulation that enables the explicit and quantitative control of structural topology \citep{zuo2023tunnels}.\par

Although topology control has been extensively studied in topology optimization, several issues remain to be further investigated. Specifically, implicit methods cannot quantitatively control the structural topology, whereas explicit methods can, albeit with certain limitations. Table~\ref{tab:explicit-stc} presents a comparative analysis of existing explicit methods for structural topology control. These limitations can be broadly categorized into the following four aspects:
\begin{enumerate}
\item[(1)] Some methods remain restricted to 2D optimization. Among those extended to 3D, most control only the number of cavities and leave tunnels unaddressed.
\item[(2)] Although a few methods can control the numbers of both tunnels and cavities in 3D, they suffer from other drawbacks: the topology-preserving method \citep{heThinningAlgorithmBased2022} requires a prescribed topological configuration, while the QSDV method \citep{hanInequalityConstraintMaximum2022} suffers from relatively low efficiency.
\item[(3)] Nearly all existing methods reduce the number of holes solely by filling, which may cause excessive structural changes in cases where cutting is more appropriate, e.g., a large void separated from the exterior by a thin neck.
\item[(4)] The explicit control of connected components is largely overlooked \citep{zuo2023tunnels}, even though optimized structures generally form a single connected component.
\end{enumerate}\par

\begin{table}[!htb]
  \centering
  \setlength{\extrarowheight}{3pt}
  \caption{Comparative analysis of explicit structural topology control methods.}
  \label{tab:explicit-stc}
  \begin{NiceTabular}{@{}>{\raggedright\arraybackslash}p{4.6cm}lcccc@{}}
    \toprule
     &  & \multicolumn{4}{c}{Number control} \\
    \cmidrule(lr){3-6}
    Explicit method & Framework &
    \begin{tabular}[c]{@{}c@{}}Connected\\ components\end{tabular} &
    \begin{tabular}[c]{@{}c@{}}2D\\ holes\end{tabular} &
    \begin{tabular}[c]{@{}c@{}}3D\\ cavities\end{tabular} &
    \begin{tabular}[c]{@{}c@{}}3D\\ tunnels\end{tabular} \\
    \midrule
    MMC method \citep{zhang2017mmc} & MMC & - & Yes & - & - \\
    \midrule
    Skeleton-based level-set method \citep{zhang2017skeleton} & Level-set & - & Yes & - & - \\
    \midrule
    Euler-Poincar\'e formula-based method \citep{liang2022complexity} & \Block{2-1}{Discrete\\ variable} & - & Yes & Yes & - \\
    QSDV method \citep{hanInequalityConstraintMaximum2022} & & - & Yes & Yes & Yes \\
    \midrule
    Graph and set theory method \citep{zhao2020direct, xiongNewApproachEliminating2020a} & \Block{4-1}{BESO} & - & Yes & Yes & - \\
    Topology-preserving method \citep{heThinningAlgorithmBased2022} & & - & Yes & Yes & Yes \\
    Gauss-Bonnet formula-based method \citep{han2021constraints} & & - & Yes & - & - \\
    Hole-Filling method \citep{heHolefillingBasedApproach2023} & & - & Yes & Yes & Yes \\
    \midrule
    Filtering variable design method \citep{wang2022topological} & \Block{2-1}{SIMP} & - & Yes & - & - \\
    SIMP and MMA-based method \citep{zuo2022explicit, zuo2023tunnels} & & - & Yes & Yes & Yes \\
    \bottomrule
  \end{NiceTabular}
\end{table}

In this study, we develop a persistent homology-based isogeometric topology optimization framework that enables explicit control of both 2D and 3D structural topology within a general SIMP interpolation scheme. 
Persistent homology directly quantifies the topological features of the structure, namely connected components, holes in 2D, and tunnels and cavities in 3D.
Differentiable topological objectives are constructed from the persistence diagram and readily incorporated into the optimization, with the resulting problem solved efficiently
by the method of moving asymptotes (MMA) \citep{svanbergMethodMovingAsymptotes1987}.
Unlike filling-based strategies, the proposed method allows each topological feature to be eliminated either by filling or by cutting, with the choice determined by the region of its corresponding persistence pair. 
This flexibility avoids the excessive structural changes caused by filling-only operations.
Moreover, the design domain and the density field are represented by the same NURBS function, exploiting the exact geometric representation and higher-order continuity of isogeometric analysis (IGA) compared with the finite element method (FEM), thereby improving both accuracy and computational efficiency.
In summary, the main contributions of this study are as follows: 
\begin{itemize}
    \item An explicit and unified topology control framework based on persistent homology is proposed, which accurately captures and constrains the numbers of connected components, holes in 2D, as well as tunnels and cavities in 3D.
    \item A family of differentiable persistence-based topological objectives is designed and readily incorporated into the SIMP framework, 
    with the resulting optimization problem solved efficiently using the MMA optimizer.
    \item A persistence-guided strategy is introduced to determine whether a given topological feature should be eliminated by cutting or by filling, according to the region of its corresponding persistence pair.
\end{itemize}\par

The remainder of this paper is organized as follows. Section~\ref{sec:preliminaries}
reviews the preliminaries on persistent homology and topology-aware optimization.
Section~\ref{sec:methodology} develops the proposed persistence-based topology
control method, which is validated on a series of 2D and 3D benchmarks in Section~\ref{sec:numerical_examples}. 
Finally, Section~\ref{sec:conclusion} concludes the paper.\par

\section{Preliminaries}
\label{sec:preliminaries}

\subsection{Persistent homology of cubical complex}
\label{subsec:ph_complex}

Persistent homology is a fundamental tool in topological data analysis (TDA) for inferring topological features at various scales from a discrete space that is filtered by a real-valued function \citep{chazalIntroductionTopologicalData2021, edelsbrunnerPersistentHomologySurvey2008}.
To characterize the topological information of the structure under consideration, the continuous density field is discretized, yielding a well-defined topological representation.
In both 2D and 3D, the discretized field naturally takes the form of grid-like data. 
For such representations, cubical complexes provide a particularly natural computational framework.\par

\begin{figure}[!htb]
  \centering
  \begin{subfigure}[b]{0.22\textwidth}
    \centering
    \includegraphics[height=2.3cm]{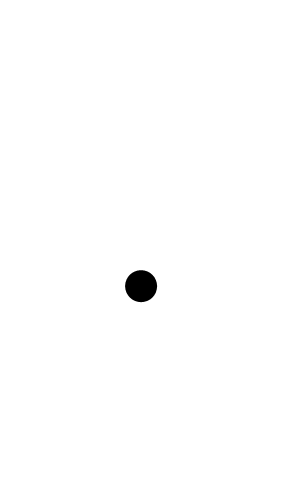}
    \caption{0-cube}
  \end{subfigure}\hfill
  \begin{subfigure}[b]{0.22\textwidth}
    \centering
    \includegraphics[height=2.3cm]{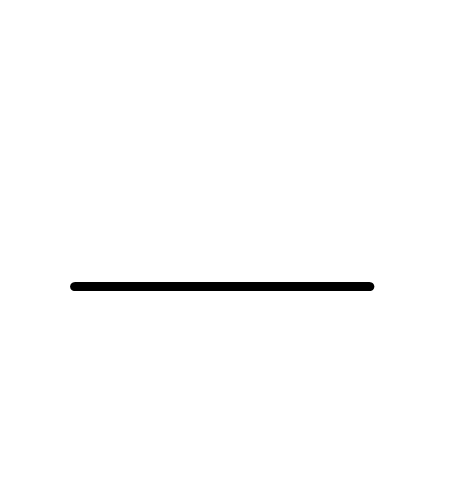}
    \caption{1-cube}
  \end{subfigure}\hfill
  \begin{subfigure}[b]{0.22\textwidth}
    \centering
    \includegraphics[height=2.3cm]{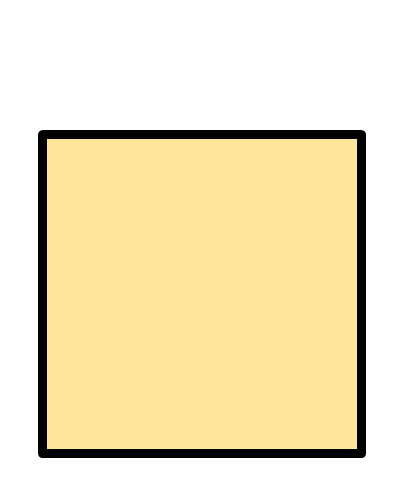}
    \caption{2-cube}
  \end{subfigure}\hfill
  \begin{subfigure}[b]{0.22\textwidth}
    \centering
    \includegraphics[height=2.3cm]{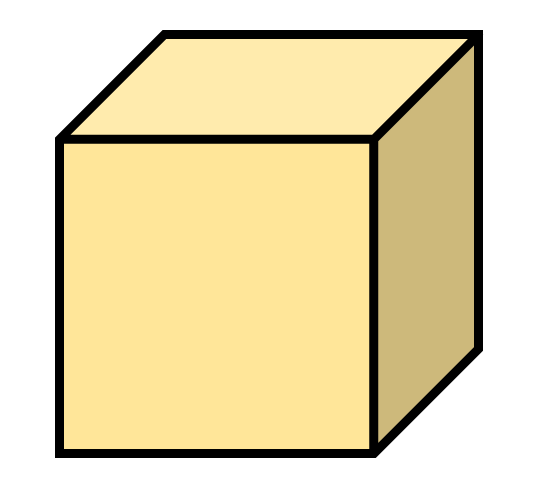}
    \caption{3-cube}
  \end{subfigure}

  \vspace{0.8em}

 \begin{subfigure}[b]{0.25\textwidth}
    \centering
    \includegraphics[width=\linewidth]{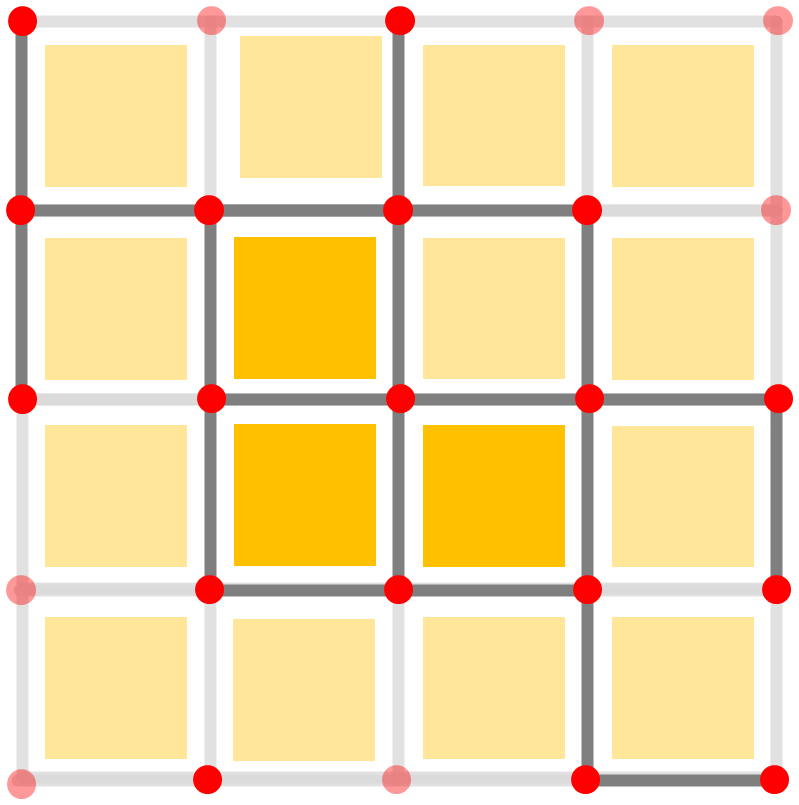}
    \caption{Cubical complex}
  \end{subfigure}\hspace{4em}
  \begin{subfigure}[b]{0.25\textwidth}
    \centering
    \includegraphics[width=\linewidth]{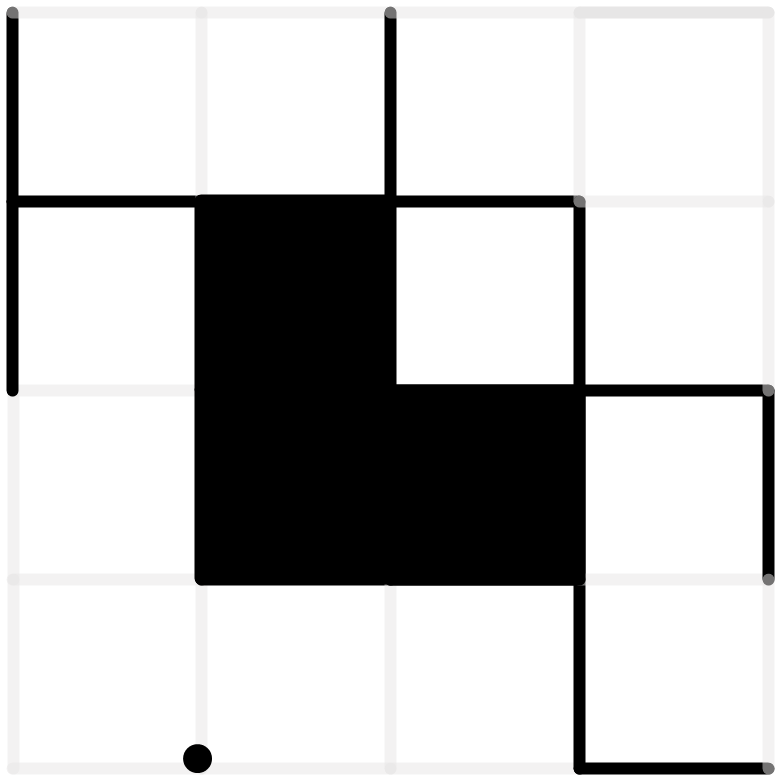}
    \caption{Underlying space}
  \end{subfigure}
  \caption{Cubical complex.
    (a)--(d) Examples of $k$-dimensional cubes for $k=0,1,2,3$.
    (e) A cubical complex with the member cubes highlighted.
    (f) The underlying space obtained from (e) by merging the cubes as a topological space on the plane.}
  \label{fig:cubical-complex}
\end{figure}

In Euclidean space, a $k$-cube $\sigma^k$ is defined as the Cartesian product of $k$ unit intervals, such as vertices (0-cubes), line segments (1-cubes), squares (2-cubes) and spatial cubes (3-cubes).
A cubical complex $\mathcal{K}$ is a collection of $k$-cubes, where the intersection of any two cubes is either empty or a common face of both.
The union of all cubes in $\mathcal{K}$ forms its underlying space $\mathcal{S}$. Fig.~\ref{fig:cubical-complex} shows examples of cubical complexes.\par

In general, the topological features of a cubical complex are classified by dimension: connected components, loops, and voids (cavities) correspond to $dim = 0$, $1$, and $2$, respectively. In the language of algebraic topology, these features are formally characterized by homology groups. For a cubical complex $\mathcal{K}$, the $d$th homology group $H_d(\mathcal{K})$ encodes the $d$-dimensional features listed above, and its rank, denoted by $\beta_d(\mathcal{K})$, is the $d$th Betti number, which equals the number of independent $d$-dimensional features in $\mathcal{K}$.
For a continuum structure, the topology is fully characterized by these Betti numbers. 
In 2D, the structure admits a single type of hole, namely loops, i.e., void regions completely surrounded by solid elements. 
In 3D, holes are distinguished into two types, namely tunnels (through-holes) and cavities (interior holes).
More precisely, $\beta_0$ is the number of connected components of the structure; $\beta_1$ gives the number of holes in 2D, whereas in 3D it is the sum of the number of tunnels and the number of independent loops contained in the cavities \citep{zuo2023tunnels}; and $\beta_2$ gives the number of cavities in 3D and vanishes identically in 2D. \par

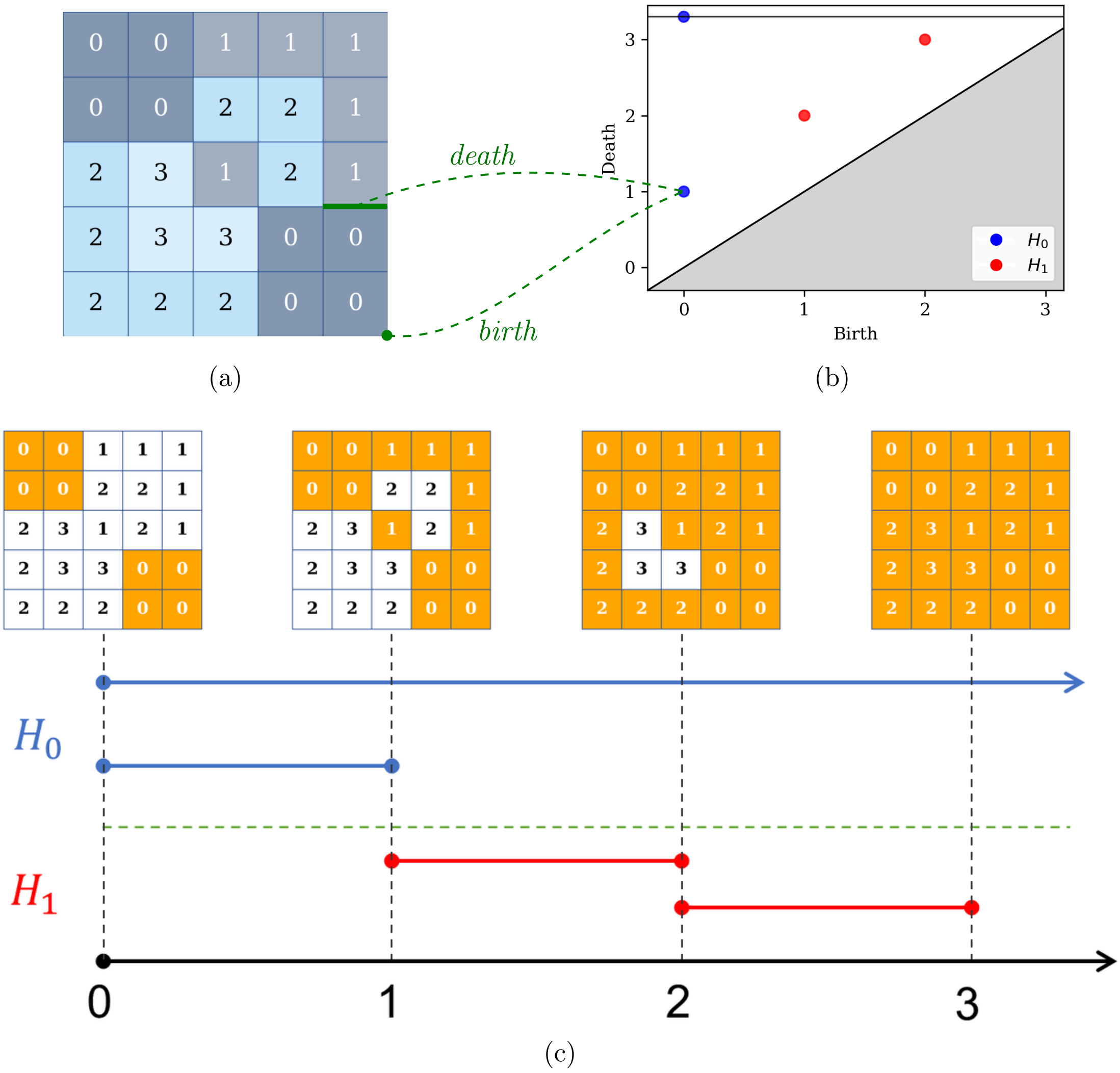
\begin{figure}[!htb]
  \centering
  \begin{tikzpicture}[baseline]
    \node[inner sep=0, anchor=north west] (a) at (0,0)
      {\includegraphics[height=4.7cm]{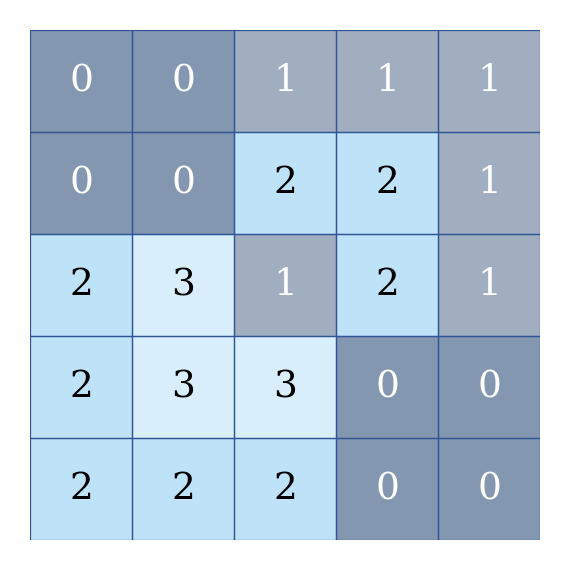}};
    \node[inner sep=0, anchor=north west] (b) at ($(a.north east)+(2.35,0)$)
      {\includegraphics[height=4.7cm]{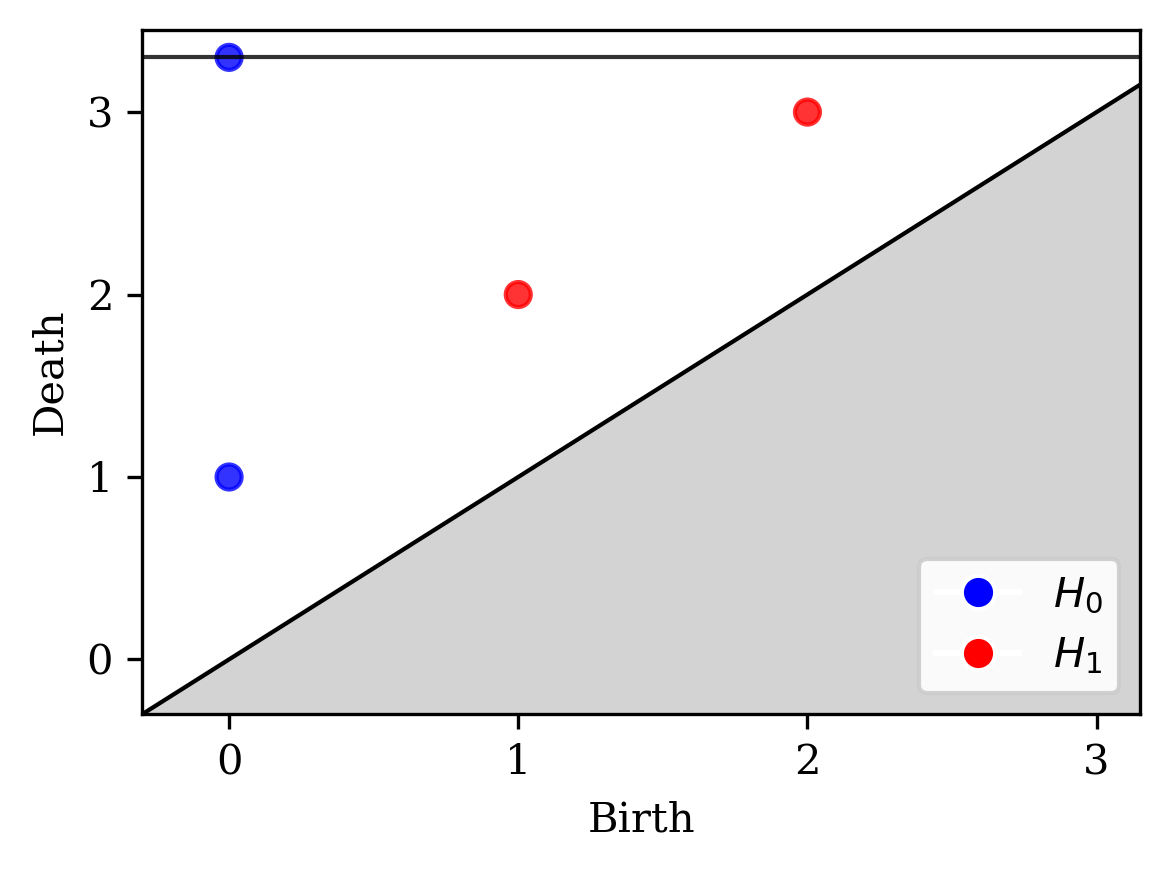}};
    \path
      let \p1=($(a.south east)-(a.north west)$),
          \p2=($(b.south east)-(b.north west)$) in
        coordinate (deathL) at ($(a.north west)+({0.768*\x1},{0.589*\y1})$)
        coordinate (deathR) at ($(a.north west)+({0.945*\x1},{0.589*\y1})$)
        coordinate (birth)  at ($(a.north west)+({0.945*\x1},{0.944*\y1})$)
        coordinate (pair)   at ($(b.north west)+({0.196*\x2},{0.548*\y2})$);
    \coordinate (death) at ($(deathL)!0.5!(deathR)$);
    \draw[green!50!black, line width=2pt] (deathL) -- (deathR);
    \fill[green!50!black] (birth) circle (2pt);
    \draw[dashed, green!50!black, line width=0.75pt]
      (death) .. controls ($(death)+(1.35,0.55)$) and ($(pair)+(-1.55,0.35)$) .. (pair)
      node[pos=0.40, above=3pt, inner sep=0.5pt, text=green!45!black, font=\small\itshape] {death};
    \draw[dashed, green!50!black, line width=0.75pt]
      (birth) .. controls ($(birth)+(1.25,-0.15)$) and ($(pair)+(-1.45,-0.55)$) .. (pair)
      node[pos=0.42, below=7pt, inner sep=0.5pt, text=green!45!black, font=\small\itshape] {birth};
    \node[anchor=north, inner sep=2pt] at ($(a.south)+(0,-0.05)$) {\footnotesize (a)};
    \node[anchor=north, inner sep=2pt] at ($(b.south)+(0,-0.05)$) {\footnotesize (b)};
  \end{tikzpicture}

  \vspace{0.7em}
  \includegraphics[width=0.92\textwidth]{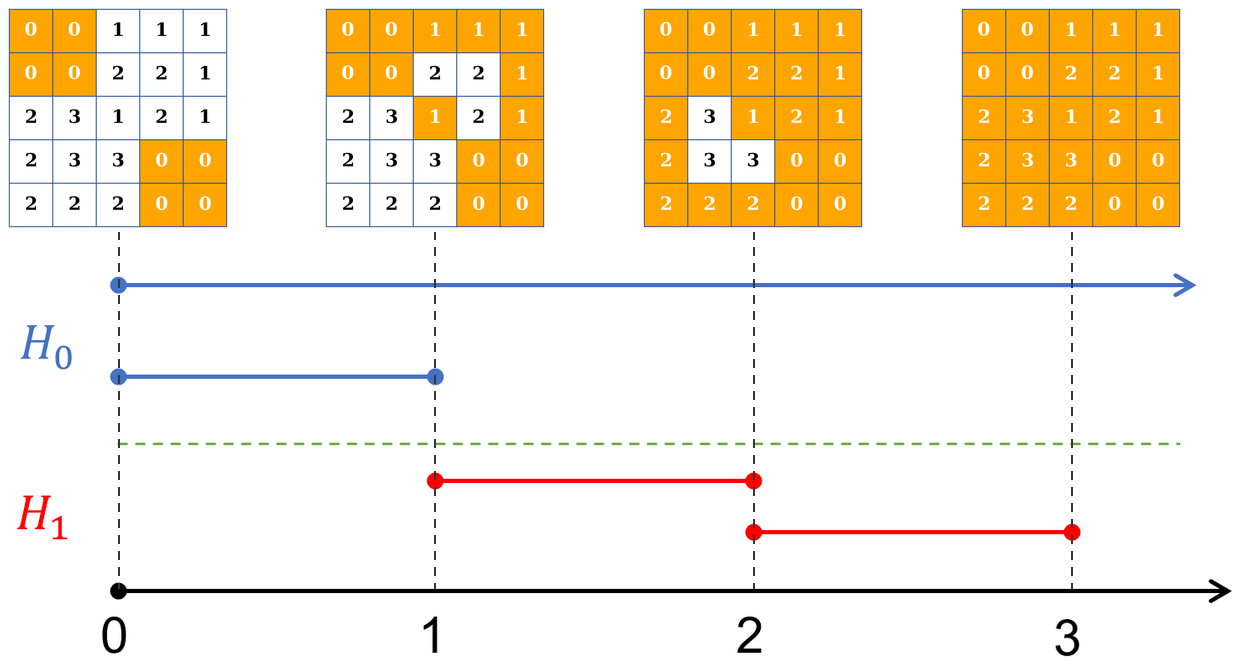}\\[-0.2em]
  {\footnotesize (c)}
  \caption{Illustration of persistent homology.
    (a) A cubical complex with height function values defined on its grid cells.
    (b) A persistence diagram records the persistence pairs.
    (c) A persistence barcode and cubical complexes (marked in orange) in the filtration.
    At time $0$, the underlying space contains two connected components, i.e., two
    0-dimensional topological features. At time $1$, the two connected components merge and one of them
    dies, corresponding to the point $(0,1)$ in
    the persistence diagram and to the interval $[0,1]$ in the persistence barcode. At the same
    time, a loop appears as a 1-dimensional topological feature and vanishes at time $2$,
    yielding the persistence pair $(1,2)$. Since there is always one connected component in the underlying space, the
    point $(0,\infty)$ exists in the persistence diagram. The green dashed line shows that a 0-dimensional persistence pair
    is born at a 0-cube and dies at a 1-cube.
    }
  \label{fig:persistence}
\end{figure}

Topological features defined on discrete structures cannot be directly optimized using discrete metrics.
To overcome this limitation, persistent homology is introduced, 
	which tracks the birth and death of topological 
	features over a series of nested complexes, known as a filtration.
    Let $f:\mathcal{K}\rightarrow \mathbb{R}$ be a real-valued function defined on the cubical complex $\mathcal{K}$. The value assigned to a $k$-cube $\sigma^k$ is taken as the maximum among all its vertices (0-cubes), i.e.,
     \begin{equation}
     f(\sigma^k) = \max_{\sigma^0 \in \sigma^k}f(\sigma^0),
     \label{eq:k_cube_value}
     \end{equation}
    which guarantees that the sublevel sets form genuine subcomplexes. For a threshold $t$, the sublevel cubical complex is defined as $\mathcal{K}_f^t = \{\sigma \in \mathcal{K} \mid f(\sigma) \leq t\}$. Varying $t$ yields the sublevel filtration $\mathcal{F}$, a nested sequence of cubical complexes
    \begin{equation}
    \mathcal{K}_f^{t_0}\subseteq\mathcal{K}_f^{t_1}\subseteq \cdots \subseteq\mathcal{K}_f^{t_n},
    \label{eq:sublevel_filtration}
    \end{equation}
    where the threshold parameter satisfies $t_{0}\leq t_{1}\leq \cdots \leq t_{n}$. 
	As the parameter $t$ increases, topological features of different dimensions are created (birth) and
    destroyed (death), so that the filtration reveals the multi-scale topological structure underlying the data. 
    Persistent homology not only calculates the homology group of each complex in the filtration but also obtains the birth and death of topological features according to well-defined algebraic rules \citep{edelsbrunner2010computational}.
    In this way, topological features are encoded as continuous variables instead of discrete Betti numbers, thereby enabling gradient-based optimization of topology.
    A $k$-dimensional feature that is born in $\mathcal{K}_f^{b}$ and dies in
    $\mathcal{K}_f^{d}$ is represented by a persistence pair $(b,d)$, where $b$ and $d$
    denote the birth and death times. Its lifespan is quantified by the
    persistence $d - b$: a high persistence indicates a long-lasting topological feature, which is potentially more important, whereas a low persistence indicates a transient feature that is more likely to be noise.
    The multiset of all such persistence pairs forms the $k$-dimensional persistence diagram (PD) (see Fig.~\ref{fig:persistence}(b)).
    Equivalently, the lifespan can be visualized as an interval $[b, d)$ along the filtration parameter axis, and the collection of all such intervals is called the $k$-dimensional persistence barcode (see Fig.~\ref{fig:persistence}(c)).\par

\subsection{Topology-aware optimization}
A PD encodes the topology information of the structure, so that optimizing the diagram offers a route to modifying the underlying topology.
This requires a correspondence between the persistence pairs and the elementary cubes of the underlying complex.
Let $f:\mathcal{K}\rightarrow \mathbb{R}$ generate a filtration $\mathcal{F}$ on the cubical complex $\mathcal{K}$.
Each homology generator gives birth at an elementary cube $\sigma$ and dies at another elementary cube $\tau$.
As illustrated in Fig.~\ref{fig:persistence}, the zero-dimensional persistence pair is born at a 0-cube and dies at a 1-cube.
The topological inverse mapping technique \citep{poulenardTopologicalFunctionOptimization2018, carriereOptimizingPersistentHomology2021} realizes the mapping from the persistence pairs to 
the elementary cubes. Given a PD, the topological inverse mapping $\pi_f: \mathrm{PD}\to \mathcal{K}^2$ is defined as
\begin{equation}
    \pi_f(b, d) = (\sigma, \tau), \quad \sigma,\tau \in \mathcal{K}.
\end{equation}
Accordingly, the persistence pair $(b, d)$ can be written as $(f(\sigma), f(\tau))$. 
Furthermore, by Eq.~\eqref{eq:k_cube_value}, there exist 0-cubes $\kappa^0_b$ and $\kappa^0_d$, corresponding to the birth and death values,
such that $f(\sigma) = f(\kappa^0_{b})$ and $f(\tau)=f(\kappa^0_{d})$.
Composing $\pi_f$ with this association maps each persistence pair to a pair of 0-cubes:
\begin{equation}
    \tilde{\pi}_f(b, d) = (\kappa^0_{b}, \kappa^0_{d}).
\end{equation}\par

After establishing the relationship between the underlying complex and persistence pairs, 
the PD can be optimized directly. Assuming that the topological objective function 
$L_{top}$ is defined on the PD, and $f$ is parameterized by $\bm{\alpha}$,
the partial derivative $\partial L_{top} / \partial \bm{\alpha}$ can be obtained through the inverse mapping $\tilde{\pi}_f$
and the differential chain rule \citep{gaoPersistentHomologydrivenOptimization2025a}.\par

In this study, the density field takes a NURBS representation, 
which serves as the real-valued function to induce a sublevel filtration. Consequently, the difficulty of persistence-based optimization lies in
designing an effective and differentiable topological objective.

\section{Methodology}
\label{sec:methodology}

This section presents the proposed persistent homology-based methodology for explicit
topology control in 2D and 3D isogeometric topology optimization.
A single sublevel filtration of the
NURBS density field yields the PDs from which the 0-, 1- and
2-dimensional topological objectives are constructed. These objectives are
incorporated into the minimum compliance problem, and their analytical
sensitivities with respect to the density control coefficients follow from the
topological inverse mapping. 
The topological constraints of different dimensions are then enforced sequentially
within a staged optimization procedure, which is summarized at the end of this
section.

\subsection{Topology control via persistent homology}

\subsubsection{Filtration induced by density field}
Let $\rho$ be a density field represented by NURBS. In 2D, it is written as
\begin{equation}
	\rho(u,v)=\sum_{i=0}^{m}\sum_{j=0}^{n}\frac{\omega_{i,j}N_{i,p}(u)N_{j,q}(v)\rho_{i,j}}{\sum_{i'=0}^{m}\sum_{j'=0}^{n}\omega_{i',j'}N_{i',p}(u)N_{j',q}(v)},	
	\label{eq:density_2D}
\end{equation}
and the 3D counterpart follows directly by including the third parametric
direction,
\begin{equation}
    \rho(u,v,w)= \sum_{i=0}^{m}\sum_{j=0}^{n}\sum_{k=0}^{l}
    \frac{\omega_{i,j,k}N_{i,p}(u)N_{j,q}(v)N_{k,r}(w)\,\rho_{i,j,k}}
         {\sum_{i'=0}^{m}\sum_{j'=0}^{n}\sum_{k'=0}^{l}\omega_{i',j',k'}N_{i',p}(u)N_{j',q}(v)N_{k',r}(w)}.
    \label{eq:density_3D}
\end{equation}

For the binary structure, the density field is thresholded as
\begin{equation}
	\tilde{\rho}=
	\begin{cases}
		1, \quad \rho\geq \bar{\rho}\\
		0, \quad \rho< \bar{\rho},
	\end{cases}
\end{equation}
	where 
	$\bar{\rho}$ is the density threshold. 
	Since $\rho \ge \bar{\rho}$
is equivalent to $-\rho \le -\bar{\rho}$, the solid phase
$\tilde{\rho}=1$ is exactly the sublevel set of the function $f=-\rho$ at the
threshold $t=-\bar{\rho}$.
Therefore, by taking the filtration function as $f=-\rho$, constructing the
cubical complex $\mathcal{K}$ of the discretized density field, and inducing the
associated sublevel filtration,
	the topological features of the structure are captured. 
For each dimension $k$, the resulting $k$-dimensional
PD is given by the multiset of persistence pairs
\begin{equation}
  P_{k} = \{(b_{i,k}, d_{i,k})\}_{i=0}^{N_{k}-1},
  \label{eq:per_pair}
\end{equation}
where $N_{k}$ denotes the total number of $k$-dimensional persistence pairs, and
$b_{i,k}$ and $d_{i,k}$ are the birth and death times of the $i$-th persistence
pair.
  
\subsubsection{Persistence-based optimization}
  
\begin{figure}[!htb]
  \centering
  \begin{tikzpicture}[
      x=1.05cm, y=1.05cm,
      >=Stealth,
      pair/.style={circle, fill=pairred, inner sep=0pt, minimum size=4.6pt},
      opt/.style={optdir, line width=0.9pt, ->},
      guide/.style={dashed, black!65, line width=0.7pt}
    ]
    \colorlet{pairred}{red!72!black}
    \colorlet{optdir}{blue!65!cyan!85!black}
    \colorlet{invalid}{black!9}
    \colorlet{regI}{blue!8}
    \colorlet{regII}{orange!14}

    \coordinate (C) at (2.45,2.45);
    \begin{scope}
      \clip (0,0) rectangle (4.9,4.9);
      \fill[invalid] (0,0) -- (4.9,0) -- (4.9,4.9) -- cycle;
      % I: above the angle bisector.
      \fill[regI] (0,4.9) -- (2.45,4.9) -- (C) -- cycle;
      % II: between the bisector and the line $y=-\bar{\rho}$.
      \fill[regII] (C) -- (0,2.45) -- (0,4.9) -- cycle;
    \end{scope}

    \draw[black!60, line width=0.7pt] (0,0) -- (4.9,4.9);
    \draw[guide] (0,2.45) -- (C);
    \draw[guide] (C) -- (2.45,4.9);
    \draw[guide] (0,4.9) -- (C);

    \draw[->, line width=0.9pt] (0,0) -- (5.25,0);
    \draw[->, line width=0.9pt] (0,0) -- (0,5.25);
    \node[below=2pt] at (2.6,0) {Birth};
    \node[rotate=90, anchor=south] at (-0.18,2.6) {Death};

    \fill (C) circle (1.15pt);
    \node[anchor=north west, font=\small, inner sep=1pt] at ($(C)+(0.06,-0.08)$)
      {$(-\bar{\rho},-\bar{\rho})$};

    \node[pair] (pI) at (1.8,3.675) {};
    \draw[opt] (pI) -- ++(0.45,0);
    \node[font=\large] at (1.7,4.35) {I};

    \node[pair] (pII) at (1.225,3.1) {};
    \draw[opt] (pII) -- ++(0,-0.45);
    \node[font=\large] at (0.55,3.2) {II};

    \node[pair] at (3.1,3.675) {};
    \node[font=\large] at (3.2,4.35) {III};

    \node[pair] at (1.225,1.8) {};
    \node[font=\large] at (0.55,1.7) {IV};

    \node[pair, anchor=center] (lg1) at (5.95,3.55) {};
    \node[anchor=west, font=\small] at ($(lg1.east)+(0.47,0)$) {Persistence pair};
    \draw[opt] (5.64,2.85) -- (6.26,2.85);
    \node[anchor=west, font=\small] at (6.42,2.85) {Optimization direction};
  \end{tikzpicture}
  \caption{PD regions divided by $y=x$, $x=-\bar{\rho}$, $y=-\bar{\rho}$ and $y=-x - 2\bar{\rho}$.
    Persistence pairs in regions I and II belong to the thresholded structure.
    The blue arrows indicate the optimization directions.}
  \label{fig:pd-regions}
\end{figure}
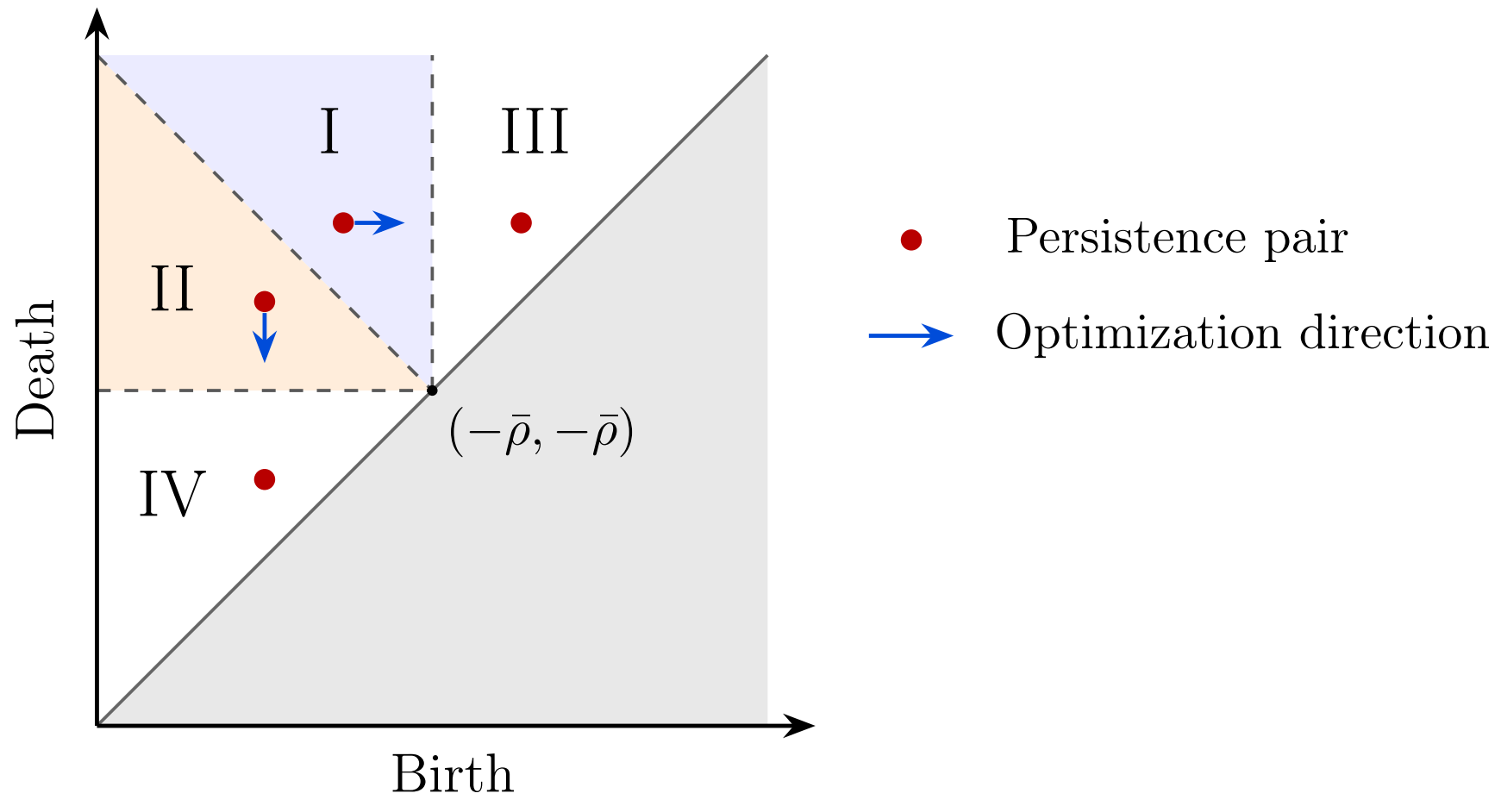

A PD shown in Fig.~\ref{fig:pd-regions} is utilized to illustrate the persistence-based optimization. 
Since the death times of the topological features are always larger than their birth times, all persistence pairs in the PD lie strictly above the diagonal $y=x$.
The diagram is further partitioned by the lines $y=x$, $x=-\bar{\rho}$, $y=-\bar{\rho}$ and $y=-x-2\bar{\rho}$ into four regions.
A topological feature is present in the current structure if and only if it is born before
$-\bar{\rho}$ and dies after $-\bar{\rho}$. Hence, only the persistence pairs in
regions I and II correspond to actual topological features.\par

\begin{figure}[!htb]
  \centering
  \setlength{\tabcolsep}{5pt}
  \providecommand{\cfpair}[2]{
    \begin{tabular}{@{}c@{}}
      \includegraphics[height=3.85cm]{#1}\\[1.5pt]
      \includegraphics[height=3.85cm]{#2}
    \end{tabular}}
  \begin{tabular}{@{}ccc@{}}
    {\small Initial} & {\small Cutting} & {\small Filling} \\[4pt]
    \cfpair{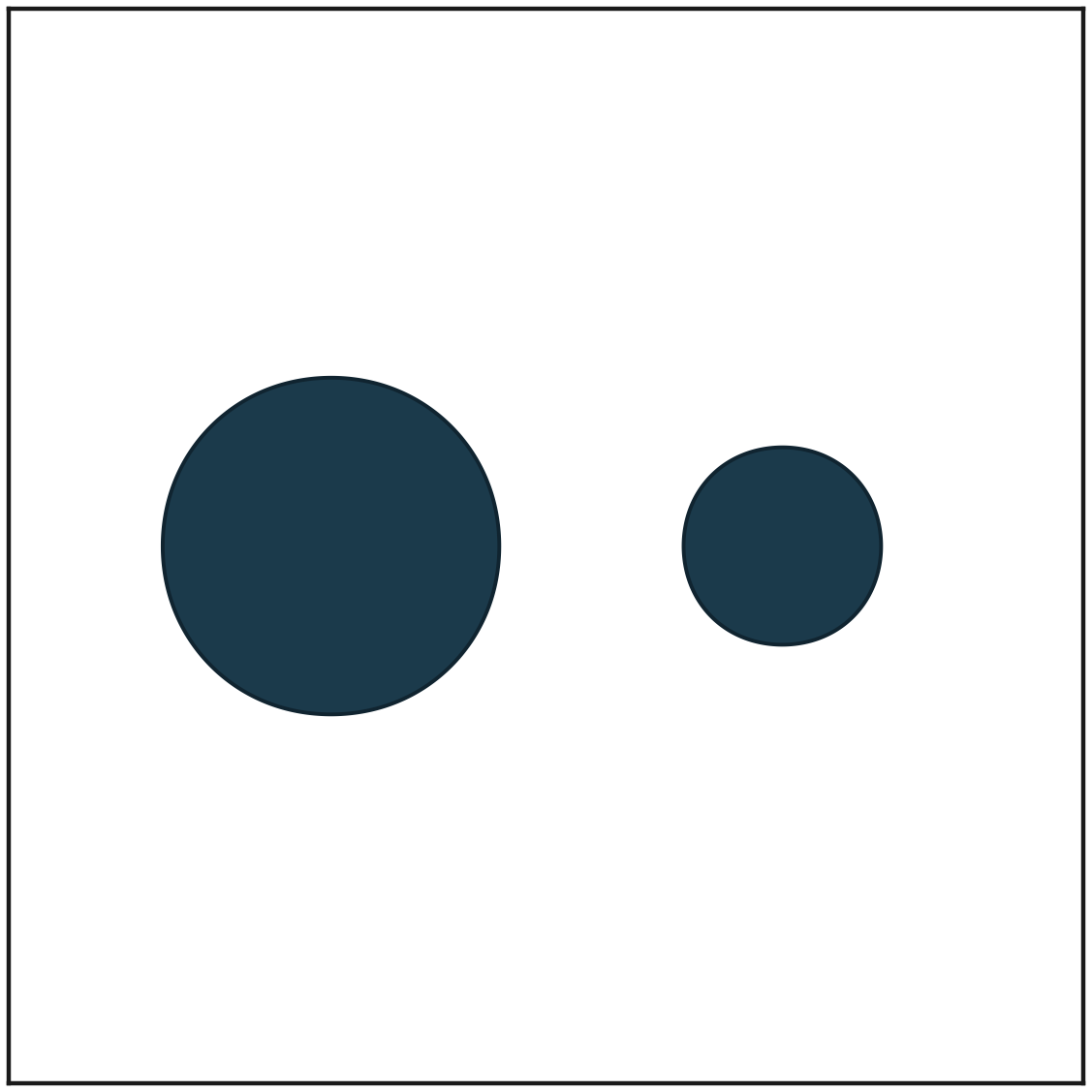}{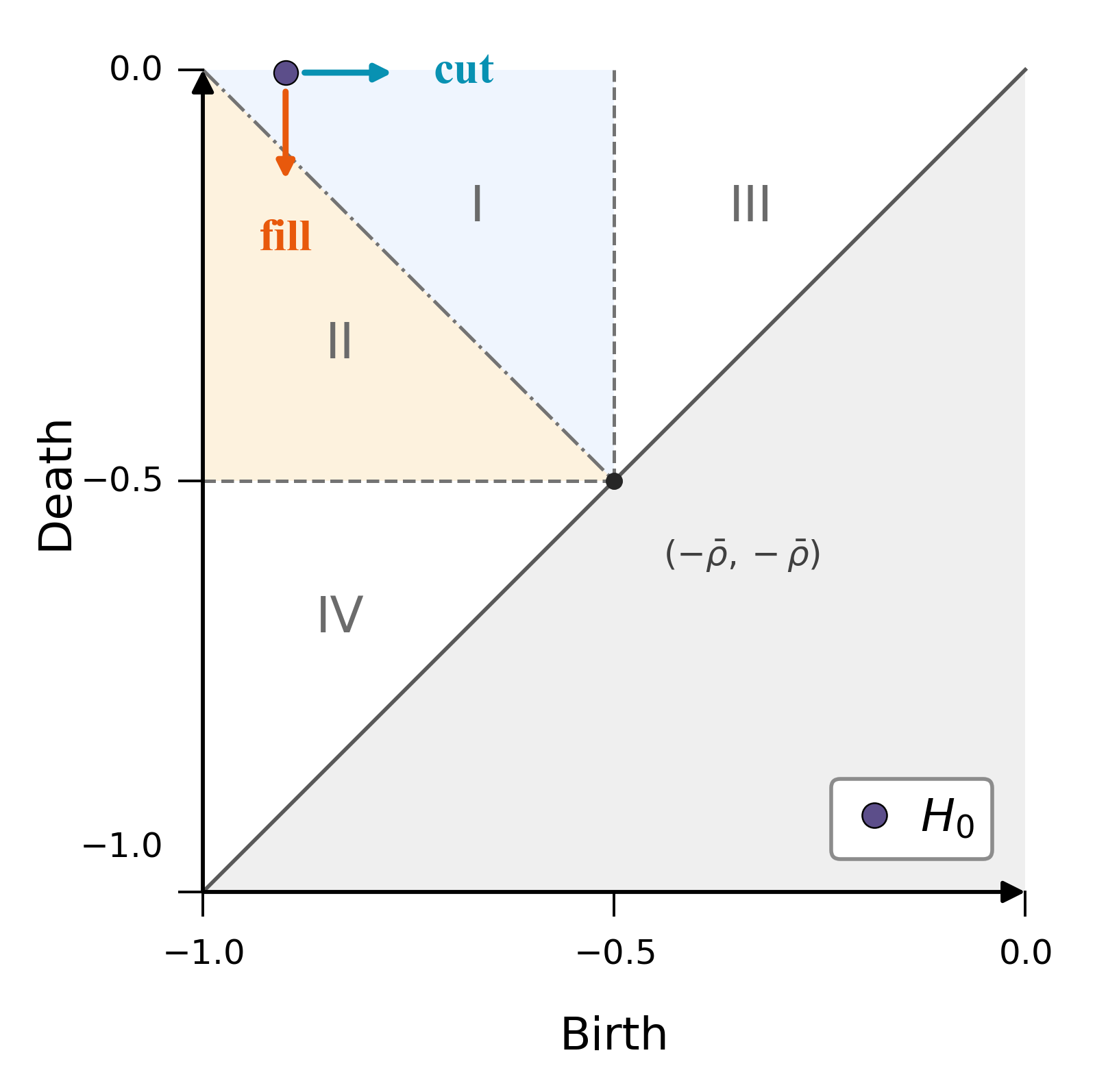}
    & \cfpair{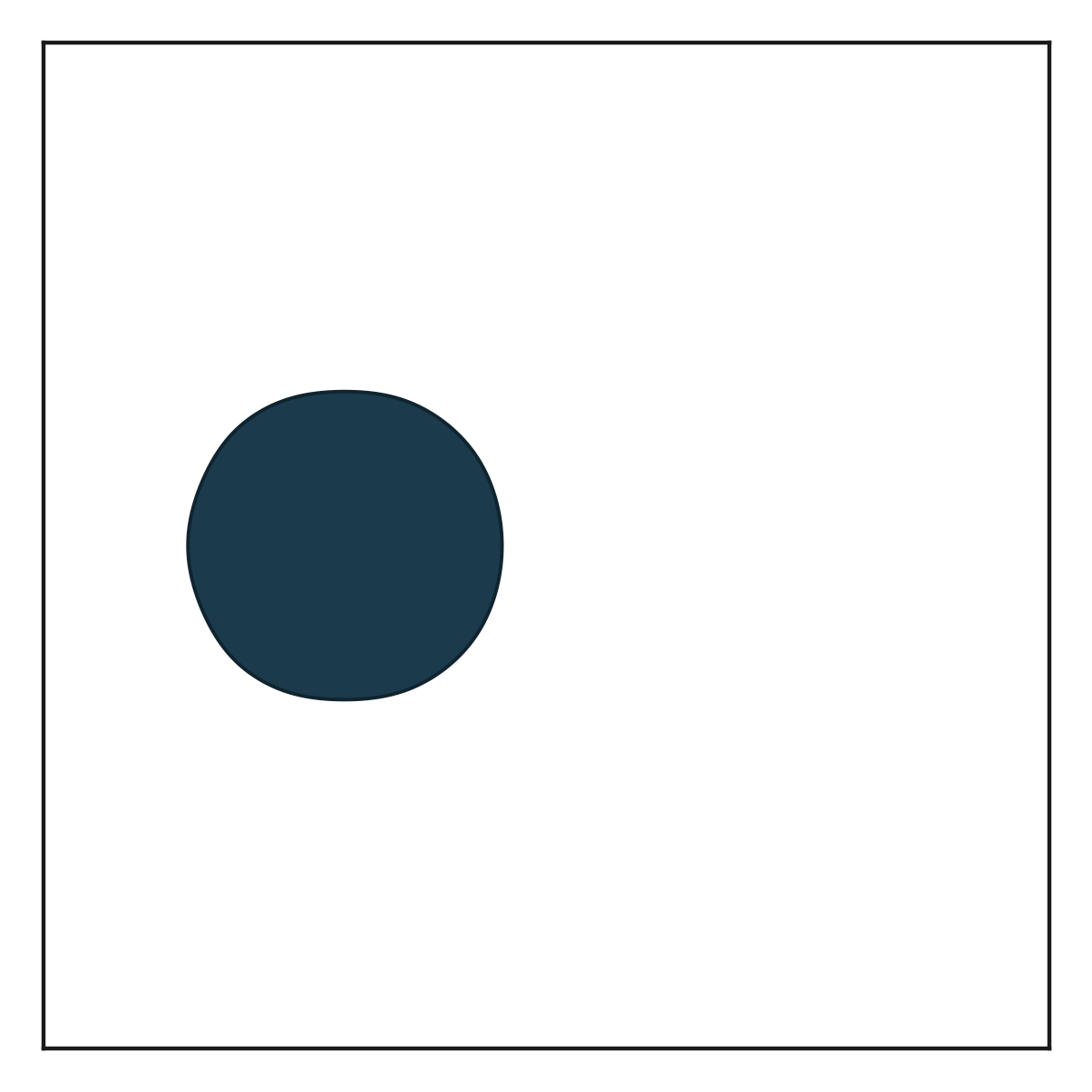}{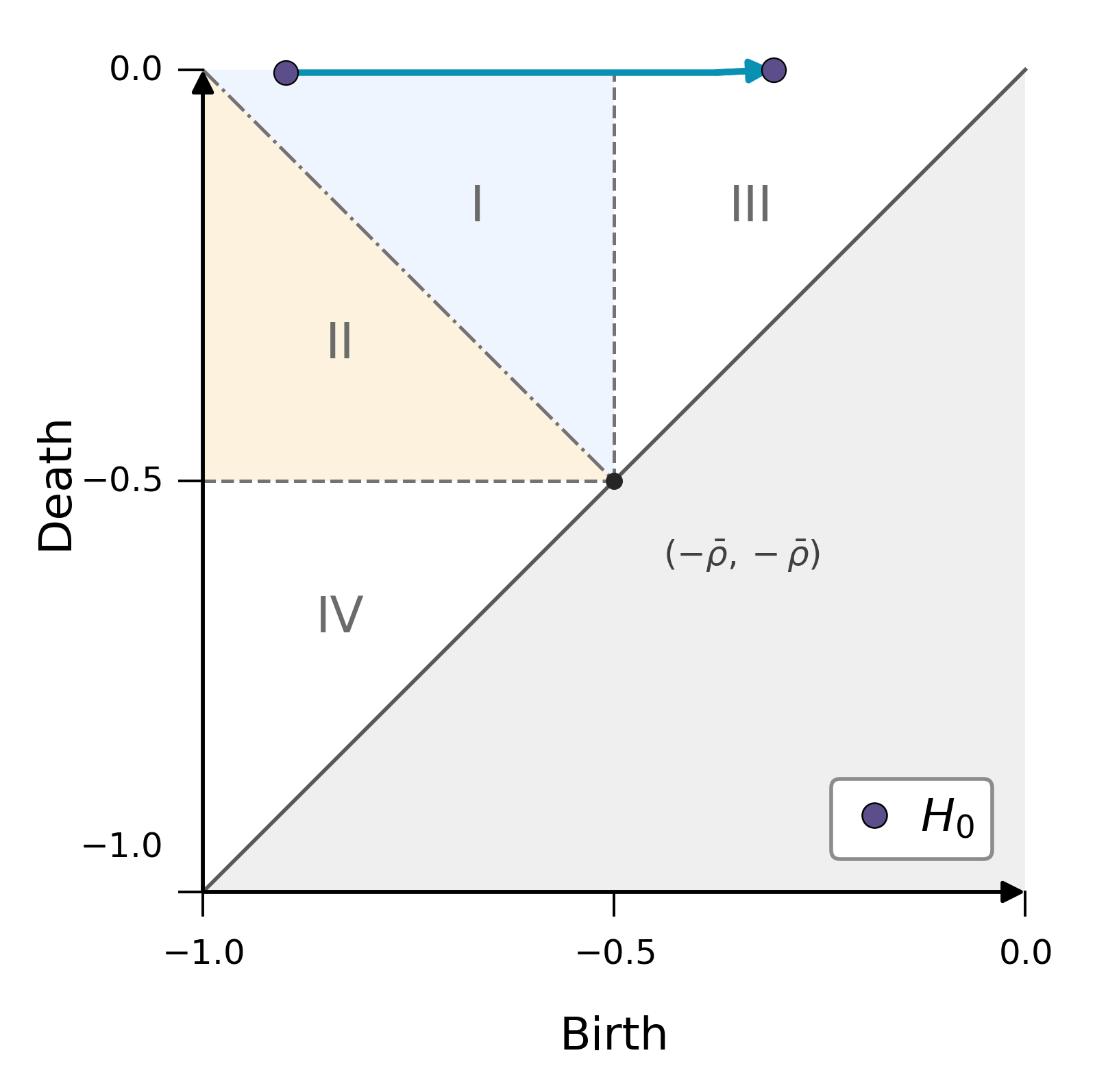}
    & \cfpair{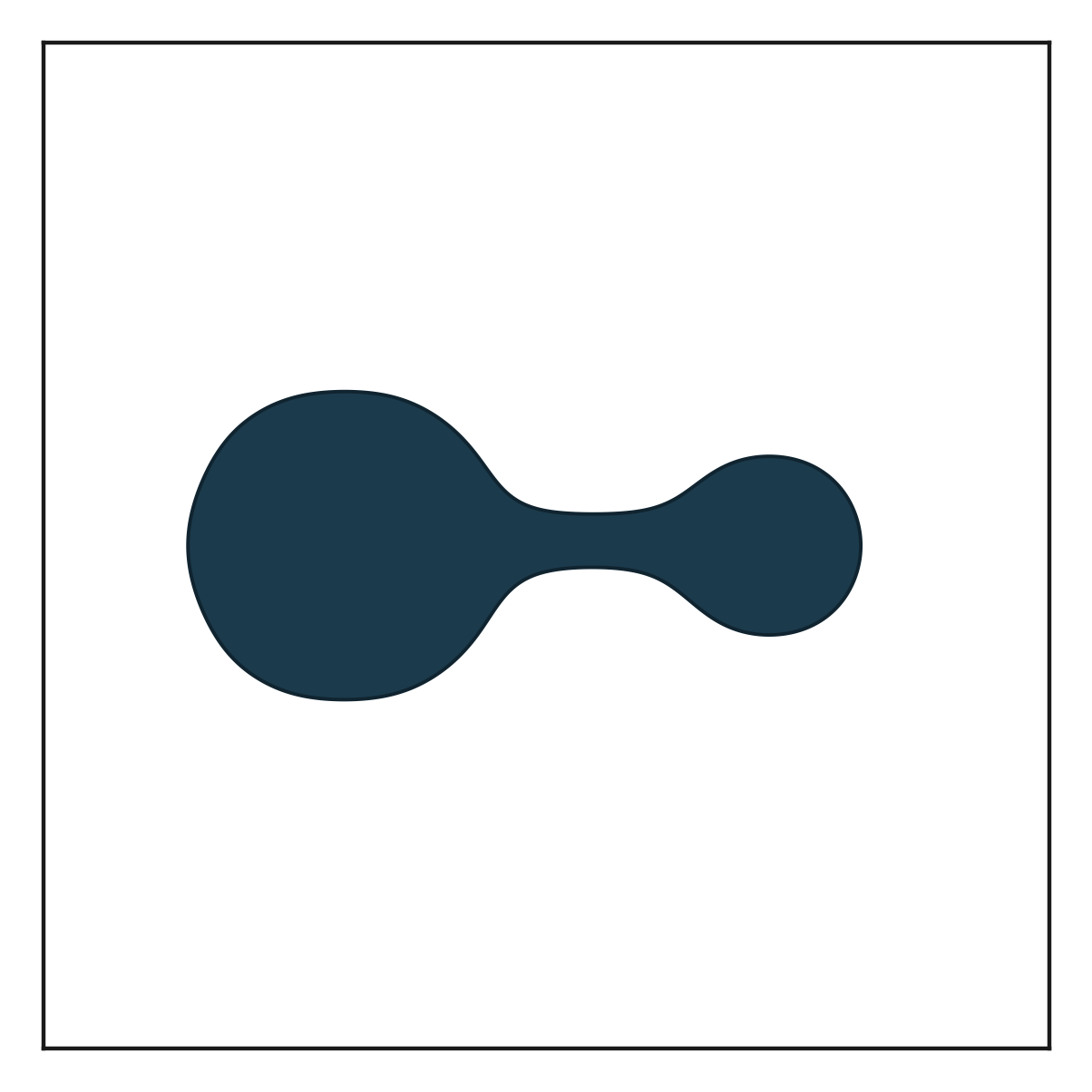}{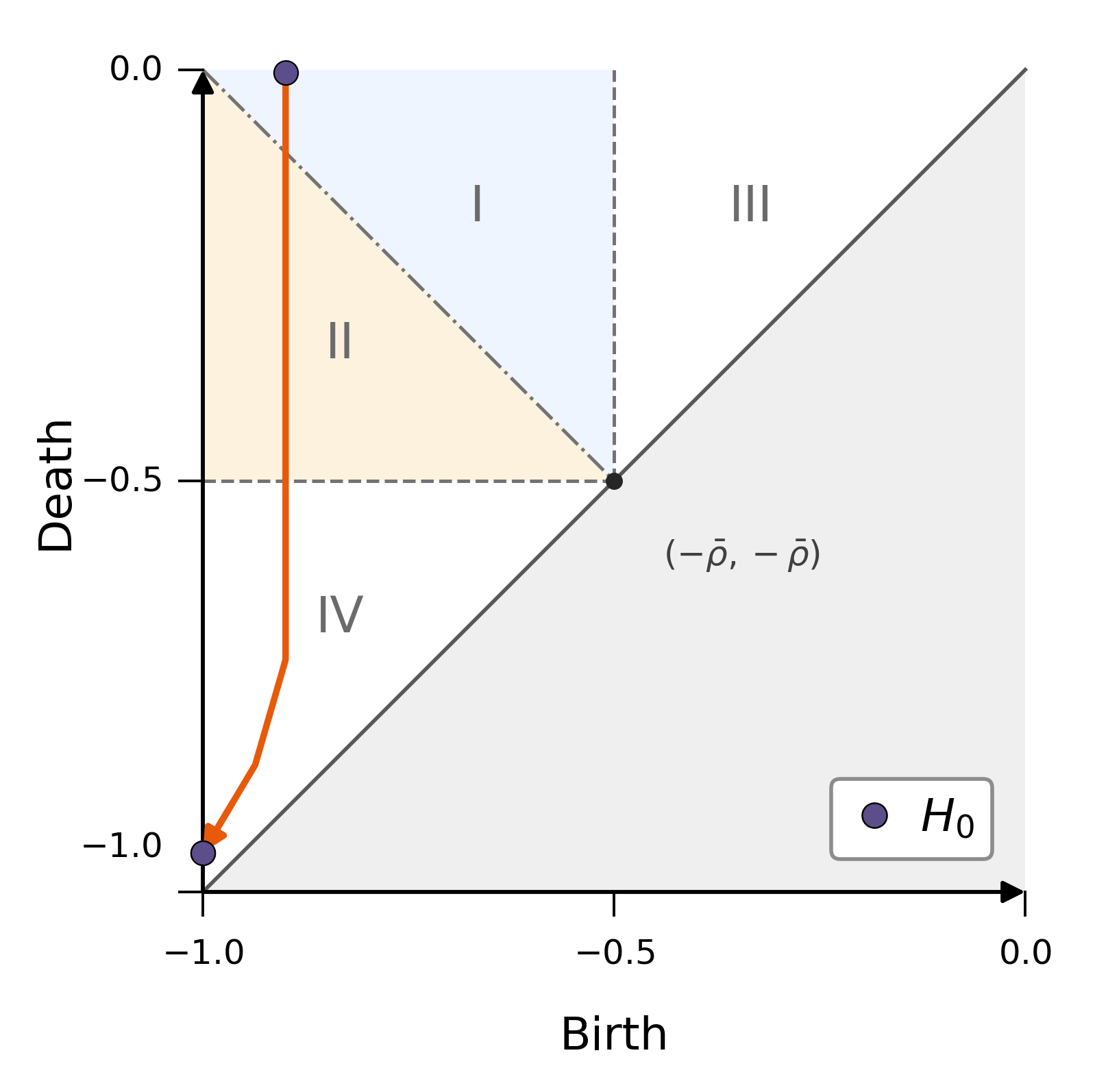} \\
  \end{tabular}
  \caption{Removal of a surplus connected component from a 2D structure by cutting (middle column)
    and filling (right column). In each column, the upper image is the thresholded 
    structure at $\bar{\rho}=0.5$ and the lower one is the corresponding 0-dimensional PD. 
    The arrows in the initial diagram mark the two
    candidate moves of the target pair, whose trajectories during optimization are
    traced by the colored curves in the remaining diagrams. 
    Cutting increases the birth time and drives the pair into region~III and thereby
    removes the surplus component, whereas filling decreases the death time and
    drives it into region~IV and thereby merges the two components into one.}
  \label{fig:cut-fill-components}
\end{figure}

\begin{figure}[!htb]
  \centering
  \setlength{\tabcolsep}{5pt}
  \providecommand{\cfpair}[2]{
    \begin{tabular}{@{}c@{}}
      \includegraphics[height=3.85cm]{#1}\\[1.5pt]
      \includegraphics[height=3.85cm]{#2}
    \end{tabular}}
  \begin{tabular}{@{}ccc@{}}
    {\small Initial} & {\small Cutting} & {\small Filling} \\[4pt]
    \cfpair{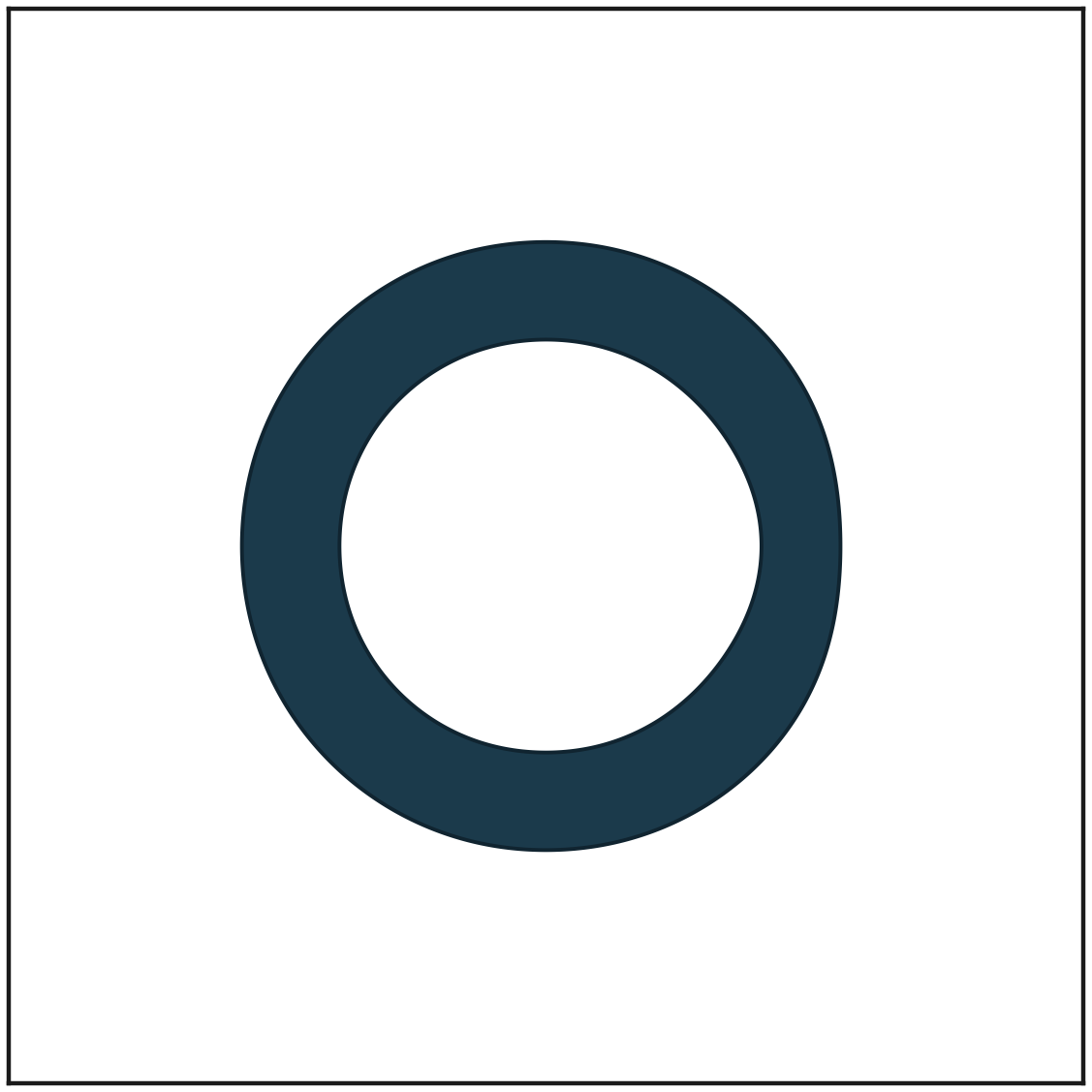}{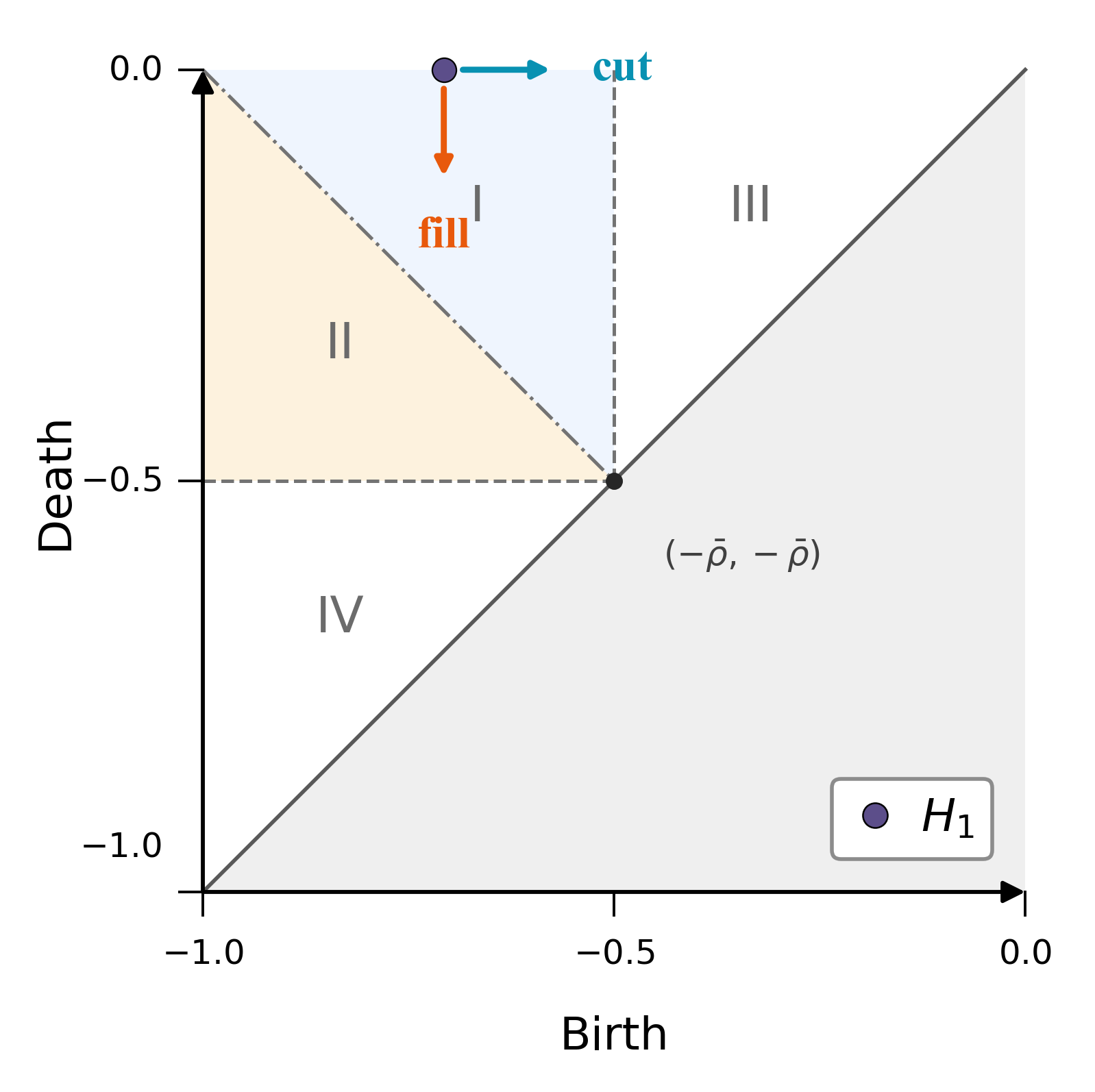}
    & \cfpair{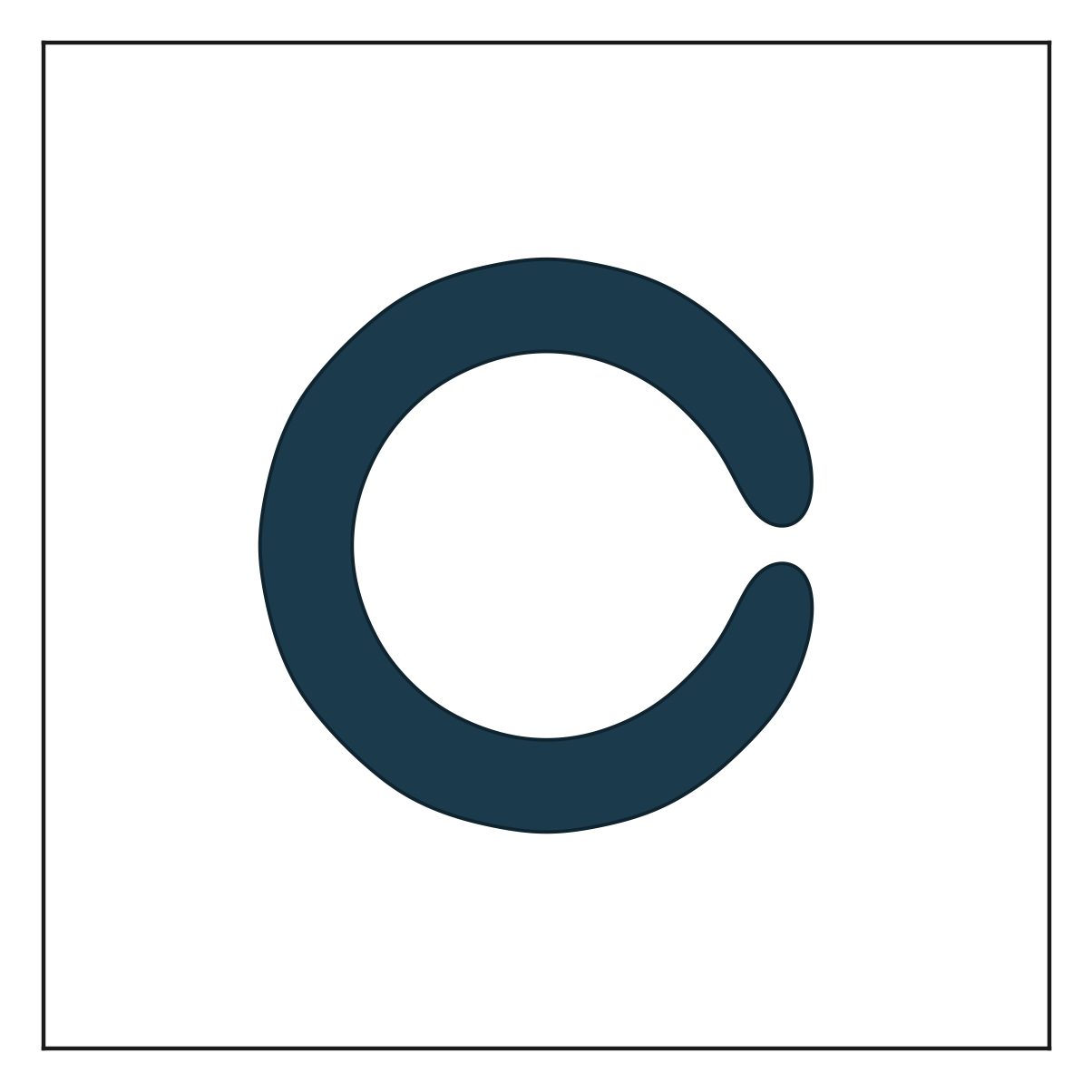}{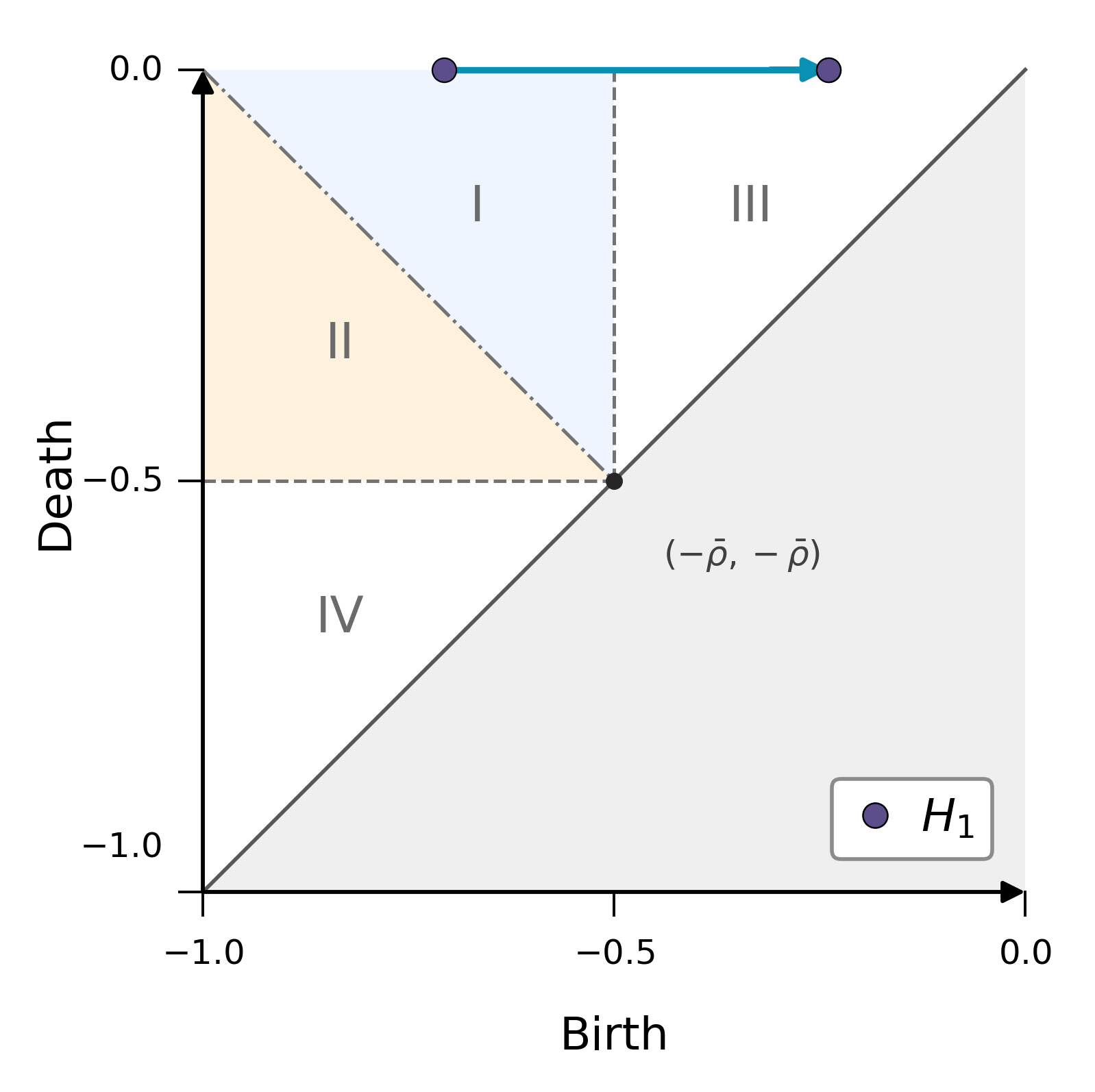}
    & \cfpair{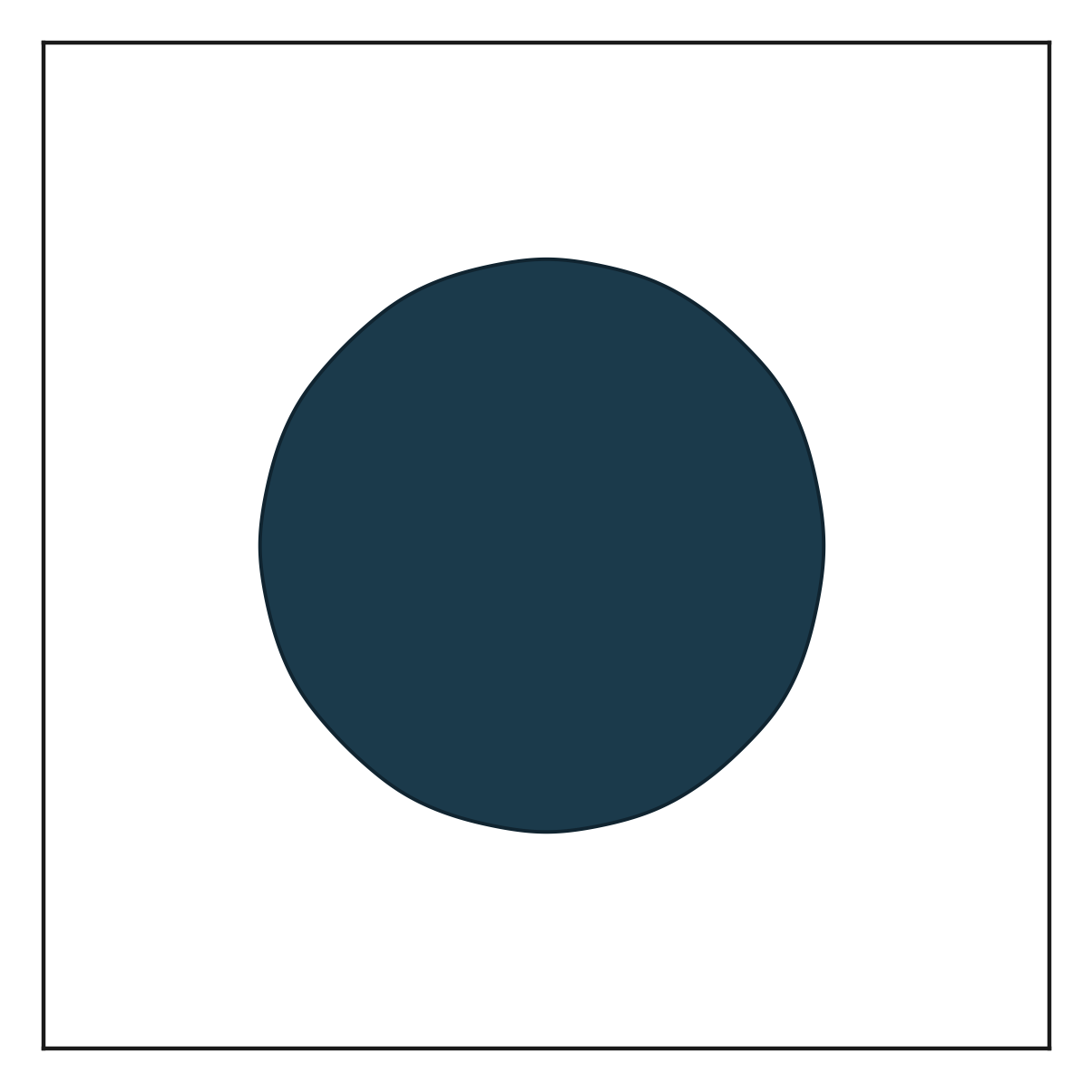}{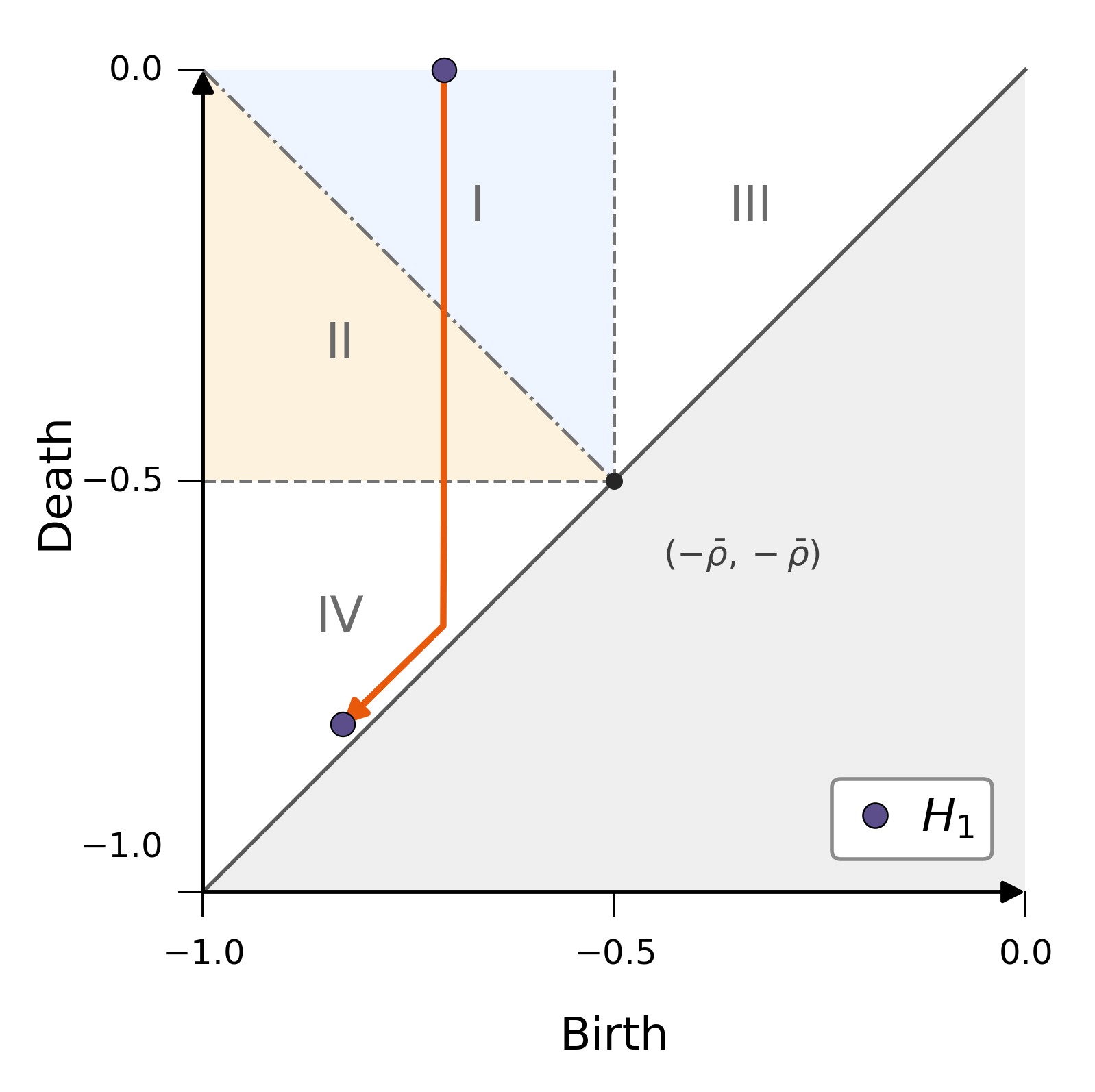} \\
  \end{tabular}
  \caption{Removal of a hole from a 2D structure by cutting (middle column) and filling (right
    column). In each column, the upper image is the thresholded structure at
    $\bar{\rho}=0.5$ and the lower one is the corresponding 1-dimensional PD.
    The arrows in the initial diagram mark the two
    candidate moves of the target pair, whose trajectories during optimization are
    traced by the colored curves in the remaining diagrams. 
    Cutting increases the birth time and drives the pair into region~III and thereby
    cuts the loop open, whereas filling decreases the death time and drives it into
    region~IV and thereby fills its interior.}
  \label{fig:cut-fill-holes}
\end{figure}

To eliminate surplus topological features, the corresponding persistence pairs should be
translocated out of regions I and II. 
Guided by a minimum-movement principle, persistence pairs in region I are moved to region III while persistence pairs in region II are moved to region IV,
so as to avoid a substantial structure change caused by a considerable movement of persistence pairs.
Furthermore, these movements admit a clear geometric interpretation. 
As illustrated in Figs.~\ref{fig:cut-fill-components} and~\ref{fig:cut-fill-holes}, moving a persistence pair to
region III amounts to maximizing its birth time, which corresponds to cutting the
associated topological feature, whereas moving a persistence pair to region IV amounts to minimizing its
death time, which corresponds to the filling operation. Despite the same topological outcome,
they can have very different impact on the geometry of the structure, thereby affecting the structural performance.
Ideally, topology control should remove as many topological features as possible while changing the geometry as little as possible \citep{zengCutFillGlobal2020}.\par

More specifically, since the filtration is induced by the sublevel sets of
$f = -\rho$, maximizing the birth time $b = f(\sigma)$ of a persistence pair can
decrease the density at the birth cube $\sigma$, i.e., at the
location where the corresponding feature is created. In this way, a connected
component is cut from the structure, a tunnel is eliminated by cutting
open its body, and a cavity is removed by cutting a tunnel that connects it to
the exterior. In contrast, minimizing the death time $d = f(\tau)$ 
increases the density at the death cube $\tau$, i.e., at the location where
the feature is destroyed. Consequently, a connected component disappears by
merging with another component, whereas a hole is removed by filling its
body: in 2D this is a loop, and in 3D either a tunnel or a cavity.\par

\subsubsection{Zero-dimensional topological objective function}
The presence of extra connected components may affect manufacturability. 
A 0-dimensional topological objective is therefore introduced to enforce
structural connectivity. Since the filtration is induced by $f=-\rho$, the
sublevel set at $t=+\infty$ is the entire design domain $\Omega$, which is
connected. Consequently, the 0-dimensional PD always contains
one essential pair $(b,+\infty)$ with $b\in\mathbb{R}$, representing the
component that never disappears. To ensure that no additional connected component
exists in the structure, all \emph{finite}
persistence pairs located in regions~I and~II of the 0-dimensional PD must be
eliminated, while the essential pair $(b,+\infty)$ is retained.\par

\begin{figure}[!htb]
  \centering
  \begin{subfigure}[b]{0.35\textwidth}
    \centering
    \includegraphics[width=\textwidth]{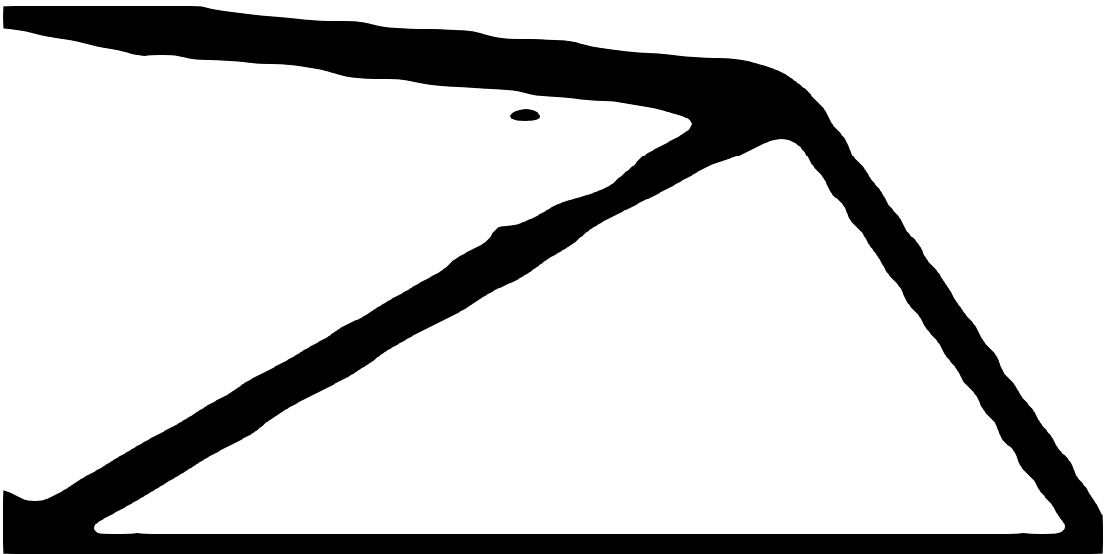}
    \caption{Before filling.}
    \label{fig:short-h0-filling-a}
  \end{subfigure}
  \hspace{3em}
  \begin{subfigure}[b]{0.35\textwidth}
    \centering
    \includegraphics[width=\textwidth]{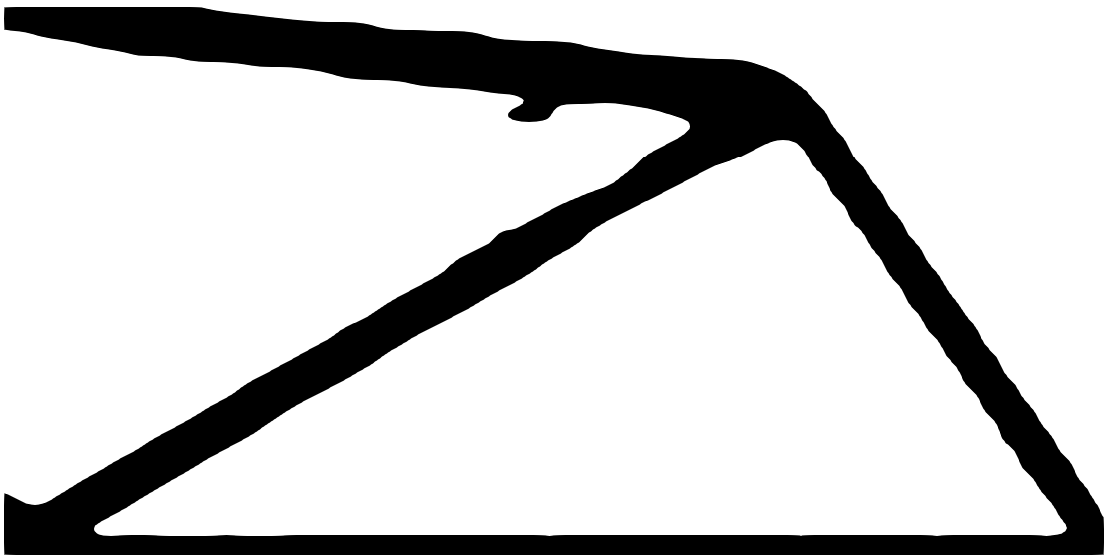}
    \caption{After filling.}
    \label{fig:short-h0-filling-b}
  \end{subfigure}
  \caption{Elimination of an extra connected component in a 2D short beam by
    filling.}
  \label{fig:short-h0-filling}
\end{figure}

In practical topology optimization, however, eliminating such a pair by
minimizing its death time, i.e., by connecting the corresponding connected
component to another one through added material, tends to generate redundant
branches that do not lie on the critical load path, which necessitates
additional post-processing, as illustrated in Fig.~\ref{fig:short-h0-filling}. Therefore, extra connected
components are removed exclusively by the cutting strategy, regardless of the
region in which the corresponding persistence pair lies.\par

The 0-dimensional persistence pairs in regions~I and~II denoted by $P_{0}^{I,II}$ can be extracted from $P_{0}$ defined in Eq.~\eqref{eq:per_pair}:
\begin{equation}
		P_{0}^{I,II} =\{(b_{i,0},d_{i,0})\}_{i=0}^{N_{0}^{I,II}-1},
\end{equation}
where $N_0^{I,II}$ equals the total number of persistence pairs within regions~I and~II except the essential pair.
The required relocation of 0-dimensional persistence pairs from regions~I and~II to region~III,
is achieved by defining the objective function as
\begin{equation}
		C_{top}^{0}=-\sum_{i=0}^{N_{0}^{I,II}-1}(b_{i,0}^{I,II} + \bar{\rho}).
		\label{eq:zero_optimization}
\end{equation}\par

\subsubsection{One-dimensional topological objective function}
To control the structural complexity, we impose an upper bound $\bar{N}_{1}$
on the number of holes in 2D structures and tunnels in 3D structures.
The excess holes or tunnels are then removed in an order determined by their
geometric significance, which is characterized as follows.\par

For a 2D structure, we first employ the topological inverse mapping to locate
the death positions of the 1-dimensional topological features in the parametric
domain. Starting from each death position, a breadth-first search (BFS) is
performed over the discretized density field to compute the hole area,
approximated by the number of discrete elements it contains.
For a 3D structure, this strategy is no longer applicable: a tunnel is a through-hole 
and therefore bounds no region that a BFS
can fill. An alternative is to adopt a persistent 1-cycle
\citep{deyPersistent1CyclesDefinition2019} as the representative cycle, such as
the volume-optimal cycle \citep{obayashiVolumeOptimalCycleTightest2018}, and to
measure the area enclosed by that tunnel. For simplicity and computational
efficiency, however, we instead sort the tunnels by the persistence
$d_{i,1}-b_{i,1}$ of the corresponding 1-dimensional persistence pair.\par

The 1-dimensional persistence pairs in regions~I and~II are denoted by
\begin{equation}
  P_{1}^{I,II} = \{(b_{i,1}, d_{i,1})\}_{i=0}^{N_{1}^{I,II}-1},
\end{equation}
where $N_{1}^{I,II}$ is the total number of 1-dimensional persistence pairs in
regions~I and~II, i.e., the total number of holes (2D) or tunnels (3D).
These pairs are sorted in ascending order of the hole area for 2D structures and
of the persistence $d_{i,1}-b_{i,1}$ for 3D structures, respectively.
When $N_{1}^{I,II} > \bar{N}_{1}$, the first $N_{1}^{I,II}-\bar{N}_{1}$ pairs
are eliminated. Grouping them according to the region they belong to yields
\begin{equation}
  \begin{aligned}
    P_{1}^{I,\mathrm{rem}}  &= \{(b_{i,1}^{I},  d_{i,1}^{I} )\}_{i=0}^{N_{1}^{I,\mathrm{rem}}-1},\\
    P_{1}^{II,\mathrm{rem}} &= \{(b_{i,1}^{II}, d_{i,1}^{II})\}_{i=0}^{N_{1}^{II,\mathrm{rem}}-1},
  \end{aligned}
\end{equation}
where $N_{1}^{I,\mathrm{rem}}$ and $N_{1}^{II,\mathrm{rem}}$ denote the numbers of
persistence pairs to be removed from regions~I and~II, respectively, satisfying
$N_{1}^{I,\mathrm{rem}} + N_{1}^{II,\mathrm{rem}} = N_{1}^{I,II} - \bar{N}_{1}$.\par

When $N_{1}^{I,II} \leq \bar{N}_{1}$, the density of the existing holes or
tunnels is frozen so as to prevent the features that are intended to be
preserved from being removed. Specifically, the interior and the boundary of each
hole or tunnel are identified through the topological inverse mapping, and the
density is decreased in the interior and increased on the boundary, respectively.
The corresponding persistence pairs are collected as
\begin{equation}
  P_{1}^{\mathrm{pres}} = \{(b_{i,1}^{I,II}, d_{i,1}^{I,II})\}_{i=0}^{N_1^{I,II}-1}.
\end{equation}

In addition, holes or tunnels removed in one iteration tend to reappear at the
same locations in subsequent iterations---most notably when a hole is repeatedly
opened and closed at the same location, or when its interior is alternately
filled in and emptied out. To suppress this oscillatory behavior, the
persistence pairs in regions~III and~IV are optimized in the same manner and are
grouped as
\begin{equation}
  \begin{aligned}
    P_{1}^{III} &= \{(b_{i,1}^{III}, d_{i,1}^{III})\}_{i=0}^{N_{1}^{III}-1},\\
    P_{1}^{IV}  &= \{(b_{i,1}^{IV},  d_{i,1}^{IV} )\}_{i=0}^{N_{1}^{IV}-1},
  \end{aligned}
\end{equation}
where $N_{1}^{III}$ and $N_{1}^{IV}$ denote the numbers of persistence pairs in
regions~III and~IV, respectively.\par

Assembling the above contributions, the 1-dimensional topological objective
function is defined as
\begin{equation}
  \begin{aligned}
    C_{\mathrm{top}}^{1}
    ={}& -\sum_{i=0}^{N_{1}^{I,\mathrm{rem}}-1}  (b_{i,1}^{I} + \bar{\rho})
        + \sum_{i=0}^{N_{1}^{II,\mathrm{rem}}-1} (d_{i,1}^{II} + \bar{\rho})\\
      &+ \sum_{i=0}^{N_1^{I,II}-1}
         \big[(b_{i,1}^{I,II} + \rho_{\max}) - (d_{i,1}^{I,II} + \rho_{\min})\big]\\
      &- \sum_{i=0}^{N_{1}^{III}-1} (b_{i,1}^{III} + \rho_{\min})
        + \sum_{i=0}^{N_{1}^{IV}-1}  (d_{i,1}^{IV} + \rho_{\max}),
  \end{aligned}
  \label{eq:one_optimization}
\end{equation}
where $\rho_{\min}$ and $\rho_{\max}$ are the lower and upper bounds of the
density field.

\subsubsection{Two-dimensional topological objective function}
For a 3D structure, the presence of existing cavities affects the evaluation of
the number of tunnels, so that the latter cannot be measured directly from the
1-dimensional Betti number $\beta_{1}$, as discussed in
Section~\ref{subsec:ph_complex}. On the other hand, as noted in
Section~\ref{sec:introduction}, cavities in a 3D structure are generally
undesirable, since they are difficult and even impossible to manufacture
by additive manufacturing. The aim is therefore to eliminate all cavities in the
3D structure before the tunnels are controlled.\par

The elimination order again follows from the geometric significance of the
features. The death positions of the 2-dimensional topological features are
first located by the topological inverse mapping; starting from each of them, a
BFS is performed over the discretized density field to compute the cavity volume,
approximated by the number of discrete elements it contains. Unlike a tunnel, a
cavity bounds a closed region, so this strategy remains applicable in 3D. The
2-dimensional persistence pairs in regions~I and~II are denoted by
\begin{equation}
  P_{2}^{I,II} = \{(b_{i,2}, d_{i,2})\}_{i=0}^{N_{2}^{I,II}-1},
\end{equation}
where $N_{2}^{I,II}$ is the total number of such pairs, i.e., the total number of
cavities. Since no cavity is to be preserved, all of them are removed, in
ascending order of the cavity volume so that the smaller ones are handled first. Grouping them
according to the region they belong to yields
\begin{equation}
  \begin{aligned}
    P_{2}^{I,\mathrm{rem}}  &= \{(b_{i,2}^{I},  d_{i,2}^{I} )\}_{i=0}^{N_{2}^{I,\mathrm{rem}}-1},\\
    P_{2}^{II,\mathrm{rem}} &= \{(b_{i,2}^{II}, d_{i,2}^{II})\}_{i=0}^{N_{2}^{II,\mathrm{rem}}-1},
  \end{aligned}
\end{equation}
where $N_{2}^{I,\mathrm{rem}}$ and $N_{2}^{II,\mathrm{rem}}$ satisfy
$N_{2}^{I,\mathrm{rem}} + N_{2}^{II,\mathrm{rem}} = N_{2}^{I,II}$.\par

Cavities removed in one iteration tend to reappear in subsequent ones as well.
To suppress this oscillatory behaviour, the pairs in regions~III and~IV are
optimized in the same manner and are grouped as
\begin{equation}
  \begin{aligned}
    P_{2}^{III} &= \{(b_{i,2}^{III}, d_{i,2}^{III})\}_{i=0}^{N_{2}^{III}-1},\\
    P_{2}^{IV}  &= \{(b_{i,2}^{IV},  d_{i,2}^{IV} )\}_{i=0}^{N_{2}^{IV}-1},
  \end{aligned}
\end{equation}
where $N_{2}^{III}$ and $N_{2}^{IV}$ denote the numbers of pairs in regions~III
and~IV.\par

Assembling the above contributions, the 2-dimensional topological objective
function is defined as
\begin{equation}
  \begin{aligned}
    C_{\mathrm{top}}^{2}
    ={}& -\sum_{i=0}^{N_{2}^{I,\mathrm{rem}}-1}  (b_{i,2}^{I} + \bar{\rho})
        + \sum_{i=0}^{N_{2}^{II,\mathrm{rem}}-1} (d_{i,2}^{II} + \bar{\rho})\\
      &- \sum_{i=0}^{N_{2}^{III}-1} (b_{i,2}^{III} + \rho_{\min})
        + \sum_{i=0}^{N_{2}^{IV}-1}  (d_{i,2}^{IV} + \rho_{\max}).
  \end{aligned}
  \label{eq:two_optimization}
\end{equation}
Compared with
Eq.~\eqref{eq:one_optimization}, the preservation term is absent, since no
2-dimensional pair is retained. Once all cavities are removed,
$\beta_{1}$ reduces exactly to the number of tunnels, so that
Eq.~\eqref{eq:one_optimization} applies directly to their control.

\subsection{Problem statement and sensitivity analysis}
\label{subsec:problem_statement_sensitivity_analysis}

The classic minimum compliance problem, augmented with the proposed
0-, 1- and 2-dimensional topological constraints, is stated as follows:
\begin{equation}
  \label{eq:optimization_problem}
  \begin{aligned}
    \text{Find:} \quad
    & \rho_{i,j}\ (2\text{D}) \ \text{or}\ \rho_{i,j,k}\ (3\text{D}),
    \ \text{with }\rho\text{ given by Eqs.~}\eqref{eq:density_2D}\text{ or }\eqref{eq:density_3D};\\
    \text{Minimize:} \quad
    & C(\rho) = \frac{1}{2}\mathbf{F}^{T}\mathbf{U}
      + C_{\mathrm{top}}(\rho);\\[2pt]
    \text{Subject to:} \quad
    & V(\rho) \leq V_{0},\\
    & \mathbf{K}\mathbf{U} = \mathbf{F},\\
    & \rho_{\min} \leq \rho \leq \rho_{\max},
  \end{aligned}
\end{equation}
where $V(\rho)$ and $V_{0}$ denote the volume fraction of the optimized structure
and its prescribed upper bound, respectively, and $\rho_{\min}$ and
$\rho_{\max}$ are the lower and upper bounds of the density field.\par

Following Eqs.~\eqref{eq:zero_optimization}, \eqref{eq:one_optimization}
and~\eqref{eq:two_optimization}, the overall topological objective is formally defined as
the weighted sum
\begin{equation}
  \label{eq:topo_opt}
  C_{\mathrm{top}}(\rho)
  = \mu_{0}\,C_{\mathrm{top}}^{0}
  + \mu_{1}\,C_{\mathrm{top}}^{1}
  + \mu_{2}\,C_{\mathrm{top}}^{2},
\end{equation}
where $\mu_{0}$, $\mu_{1}$ and $\mu_{2}$ are the weights assigned to the
0-, 1- and 2-dimensional topological objectives. In 2D no cavity exists, so that
$C_{\mathrm{top}}^{2}$ is inactive and one may set $\mu_{2}=0$; in 3D all three
terms are active, with $C_{\mathrm{top}}^{1}$ controlling the tunnels and
$C_{\mathrm{top}}^{2}$ the cavities.

\paragraph*{Sensitivity analysis.}
The optimization problem of Eq.~\eqref{eq:optimization_problem} is solved by a
gradient-based optimizer, which requires the sensitivities of the objective
function and of the constraints with respect to the design variables. For a unified
treatment of 2D and 3D cases, let $\bm{i}$ denote the multi-index $(i,j)$ in
2D and $(i,j,k)$ in 3D, and let $\bm{u}$ denote the corresponding parametric
coordinates $(u,v)$ or $(u,v,w)$, so that the design variables are the
control coefficients $\rho_{\bm{i}}$ appearing in Eq.~\eqref{eq:optimization_problem}.
In terms of these variables, the NURBS density field of
Eq.~\eqref{eq:density_2D} (or Eq.~\eqref{eq:density_3D} in 3D) is written
compactly as
\begin{equation}
  \rho(\bm{u}) = \sum_{\bm{i}} R_{\bm{i}}(\bm{u})\,\rho_{\bm{i}},
  \qquad
  R_{\bm{i}}(\bm{u}) =
  \frac{\omega_{\bm{i}}N_{\bm{i}}(\bm{u})}
       {\sum_{\bm{i}'}\omega_{\bm{i'}}N_{\bm{i'}}(\bm{u})},
  \label{eq:rho_basis}
\end{equation}
where $N_{\bm{i}}(\bm{u})=N_{i,p}(u)N_{j,q}(v)$ in 2D and
$N_{i,p}(u)N_{j,q}(v)N_{k,r}(w)$ in 3D, and $R_{\bm{i}}$ is the corresponding
rational basis function. Consequently,
\begin{equation}
  \frac{\partial \rho(\bm{u})}{\partial \rho_{\bm{i}}} = R_{\bm{i}}(\bm{u}).
  \label{eq:drho_dcontrol}
\end{equation}

Using the topological inverse mapping introduced above, which assigns to each
persistence pair $(b_{k,m},d_{k,m})$ its birth and death locations
$\bigl(\pi_b(b_{k,m}),\pi_d(d_{k,m})\bigr)$ in the parametric domain, and
recalling that the filtration is induced by $f=-\rho$, one has
\begin{equation}
  b_{k,m} = -\rho\bigl(\pi_b(b_{k,m})\bigr),
  \qquad
  d_{k,m} = -\rho\bigl(\pi_d(d_{k,m})\bigr),
\end{equation}
and thus, by Eq.~\eqref{eq:drho_dcontrol},
\begin{equation}
  \frac{\partial b_{k,m}}{\partial \rho_{\bm{i}}}
  = -R_{\bm{i}}\bigl(\pi_b(b_{k,m})\bigr),
  \qquad
  \frac{\partial d_{k,m}}{\partial \rho_{\bm{i}}}
  = -R_{\bm{i}}\bigl(\pi_d(d_{k,m})\bigr).
  \label{eq:bd_sens}
\end{equation}

The sensitivities of the topological objectives then follow from the chain rule.
For the 0-dimensional objective, in which every finite pair is removed by
maximizing its birth time irrespective of the region it belongs to, one obtains
\begin{equation}
  \frac{\partial C_{\mathrm{top}}^{0}}{\partial \rho_{\bm{i}}}
  = \sum_{(b_{k,0},d_{k,0})\in P_{0}^{I,II}}
    R_{\bm{i}}\bigl(\pi_b(b_{k,0})\bigr).
  \label{eq:sens_0d}
\end{equation}
For the 1- and 2-dimensional objectives ($m=1,2$), assembling the
contributions of the removal, preservation and regularization terms of
Eqs.~\eqref{eq:one_optimization} and~\eqref{eq:two_optimization} gives
\begin{equation}
  \begin{aligned}
    \frac{\partial C_{\mathrm{top}}^{m}}{\partial \rho_{\bm{i}}}
    ={}& \sum_{P_{m}^{I,\mathrm{rem}}} R_{\bm{i}}\bigl(\pi_b(b_{k,m})\bigr)
        - \sum_{P_{m}^{II,\mathrm{rem}}} R_{\bm{i}}\bigl(\pi_d(d_{k,m})\bigr)\\
      &+ \sum_{P_{m}^{\mathrm{pres}}}
         \Bigl[ R_{\bm{i}}\bigl(\pi_d(d_{k,m})\bigr)
              - R_{\bm{i}}\bigl(\pi_b(b_{k,m})\bigr) \Bigr]\\
      &+ \sum_{P_{m}^{III}} R_{\bm{i}}\bigl(\pi_b(b_{k,m})\bigr)
        - \sum_{P_{m}^{IV}}  R_{\bm{i}}\bigl(\pi_d(d_{k,m})\bigr),
  \end{aligned}
  \label{eq:sens_md}
\end{equation}
where the summation sets are those introduced in the preceding subsections. In
2D no cavity exists, so that the 2-dimensional objective is inactive
($\mu_{2}=0$) and only Eq.~\eqref{eq:sens_0d} together with the case $m=1$ of
Eq.~\eqref{eq:sens_md} is used.

The sensitivity of the complete objective
$C(\rho)=\frac12\mathbf{F}^{T}\mathbf{U}+\sum_{m=0}^{2}\mu_m C_{\mathrm{top}}^{m}$
follows by superposition. For the compliance term, the modified SIMP
interpolation \citep{sigmundMorphologybasedBlackWhite2007a} is adopted:
\begin{equation}
  E(\rho) = E_{\min} + \rho^{p}\,(E_0 - E_{\min}),
  \qquad \rho \in [\rho_{\min},\rho_{\max}],
  \label{eq:simp}
\end{equation}
where $E_0$ is the Young's modulus of the solid material, $E_{\min}$ is a small
positive value with $0 < E_{\min} \ll E_0$ introduced to prevent the global
stiffness matrix from becoming singular, and $p$ is the penalization factor.
Under the modified SIMP law of Eq.~\eqref{eq:simp}, the global stiffness matrix
depends on the design variables through the density field. Since the compliance
is self-adjoint, its sensitivity takes the compact form
\begin{equation}
  \begin{aligned}
  \frac{\partial}{\partial\rho_{\bm{i}}}
  \left(\frac{1}{2}\mathbf{F}^{T}\mathbf{U}\right)
  &= -\frac{1}{2}\,\mathbf{U}^{T}
    \frac{\partial\mathbf{K}}{\partial\rho_{\bm{i}}}\mathbf{U}\\
  &= -\frac{1}{2}\,\mathbf{U}^{T}
      \biggl[\int_{\Omega}
    p\,[\rho(\bm{u})]^{p-1}\,(E_0-E_{\min})R_{\bm{i}}(\bm{u})\,
    \mathbf{B}^{T}(\bm{u})\,\mathbf{D}^{0}\,
    \mathbf{B}(\bm{u})\,d\Omega\biggr]\mathbf{U},  
  \end{aligned}
  \label{eq:sens_compliance}
\end{equation}
where $\mathbf{B}^{T}(\bm{u})$ denotes the strain--displacement matrix and $\mathbf{D}^{0}$ is the
constitutive matrix for unit Young's modulus.\par
Hence,
\begin{equation}
  \frac{\partial C(\rho)}{\partial\rho_{\bm{i}}}
  = -\frac12\,\mathbf{U}^{T}
    \frac{\partial\mathbf{K}}{\partial\rho_{\bm{i}}}\,\mathbf{U}
    + \sum_{m=0}^{2}\mu_m\,
      \frac{\partial C_{\mathrm{top}}^{m}}{\partial\rho_{\bm{i}}}.
  \label{eq:sens_total}
\end{equation}

Finally, taking the volume fraction as
$V(\rho)=\frac{1}{|\Omega|}\int_{\Omega}\rho\,d\Omega$, Eqs.~\eqref{eq:rho_basis}
and~\eqref{eq:drho_dcontrol} give
\begin{equation}
  \frac{\partial V(\rho)}{\partial \rho_{\bm{i}}}
  = \frac{1}{|\Omega|}\int_{\Omega}
    \frac{\partial \rho(\bm{u})}{\partial \rho_{\bm{i}}}\,d\Omega
  = \frac{1}{|\Omega|}\int_{\Omega} R_{\bm{i}}(\bm{u})\,d\Omega,
  \label{eq:sens_vol}
\end{equation}
which is independent of the design variables and can be evaluated once in a
preprocessing step.\par

\subsection{Staged optimization}
\label{subsec:overall_flow}

\begin{figure}[!htb]
  \centering
  \begingroup
  \everymath{\textstyle}
  \resizebox{\linewidth}{!}{
  \begin{tikzpicture}[
      font=\footnotesize,
      node distance=0.6cm and 0.9cm,
      lbl/.style={font=\scriptsize, inner sep=1.2pt, text=black!75},
      box/.style={draw, rounded corners=2pt, align=center, inner xsep=5pt,
                  minimum height=0.62cm, text height=1.8ex, text depth=0.6ex,
                  line width=0.55pt},
      io/.style={box, fill=green!7, draw=green!45!black},
      proc/.style={box, fill=blue!5, draw=blue!45!black},
      dec/.style={diamond, sharp corners, align=center, aspect=4,
                  inner xsep=1.5pt, inner ysep=1pt,
                  fill=orange!9, draw=orange!70!black,
                  line width=0.55pt},
      act/.style={box, fill=violet!6, draw=violet!55!black},
      keep/.style={box, fill=cyan!7, draw=cyan!50!black},
      arr/.style={-{Stealth[length=1.8mm,width=1.4mm]},
                  draw=black!70, line width=0.55pt},
      rail/.style={draw=black!70, line width=0.55pt},
    ]
    \node[io] (init) {Initialize design domain and parameters};
    \node[proc] (fea) [right=of init] {Solve $\mathbf{K}\mathbf{U}=\mathbf{F}$};
    \node[proc] (sens) [right=of fea] {Evaluate compliance sensitivity};
    \node[proc] (filt) [below=of sens] {Construct sublevel filtration of $f=-\rho$};
    \node[proc] (num) at (init |- filt) {Extract $N_{0}$, $N_{1}$, and $N_{2}$};
    \node[proc] (pd) at ($(num.east)!0.5!(filt.west)$) {Compute PDs};

    \node[dec] (chk2) [below=of num] {$N_{2}>0$?};
    \node[dec] (chk1) [below=of chk2] {$N_{1}>\bar{N}_{1}$?};
    \node[keep] (frz) [below=of chk1]
      {Activate preservation and regularization terms of $C_{\mathrm{top}}^{1}$ for retained features};
    \node[dec] (chk0) [below=of frz] {$N_{0}>1$?};
    \node[proc] (asm) [below=of chk0] {Assemble sensitivities of active terms};
    \node[proc] (mma) [below=of asm] {Update density control coefficients $\rho_{\bm{i}}$ by MMA};
    \node[dec] (conv) [below=of mma]
      {Converged or maximum number of iterations reached?};

    \node[act, anchor=east] (top1) at (filt.east |- chk1)
      {Activate removal and regularization terms of $C_{\mathrm{top}}^{1}$};
    \path let \p1=($(top1.east)-(top1.west)$) in
      node[act, minimum width=\x1, anchor=east] (top2) at (filt.east |- chk2)
        {Activate $C_{\mathrm{top}}^{2}$}
      node[act, minimum width=\x1, anchor=east] (top0) at (filt.east |- chk0)
        {Activate $C_{\mathrm{top}}^{0}$};
    \node[io, anchor=west] (out) at ($(conv.east)+(0.9,0)$)
      {Output optimized density field};

    \draw[arr] (init) -- (fea);
    \draw[arr] (fea)  -- (sens);
    \draw[arr] (sens) -- (filt);
    \draw[arr] (filt) -- (pd);
    \draw[arr] (pd)   -- (num);
    \draw[arr] (num)  -- (chk2);
    \draw[arr] (chk2) -- node[lbl, right] {no} (chk1);
    \draw[arr] (chk1) -- node[lbl, right] {no} (frz);
    \draw[arr] (frz)  -- (chk0);
    \draw[arr] (chk0) -- node[lbl, right] {no} (asm);
    \draw[arr] (asm)  -- (mma);
    \draw[arr] (mma)  -- (conv);
    \draw[arr] (chk2) -- node[lbl, above] {yes} (top2);
    \draw[arr] (chk1) -- node[lbl, above] {yes} (top1);
    \draw[arr] (chk0) -- node[lbl, above] {yes} (top0);
    \draw[arr] (conv) -- node[lbl, above] {yes} (out);

    \coordinate (rx) at ($(top1.east)+(0.35,0)$);
    \draw[rail] (top2.east) -- (top2.east -| rx);
    \draw[rail] (top1.east) -- (top1.east -| rx);
    \draw[rail] (top0.east) -- (top0.east -| rx);
    \draw[arr] (top2.east -| rx) -- (rx |- asm.east) -- (asm.east);

    \coordinate (lx) at ($(frz.west)+(-0.4,0)$);
    \coordinate (ltop) at ($(lx |- init.north)+(0,0.35)$);
    \draw[arr] (conv.west) -| node[lbl, above, pos=0.25] {no} (ltop) -| (fea.north);
  \end{tikzpicture}}
  \endgroup
  \caption{Flowchart of the proposed persistent homology-based isogeometric
  topology optimization framework.}
  \label{fig:flowchart}
\end{figure}
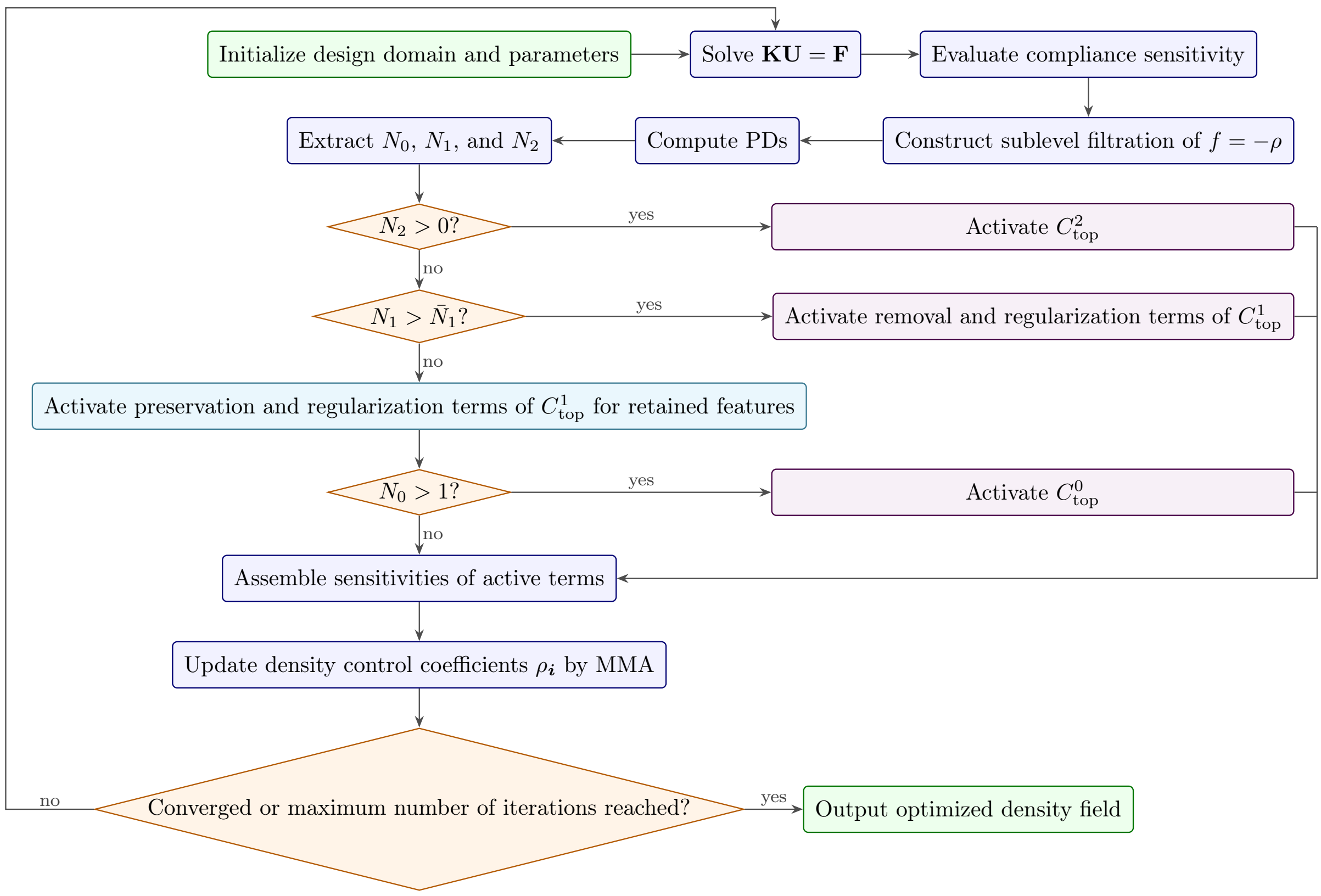

The preceding sections have introduced the three ingredients of the proposed
framework: the 0-dimensional objective $C_{\mathrm{top}}^{0}$, which removes
redundant connected components by cutting; the 1-dimensional objective
$C_{\mathrm{top}}^{1}$, which controls the number of holes in 2D, or of tunnels
in 3D, and preserves the ones to be retained; and the 2-dimensional objective
$C_{\mathrm{top}}^{2}$, which eliminates the cavities and is active in 3D only.
This subsection summarizes how these ingredients are assembled into the overall
optimization procedure, which is illustrated in Fig.~\ref{fig:flowchart}.\par

\paragraph*{Topological measures from a single filtration.}
At each iteration, the discretized density field is used to induce a single
sublevel filtration, namely the filtration of the function $f=-\rho$, whose
sublevel set at the threshold $t=-\bar{\rho}$ coincides with the solid phase
$\rho\ge \bar{\rho}$. From this filtration, the PDs in
dimensions $0$, $1$, and $2$ (in 3D) are computed simultaneously. The
0-dimensional diagram yields the number of connected components $N_{0}$ of the
current structure, i.e., the essential pair together with the finite 0-dimensional
pairs located in regions~I and~II. The 1-dimensional diagram yields the number of
holes in 2D, or of tunnels in 3D, denoted by $N_{1}$, that is, the number of
1-dimensional pairs lying in regions~I and~II. In 3D, the 2-dimensional diagram
additionally yields the number of cavities $N_{2}$, i.e., the number of
2-dimensional pairs in regions~I and~II.\par

\paragraph*{Motivation for a staged activation.}
The topological constraints of the different dimensions are not enforced
simultaneously, for two reasons.
\begin{itemize}
    \item \emph{Conflict between the 0- and 1-dimensional controls.} When a hole
    in 2D, or a tunnel in 3D, is eliminated by filling, the added material
    generally forms a new connected component before it merges with the existing
    structure. The 0-dimensional control would then remove this newly created
    component and thereby undo the filling, so that the two objectives compete
    with each other and the hole can never be closed.
    \item \emph{Cavities obscure the tunnel count in 3D.} The presence of
    cavities affects the evaluation of the number of tunnels, since $\beta_{1}$
    no longer equals the number of tunnels when cavities are present.
\end{itemize}
Consequently, the 2-dimensional control is applied first, the 1-dimensional one
next, and the 0-dimensional one last.\par

\paragraph*{Staged assembly of the objective and its sensitivity.}
The augmented objective is assembled in three stages, and only the terms of the
current stage are activated:
\begin{enumerate}
    \item \textbf{Stage 1: cavities (3D only).} If $N_{2}>0$, $C_{\mathrm{top}}^{2}$
    is activated to remove all cavities, which are sorted in ascending order of
    their BFS-approximated volume.
    In addition, the two regularization terms acting on the persistence pairs in
    regions~III and~IV are always included so as
    to suppress the oscillatory reappearance of cavities. Once $N_{2}=0$, the
    structure is treated as in 2D and this stage terminates.
    \item \textbf{Stage 2: holes or tunnels.} If $N_{1}>\bar{N}_{1}$,
    $C_{\mathrm{top}}^{1}$ is activated to remove the $N_{1}-\bar{N}_{1}$ least
    significant features, ordered by the BFS-approximated area in 2D and by the
    persistence $d_{i,1}-b_{i,1}$ in 3D. Once $N_{1}\leq\bar{N}_{1}$ holds, the
    preservation term of $C_{\mathrm{top}}^{1}$ is instead used to freeze the
    density of the existing holes or tunnels, so that they are not removed while
    the remaining material is further optimized for compliance.
     In both subcases, i.e.\ regardless of
    whether the removal or the preservation term is active, the two
    regularization terms for regions~III and~IV are always included.
    \item \textbf{Stage 3: connected components.} Only after the bound on the
    number of holes or tunnels is met, and if connectivity is to be enforced
    while $N_{0}>1$, is the 0-dimensional objective $C_{\mathrm{top}}^{0}$
    activated to remove the extra connected components by cutting. 
\end{enumerate}\par

In every stage, the sensitivities of the active topological terms, derived in
Section~\ref{subsec:problem_statement_sensitivity_analysis}, are added to that of the compliance to form the objective gradient.
The design variables are then updated by the MMA optimizer, and the
iteration proceeds until the maximum number of iterations is reached or the
convergence criterion is met.\par

\section{Numerical examples}
\label{sec:numerical_examples}

In this section, several benchmark examples are presented to validate the
effectiveness of the proposed method. Unless otherwise stated, all problems are
solved by the MMA optimizer, the density field is represented by NURBS, and a
single sublevel filtration of $f=-\rho$ is induced at each iteration to extract
the topological features of the structure, while the topological constraints are
activated in the staged manner described in Section~\ref{subsec:overall_flow}.
The material is modeled by the modified SIMP law of Eq.~\eqref{eq:simp}, and all the
parameter settings used in the following numerical examples are listed in 
Table~\ref{tab:parameters}.\par

Two numerical techniques are applied to the density field. First, a smoothing mechanism
based on the Shepard function is applied to the density control coefficients so as to
ensure sufficient smoothness of the density field~\citep{gaoIgaTopImplementationTopology2021}. Second, a continuation-based
Heaviside projection is employed to obtain binary designs, in which the
projection steepness parameter $\beta$ is doubled at fixed update intervals~\citep{guestAchievingMinimumLength2004}.\par

Since the structural topology evolves rapidly during the early iterations, the topology
control is switched on only after the structure has stabilized rather than from
the very first iteration. The optimization terminates either when the maximum
number of iterations is reached or when the design variables have converged.
For the computation of persistent homology, the continuous density field is
discretized on a uniform grid over the parametric domain, which is taken as
$[0,1]^{2}$ in 2D and $[0,1]^{3}$ in 3D, with a resolution fine enough to guarantee that the topological
information extracted from the discretized field remains consistent with that of
the continuous field.\par

\begin{table}[!htb]
  \centering
  \caption{Default parameter settings of the numerical examples.}
  \label{tab:parameters}
  \begin{tabular}{@{}ll@{}}
    \hline
    Parameter                                        & Value \\
    \hline
    Young's modulus of the solid phase, $E_0$        & $1.0$ \\
    Young's modulus of the void phase, $E_{\min}$    & $10^{-9}$ \\
    Penalization factor, $p$                         & $3$ \\
    Poisson's ratio, $\nu$                           & $0.3$ \\
    Bounds of the density field, $[\rho_{\min},\rho_{\max}]$ & $[0.0,\,1.0]$ \\
    Initial density control coefficients             & 1.0 ($\rho_{\max}$)\\
    Shepard function's radius parameter, $r_m$       & 2.0\\
    Heaviside function's threshold parameter, $\eta$ & 0.5\\
    Heaviside function's steepness parameter, $\beta$ & $1 \to 8$, doubled every $50$ iterations\\
    Maximum number of iterations                     & $400$ \\
    Convergence tolerance on the design variables    & $10^{-2}$ \\
    Activation iteration of the topology control     & $30$--$180$ (example-dependent) \\
    Density threshold, $\bar{\rho}$                  & $0.5$ \\
    Resolution for persistent homology (2D)          & $200 \times 200$ \\
    Resolution for persistent homology (3D)          & $50 \times 50 \times 50$ \\
    \hline
  \end{tabular}
\end{table}

\subsection{2D numerical examples}

\subsubsection{Short beam}
\label{subsec:short_beam_2D}

\begin{figure}[!htb]
  \centering
  \begin{subfigure}[t]{0.30\textwidth}
    \centering
    \includegraphics[width=\linewidth]{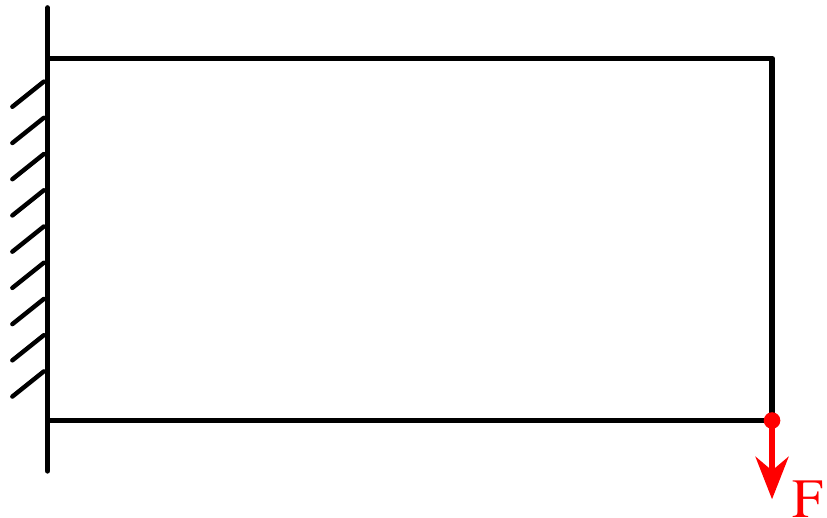}
    \caption{}
    \label{fig:short-beam-model}
  \end{subfigure}
  \hfill
  \begin{subfigure}[t]{0.36\textwidth}
    \centering
    \includegraphics[width=\linewidth]{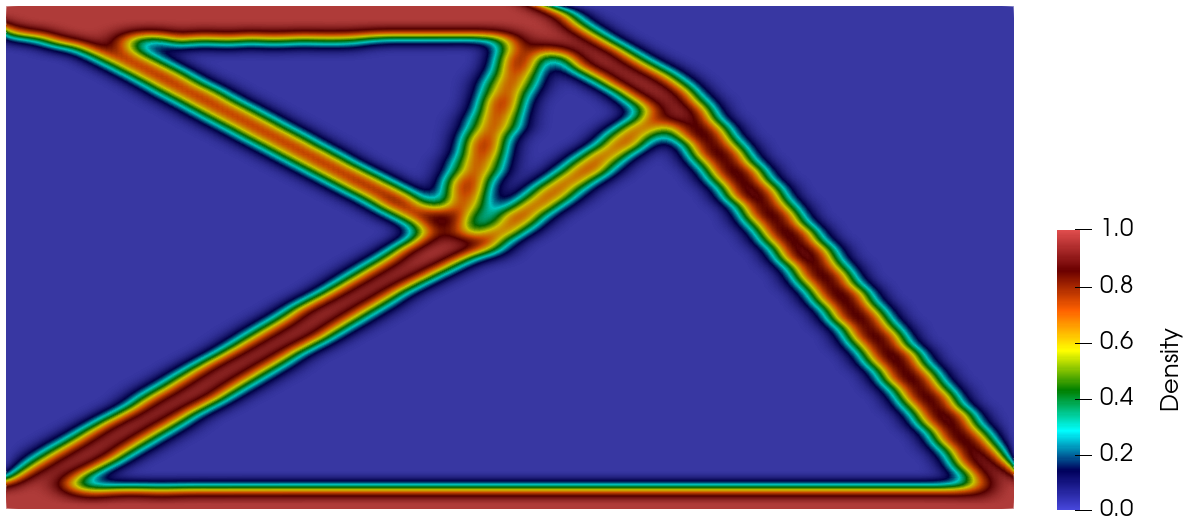}
    \caption{}
    \label{fig:short-beam-density}
  \end{subfigure}
  \hfill
  \begin{subfigure}[t]{0.30\textwidth}
    \centering
    \includegraphics[width=\linewidth]{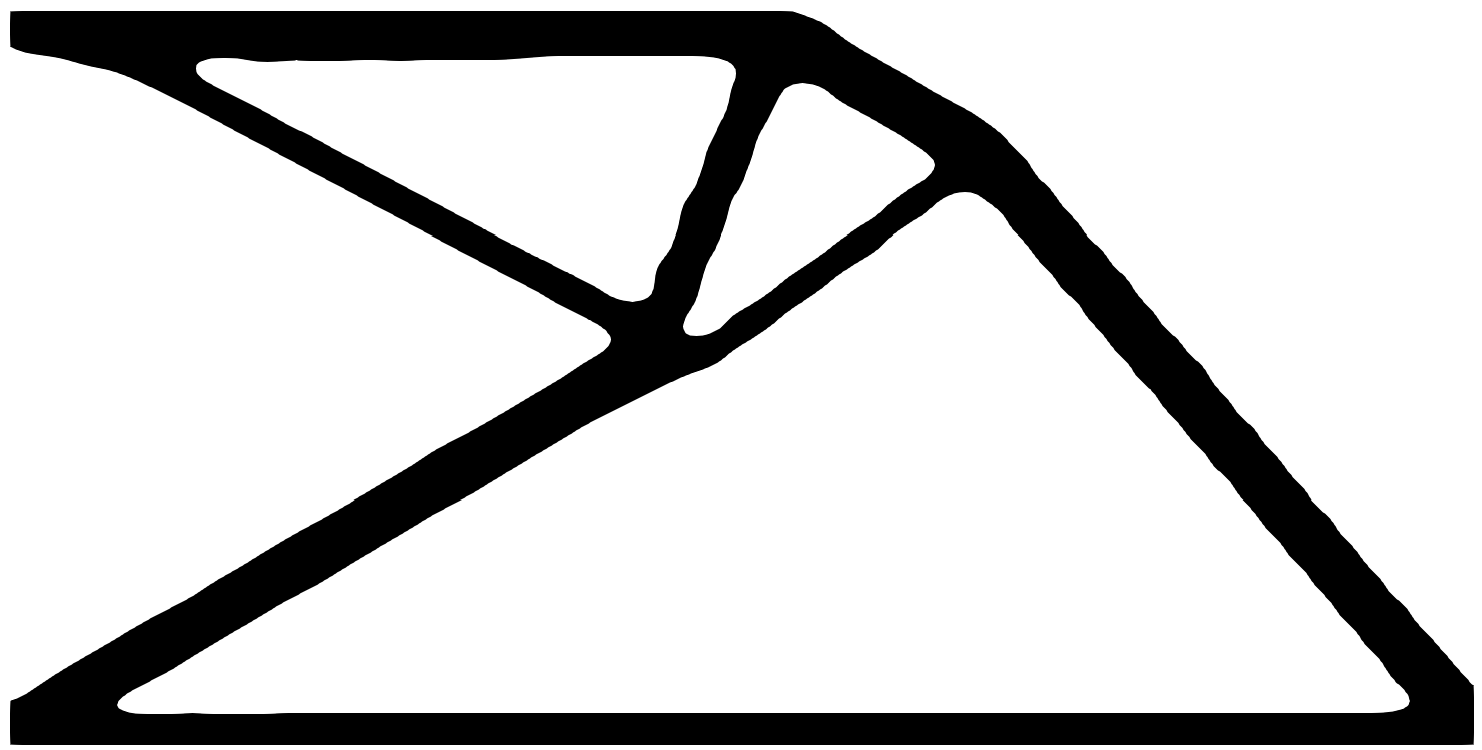}
    \caption{}
    \label{fig:short-beam-topology}
  \end{subfigure}
  \caption{(a) The short beam model.
   (b) Optimized density field without topology control.
   (c) The corresponding binary structure with $N_1=3$.}
  \label{fig:short-beam-nocontrol}
  \end{figure}

The classic short beam benchmark is considered to examine the capability of the
proposed method for explicit topology control. The design domain is the rectangle
$[0,1]\times[0,0.5]$, represented by a bicubic NURBS surface with
$61\times 61$ control points. As illustrated in Fig.~\ref{fig:short-beam-model}, the left edge is fully clamped and a unit downward
force is applied at the lower-right corner. Meanwhile, the maximum volume fraction is set to $V_{0}=0.2$.
Solving the minimum compliance problem without topological
constraint yields the baseline structure shown in
Figs.~\ref{fig:short-beam-density} and~\ref{fig:short-beam-topology}, which contains three holes.\par

\begin{figure}[!htb]
  \centering
  \begin{subfigure}[b]{0.40\textwidth}
    \centering
    \includegraphics[width=\textwidth]{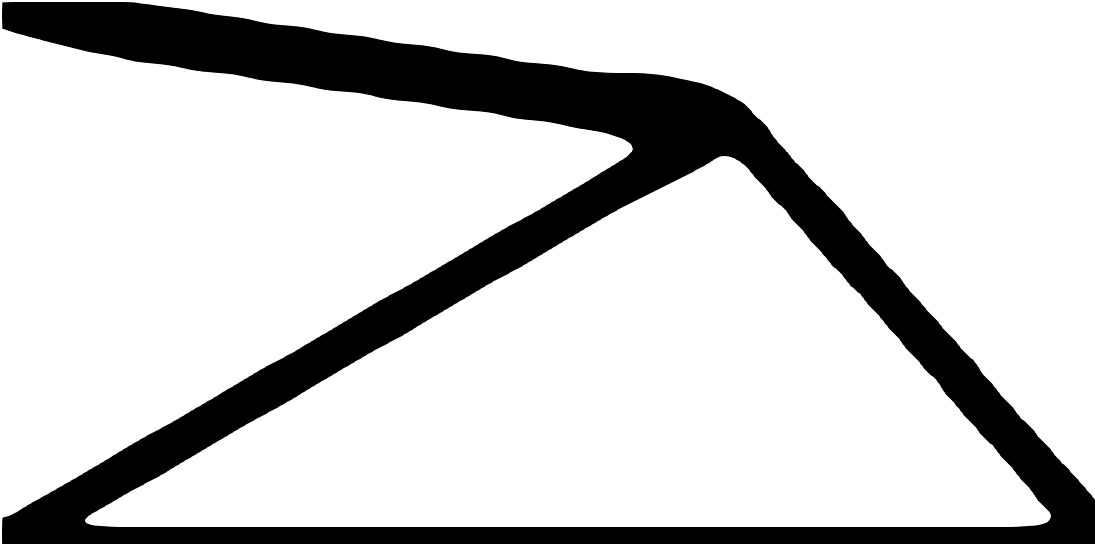}
    \caption{$\bar{N}_{1}=1$}
    \label{fig:short-beam-h1-1}
  \end{subfigure}
  \hfill
  \begin{subfigure}[b]{0.40\textwidth}
    \centering
    \includegraphics[width=\textwidth]{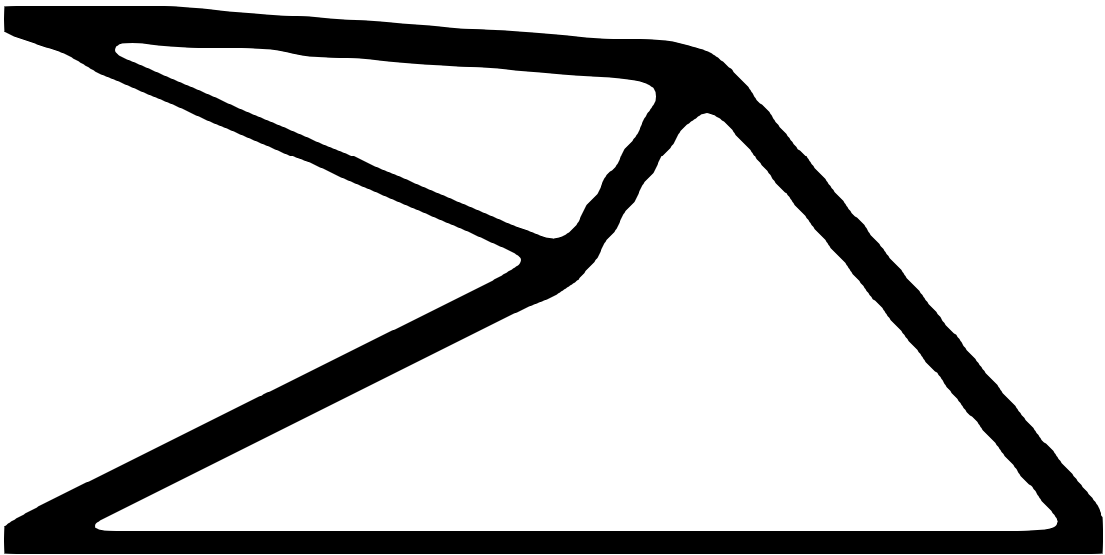}
    \caption{$\bar{N}_{1}=2$}
    \label{fig:short-beam-h1-2}
  \end{subfigure}
  \caption{Optimized short beam structures for $\bar{N}_{1}=1,2$, respectively.}
  \label{fig:short-beam-control}
\end{figure}

We impose topological constraints on the maximum number of holes $\bar{N}_{1}$ in the short beam 
and guarantee the structural connectivity at the same time.
The optimized structures for $\bar{N}_{1}=1$ and $\bar{N}_{1}=2$ are shown in
Fig.~\ref{fig:short-beam-control}, and the quantitative comparison is reported
in Table~\ref{tab:short-beam}.
As the admissible number of holes decreases, the compliance increases gradually.
This trend is expected and intuitive, since fewer admissible holes reduce the design freedom and
restrict the ability of the optimizer to distribute material in a compliance-optimal manner.
Nevertheless, the compliance variation remains moderate, which
indicates that the proposed topology control strategy achieves 
structural regularization with only a limited performance sacrifice.
In addition, the optimized structures exhibit smooth geometric transitions 
without artificial discontinuities.
The structural connectivity is consistently preserved across all cases, 
confirming the effectiveness of the zero-dimensional topology control.\par

\begin{table}[!htb]
  \centering
  \caption{Comparison of the compliance values and of the resulting numbers of holes ($N_{1}$) for the 2D
short beam under different prescribed bounds $\bar{N}_{1}$.}
  \label{tab:short-beam}
  \small
  \begin{tabular}{@{}lccc@{}}
    \toprule
    Quantity & $\bar{N}_{1}=1$ & $\bar{N}_{1}=2$ & Baseline\\
    \midrule
    Weights $(\mu_{0},\mu_{1})$ & $(100,100)$ & $(100,100)$ & $(0,0)$\\
    Compliance & $179.97$ & $178.99$ & $178.94$\\
    $N_1$ & 1 & 2 & 3\\
    \bottomrule
  \end{tabular}
\end{table}

\begin{figure}[!htb]
  \centering
  \begin{subfigure}[b]{0.33\textwidth}
    \centering
    \includegraphics[width=\linewidth]{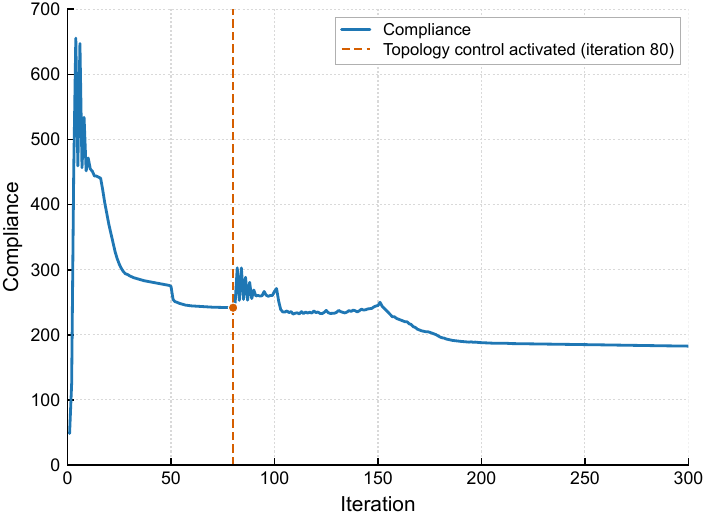}
    \caption{Compliance vs. iteration}
    \label{fig:short-beam-evolution-a}
  \end{subfigure}\hfill
  \begin{subfigure}[b]{0.30\textwidth}
    \centering
    \includegraphics[width=\linewidth]{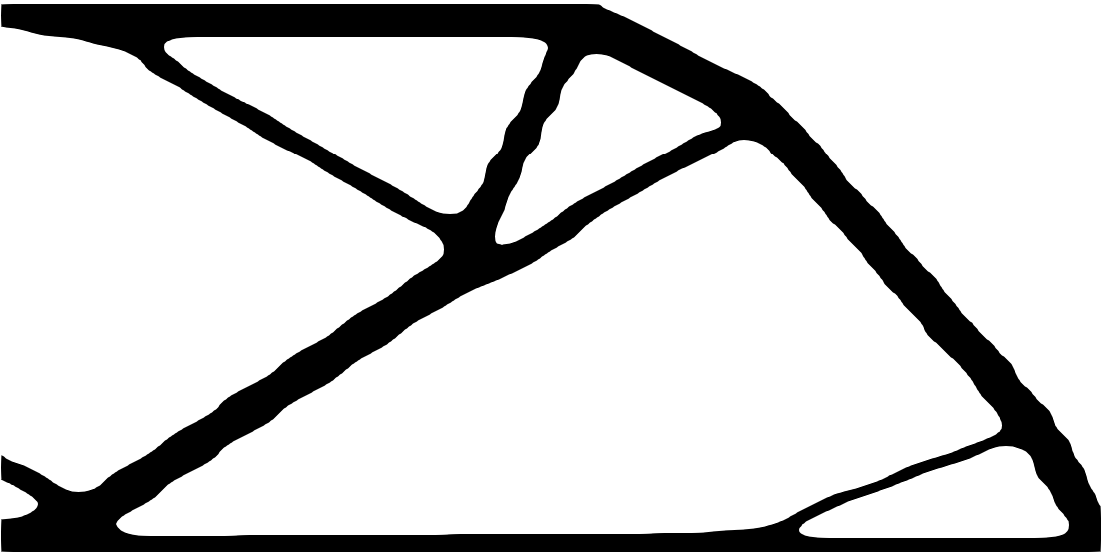}
    \caption{iter${}={}$36: $N_{1}=4$}
    \label{fig:short-beam-evolution-b}
  \end{subfigure}\hfill
  \begin{subfigure}[b]{0.30\textwidth}
    \centering
    \includegraphics[width=\linewidth]{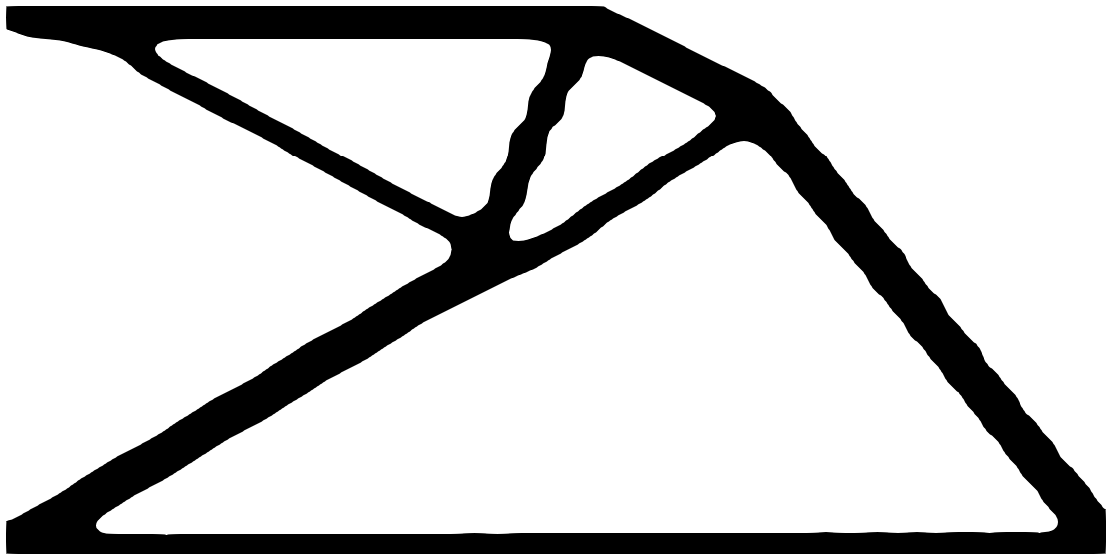}
    \caption{iter${}={}$80: $N_{1}=3$}
    \label{fig:short-beam-evolution-c}
  \end{subfigure}\\[8pt]
  \begin{subfigure}[b]{0.30\textwidth}
    \centering
    \includegraphics[width=\linewidth]{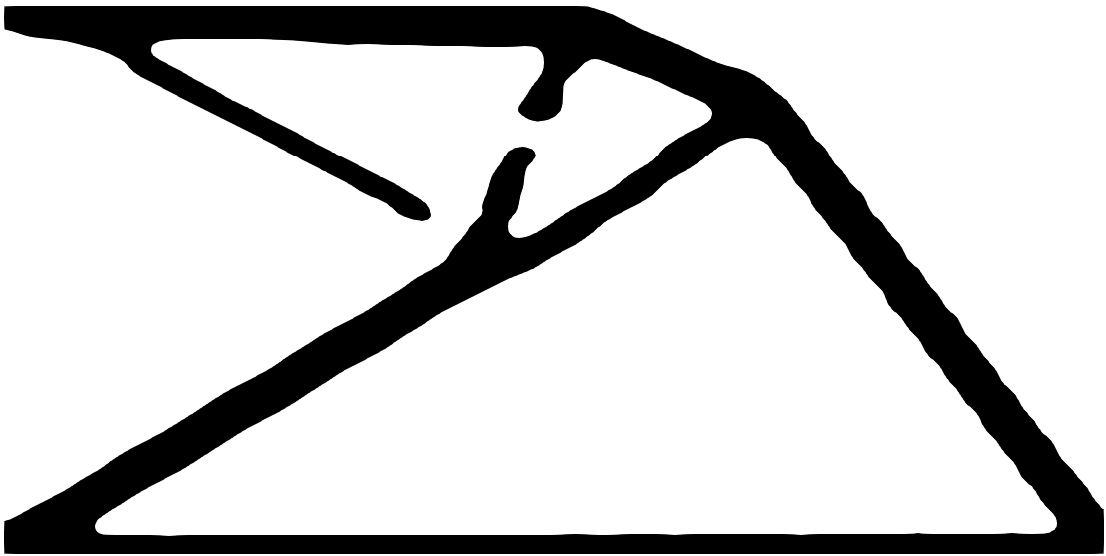}
    \caption{iter${}={}$90: $N_{1}=1$}
    \label{fig:short-beam-evolution-d}
  \end{subfigure}\hfill
  \begin{subfigure}[b]{0.30\textwidth}
    \centering
    \includegraphics[width=\linewidth]{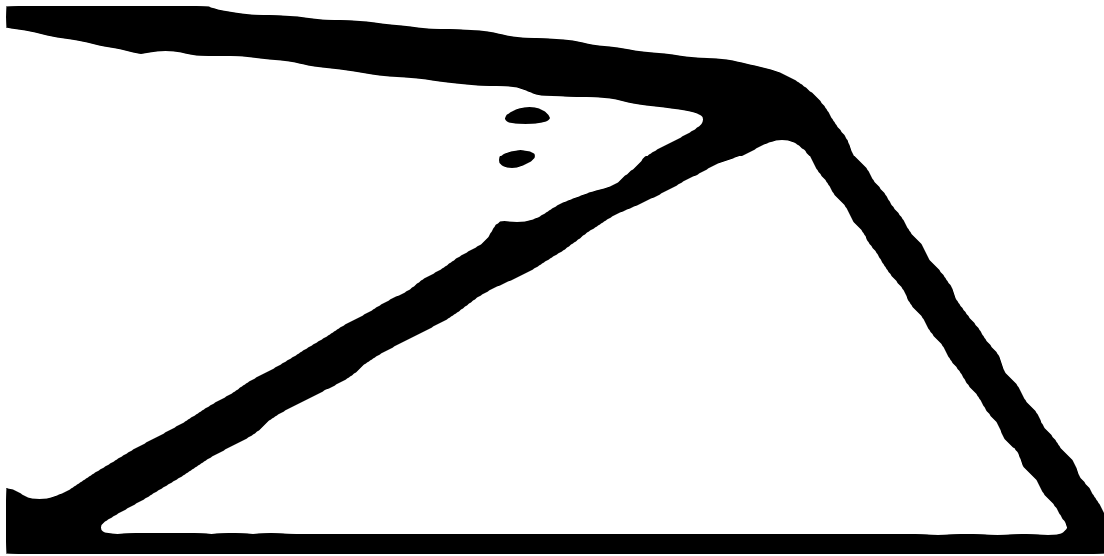}
    \caption{iter${}={}$178: $N_{1}=1$}
    \label{fig:short-beam-evolution-e}
  \end{subfigure}\hfill
  \begin{subfigure}[b]{0.30\textwidth}
    \centering
    \includegraphics[width=\linewidth]{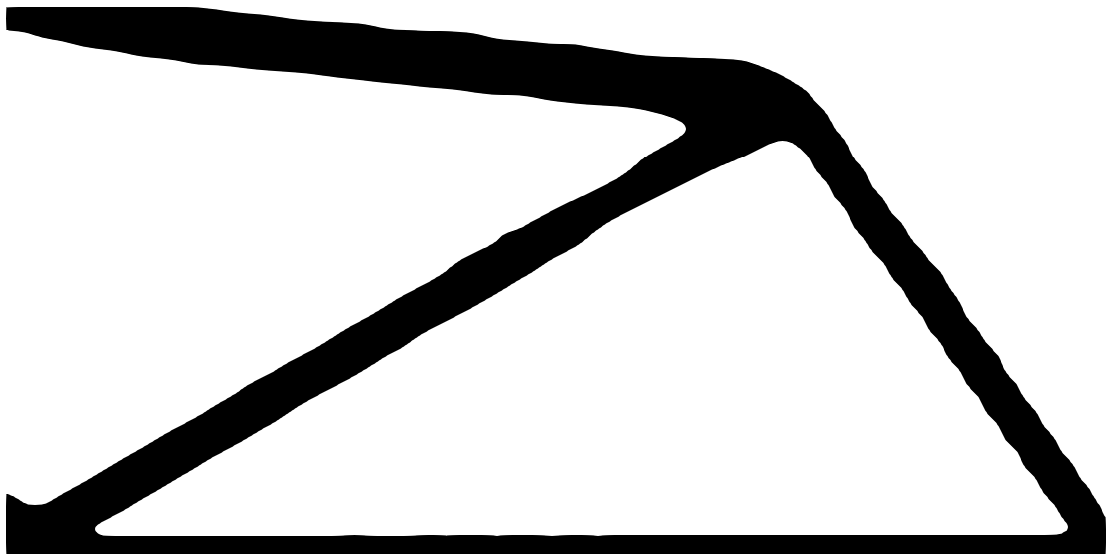}
    \caption{iter${}={}$190: $N_{1}=1$}
    \label{fig:short-beam-evolution-f}
  \end{subfigure}
  \caption{Evolution of the topology for the short beam when $\bar{N}_{1}=1$. (a) compliance history; (b)--(f) binary structures at representative
    iterations.}
  \label{fig:short-beam-evolution}
\end{figure}

The evolution of the compliance and the topology for
$\bar{N}_{1}=1$ is illustrated in Fig.~\ref{fig:short-beam-evolution}. 
The iteration history in Fig.~\ref{fig:short-beam-evolution-a} starts with a peak,
which stems from the initialization of the density field and the SIMP penalization of the gray
field at the initial Heaviside steepness $\beta=1$.
Figs.~\ref{fig:short-beam-evolution-b}--\ref{fig:short-beam-evolution-f} depict the structural evolution 
and demonstrate the effectiveness of the proposed topological objective function.
At the 80th iteration (Fig.~\ref{fig:short-beam-evolution-c}) the structure contains three holes and violates the topological constraint.
Once the topology control is activated, the number of holes is reduced to one by
the 90th iteration through cutting along the hole boundaries
(Fig.~\ref{fig:short-beam-evolution-d}), thereby satisfying the constraint and
confirming that the persistence-based objective effectively suppresses redundant
topological features.
Subsequently, at the 178th iteration, the structure contains two extra connected components, which are
eliminated in the following iterations.
Finally, by the $190$th iteration the excess material is completely removed 
	and the structure reaches a stable topology.
Throughout the late stage of the optimization the compliance remains nearly
constant, which shows that the proposed topology control mechanism does not compromise the
convergence stability. \par

\subsubsection{L-shaped beam}
\label{subsec:l_beam_2D}

\begin{figure}[!htb]
  \centering
  \begin{subfigure}[t]{0.30\textwidth}
    \centering
    \includegraphics[width=\linewidth]{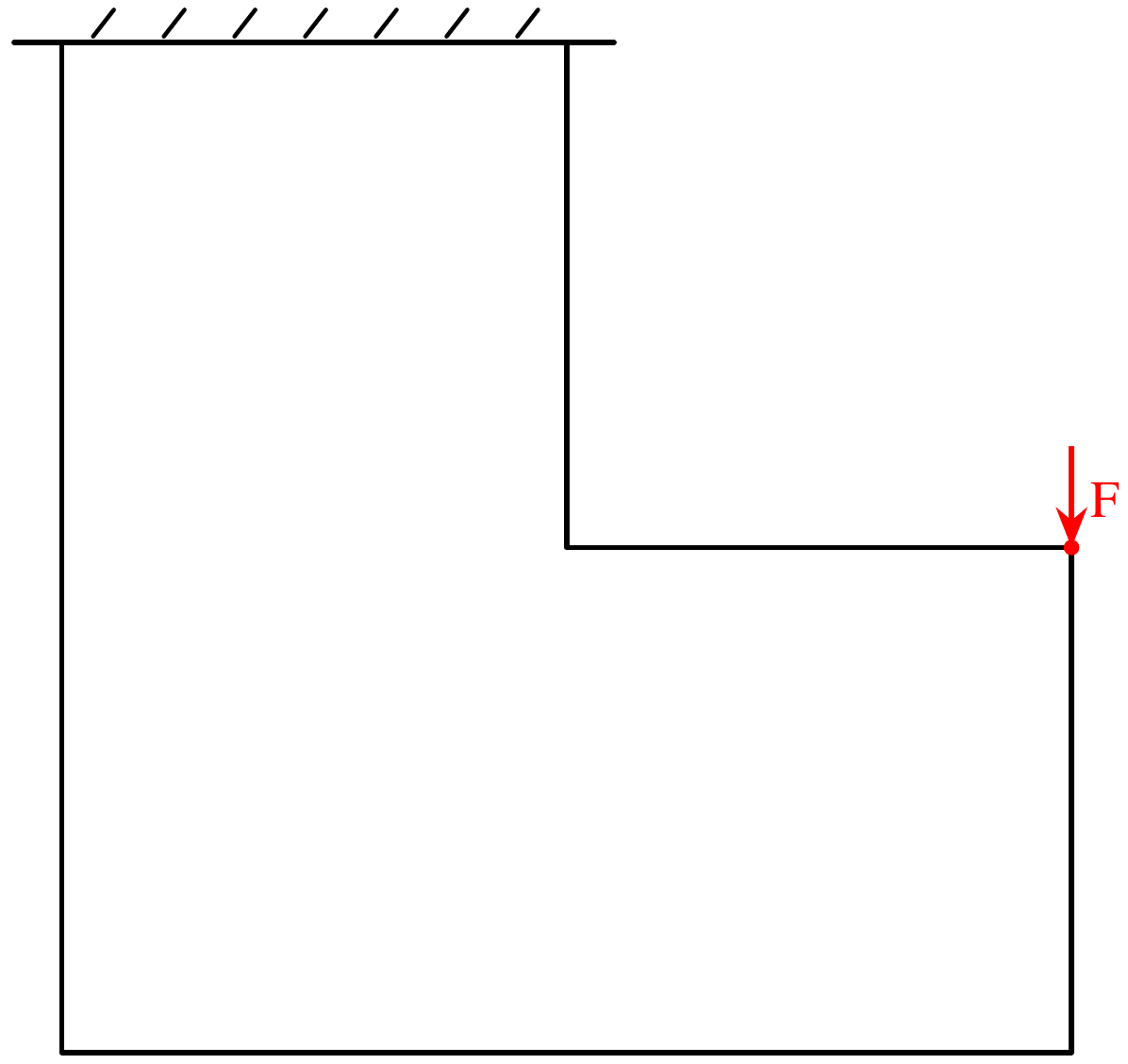}
    \caption{}
    \label{fig:lbeam-model}
  \end{subfigure}
  \hfill
  \begin{subfigure}[t]{0.36\textwidth}
    \centering
    \includegraphics[width=\linewidth]{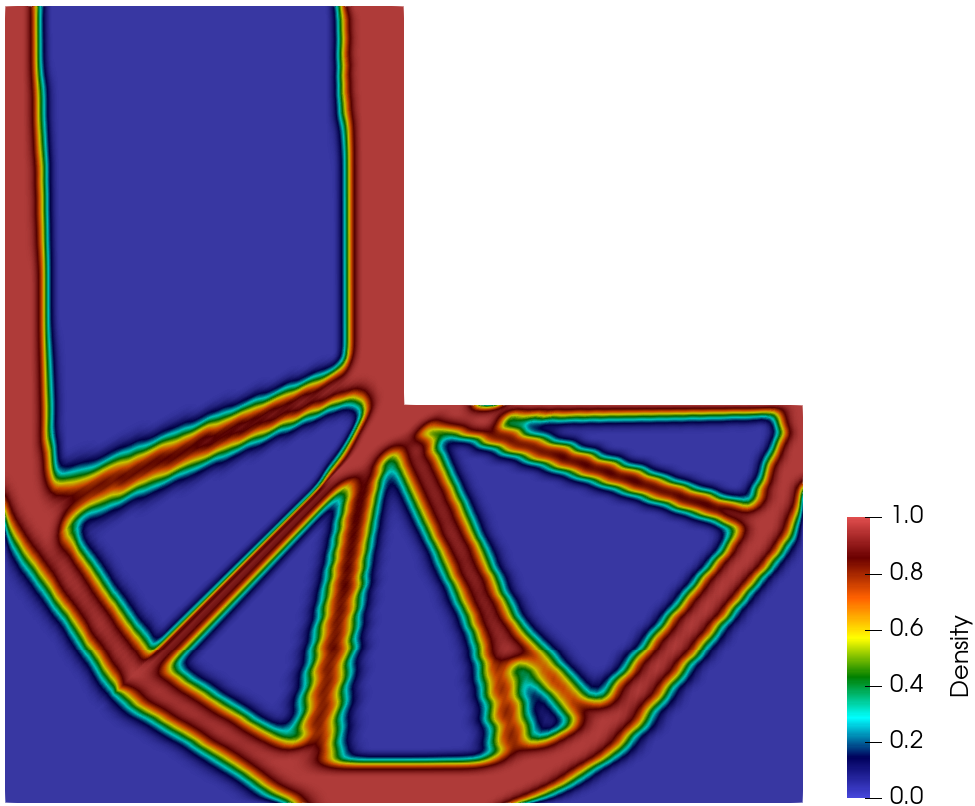}
    \caption{}
    \label{fig:lbeam-density}
  \end{subfigure}
  \hfill
  \begin{subfigure}[t]{0.30\textwidth}
    \centering
    \includegraphics[width=\linewidth]{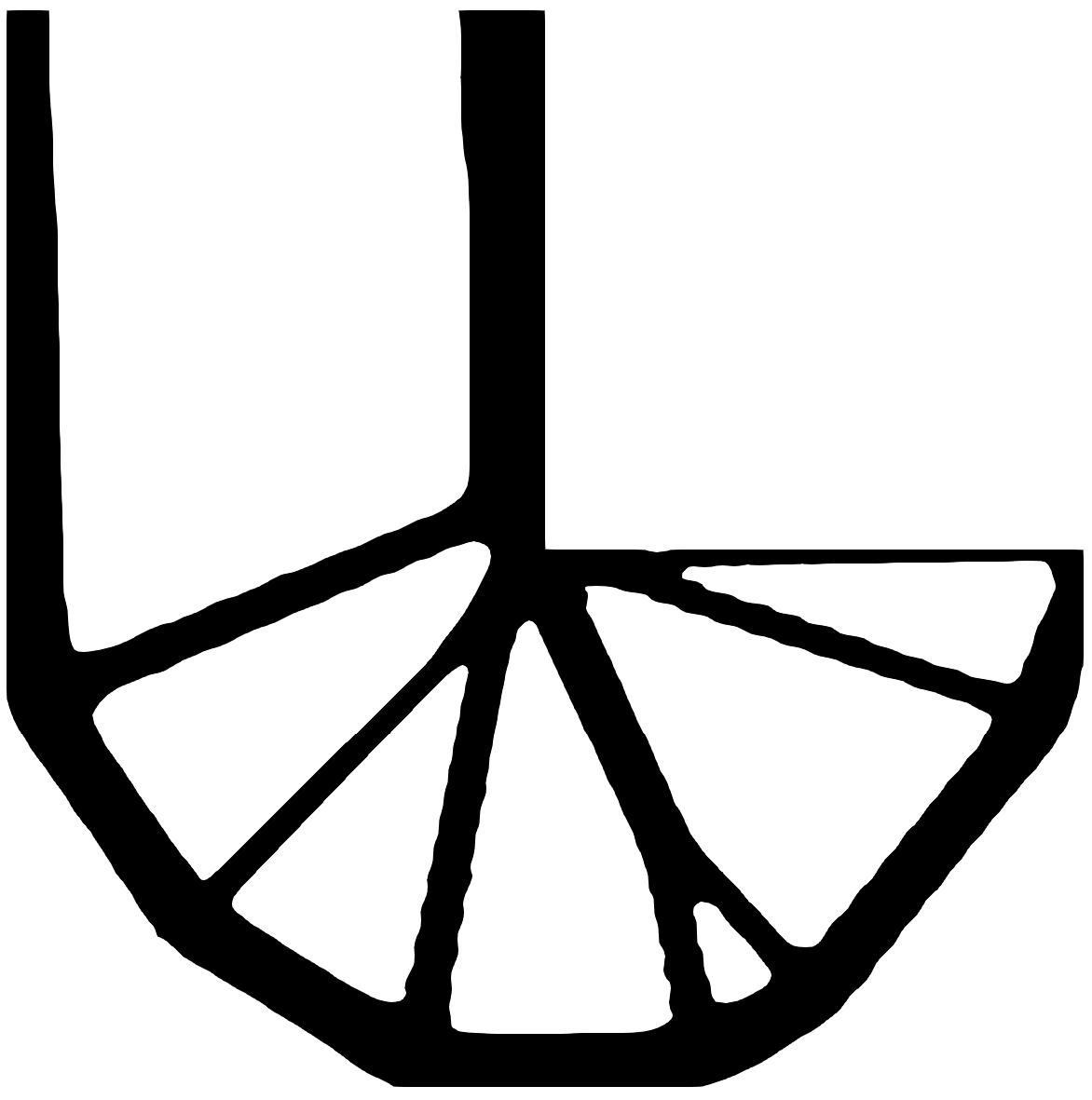}
    \caption{}
    \label{fig:lbeam-topology}
  \end{subfigure}
  \caption{(a) The L-shaped beam model.
   (b) Optimized density field without topology control.
   (c) The corresponding binary structure with $N_1=6$.}
    \label{fig:lbeam-nocontrol}
  \end{figure}

Fig.~\ref{fig:lbeam-nocontrol} presents the L-shaped beam benchmark.
The design domain is the square $[0,10]^{2}$ with the
corner region $[5,10]\times[5,10]$ removed, represented by a biquadratic NURBS
surface with $52\times 103$ control points. The top
edge is fully clamped and a unit downward point load is applied at the free end
$(10,5)$, as illustrated in Fig.~\ref{fig:lbeam-model}. The minimum compliance problem without topological constraint is solved with the maximum volume fraction set to $V_{0}=0.3$.
As shown in Figs.~\ref{fig:lbeam-density} and~\ref{fig:lbeam-topology}, there exist six holes in the optimized structure.\par

\begin{figure}[!htb]
  \centering
  \begin{subfigure}[b]{0.30\textwidth}
    \centering
    \includegraphics[width=\textwidth]{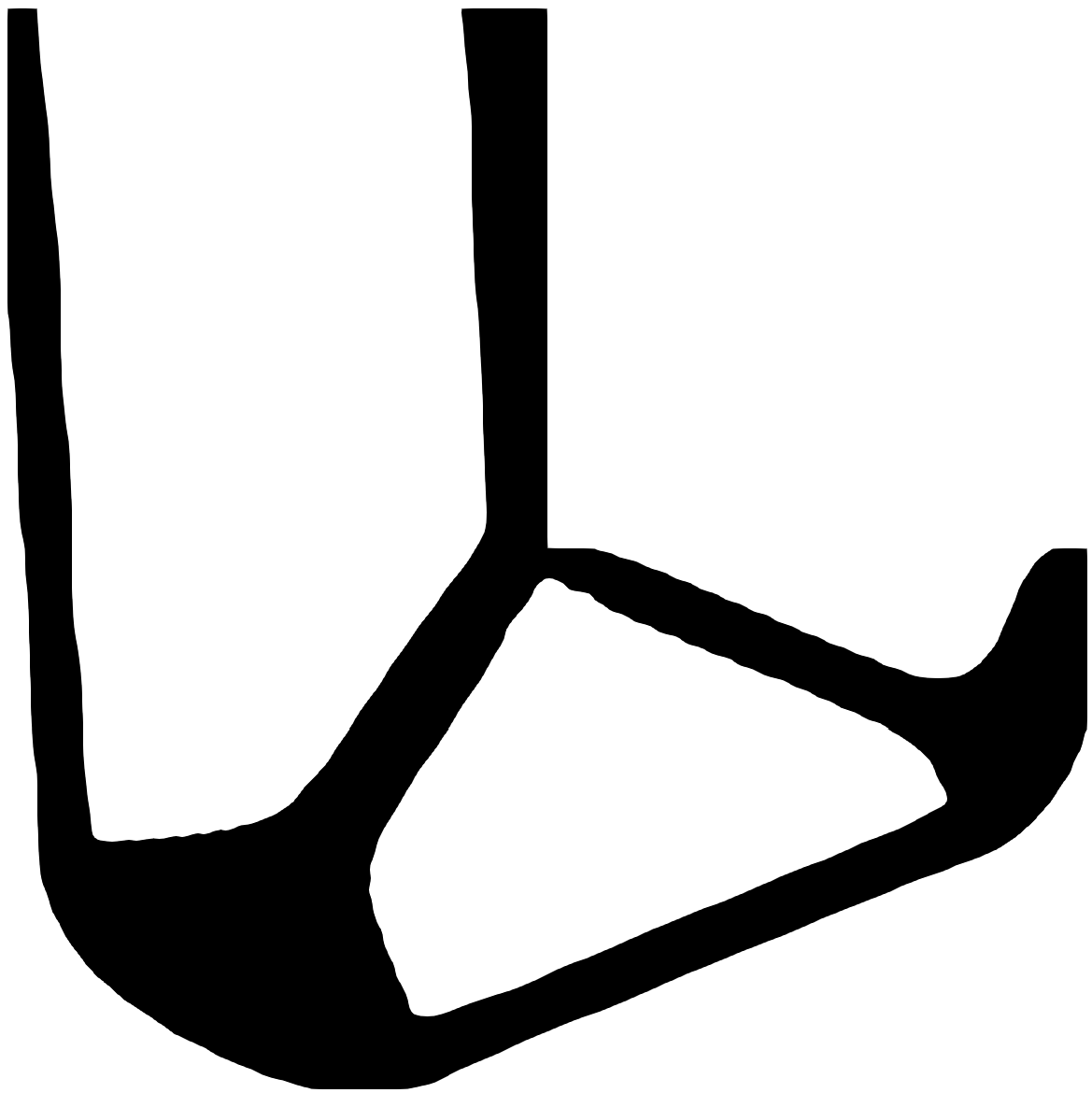}
    \caption{$\bar{N}_{1}=1$}
    \label{fig:lbeam-h1-1}
  \end{subfigure}
  \hfill
  \begin{subfigure}[b]{0.30\textwidth}
    \centering
    \includegraphics[width=\textwidth]{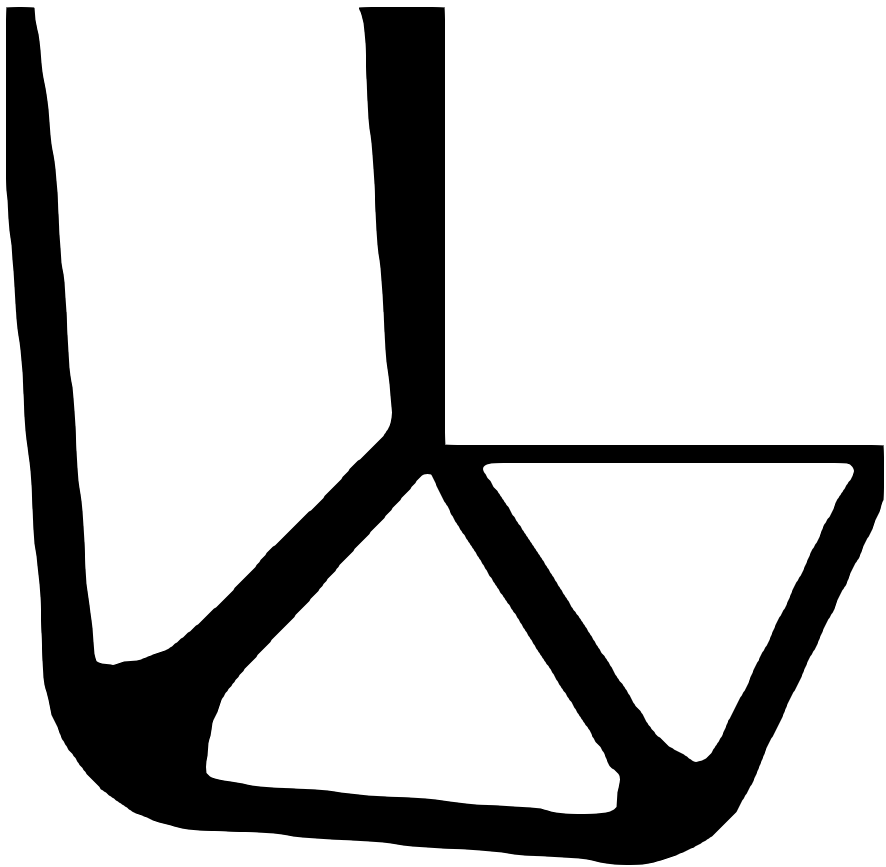}
    \caption{$\bar{N}_{1}=2$}
    \label{fig:lbeam-h1-2}
  \end{subfigure}
  \hfill
  \begin{subfigure}[b]{0.30\textwidth}
    \centering
    \includegraphics[width=\textwidth]{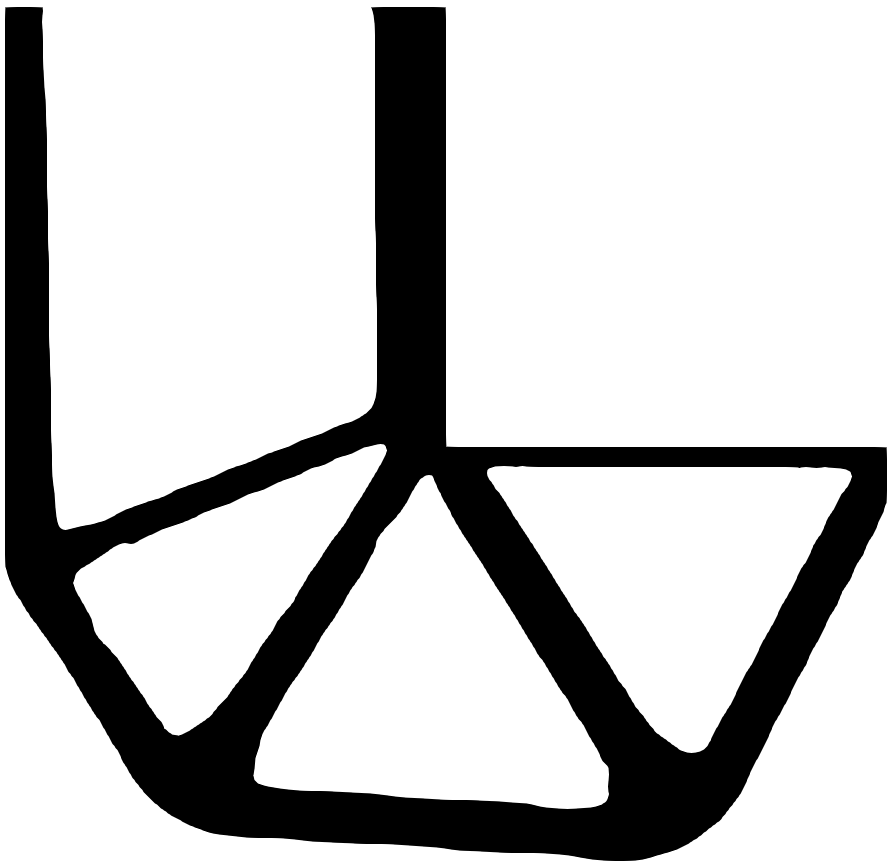}
    \caption{$\bar{N}_{1}=3$}
    \label{fig:lbeam-h1-3}
  \end{subfigure}\\[8pt]
  \begin{subfigure}[b]{0.30\textwidth}
    \centering
    \includegraphics[width=\textwidth]{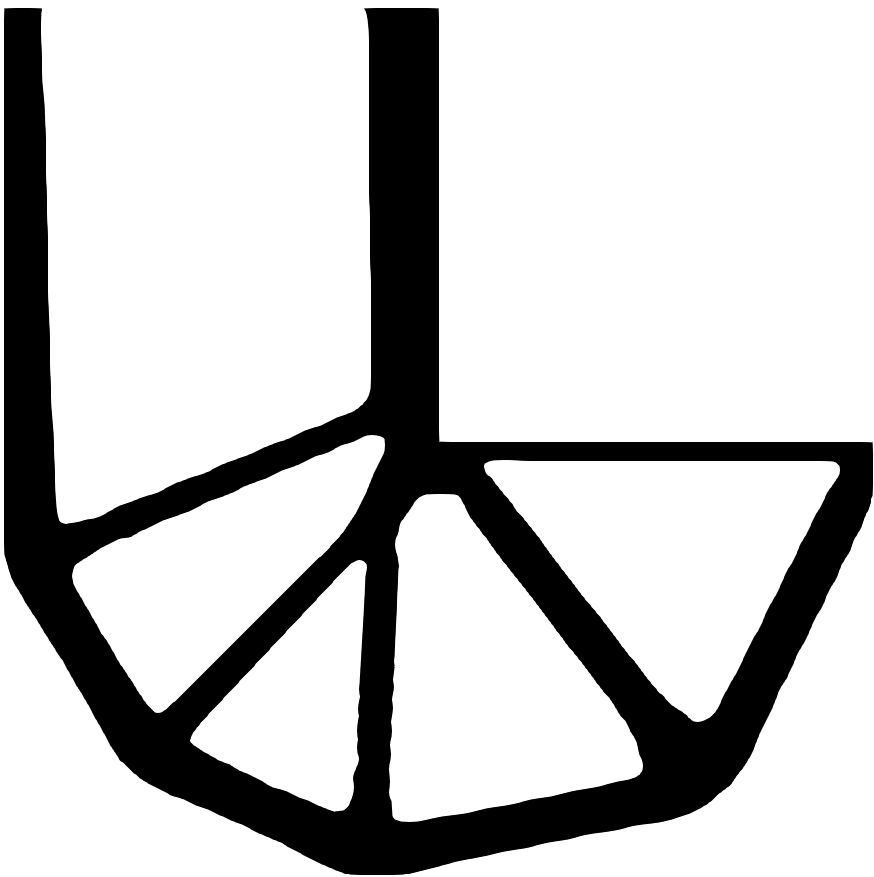}
    \caption{$\bar{N}_{1}=4$}
    \label{fig:lbeam-h1-4}
  \end{subfigure}
  \hspace{2em}
  \begin{subfigure}[b]{0.30\textwidth}
    \centering
    \includegraphics[width=\textwidth]{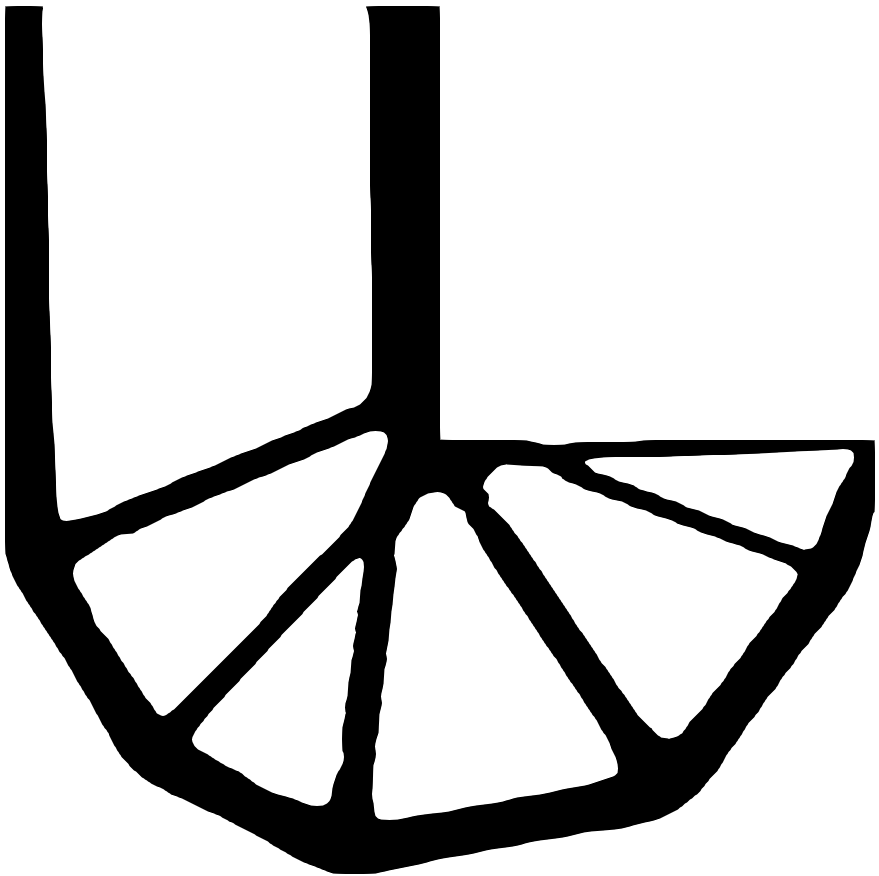}
    \caption{$\bar{N}_{1}=5$}
    \label{fig:lbeam-h1-5}
  \end{subfigure}
  \caption{Optimized L-shaped beam structures for $\bar{N}_{1}=1,2,3,4,5$, respectively.}
  \label{fig:lbeam-control}
\end{figure}

Topology control is then imposed by prescribing an upper bound on
the number of holes, and the optimized structures for $\bar{N}_{1}=1,2,3,4,5$
are shown in Fig.~\ref{fig:lbeam-control}. In every case the resulting number of
holes coincides with $\bar{N}_{1}$, confirming that the
proposed method regulates the topology over a wide range of bounds instead of
merely driving the structure towards a hole-free design.
As reported in Table~\ref{tab:lbeam}, the
compliance decrease moderately as $\bar{N_1}$ increases and approaches the baseline
value obtained without topology control, which exhibits a clear 
trade-off between mechanical performance and topological complexity.\par

\begin{table}[!htb]
  \centering
  \caption{Comparison of the compliance values and of the resulting numbers of holes ($N_{1}$) for the 2D
L-shaped beam under different prescribed bounds $\bar{N}_{1}$.}
  \label{tab:lbeam}
  \small
  \setlength{\tabcolsep}{5pt}
  \begin{tabular}{@{}lcccccc@{}}
    \toprule
    Quantity & $\bar{N}_{1}=1$ & $\bar{N}_{1}=2$ & $\bar{N}_{1}=3$ & $\bar{N}_{1}=4$ & $\bar{N}_{1}=5$ & Baseline\\
    \midrule
    Weights $(\mu_{0},\mu_{1})$ & $(100,100)$ & $(100,100)$ & $(100,100)$ & $(100,100)$ & $(100,100)$ & (0, 0)\\
    Compliance & $149.82$ & $139.14$ & $132.35$ & $130.31$ & $130.14$ & 129.59\\
    $N_1$ & 1 & 2 & 3 & 4 & 5 & 6\\
    \bottomrule
  \end{tabular}
\end{table}

\begin{figure}[!htb]
  \centering
  \begin{subfigure}[b]{0.33\textwidth}
    \centering
    \includegraphics[width=\linewidth]{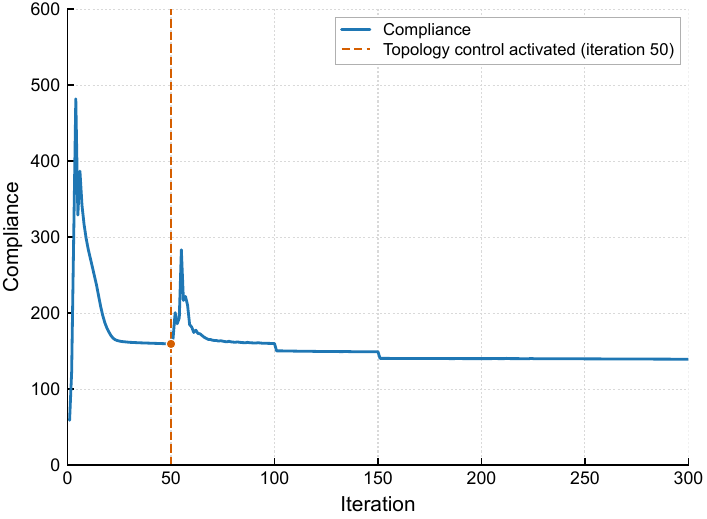}
    \caption{Compliance vs. iteration}
    \label{fig:l-beam-evolution-a}
  \end{subfigure}\hfill
  \begin{subfigure}[b]{0.30\textwidth}
    \centering
    \includegraphics[width=\linewidth]{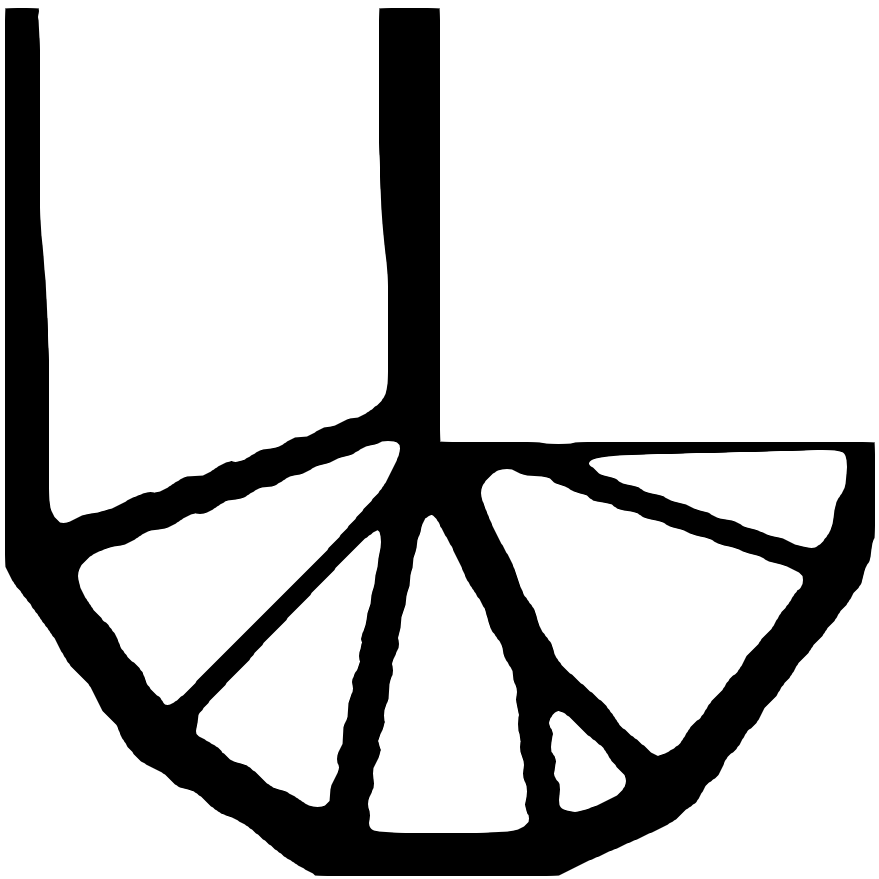}
    \caption{iter${}={}$50: $N_{1}=6$}
    \label{fig:l-beam-evolution-b}
  \end{subfigure}\hfill
  \begin{subfigure}[b]{0.30\textwidth}
    \centering
    \includegraphics[width=\linewidth]{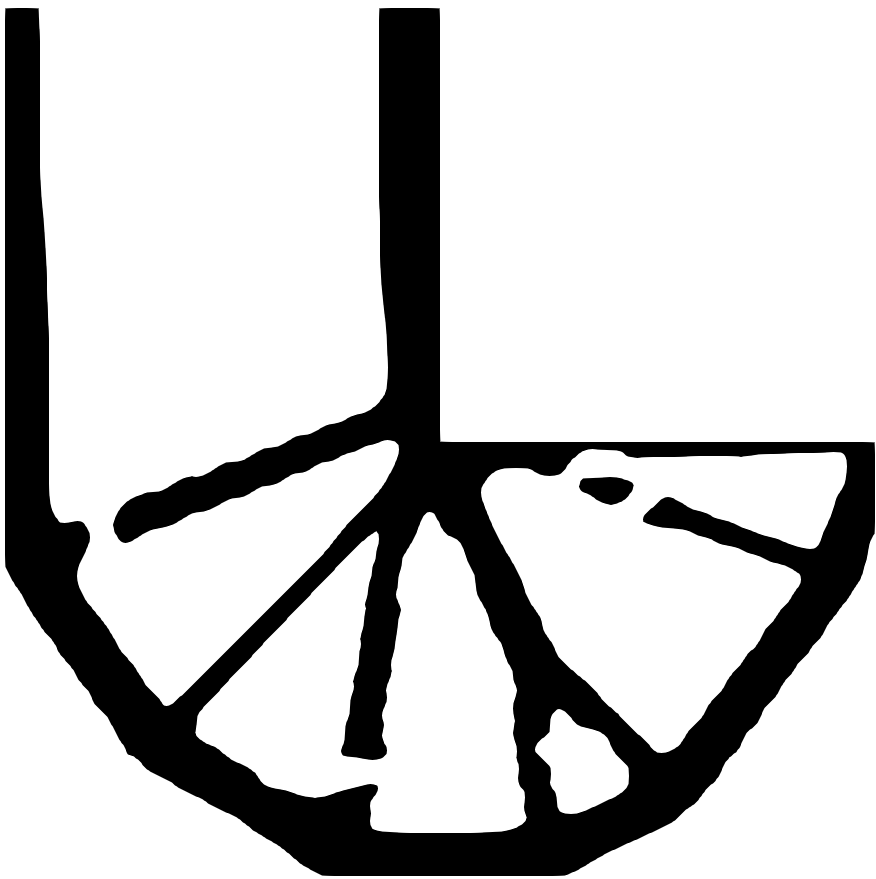}
    \caption{iter${}={}$52: $N_{1}=3$}
    \label{fig:l-beam-evolution-c}
  \end{subfigure}\\[8pt]
  \begin{subfigure}[b]{0.30\textwidth}
    \centering
    \includegraphics[width=\linewidth]{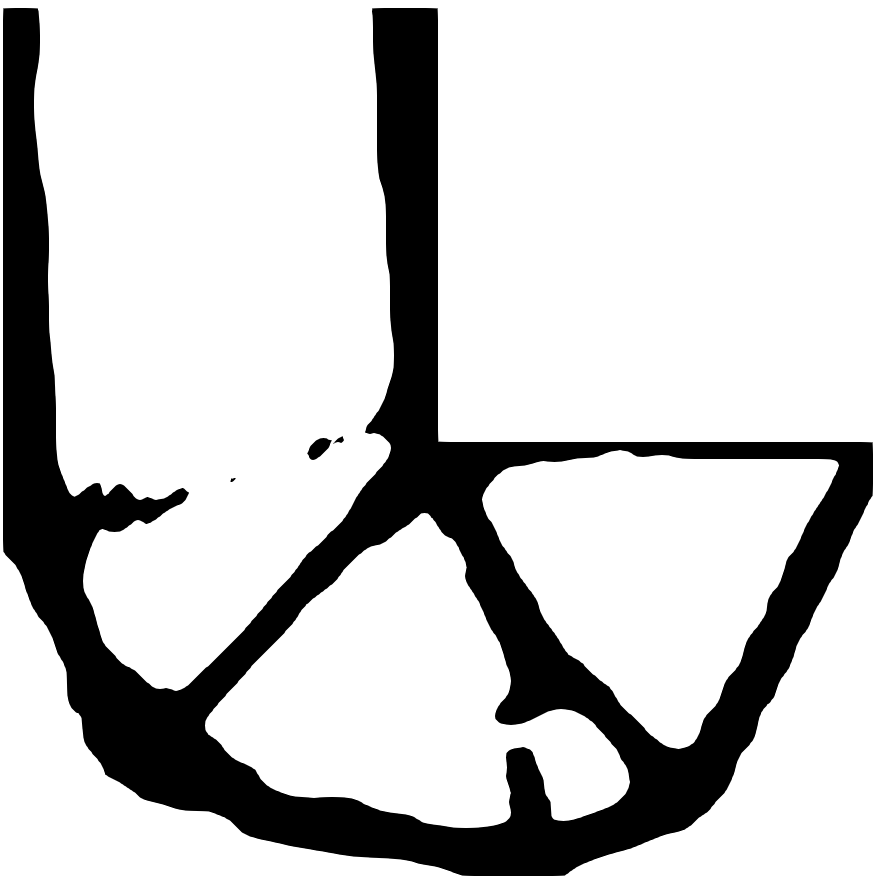}
    \caption{iter${}={}$56: $N_{1}=2$}
    \label{fig:l-beam-evolution-d}
  \end{subfigure}\hfill
  \begin{subfigure}[b]{0.30\textwidth}
    \centering
    \includegraphics[width=\linewidth]{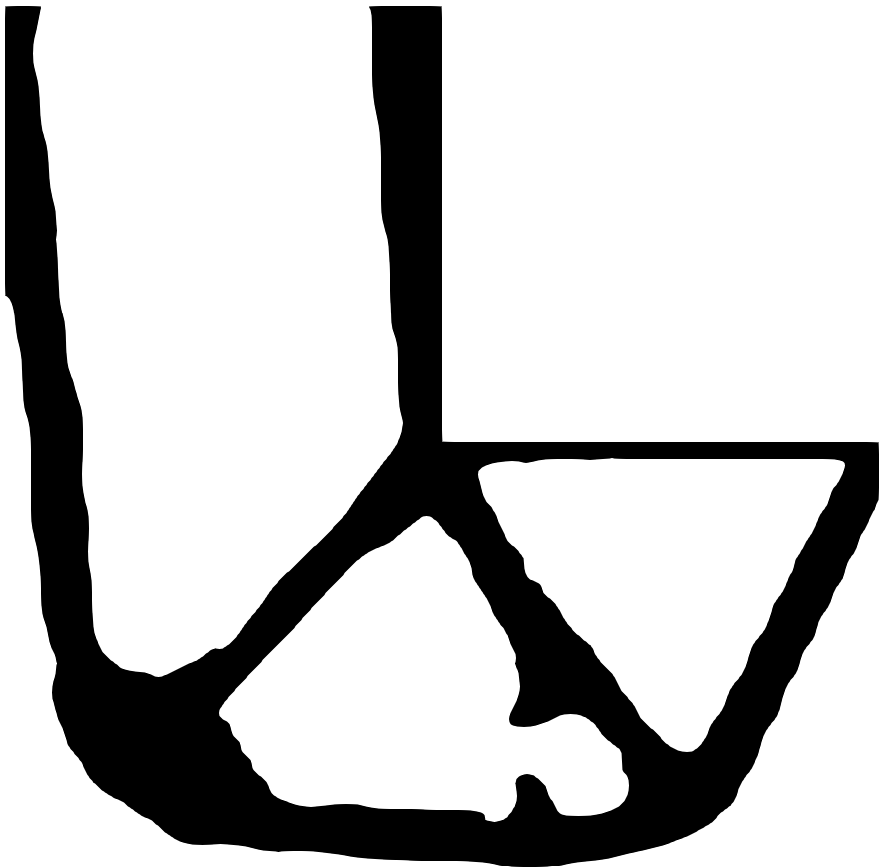}
    \caption{iter${}={}$60: $N_{1}=2$}
    \label{fig:l-beam-evolution-e}
  \end{subfigure}\hfill
  \begin{subfigure}[b]{0.30\textwidth}
    \centering
    \includegraphics[width=\linewidth]{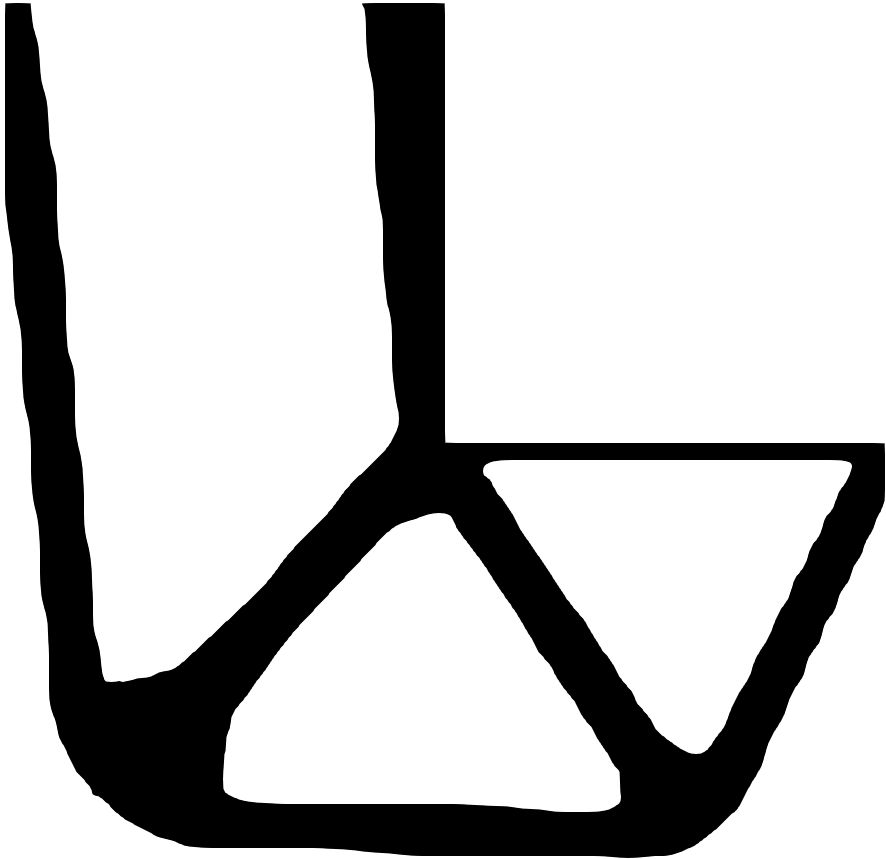}
    \caption{iter${}={}$120: $N_{1}=2$}
    \label{fig:l-beam-evolution-f}
  \end{subfigure}
  \caption{Evolution of the topology for the L-shaped beam when $\bar{N}_{1}=2$. (a) compliance history; (b)--(f) binary structures at representative
    iterations.}
  \label{fig:l-beam-evolution}
\end{figure}

Fig.~\ref{fig:l-beam-evolution} shows the compliance history together with the
topology evolution for the case $\bar{N}_{1}=2$. The topology control is activated
at iteration~50, at which point the structure contains six holes and therefore
exceeds the prescribed bound. The redundant holes are eliminated by cutting along
their boundaries, so that the number of holes falls to three within two iterations
and reaches the prescribed value of two by iteration~56, as shown in
Figs.~\ref{fig:l-beam-evolution-b}--\ref{fig:l-beam-evolution-d}. This stage is
accompanied by the appearance of transient extra connected components and by a
temporary increase in the compliance, and these components are subsequently
removed by the 0-dimensional topological objective. Once the prescribed topology
is attained at iteration~60, the remaining iterations are devoted to compliance
minimization, as shown by the nearly unchanged structures at iterations~60 and~120
(Figs.~\ref{fig:l-beam-evolution-e} and~\ref{fig:l-beam-evolution-f}).\par

\subsubsection{Quarter annulus}
\begin{figure}[!htb]
  \centering
  \begin{subfigure}[t]{0.30\textwidth}
    \centering
    \includegraphics[width=\linewidth]{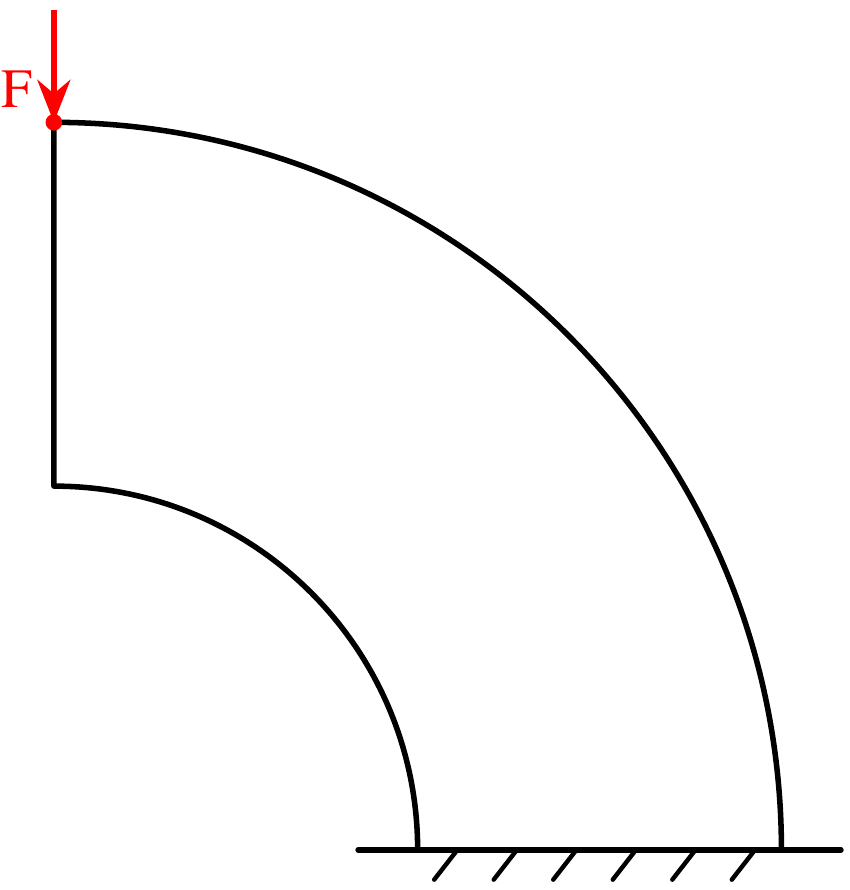}
    \caption{}
    \label{fig:annulus-model}
  \end{subfigure}
  \hfill
  \begin{subfigure}[t]{0.36\textwidth}
    \centering
    \includegraphics[width=\linewidth]{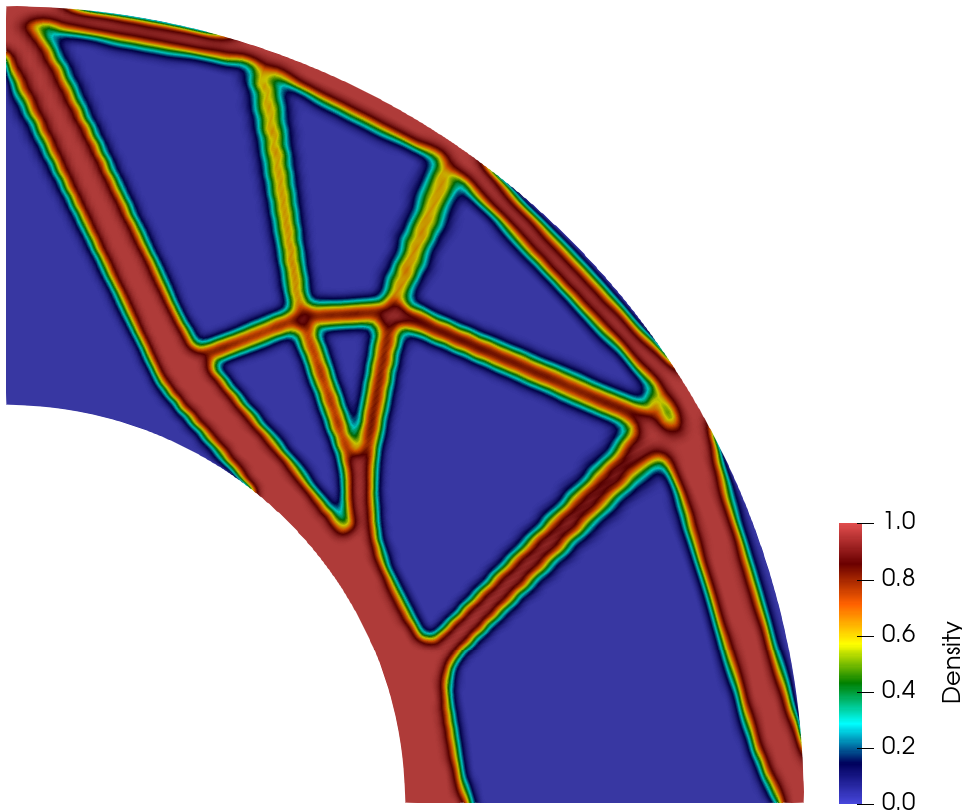}
    \caption{}
    \label{fig:annulus-density}
  \end{subfigure}
  \hfill
  \begin{subfigure}[t]{0.30\textwidth}
    \centering
    \includegraphics[width=\linewidth]{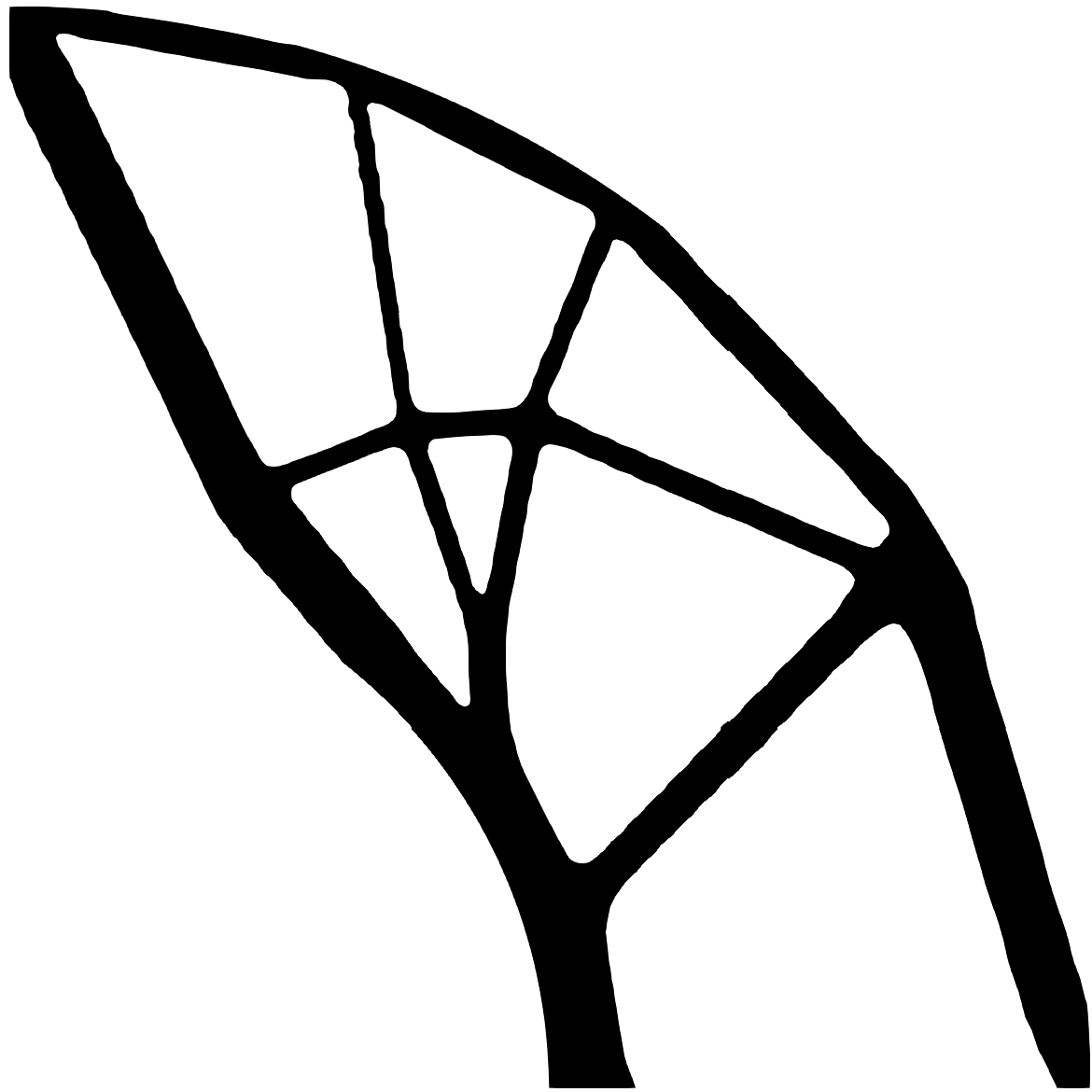}
    \caption{}
    \label{fig:annulus-topology}
  \end{subfigure}
  \caption{(a) The quarter annulus model.
   (b) Optimized density field without topology control.
   (c) The corresponding binary structure with $N_1=6$.}
  \label{fig:annulus-nocontrol}
  \end{figure}

In this example, we consider the classic quarter annulus benchmark, whose inner and outer radii are $5$ and $10$, respectively. 
As illustrated in Fig.~\ref{fig:annulus-model}, the bottom edge of the design domain is fully clamped, 
while a unit downward load is applied at the top of the outer arc.
This domain is represented exactly by a biquadratic NURBS surface with $102\times 52$ control points. 
With this discretization in place, the optimization is performed subject to a volume constraint, where the maximum volume fraction is set to $V_{0}=0.3$.
When no topology control is imposed, the optimized structure exhibits six holes, as shown in Figs.~\ref{fig:annulus-density} and \ref{fig:annulus-topology}.\par

\begin{figure}[!htb]
  \centering
  \begin{subfigure}[b]{0.30\textwidth}
    \centering
    \includegraphics[width=\textwidth]{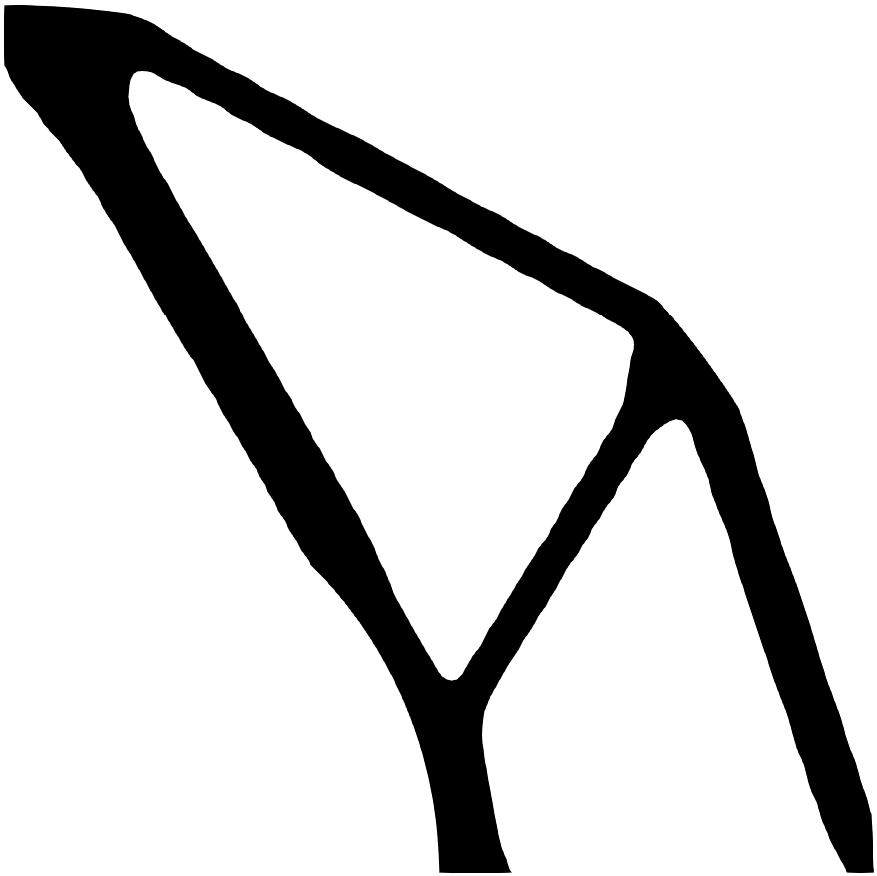}
    \caption{$\bar{N}_{1}=1$}
    \label{fig:annulus-h1-1}
  \end{subfigure}
  \hfill
  \begin{subfigure}[b]{0.30\textwidth}
    \centering
    \includegraphics[width=\textwidth]{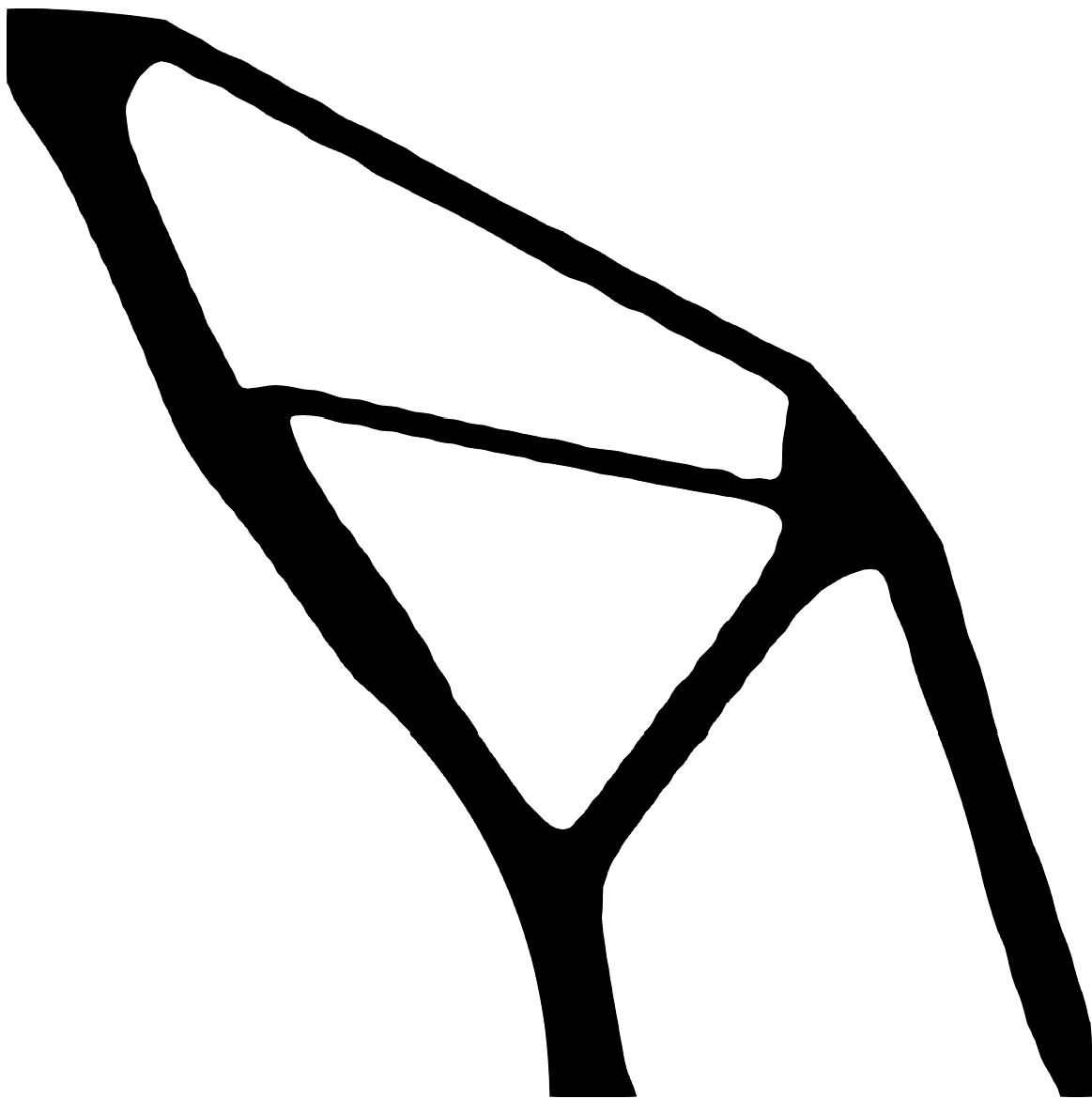}
    \caption{$\bar{N}_{1}=2$}
    \label{fig:annulus-h1-2}
  \end{subfigure}
  \hfill
  \begin{subfigure}[b]{0.30\textwidth}
    \centering
    \includegraphics[width=\textwidth]{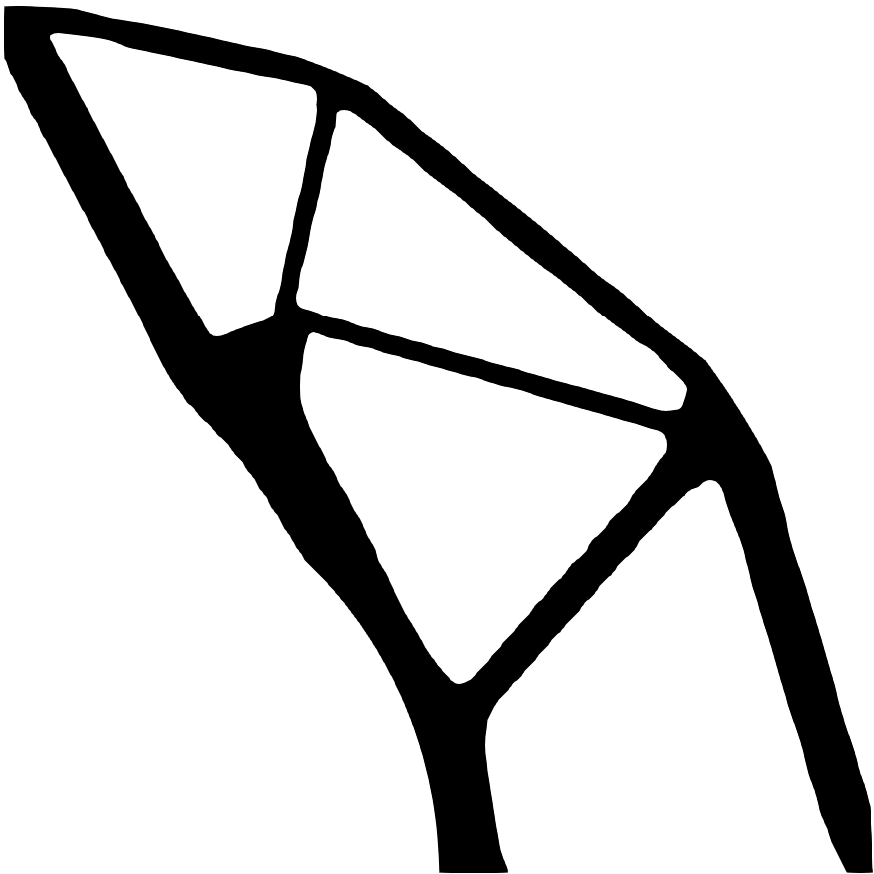}
    \caption{$\bar{N}_{1}=3$}
    \label{fig:annulus-h1-3}
  \end{subfigure}\\[8pt]
  \begin{subfigure}[b]{0.30\textwidth}
    \centering
    \includegraphics[width=\textwidth]{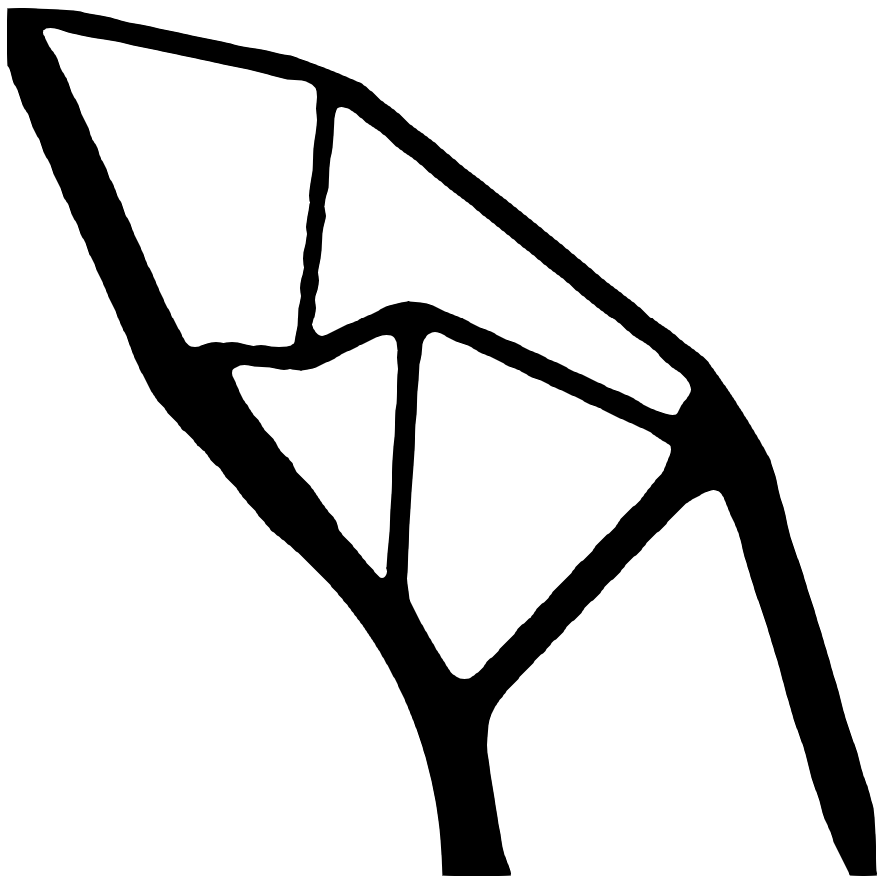}
    \caption{$\bar{N}_{1}=4$}
    \label{fig:annulus-h1-4}
  \end{subfigure}
  \hspace{2em}
  \begin{subfigure}[b]{0.30\textwidth}
    \centering
    \includegraphics[width=\textwidth]{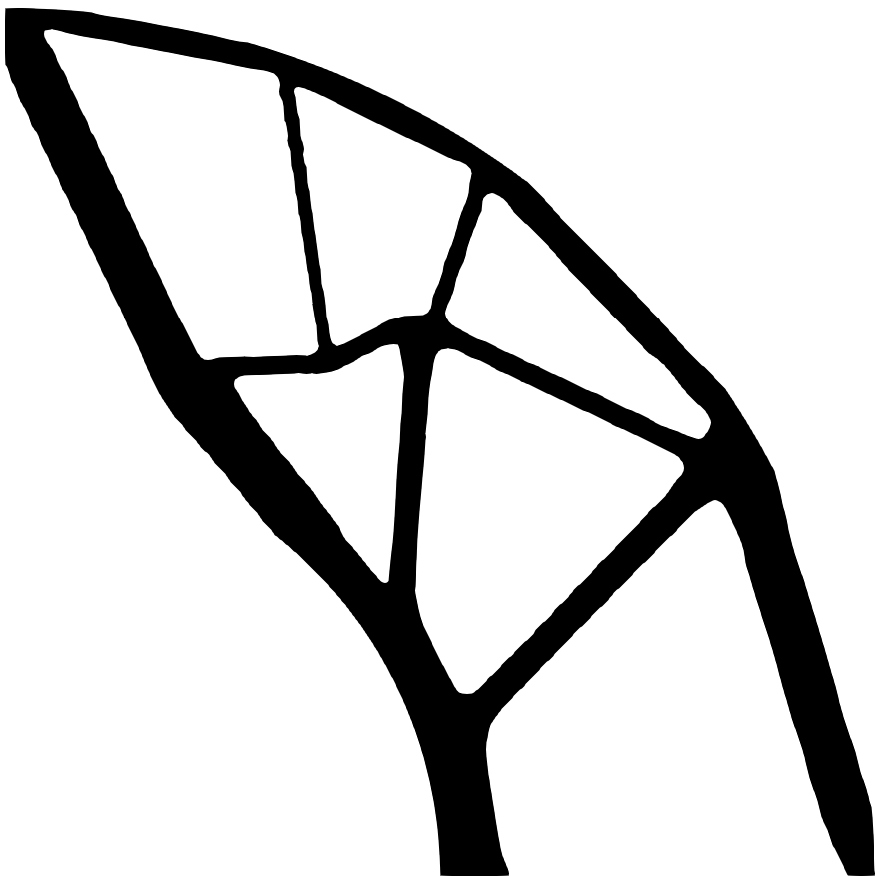}
    \caption{$\bar{N}_{1}=5$}
    \label{fig:annulus-h1-5}
  \end{subfigure}
  \caption{Optimized quarter annulus structures for $\bar{N}_{1}=1,2,3,4,5$, respectively.}
  \label{fig:annulus-control}
\end{figure}

Table~\ref{tab:annulus} and Fig.~\ref{fig:annulus-control} summarize the
topology-controlled results for $\bar{N}_{1}=1,2,3,4,5$. Overall, a smaller
admissible number of holes leads to a moderate increase of the compliance, which
is consistent with the reduced design freedom under a tighter topological
constraint. Unlike the preceding examples, however, the compliance does not vary
monotonically with $\bar{N}_{1}$, since the compliance obtained for
$\bar{N}_{1}=2$ is slightly higher than that for $\bar{N}_{1}=1$. This suggests
that different bounds on the number of holes lead the optimizer towards different
local optima, so that the mechanical cost of the constraint is not necessarily a
monotone function of $\bar{N}_{1}$.\par

\begin{table}[!htb]
  \centering
  \caption{Comparison of the compliance values and of the resulting numbers of holes ($N_{1}$) for the 2D
quarter annulus under different prescribed bounds $\bar{N}_{1}$.}
  \label{tab:annulus}
  \small
  \setlength{\tabcolsep}{5pt}
  \begin{tabular}{@{}lcccccc@{}}
    \toprule
    Quantity & $\bar{N}_{1}=1$ & $\bar{N}_{1}=2$ & $\bar{N}_{1}=3$ & $\bar{N}_{1}=4$ & $\bar{N}_{1}=5$ & Baseline\\
    \midrule
    Weights $(\mu_{0},\mu_{1})$ & $(100,200)$ & $(100,200)$ & $(100,200)$ & $(100,200)$ & $(100,200)$ & $(0, 0)$\\
    Compliance & $104.08$ & $104.64$ & $101.76$ & $100.33$ & $98.50$ & $98.16$\\
    $N_1$ & 1 & 2 & 3 & 4 & 5 & 6\\
    \bottomrule
  \end{tabular}
\end{table}

\begin{figure}[!htb]
  \centering
  \begin{subfigure}[b]{0.33\textwidth}
    \centering
    \includegraphics[width=\linewidth]{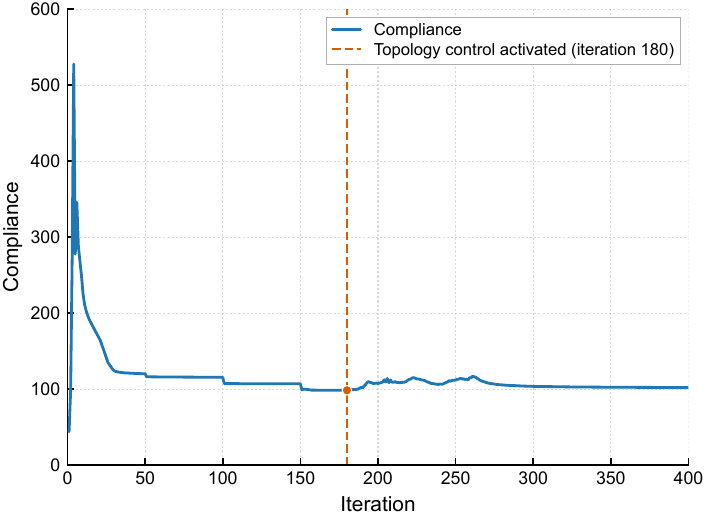}
    \caption{Compliance vs. iteration}
    \label{fig:annulus-evolution-a}
  \end{subfigure}\hfill
  \begin{subfigure}[b]{0.30\textwidth}
    \centering
    \includegraphics[width=\linewidth]{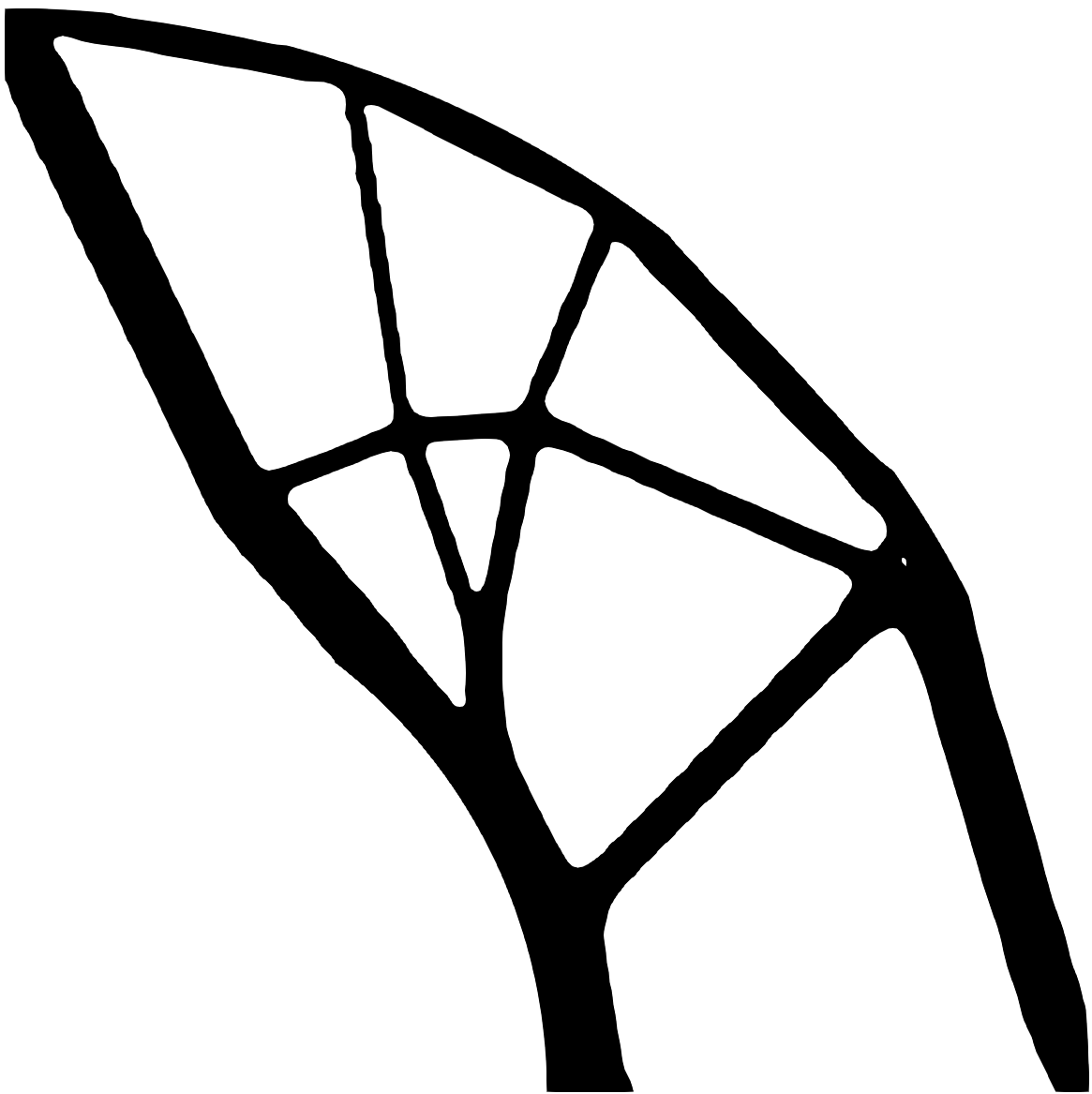}
    \caption{iter${}={}$180: $N_{1}=7$}
    \label{fig:annulus-evolution-b}
  \end{subfigure}\hfill
  \begin{subfigure}[b]{0.30\textwidth}
    \centering
    \includegraphics[width=\linewidth]{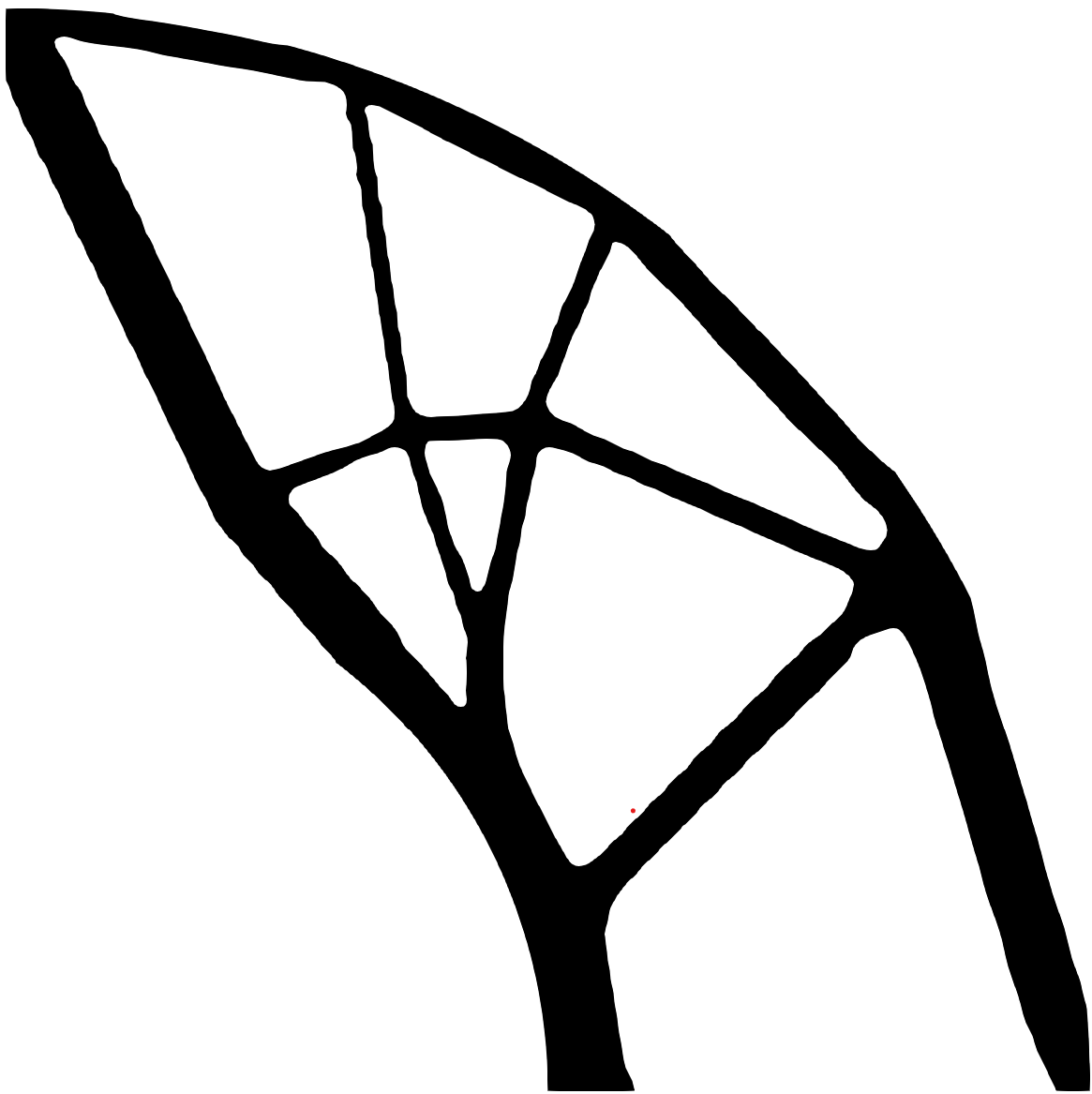}
    \caption{iter${}={}$181: $N_{1}=6$}
    \label{fig:annulus-evolution-c}
  \end{subfigure}\\[8pt]
  \begin{subfigure}[b]{0.30\textwidth}
    \centering
    \includegraphics[width=\linewidth]{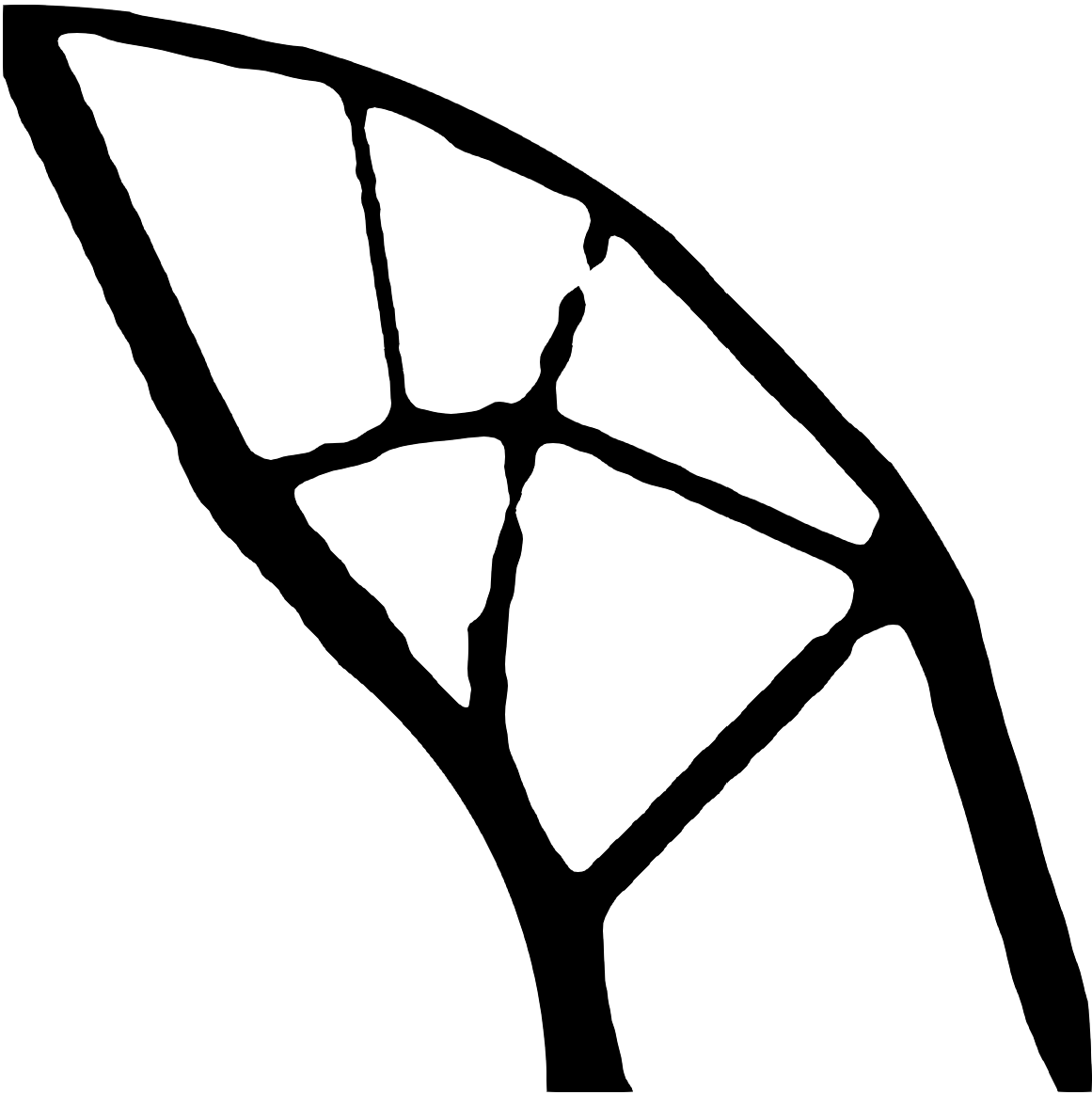}
    \caption{iter${}={}$201: $N_{1}=4$}
    \label{fig:annulus-evolution-d}
  \end{subfigure}\hfill
  \begin{subfigure}[b]{0.30\textwidth}
    \centering
    \includegraphics[width=\linewidth]{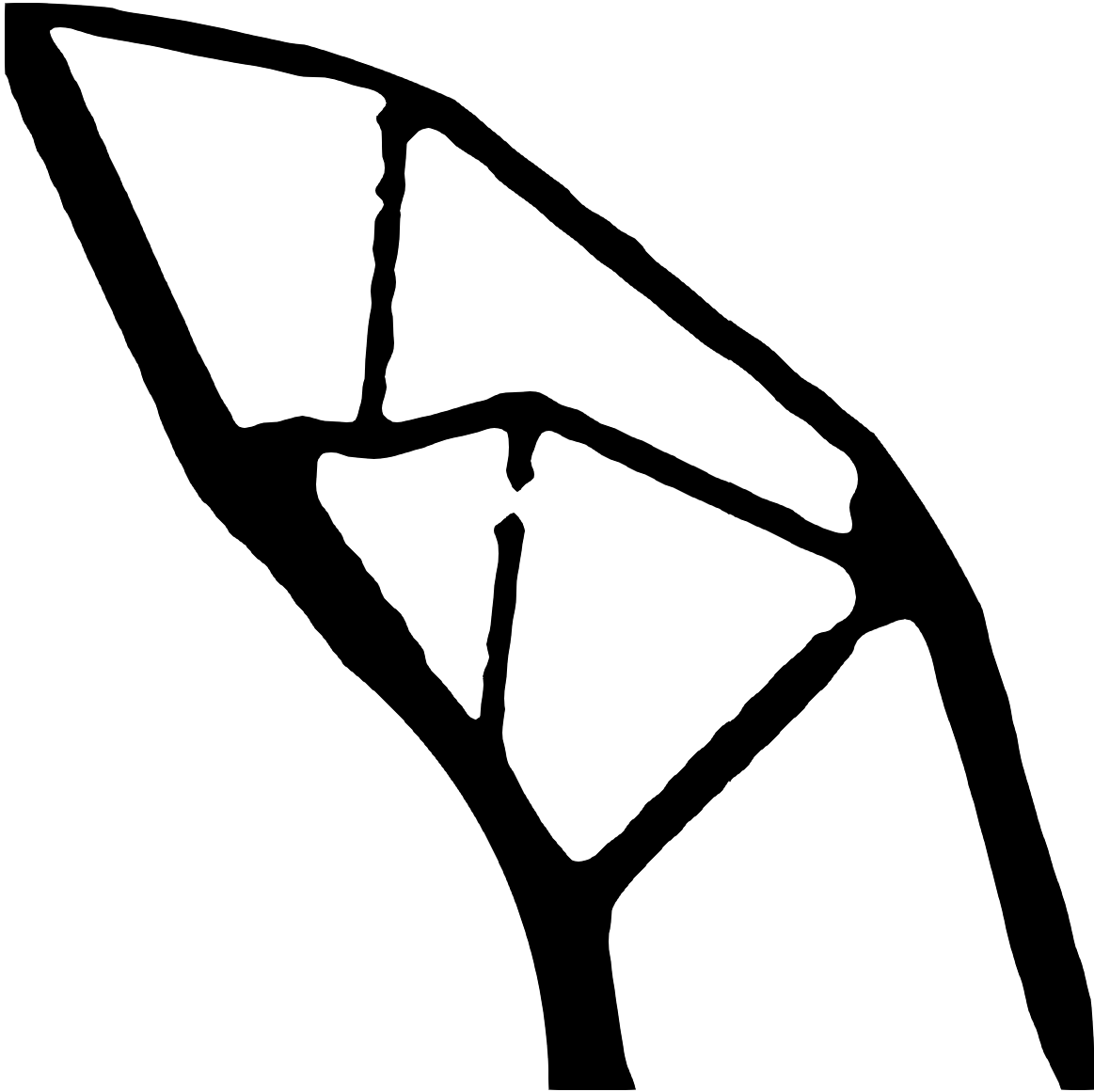}
    \caption{iter${}={}$245: $N_{1}=3$}
    \label{fig:annulus-evolution-e}
  \end{subfigure}\hfill
  \begin{subfigure}[b]{0.30\textwidth}
    \centering
    \includegraphics[width=\linewidth]{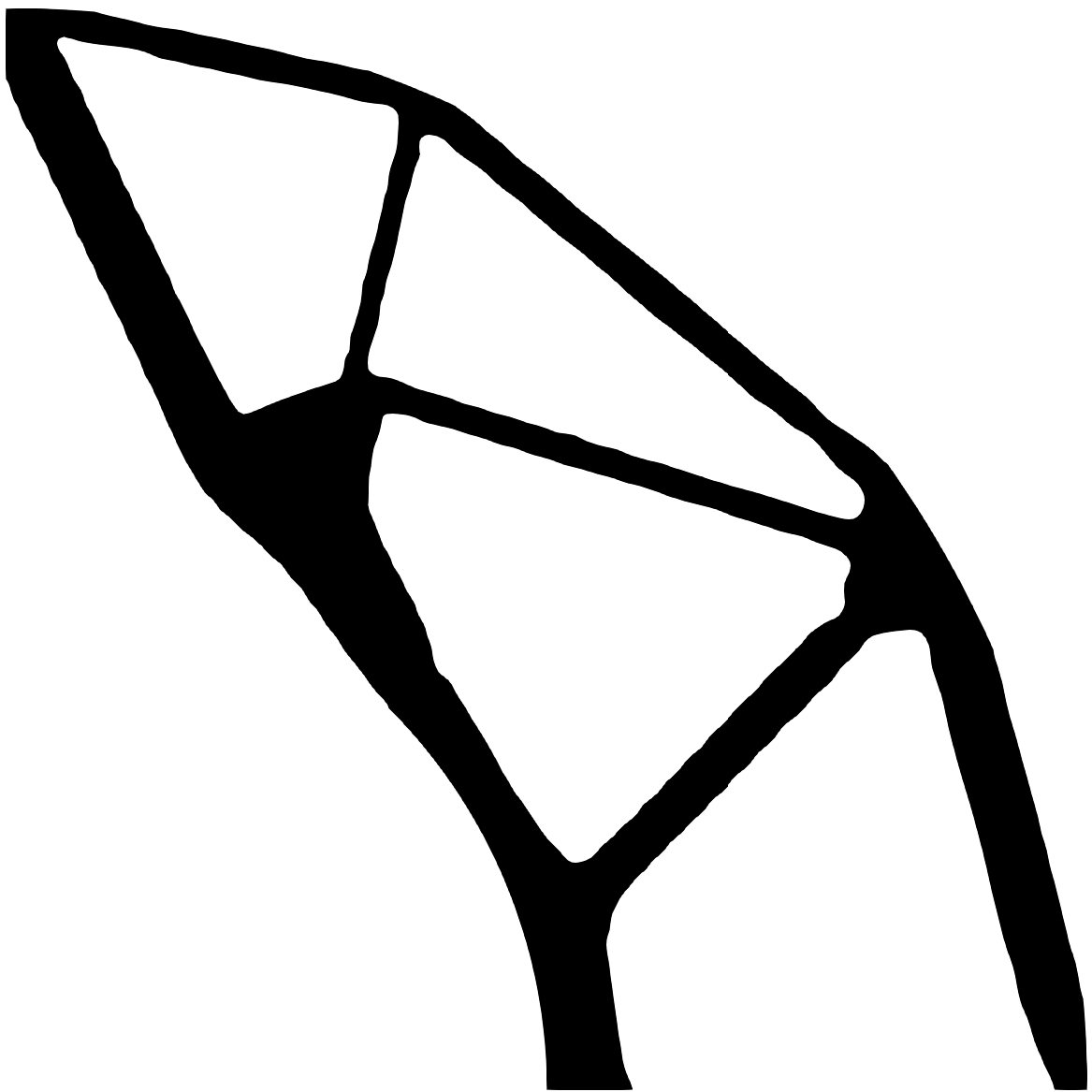}
    \caption{iter${}={}$290: $N_{1}=3$}
    \label{fig:annulus-evolution-f}
  \end{subfigure}
  \caption{Evolution of the topology for the quarter annulus when $\bar{N}_{1}=3$. (a) compliance history; (b)--(f) binary structures at representative
    iterations.}
  \label{fig:annulus-evolution}
\end{figure}

Fig.~\ref{fig:annulus-evolution} shows the compliance history together with the
topology evolution for $\bar{N}_{1}=3$. The topology control is switched on at
iteration~180, where the structure still has seven holes, four more than allowed.
Since the holes to be removed are ordered by their area, 
the smallest one is taken first and filled already in the very next
iteration, in agreement with the proposed removal strategy.
The remaining redundant holes are then eliminated gradually,
so that the number of holes drops to five and four and finally attains the prescribed
value of three around iteration~245, as illustrated in
Figs.~\ref{fig:annulus-evolution-b}--\ref{fig:annulus-evolution-e}.
Subsequently, the optimizer concentrates on
redistributing the remaining material: the ineffective branches are gradually
pruned away, and by iteration~290 the layout has settled into a smooth and stable
configuration (Fig.~\ref{fig:annulus-evolution-f}).\par

\subsubsection{Michell structure}

\begin{figure}[!htb]
  \centering
  \begin{subfigure}[t]{0.30\textwidth}
    \centering
    \includegraphics[width=\linewidth]{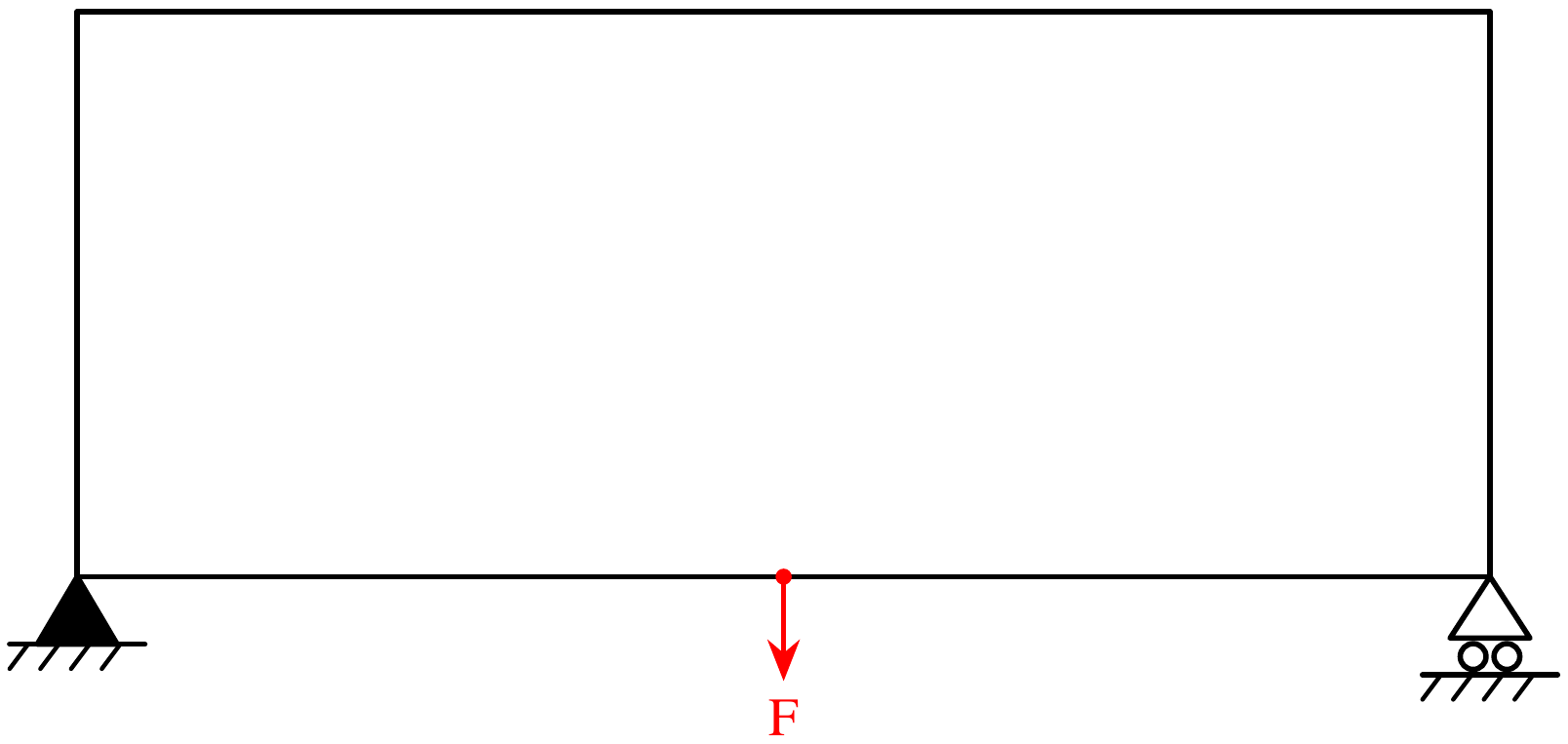}
    \caption{}
    \label{fig:michell-model}
  \end{subfigure}
  \hfill
  \begin{subfigure}[t]{0.36\textwidth}
    \centering
    \includegraphics[width=\linewidth]{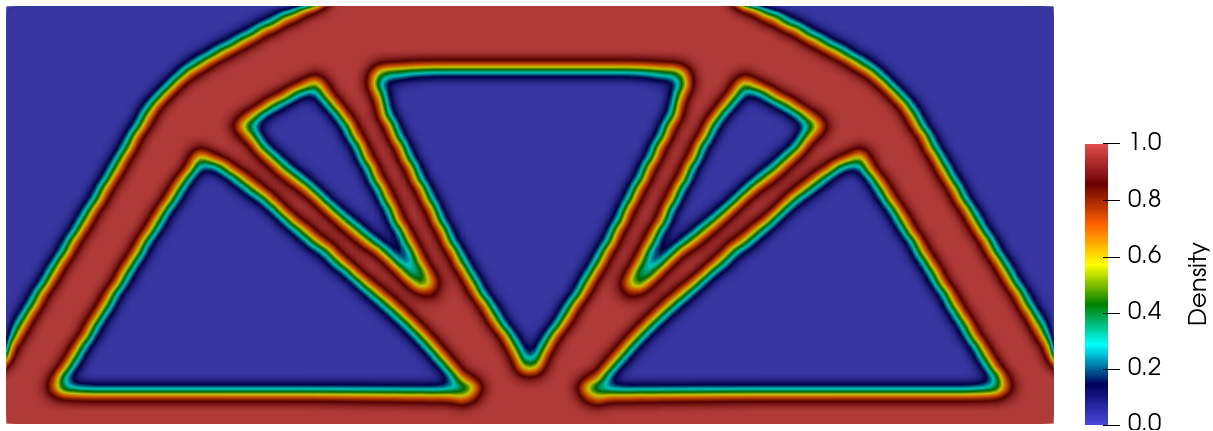}
    \caption{}
    \label{fig:michell-density}
  \end{subfigure}
  \hfill
  \begin{subfigure}[t]{0.30\textwidth}
    \centering
    \includegraphics[width=\linewidth]{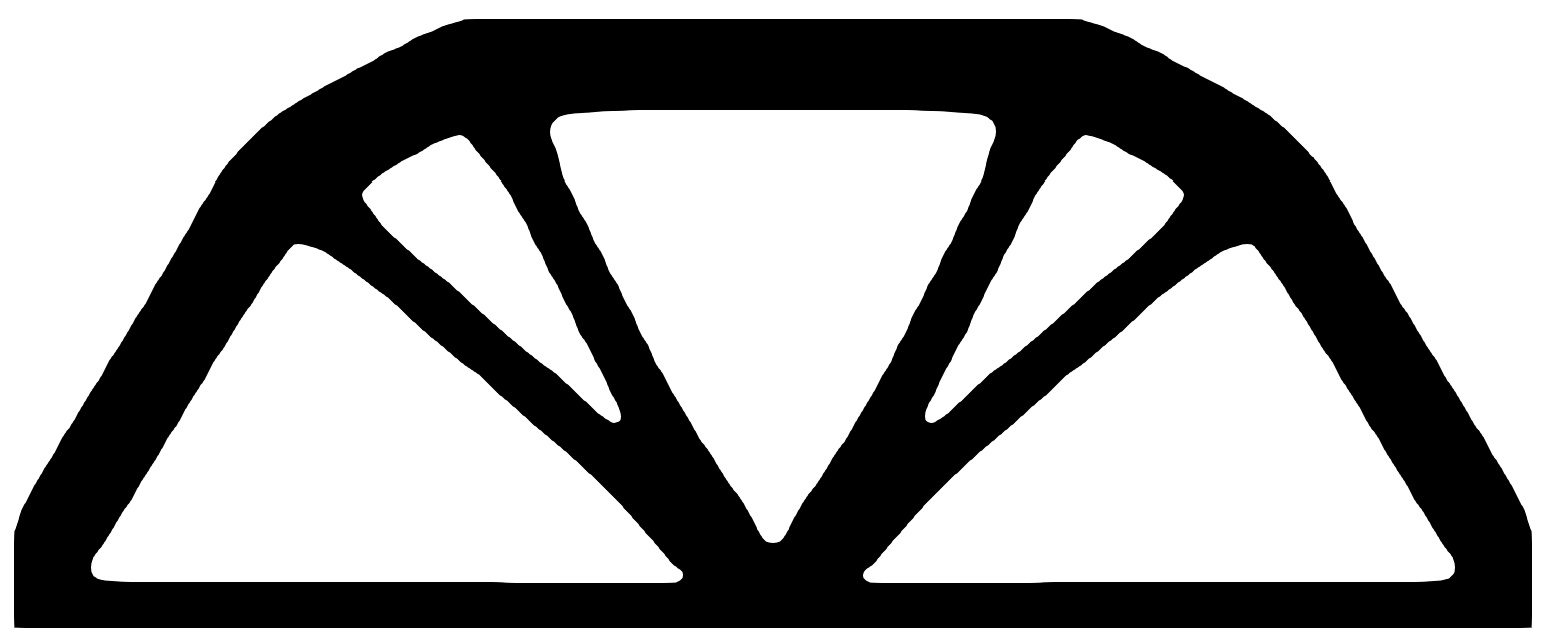}
    \caption{}
    \label{fig:michell-topology}
  \end{subfigure}
  \caption{(a) The Michell structure model.
   (b) Optimized density field without topology control.
   (c) The corresponding binary structure with $N_1=5$.}
  \label{fig:michell-nocontrol}
  \end{figure}

In this subsection, the Michell structure benchmark is considered. The rectangular
design domain $[0,10]\times[0,4]$ is represented by a biquadratic NURBS surface
with $102\times 42$ control points. Michell-type boundary conditions are adopted, namely
a pin support at the lower-left corner, a roller support at the lower-right corner,
and a unit downward load applied at the midpoint of the bottom edge, as
illustrated in Fig.~\ref{fig:michell-model}. The maximum volume fraction is set to
$V_{0}=0.4$. As shown in Figs.~\ref{fig:michell-density} and
\ref{fig:michell-topology}, the optimized structure without topology control
contains five holes.

\begin{figure}[!htb]
  \centering
  \begin{subfigure}[b]{0.40\textwidth}
    \centering
    \includegraphics[width=\textwidth]{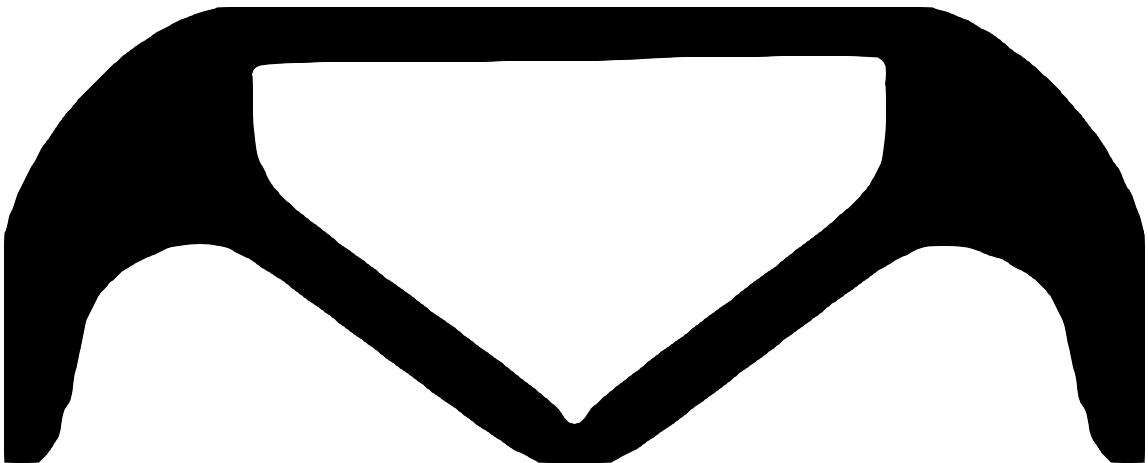}
    \caption{$\bar{N}_{1}=1$}
    \label{fig:michell-h1-1}
  \end{subfigure}
  \hfill
  \begin{subfigure}[b]{0.40\textwidth}
    \centering
    \includegraphics[width=\textwidth]{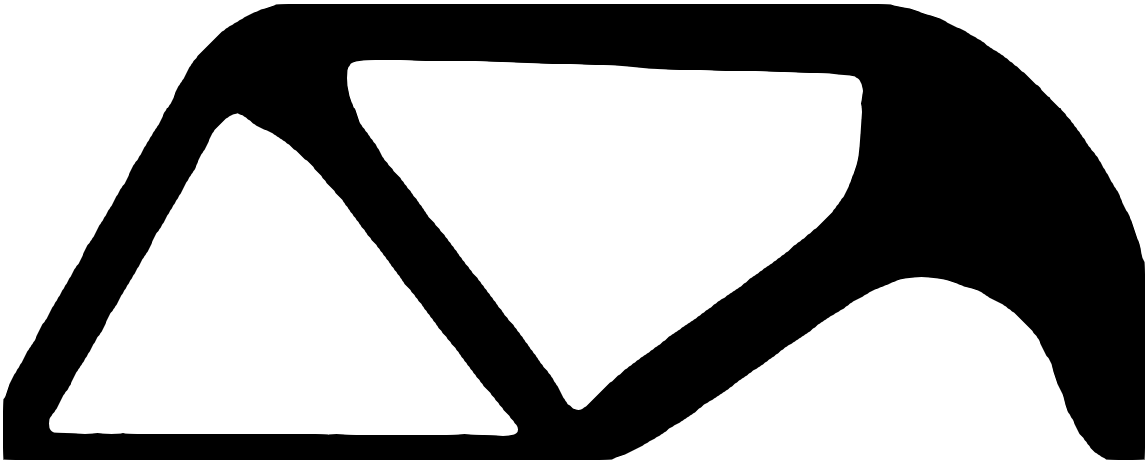}
    \caption{$\bar{N}_{1}=2$}
    \label{fig:michell-h1-2}
  \end{subfigure}\\[8pt]
  \begin{subfigure}[b]{0.40\textwidth}
    \centering
    \includegraphics[width=\textwidth]{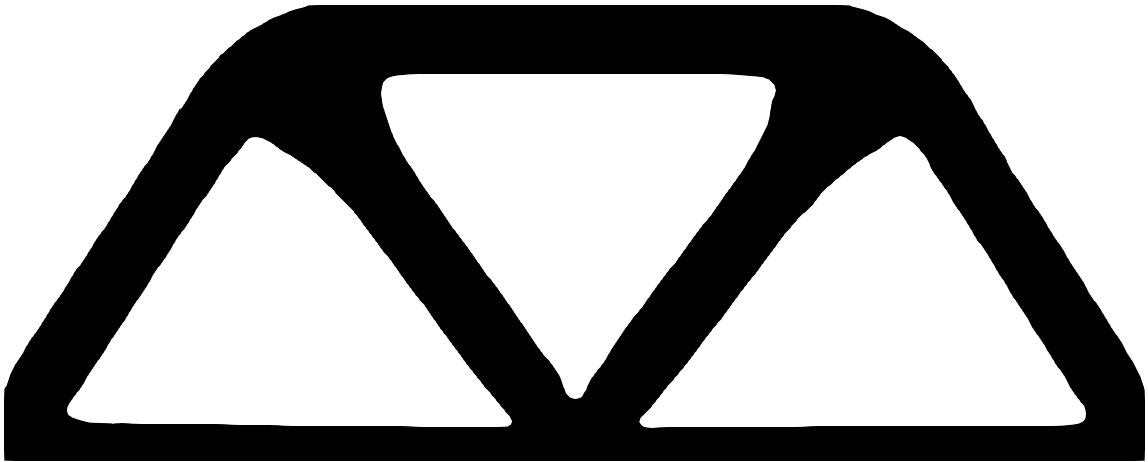}
    \caption{$\bar{N}_{1}=3$}
    \label{fig:michell-h1-3}
  \end{subfigure}
  \hfill
  \begin{subfigure}[b]{0.40\textwidth}
    \centering
    \includegraphics[width=\textwidth]{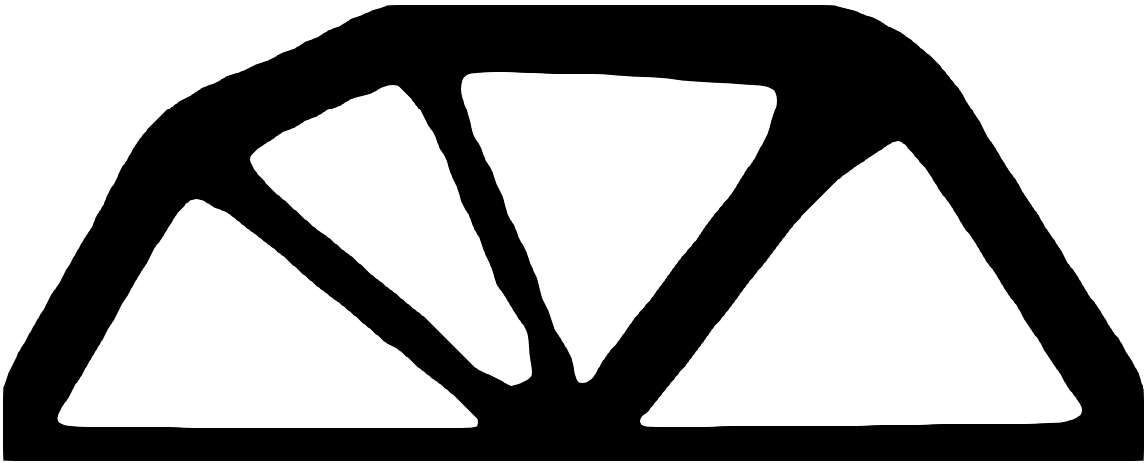}
    \caption{$\bar{N}_{1}=4$}
    \label{fig:michell-h1-4}
  \end{subfigure}
  \caption{Optimized Michell structures for $\bar{N}_{1}=1,2,3,4$, respectively.}
  \label{fig:michell-control}
\end{figure}

Table~\ref{tab:michell} and Fig.~\ref{fig:michell-control} summarize the
topology-controlled results for $\bar{N}_{1}=1,2,3,4$.
A similar performance–complexity trade-off is observed, namely that restricting the admissible number of
holes raises the compliance. A further observation concerns symmetry.
Although the unconstrained solution is indeed a symmetric layout (Fig.~\ref{fig:michell-topology}),
the optimized designs may become non-symmetric once topology control is imposed.
This behavior is expected, because the topological constraint restricts the admissible hole configurations and therefore
can break the symmetry even when the loading and boundary conditions are symmetric. 
The resulting designs remain mechanically efficient and satisfy the prescribed topology.\par

\begin{table}[!htb]
  \centering
  \caption{Comparison of the compliance values and of the resulting numbers of holes ($N_{1}$) for the 2D
Michell structure under different prescribed bounds $\bar{N}_{1}$.}
  \label{tab:michell}
  \small
  \begin{tabular}{@{}lccccc@{}}
    \toprule
    Quantity & $\bar{N}_{1}=1$ & $\bar{N}_{1}=2$ & $\bar{N}_{1}=3$ & $\bar{N}_{1}=4$ & Baseline\\
    \midrule
    Weights $(\mu_{0},\mu_{1})$ & $(100,100)$ & $(100,100)$ & $(100,100)$ & $(100,100)$ & $(0, 0)$\\
    Compliance & $32.31$ & $26.50$ & $22.42$ & $22.11$ & $21.75$\\
    $N_1$ & 1 & 2 & 3 & 4 & 5\\
    \bottomrule
  \end{tabular}
\end{table}

\begin{figure}[!htb]
  \centering
  \begin{subfigure}[b]{0.33\textwidth}
    \centering
    \includegraphics[width=\linewidth]{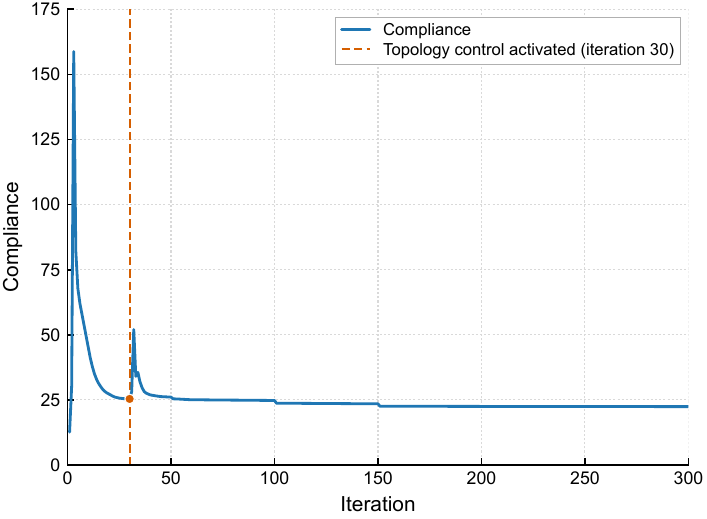}
    \caption{Compliance vs. iteration}
    \label{fig:michell-evolution-a}
  \end{subfigure}\hfill
  \begin{subfigure}[b]{0.30\textwidth}
    \centering
    \includegraphics[width=\linewidth]{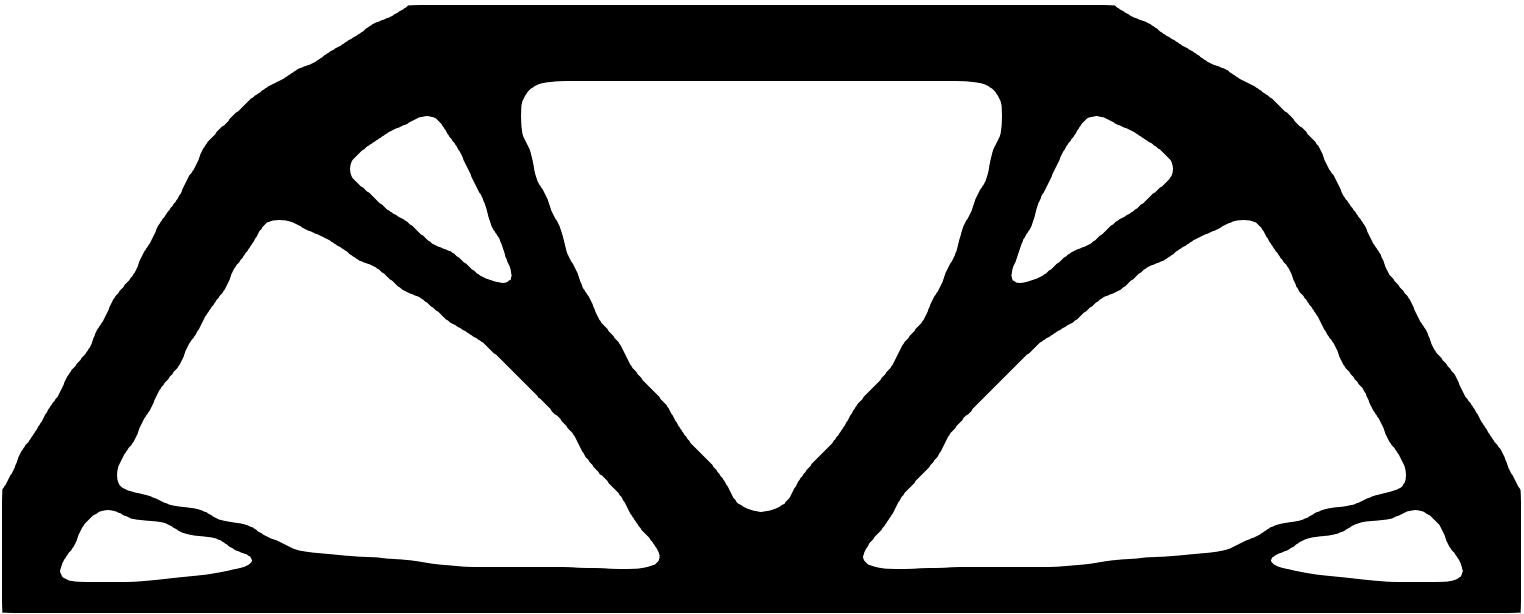}
    \caption{iter${}={}$30: $N_{1}=7$}
    \label{fig:michell-evolution-b}
  \end{subfigure}\hfill
  \begin{subfigure}[b]{0.30\textwidth}
    \centering
    \includegraphics[width=\linewidth]{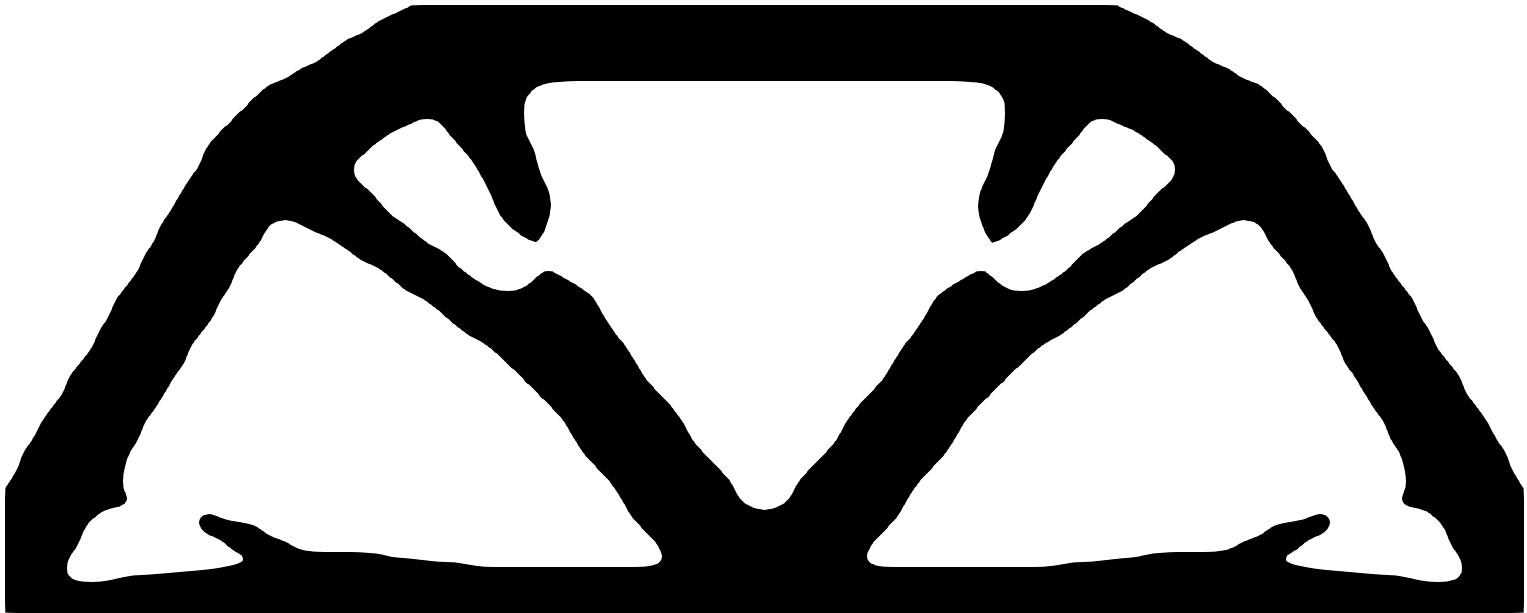}
    \caption{iter${}={}$31: $N_{1}=3$}
    \label{fig:michell-evolution-c}
  \end{subfigure}\\[8pt]
  \begin{subfigure}[b]{0.30\textwidth}
    \centering
    \includegraphics[width=\linewidth]{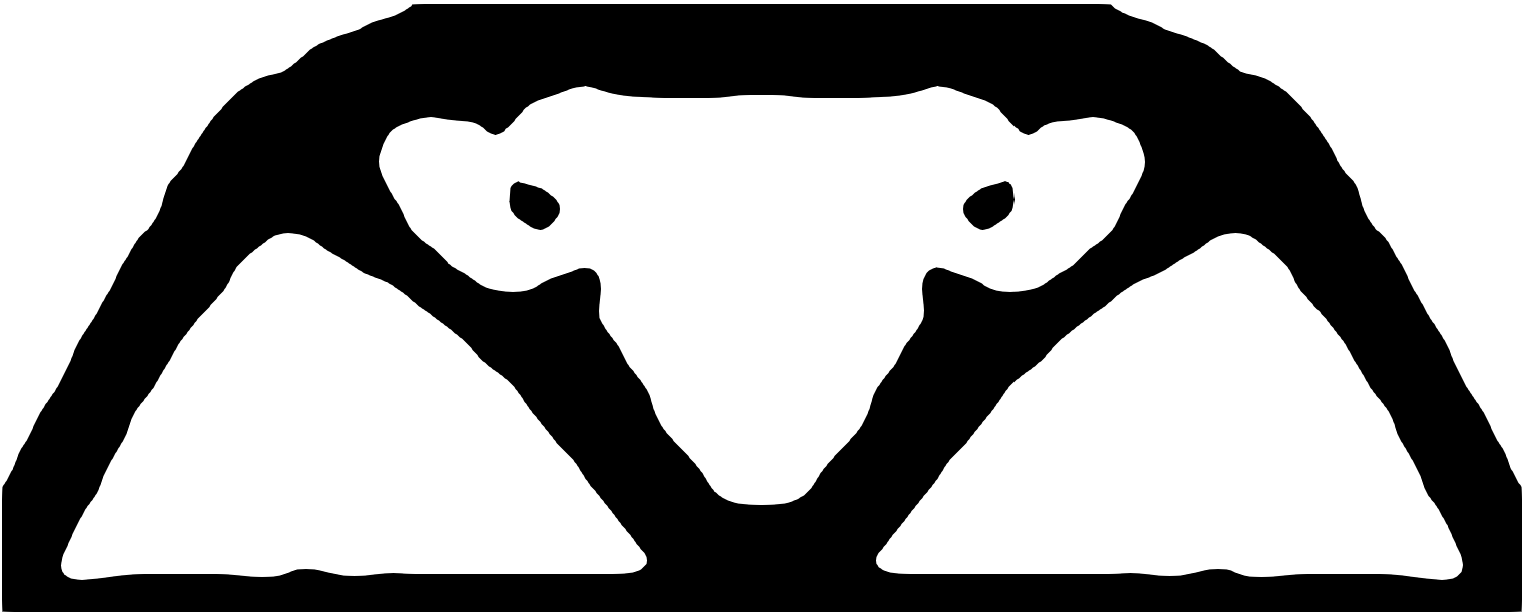}
    \caption{iter${}={}$33: $N_{1}=3$}
    \label{fig:michell-evolution-d}
  \end{subfigure}\hfill
  \begin{subfigure}[b]{0.30\textwidth}
    \centering
    \includegraphics[width=\linewidth]{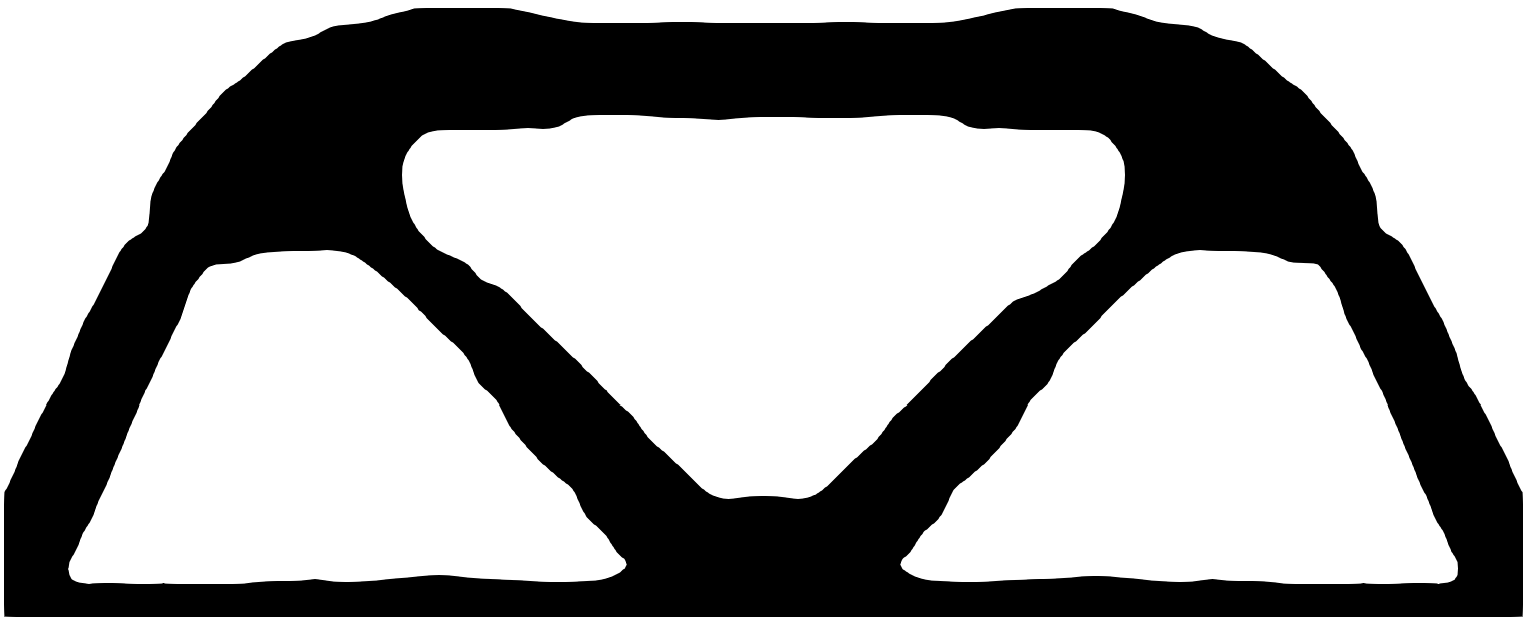}
    \caption{iter${}={}$36: $N_{1}=3$}
    \label{fig:michell-evolution-e}
  \end{subfigure}\hfill
  \begin{subfigure}[b]{0.30\textwidth}
    \centering
    \includegraphics[width=\linewidth]{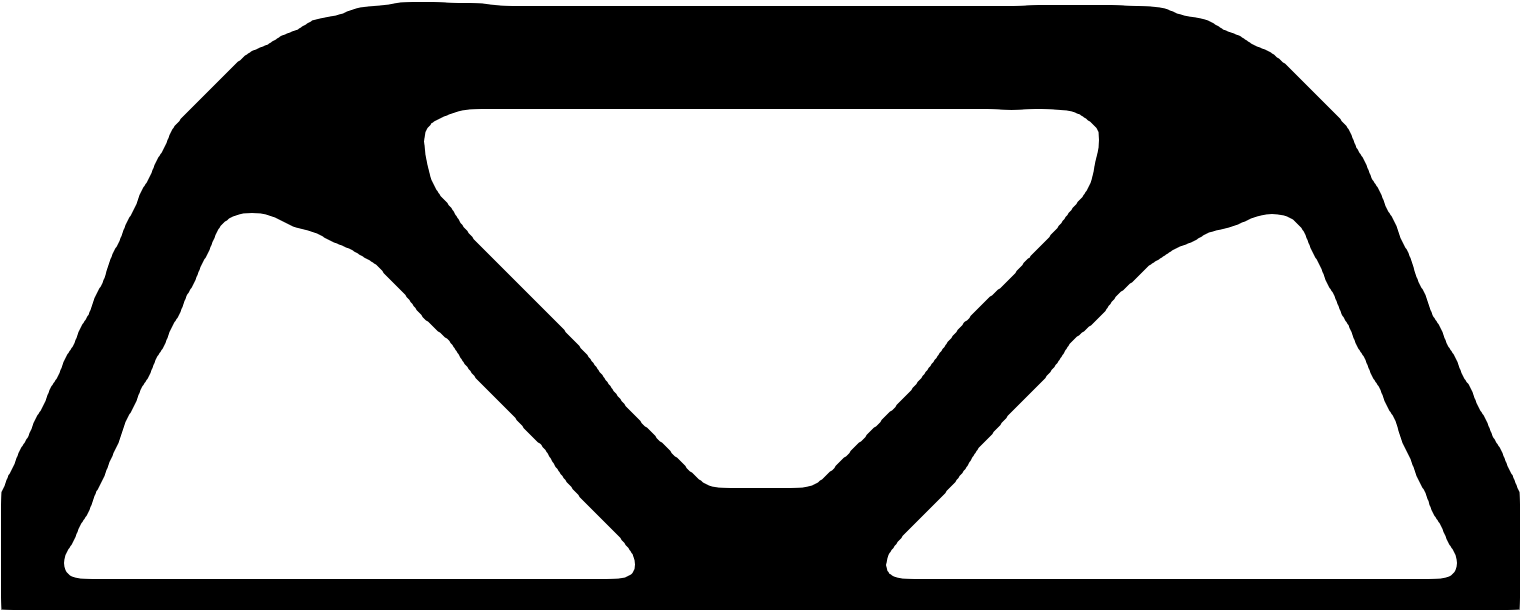}
    \caption{iter${}={}$45: $N_{1}=3$}
    \label{fig:michell-evolution-f}
  \end{subfigure}
  \caption{Evolution of the topology for the Michell structure when $\bar{N}_{1}=3$. (a) compliance history; (b)--(f) binary structures at representative
    iterations.}
  \label{fig:michell-evolution}
\end{figure}

Fig. \ref{fig:michell-evolution} illustrates the optimization process for $\bar{N}_{1}=3$.
The topology control is activated at iteration~30, at which point the structure contains seven holes and therefore
exceeds the prescribed bound. The response is extremely fast: the number of holes drops from seven to three within a
single iteration, so that the prescribed topology is already attained at
iteration~31 (Figs.~\ref{fig:michell-evolution-b}--\ref{fig:michell-evolution-c}). 
This abrupt change is accompanied by a sharp rise in compliance, followed by the
appearance of extra connected components at iteration~33, which are eliminated immediately (Figs.~\ref{fig:michell-evolution-d}--\ref{fig:michell-evolution-e}).
As the load path is re-established, the compliance decreases again and stabilizes,
while the topology remains stable from iteration~45 onwards (Fig.~\ref{fig:michell-evolution-f}).\par

\subsubsection{Cantilever beam}

\begin{figure}[!htb]
  \centering
  \begin{subfigure}[t]{0.30\textwidth}
    \centering
    \includegraphics[width=\linewidth]{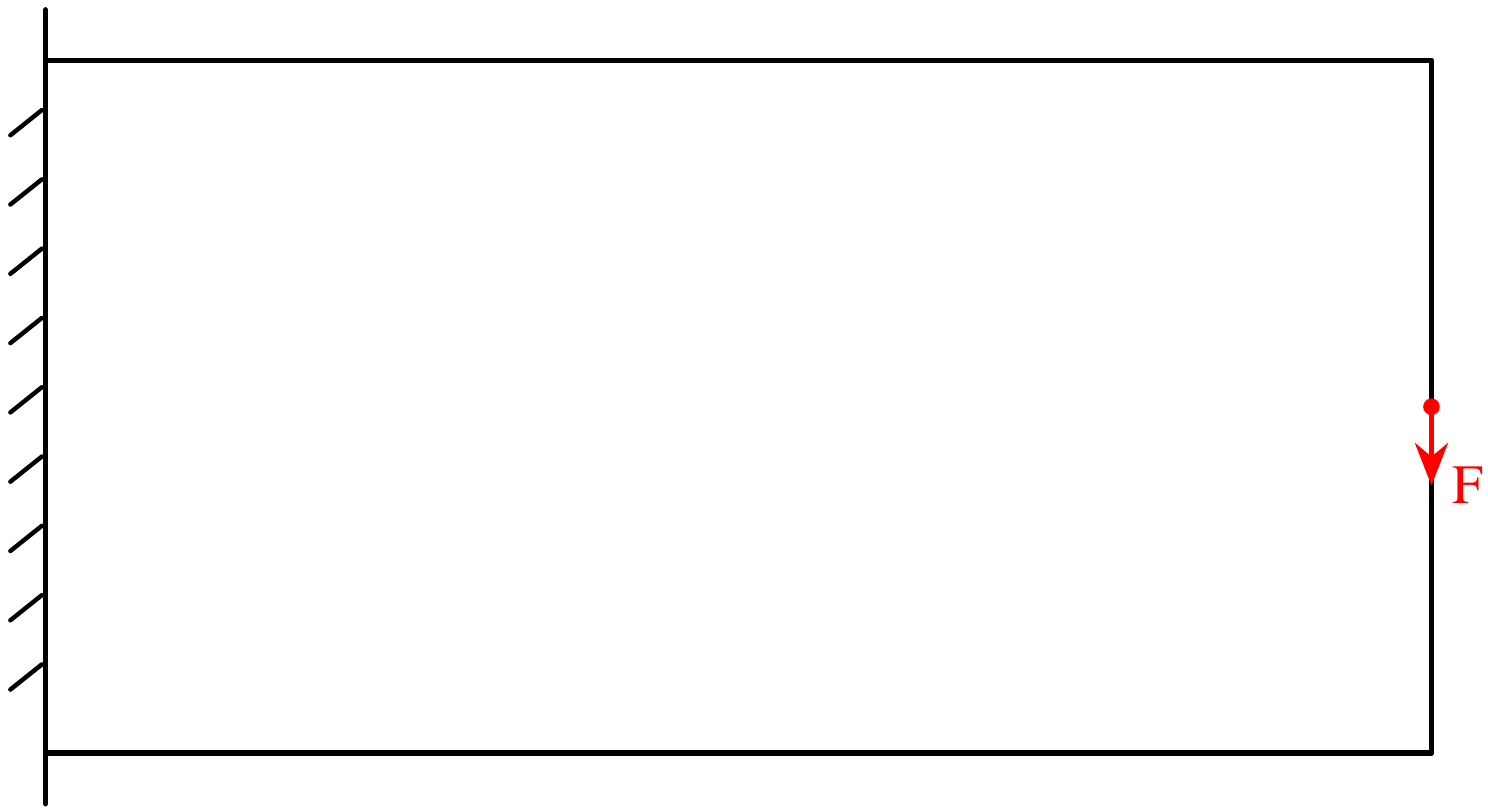}
    \caption{}
    \label{fig:cantilever-model}
  \end{subfigure}
  \hfill
  \begin{subfigure}[t]{0.36\textwidth}
    \centering
    \includegraphics[width=\linewidth]{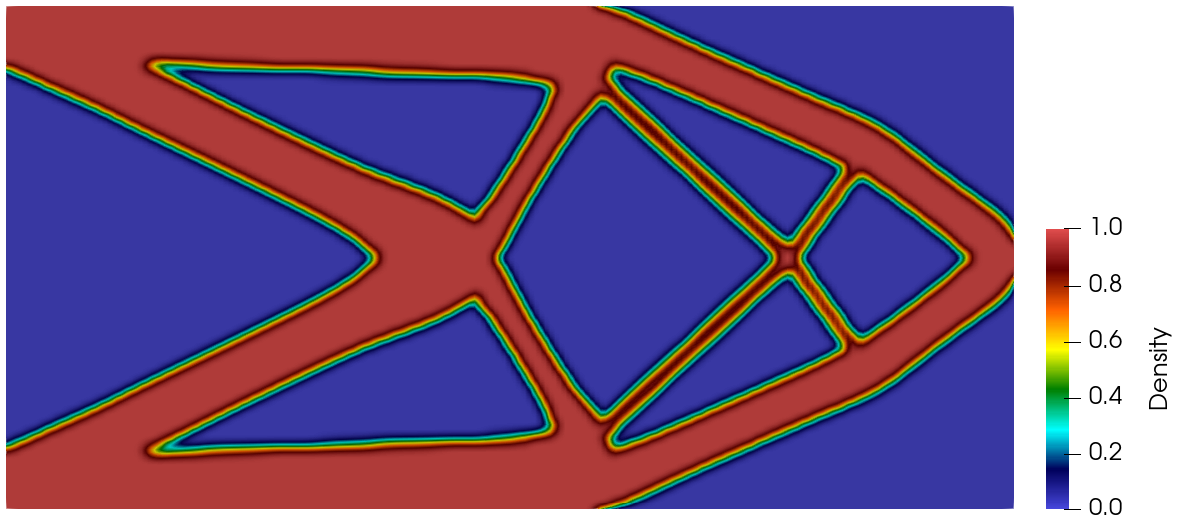}
    \caption{}
    \label{fig:cantilever-density}
  \end{subfigure}
  \hfill
  \begin{subfigure}[t]{0.30\textwidth}
    \centering
    \includegraphics[width=\linewidth]{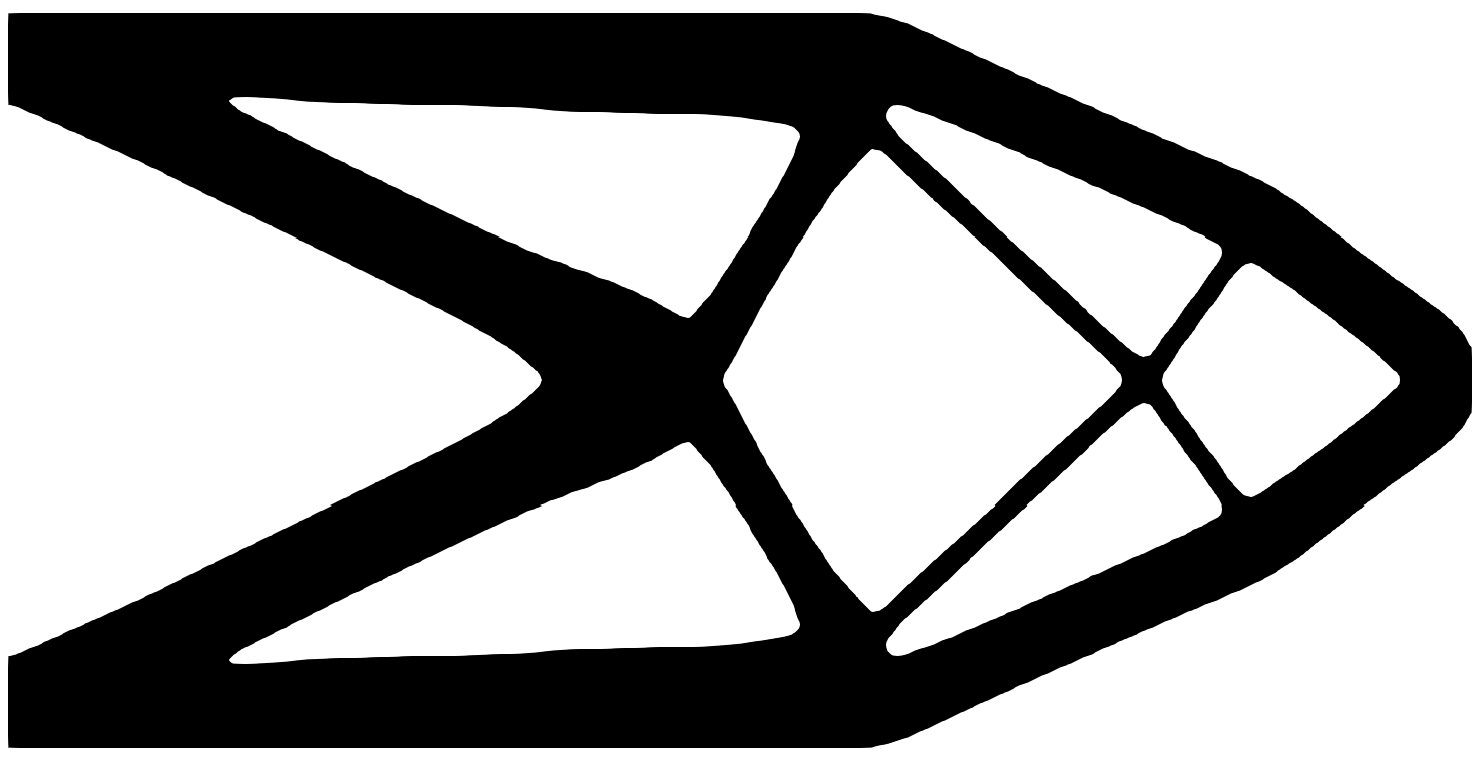}
    \caption{}
    \label{fig:cantilever-topology}
  \end{subfigure}
  \caption{(a) The cantilever beam model.
   (b) Optimized density field without topology control.
   (c) The corresponding binary structure with $N_1=6$.}
  \label{fig:cantilever-nocontrol}
  \end{figure}

In this subsection, the cantilever beam is taken as the final 2D benchmark. The
rectangular design domain $[0,10]\times[0,5]$ is represented by a biquadratic NURBS
surface with $162\times 82$ control points. The boundary conditions are illustrated
in Fig.~\ref{fig:cantilever-model}, where the left edge is fully clamped and a unit
downward point load is applied at the midpoint of the right free edge. The support
and the loading are symmetric about the horizontal midline of the beam, and the
maximum volume fraction is set to $V_{0}=0.4$. Without topology control, the optimized
layout is symmetric and contains six holes, as shown in
Figs.~\ref{fig:cantilever-density} and~\ref{fig:cantilever-topology}.\par

\begin{figure}[!htb]
  \centering
  \begin{subfigure}[b]{0.30\textwidth}
    \centering
    \includegraphics[width=\textwidth]{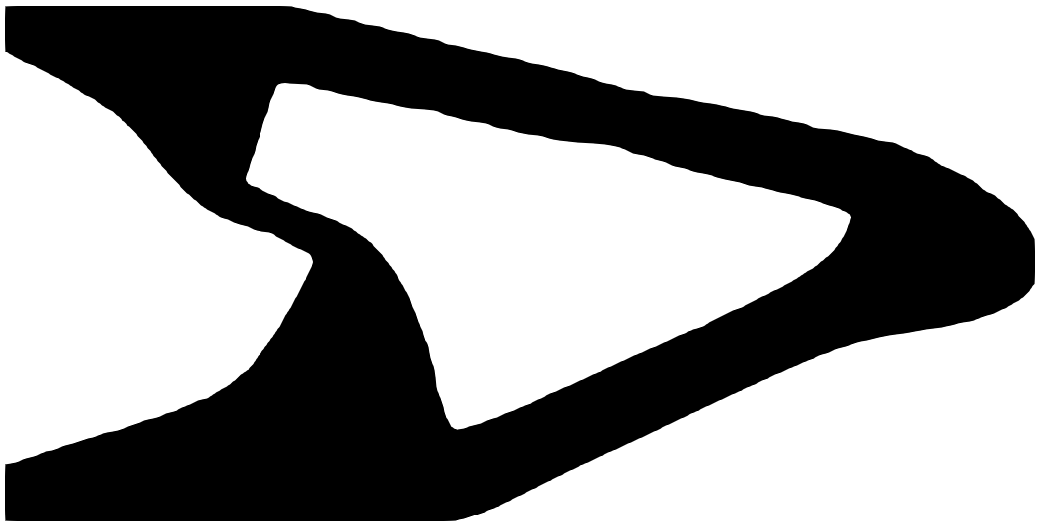}
    \caption{$\bar{N}_{1}=1$}
    \label{fig:cantilever-h1-1}
  \end{subfigure}
  \hfill
  \begin{subfigure}[b]{0.30\textwidth}
    \centering
    \includegraphics[width=\textwidth]{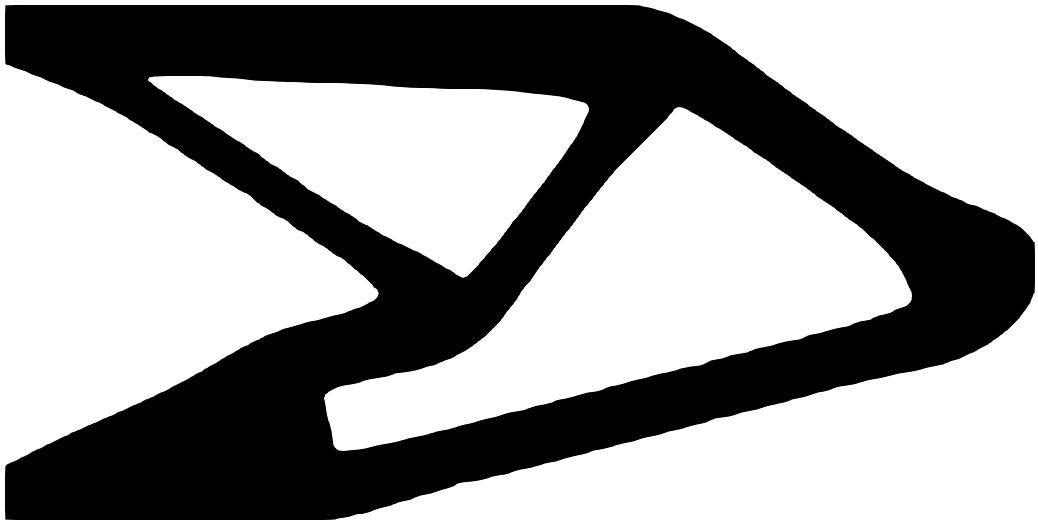}
    \caption{$\bar{N}_{1}=2$}
    \label{fig:cantilever-h1-2}
  \end{subfigure}
  \hfill
  \begin{subfigure}[b]{0.30\textwidth}
    \centering
    \includegraphics[width=\textwidth]{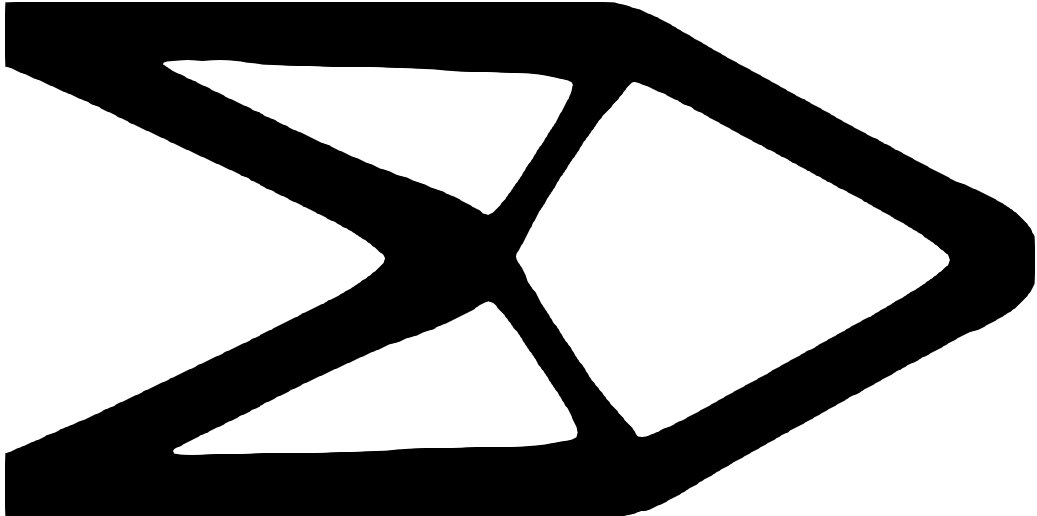}
    \caption{$\bar{N}_{1}=3$}
    \label{fig:cantilever-h1-3}
  \end{subfigure}\\[8pt]
  \begin{subfigure}[b]{0.30\textwidth}
    \centering
    \includegraphics[width=\textwidth]{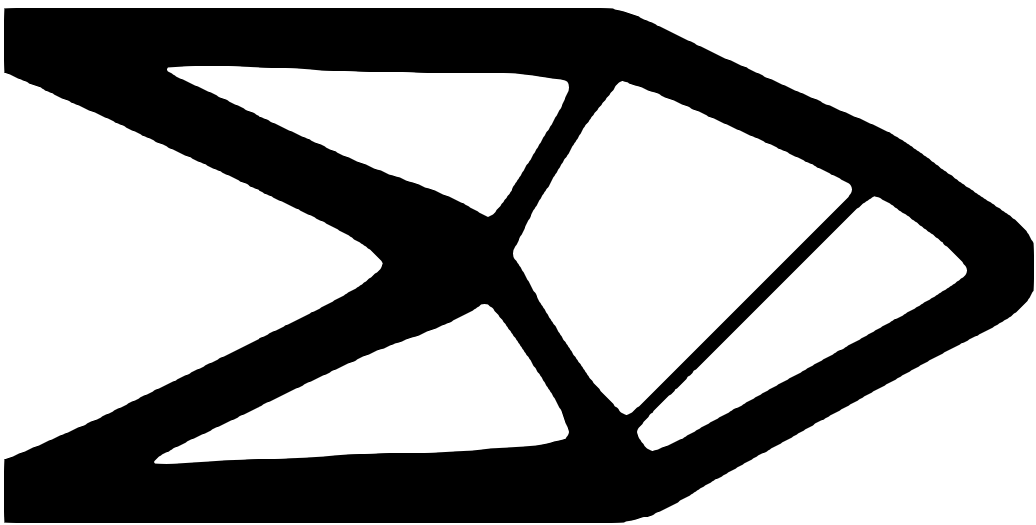}
    \caption{$\bar{N}_{1}=4$}
    \label{fig:cantilever-h1-4}
  \end{subfigure}
  \hspace{2em}
  \begin{subfigure}[b]{0.30\textwidth}
    \centering
    \includegraphics[width=\textwidth]{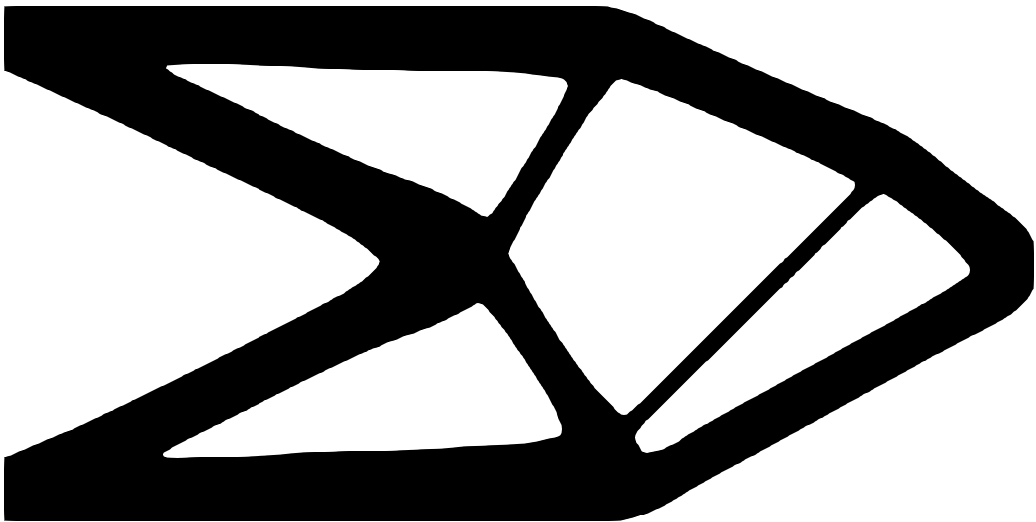}
    \caption{$\bar{N}_{1}=5$}
    \label{fig:cantilever-h1-5}
  \end{subfigure}
  \caption{Optimized cantilever beam structures for $\bar{N}_{1}=1,2,3,4,5$, respectively.}
  \label{fig:cantilever-control}
\end{figure}

The topology-controlled results for $\bar{N}_{1}=1,2,3,4,5$ are listed in Table~\ref{tab:cantilever} and Figure~\ref{fig:cantilever-control}.
A practical observation concerns the choice of the topological weights: while $\mu_{0}$ is kept fixed, a tighter bound
on the number of holes requires a larger weight $\mu_{1}$ assigned to the
1-dimensional objective. This is because when the prescribed topology deviates
more significantly from the unconstrained optimum, a stronger topological
driving force is needed to overcome the compliance-minimizing tendency of the
optimizer. The results thus confirm that the proposed persistence-based topological objective
provides an adjustable mechanism for balancing structural performance and topological complexity.\par

\begin{table}[!htb]
  \centering
  \caption{Comparison of the compliance values and of the resulting numbers of holes ($N_{1}$) for the 2D
cantilever beam under different prescribed bounds $\bar{N}_{1}$.}
  \label{tab:cantilever}
  \small
  \setlength{\tabcolsep}{4pt}
  \begin{tabular}{@{}lcccccc@{}}
    \toprule
    Quantity & $\bar{N}_{1}=1$ & $\bar{N}_{1}=2$ & $\bar{N}_{1}=3$ & $\bar{N}_{1}=4$ & $\bar{N}_{1}=5$ & Baseline\\
    \midrule
    Weights $(\mu_{0},\mu_{1})$ & $(100,200)$ & $(100,200)$ & $(100,50)$ & $(100,10)$ & $(100,5)$ & $(0, 0)$\\
    Compliance & $88.62$ & $79.93$ & $75.24$ & $75.23$ & $75.21$ & $75.08$\\
    $N_1$ & 1 & 2 & 3 & 4 & 4 & 6\\
    \bottomrule
  \end{tabular}
\end{table}

\begin{figure}[!htb]
  \centering
  \begin{subfigure}[b]{0.33\textwidth}
    \centering
    \includegraphics[width=\linewidth]{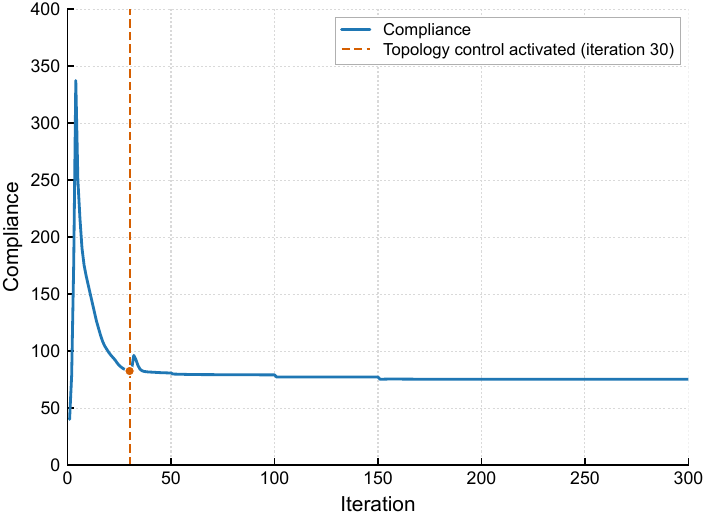}
    \caption{Compliance vs. iteration}
    \label{fig:cantilever-evolution-a}
  \end{subfigure}\hfill
  \begin{subfigure}[b]{0.30\textwidth}
    \centering
    \includegraphics[width=\linewidth]{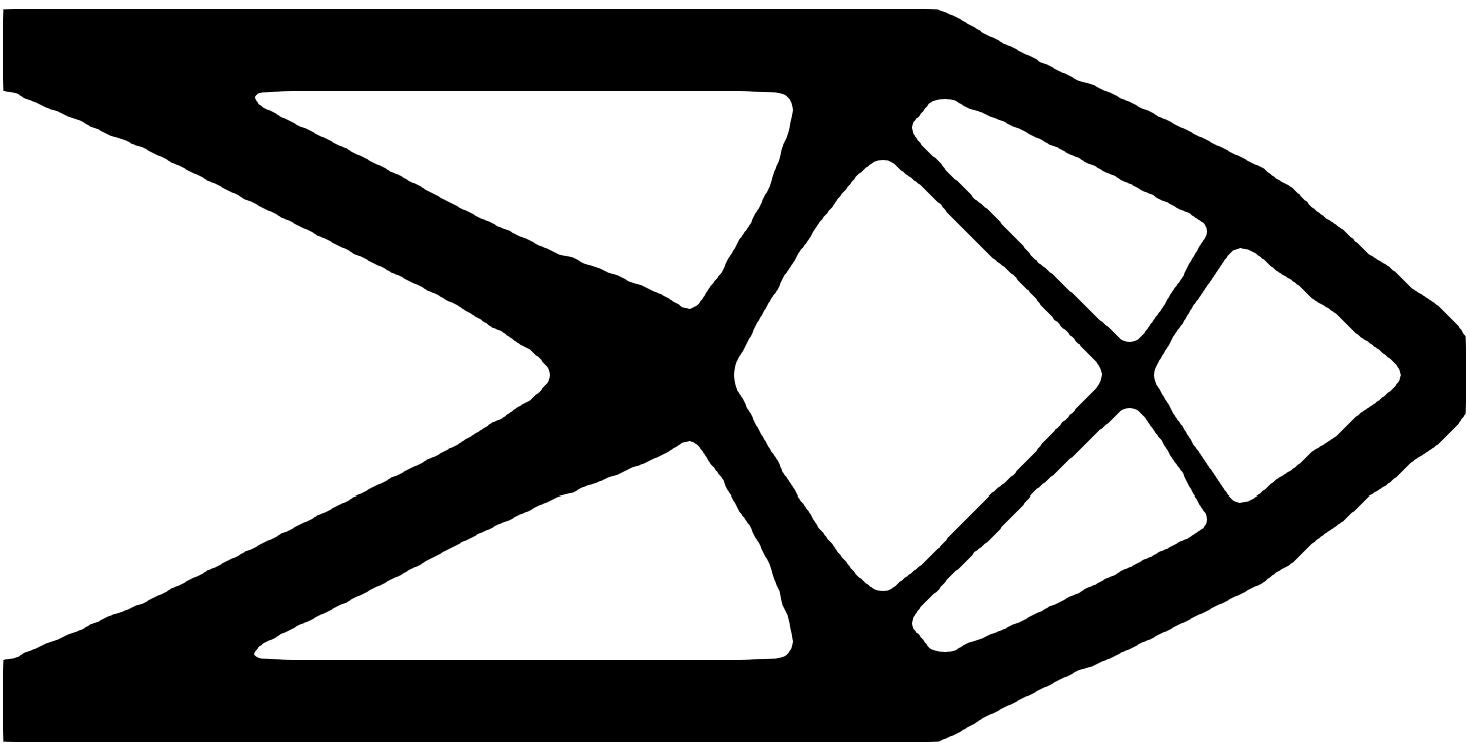}
    \caption{iter${}={}$30: $N_{1}=6$}
    \label{fig:cantilever-evolution-b}
  \end{subfigure}\hfill
  \begin{subfigure}[b]{0.30\textwidth}
    \centering
    \includegraphics[width=\linewidth]{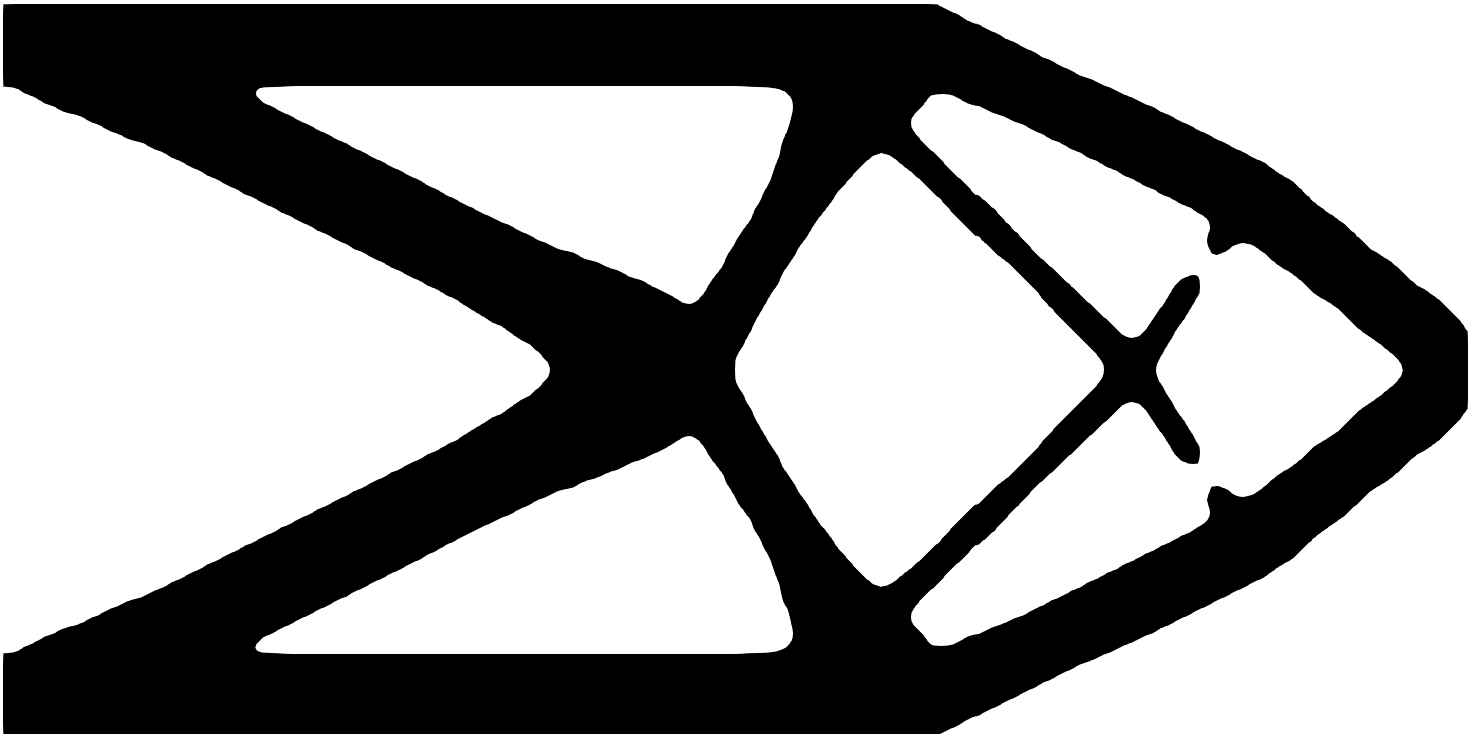}
    \caption{iter${}={}$31: $N_{1}=4$}
    \label{fig:cantilever-evolution-c}
  \end{subfigure}\\[8pt]
  \begin{subfigure}[b]{0.30\textwidth}
    \centering
    \includegraphics[width=\linewidth]{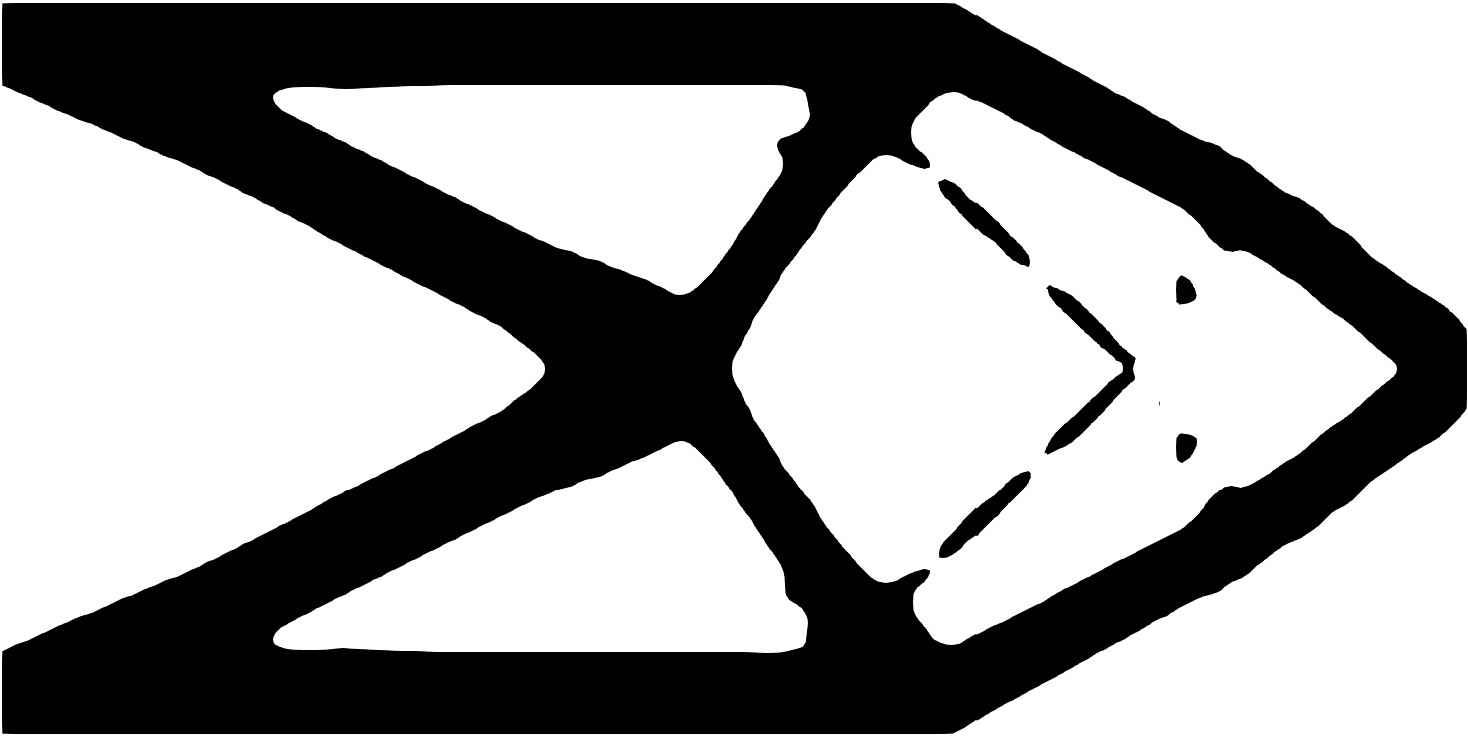}
    \caption{iter${}={}$33: $N_{1}=3$}
    \label{fig:cantilever-evolution-d}
  \end{subfigure}\hfill
  \begin{subfigure}[b]{0.30\textwidth}
    \centering
    \includegraphics[width=\linewidth]{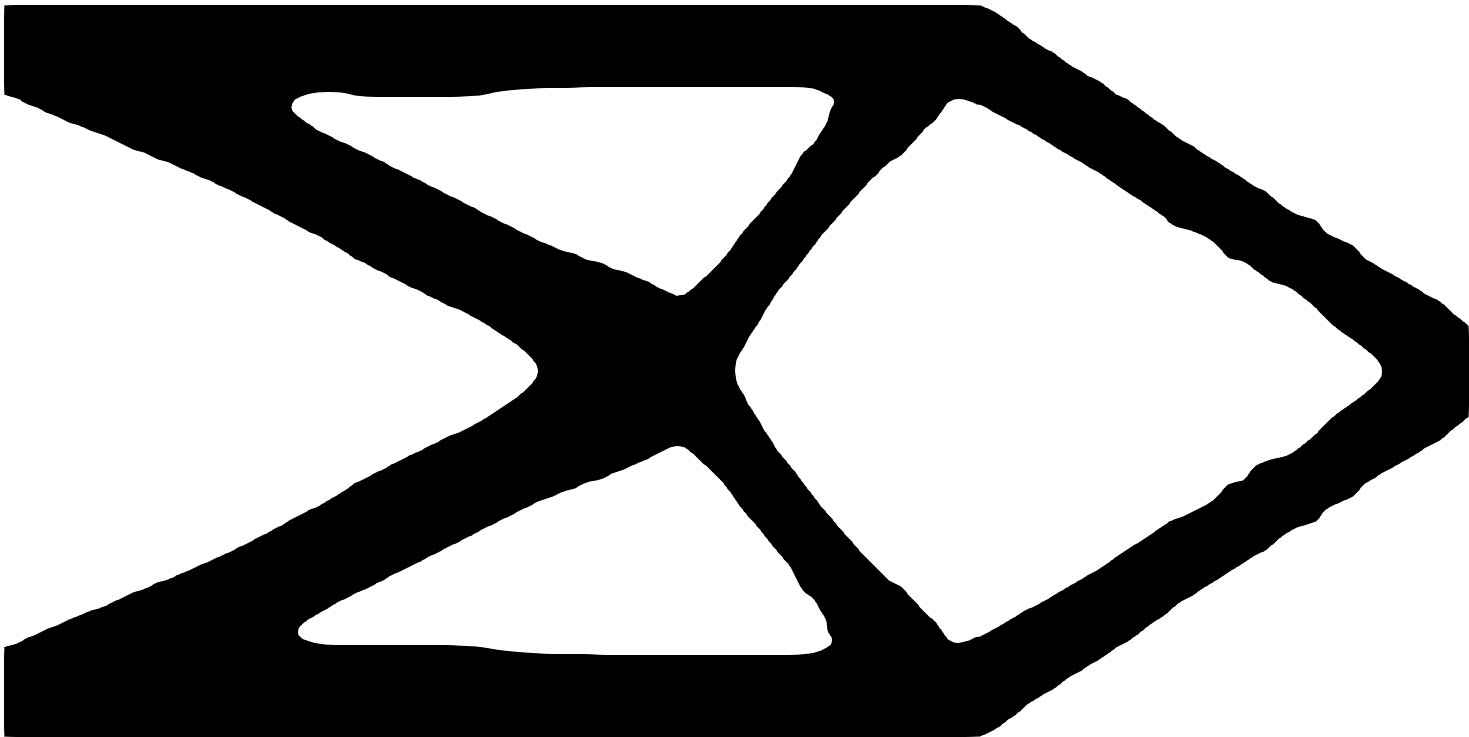}
    \caption{iter${}={}$35: $N_{1}=3$}
    \label{fig:cantilever-evolution-e}
  \end{subfigure}\hfill
  \begin{subfigure}[b]{0.30\textwidth}
    \centering
    \includegraphics[width=\linewidth]{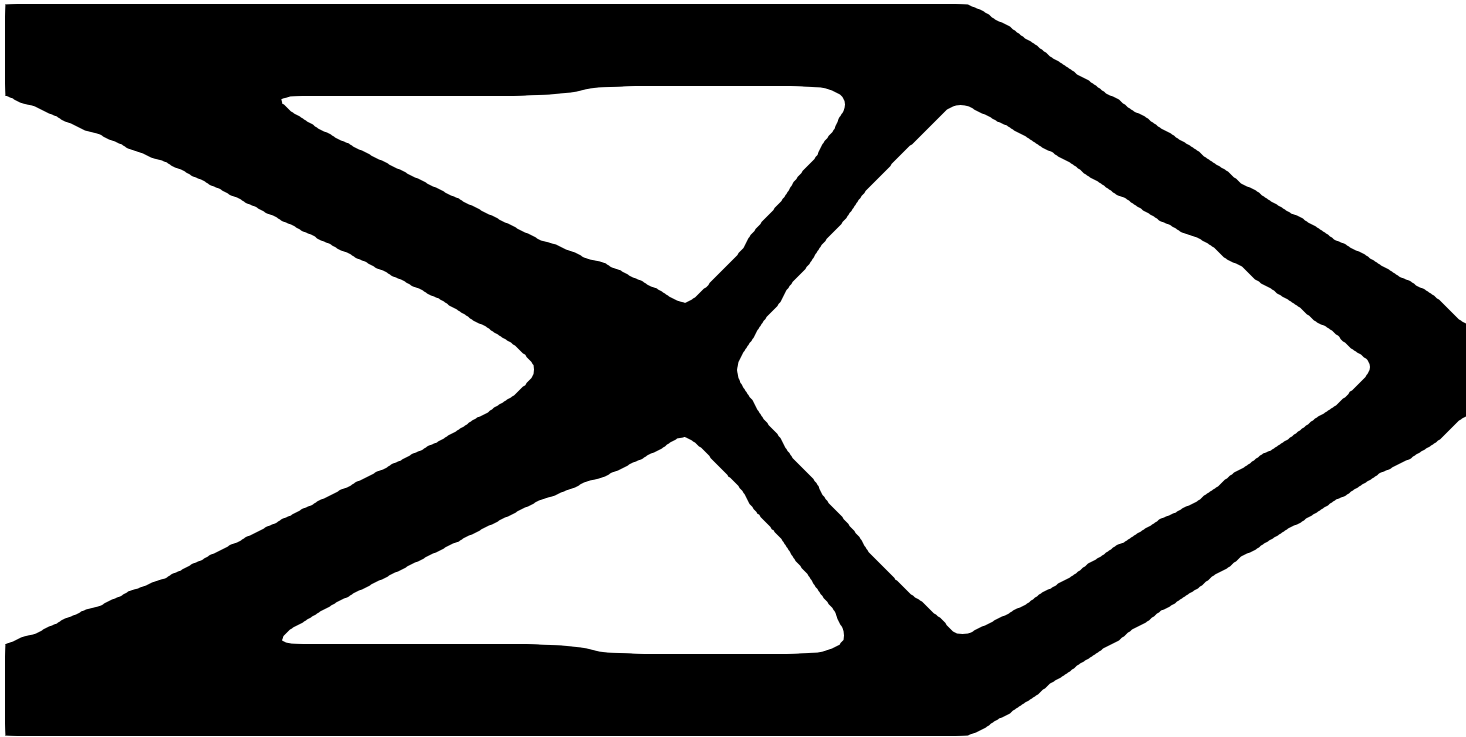}
    \caption{iter${}={}$45: $N_{1}=3$}
    \label{fig:cantilever-evolution-f}
  \end{subfigure}
  \caption{Evolution of the topology for the cantilever beam when $\bar{N}_{1}=3$. (a) compliance history; (b)--(f) binary structures at representative
    iterations.}
  \label{fig:cantilever-evolution}
\end{figure}

The optimization process for $\bar{N}_{1}=3$ is illustrated in
Fig.~\ref{fig:cantilever-evolution}. Once the topology control is activated at
iteration~30, the number of holes drops from six to four within a single iteration
and reaches the prescribed value of three within only
three iterations (Figs.~\ref{fig:cantilever-evolution-b}--
\ref{fig:cantilever-evolution-d}). The extra connected components that appear
during this rapid removal are in turn cleared within two further iterations, as shown in
Figs.~\ref{fig:cantilever-evolution-d} and \ref{fig:cantilever-evolution-e}. This
fast correction demonstrates the strong capability of the proposed topological
objective to enforce the prescribed topology.\par

\subsection{3D numerical examples}
\subsubsection{Cantilever beam}

\begin{figure}[!htb]
  \centering
  \begin{subfigure}[t]{0.40\textwidth}
    \centering
    \includegraphics[width=\linewidth]{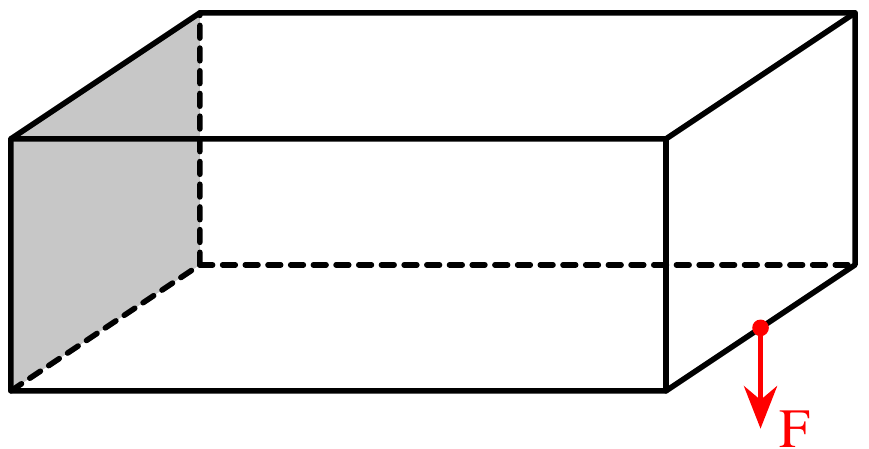}
    \caption{}
    \label{fig:cantilever3d-model}
  \end{subfigure}
  \hfill
  \begin{subfigure}[t]{0.42\textwidth}
    \centering
    \includegraphics[width=\linewidth]{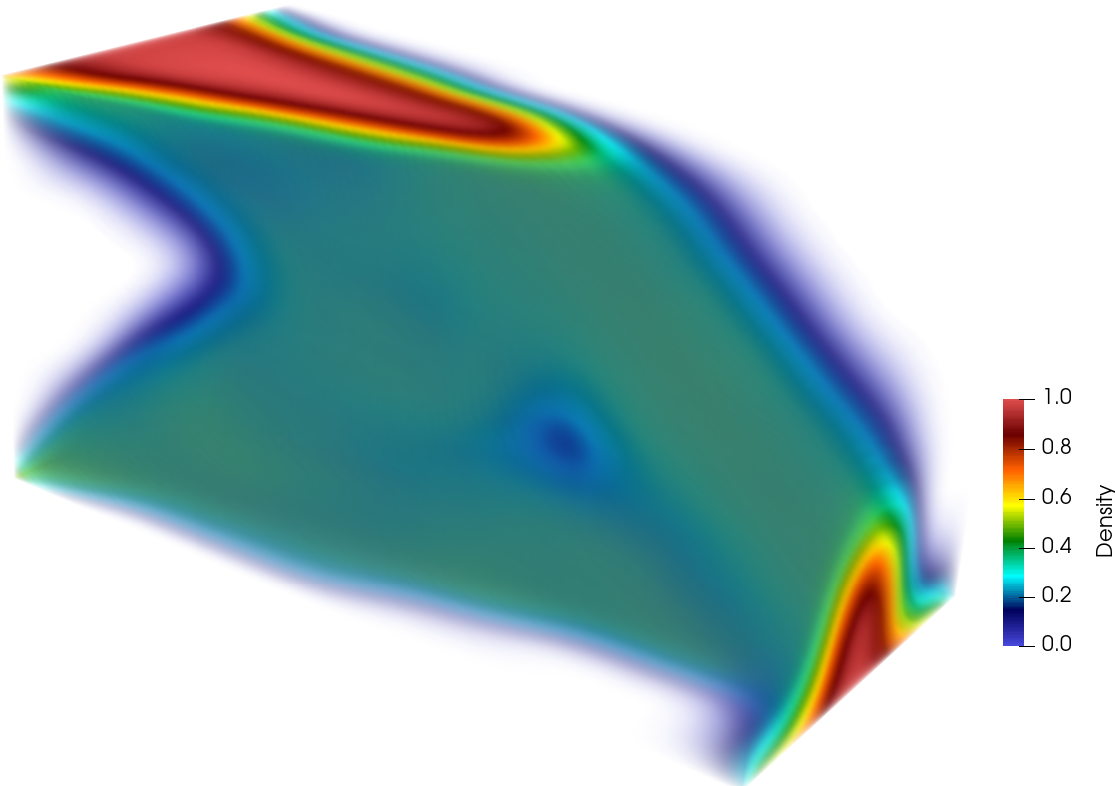}
    \caption{}
    \label{fig:cantilever3d-density}
  \end{subfigure}\\[8pt]
  \begin{subfigure}[t]{0.40\textwidth}
    \centering
    \includegraphics[width=\linewidth]{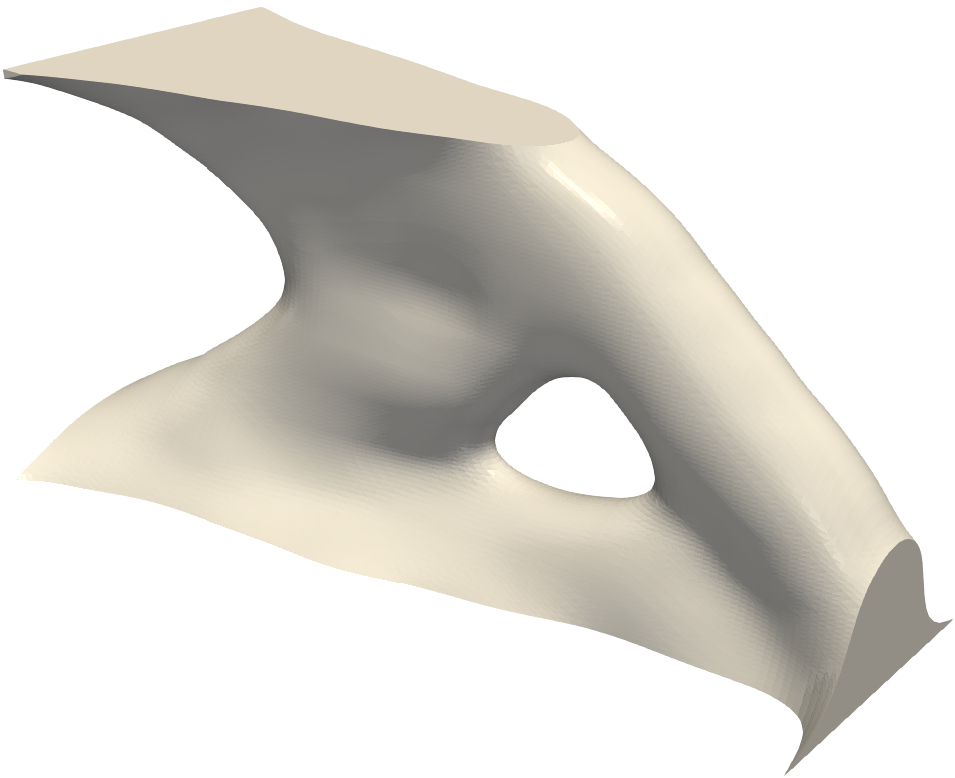}
    \caption{}
    \label{fig:cantilever3d-topology}
  \end{subfigure}
  \hfill
  \begin{subfigure}[t]{0.40\textwidth}
    \centering
    \includegraphics[width=\linewidth]{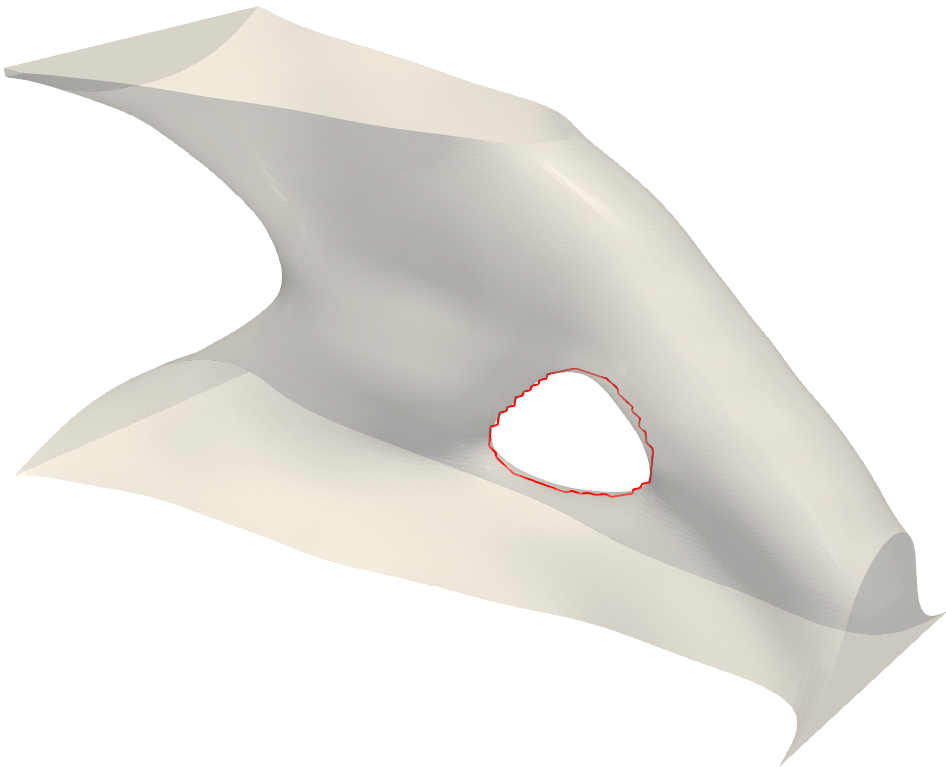}
    \caption{}
    \label{fig:cantilever3d-topology-transparent}
  \end{subfigure}
  \caption{(a) The 3D cantilever beam model.
   (b) Optimized density field without topology control.
   (c) The corresponding binary structure with $N_1=1, N_2=0$.
   (d) Semi-transparent view of (c) for visualization, with the tunnel marked by a red loop.}
\end{figure}

The first 3D example considers the cantilever beam benchmark, whose design domain is the cuboid
$[0,4]\times[0,2]\times[0,2]$, represented by a triquadratic NURBS volume with
$21\times 11\times 11$ control points. As illustrated in
Fig.~\ref{fig:cantilever3d-model}, the left face is fully clamped and a unit
downward point load is applied at the midpoint of the bottom edge on the free
end. The maximum volume fraction is set to $V_{0}=0.2$. 
For comparison, the minimum compliance problem without topology control is solved first,
yielding the baseline structure shown in
Figs.~\ref{fig:cantilever3d-density}--\ref{fig:cantilever3d-topology-transparent},
which contains one tunnel and no cavities.\par

\begin{figure}[!htb]
  \centering
  \begin{subfigure}[b]{0.40\textwidth}
    \centering
    \includegraphics[width=\textwidth]{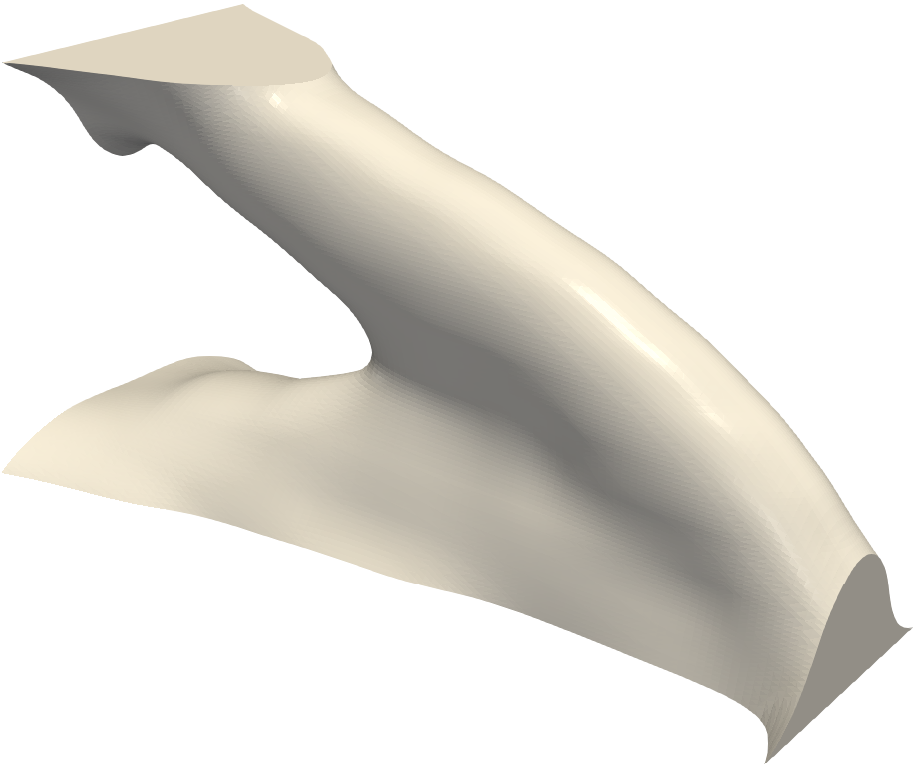}
    \caption{Opaque view}                                                     
    \label{fig:cantilever3d-h1-0}
  \end{subfigure}\hfill
  \begin{subfigure}[b]{0.40\textwidth}
    \centering
    \includegraphics[width=\textwidth]{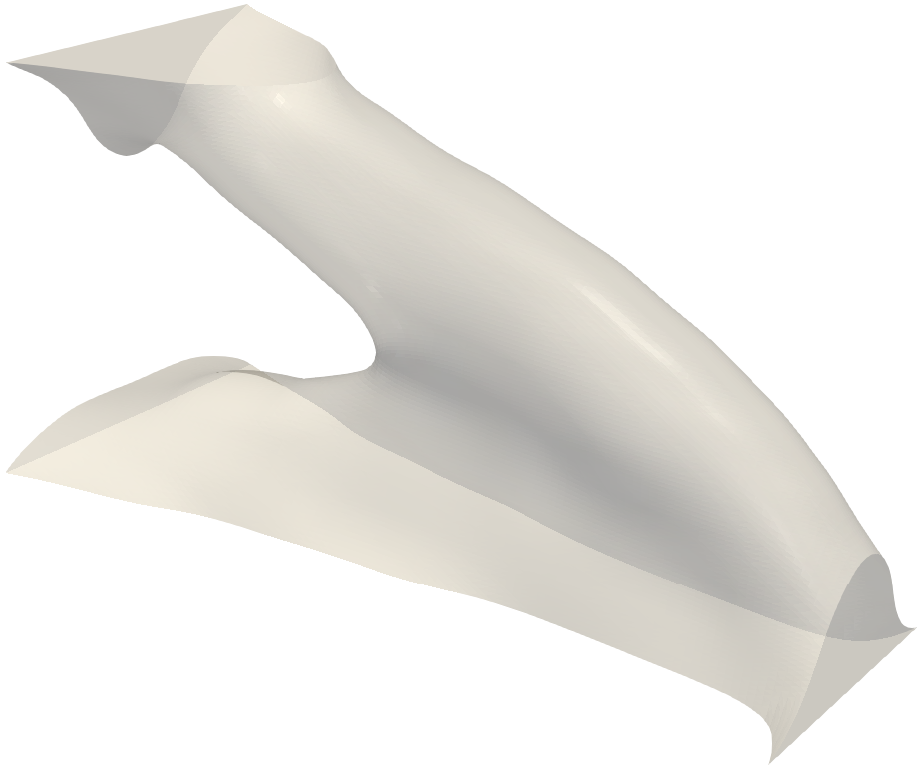}
    \caption{Semi-transparent view}                                                     
    \label{fig:cantilever3d-h1-0-transparent}
  \end{subfigure}
  \caption{Optimized 3D cantilever beam structure for $\bar{N}_{1}=0, \bar{N}_2=0$.}
  \label{fig:cantilever3d-control}
\end{figure}

The same problem is then solved with topology control under the prescribed bounds
$\bar{N}_{1}=0$ and $\bar{N}_{2}=0$.
The optimized structure is shown in Fig.~\ref{fig:cantilever3d-control}, where
filling removes the single tunnel and yields a
smooth final structure. However, as reported in Table~\ref{tab:cantilever3d}, the compliance of the resulting design is higher than
that of the baseline, suggesting that the topological constraint is enforced at
the expense of mechanical performance.\par

\begin{table}[!htb]
  \centering
  \caption{Comparison of the compliance values and of the resulting numbers of tunnels
($N_{1}$) and cavities ($N_{2}$) for the 3D cantilever beam under different
prescribed bounds $\bar{N}_{1}$, with $\bar{N}_{2}=0$.}
  \label{tab:cantilever3d}
  \small
  \begin{tabular}{@{}lcc@{}}
    \toprule
    Quantity & $\bar{N}_{1}=0$ & Baseline\\
    \midrule
    Weights $(\mu_{0},\mu_{1},\mu_{2})$ & $(100,100,100)$ & $(0, 0, 0)$\\
    Compliance & 120.31 & $111.38$\\
    $(N_{1},N_{2})$ & $(0, 0)$ & $(1, 0)$\\
    \bottomrule
  \end{tabular}
\end{table}

\subsubsection{MBB beam}

\begin{figure}[!htb]
  \centering
  \begin{subfigure}[t]{0.40\textwidth}
    \centering
    \includegraphics[width=\linewidth]{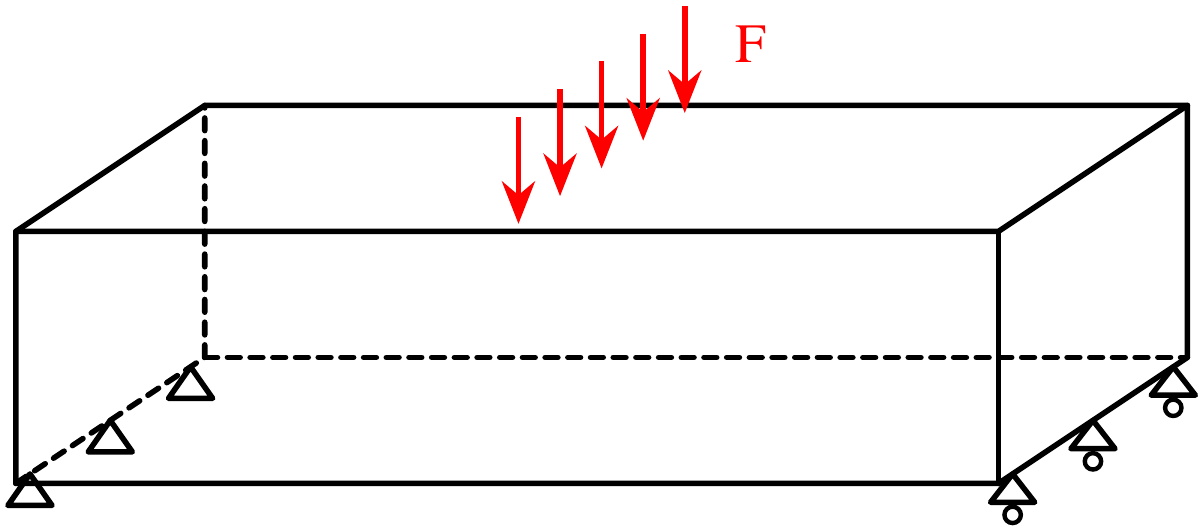}
    \caption{}
    \label{fig:mbb3d-model}
  \end{subfigure}
  \hfill
  \begin{subfigure}[t]{0.42\textwidth}
    \centering
    \includegraphics[width=\linewidth]{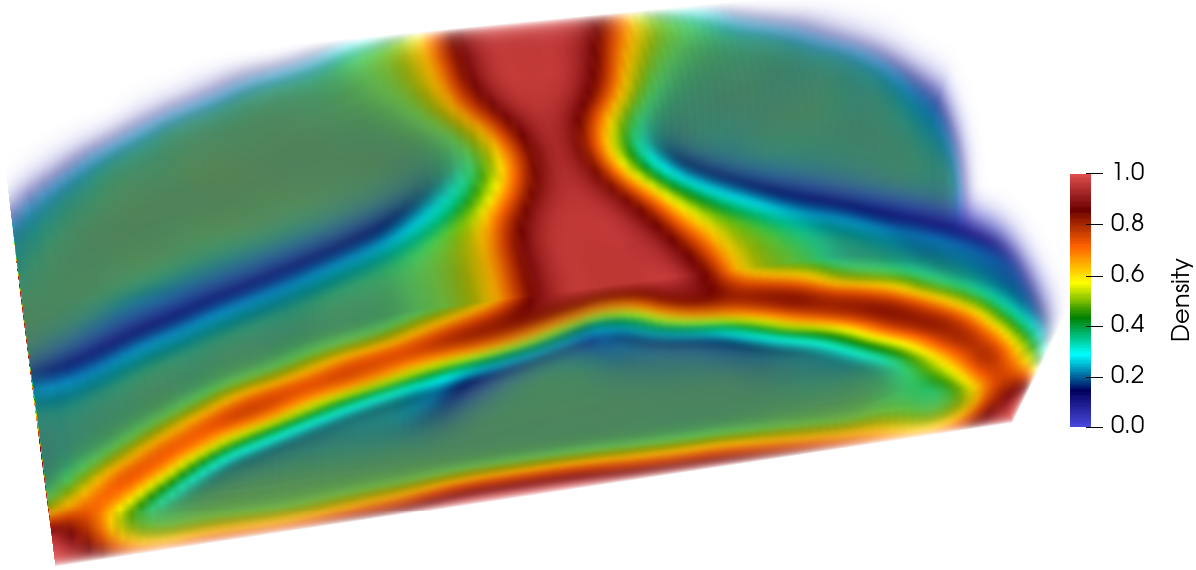}
    \caption{}
    \label{fig:mbb3d-density}
  \end{subfigure}\\[8pt]
  \begin{subfigure}[t]{0.40\textwidth}
    \centering
    \includegraphics[width=\linewidth]{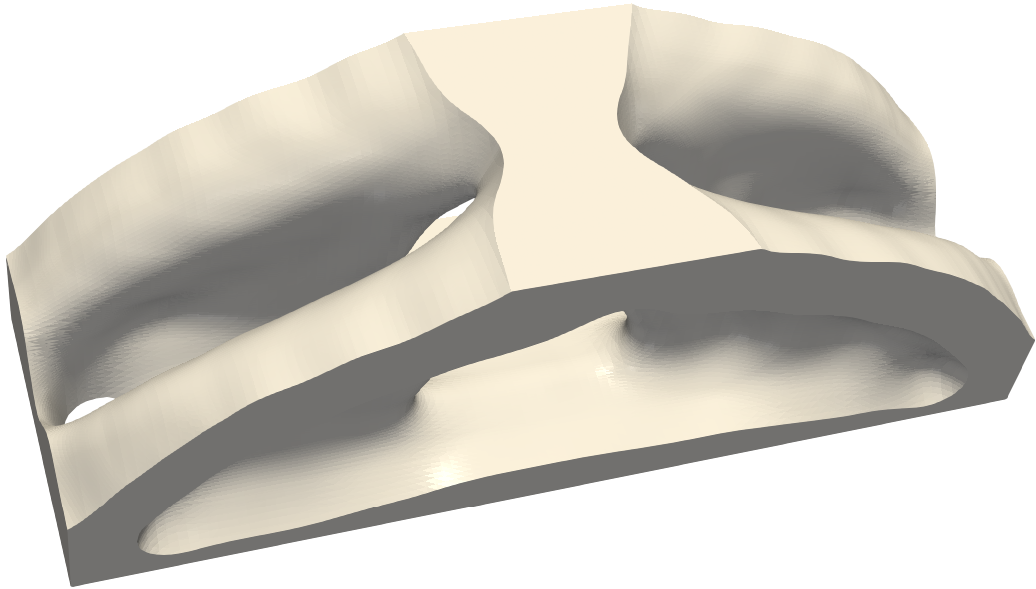}
    \caption{}
    \label{fig:mbb3d-topology}
  \end{subfigure}
    \hfill
    \begin{subfigure}[t]{0.40\textwidth}
    \centering
    \includegraphics[width=\linewidth]{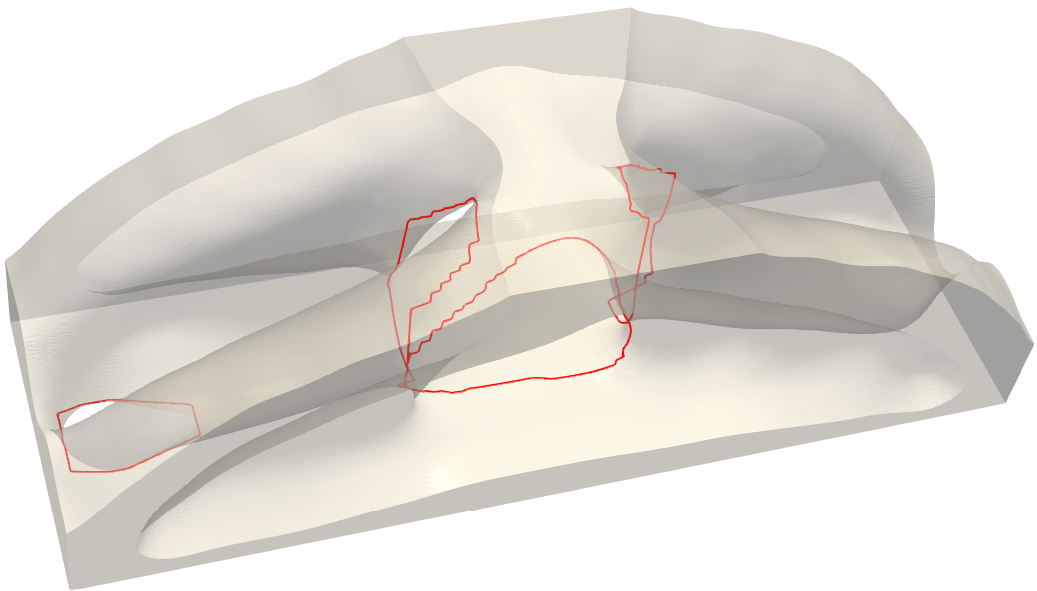}
    \caption{}
    \label{fig:mbb3d-topology-transparent}
  \end{subfigure}
  \caption{(a) The 3D MBB beam model.
   (b) Optimized density field without topology control.
   (c) The corresponding binary structure with $N_1=4, N_2=0$.
   (d) Semi-transparent view of (c), with the tunnels marked by red loops.}
\end{figure}

\begin{figure}[!htb]
  \centering
  \begin{subfigure}[b]{0.40\textwidth}
    \centering
    \includegraphics[width=\textwidth]{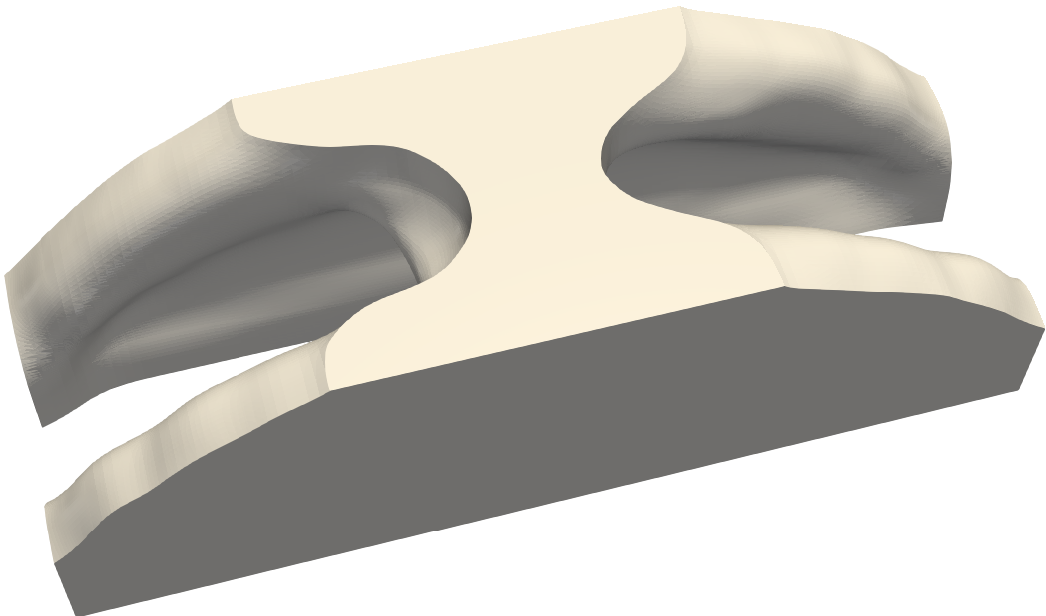}
    \caption{Opaque view}
    \label{fig:mbb3d-h1-0-topology}
  \end{subfigure}
  \hfill
  \begin{subfigure}[b]{0.40\textwidth}
    \centering
    \includegraphics[width=\textwidth]{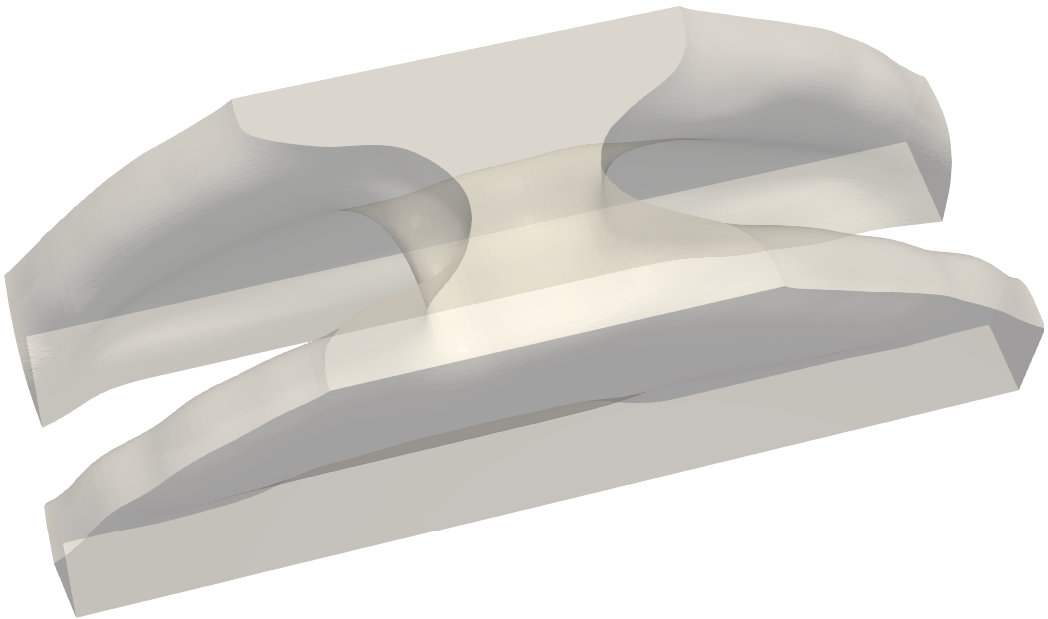}
    \caption{Semi-transparent view}
    \label{fig:mbb3d-h1-0-transparent}
  \end{subfigure}
  \caption{Optimized 3D MBB beam structure for $\bar{N}_{1}=0, \bar{N}_{2}=0$.}
  \label{fig:mbb3d-h1-0}
\end{figure}

\begin{figure}[!htb]
  \centering
  \begin{subfigure}[b]{0.40\textwidth}
    \centering
    \includegraphics[width=\textwidth]{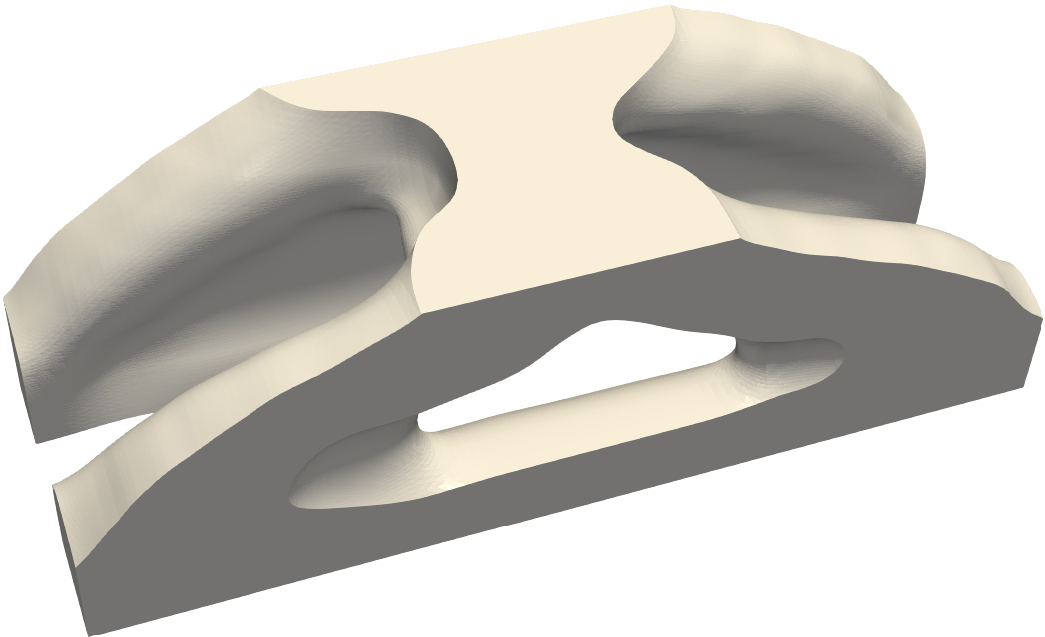}
    \caption{Opaque view}
    \label{fig:mbb3d-h1-1-topology}
  \end{subfigure}
  \hfill
  \begin{subfigure}[b]{0.40\textwidth}
    \centering
    \includegraphics[width=\textwidth]{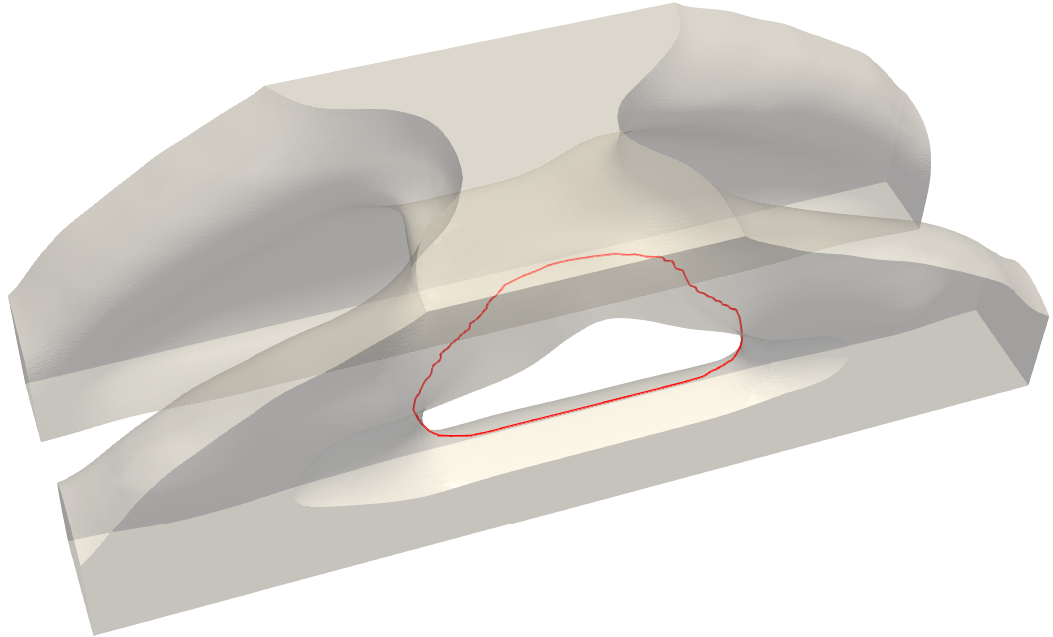}
    \caption{Semi-transparent view}
    \label{fig:mbb3d-h1-transparent}
  \end{subfigure}
  \caption{Optimized 3D MBB beam structure for $\bar{N}_{1}=1, \bar{N}_{2}=0$.}
  \label{fig:mbb3d-h1-1}
\end{figure}

\begin{figure}[!htb]
  \centering
  \begin{subfigure}[b]{0.40\textwidth}
    \centering
    \includegraphics[width=\textwidth]{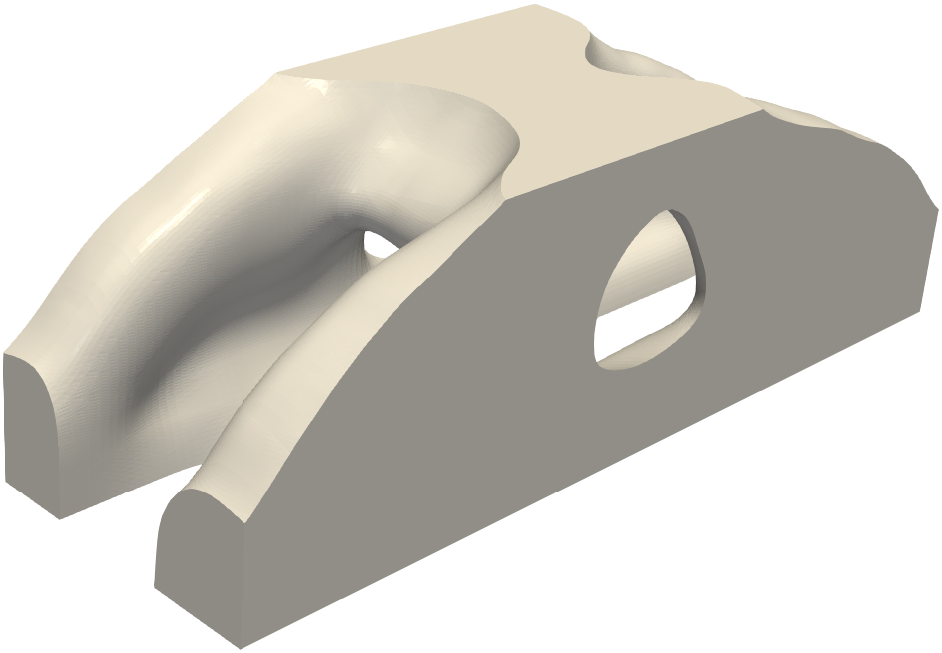}
    \caption{Opaque view}
    \label{fig:mbb3d-h1-2-topology}
  \end{subfigure}
  \hfill
  \begin{subfigure}[b]{0.40\textwidth}
    \centering
    \includegraphics[width=\textwidth]{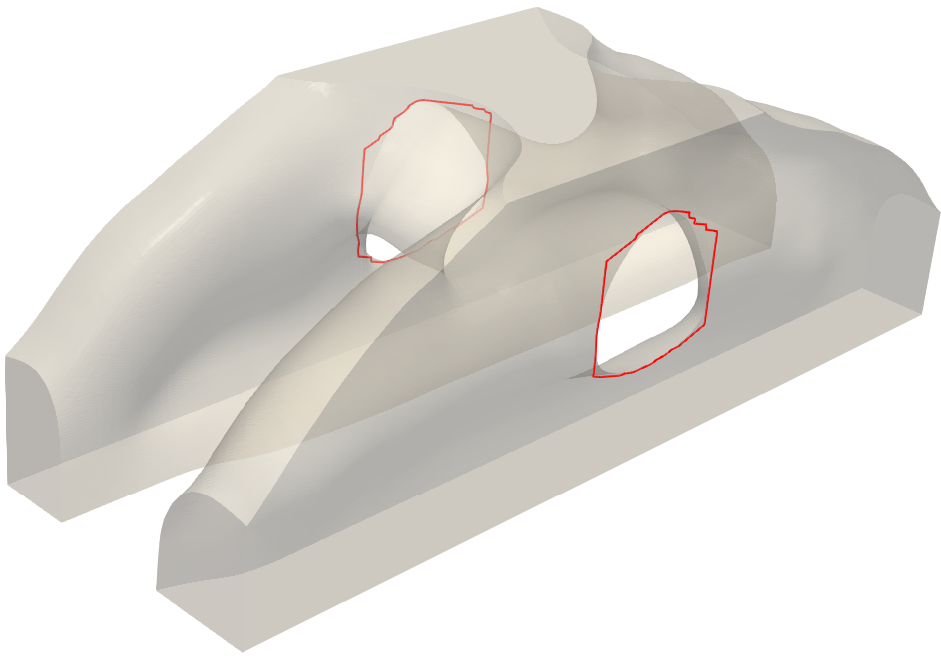}
    \caption{Semi-transparent view}
    \label{fig:mbb3d-h1-2-transparent}
  \end{subfigure}
  \caption{Optimized 3D MBB beam structure for $\bar{N}_{1}=2, \bar{N}_{2}=0$.}
  \label{fig:mbb3d-h1-2}
\end{figure}

\begin{figure}[!htb]
  \centering
  \begin{subfigure}[b]{0.40\textwidth}
    \centering
    \includegraphics[width=\textwidth]{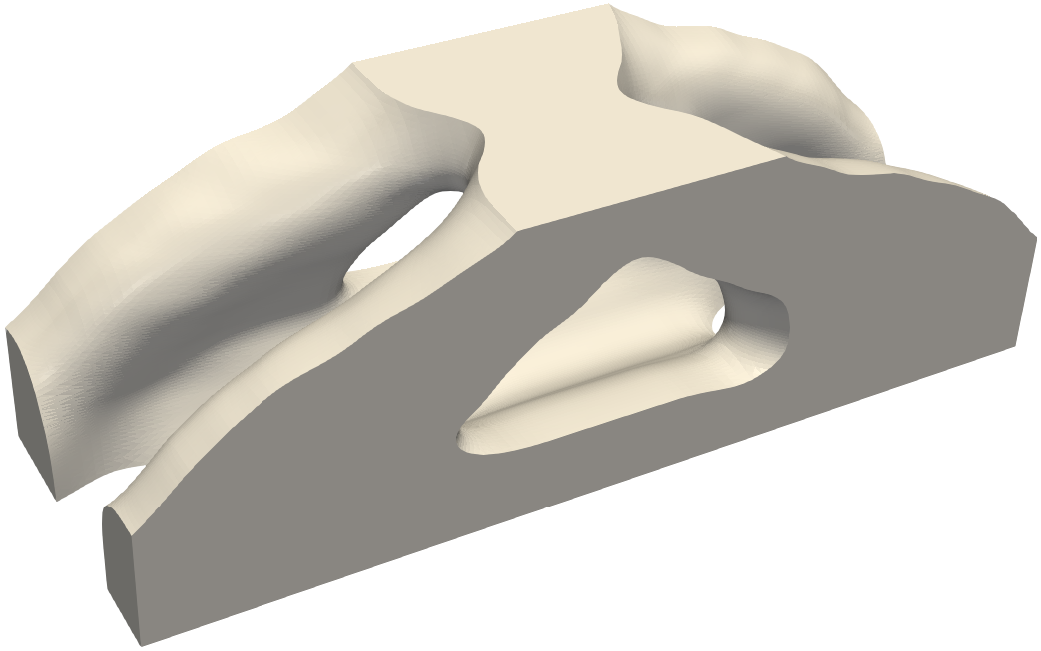}
    \caption{Opaque view}
    \label{fig:mbb3d-h1-3-topology}
  \end{subfigure}
  \hfill
  \begin{subfigure}[b]{0.40\textwidth}
    \centering
    \includegraphics[width=\textwidth]{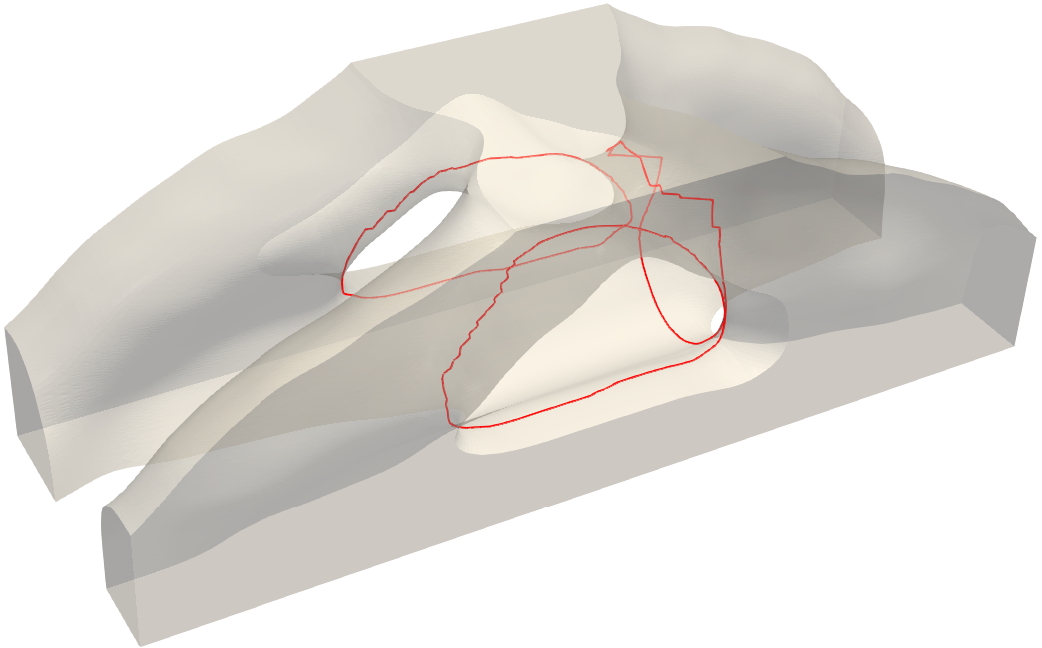}
    \caption{Semi-transparent view}
    \label{fig:mbb3d-h1-3-transparent}
  \end{subfigure}
  \caption{Optimized 3D MBB beam structure for $\bar{N}_{1}=3, \bar{N}_{2}=0$.}
  \label{fig:mbb3d-h1-3}
\end{figure}

In this subsection, we consider the 3D MBB beam benchmark occupying the cuboid
$[0,6]\times[0,2]\times[0,2]$, represented by a triquadratic NURBS volume with
$31\times 11\times 11$ control points. The boundary conditions are shown in
Fig.~\ref{fig:mbb3d-model}, with one bottom edge fully constrained and the
opposite one restrained only in the vertical direction, and a unit downward line
load is applied along the transverse centerline of the top face.
The maximum volume fraction is set to $V_{0}=0.4$.
Solving the minimum compliance problem
without topology control yields the baseline structure illustrated in
Figs.~\ref{fig:mbb3d-density}--\ref{fig:mbb3d-topology-transparent},
which contains four tunnels and no cavities.
Subsequently, the same problem is solved under topology control, with
$\bar{N}_{2}$ fixed to $0$ and $\bar{N}_{1}=0,1,2,3$. The optimized structures
are shown in Figs.~\ref{fig:mbb3d-h1-0}--\ref{fig:mbb3d-h1-3}, and the quantitative comparison is reported in Table~\ref{tab:mbb3d}.\par

\begin{table}[!htb]
  \centering
  \caption{Comparison of the compliance values and of the resulting numbers of tunnels
($N_{1}$) and cavities ($N_{2}$) for the 3D MBB beam under different
prescribed bounds $\bar{N}_{1}$, with $\bar{N}_{2} = 0$.}
  \label{tab:mbb3d}
  \small
  \setlength{\tabcolsep}{5pt}
  \begin{tabular}{@{}lccccc@{}}
    \toprule
    \multirow{2}{*}{Quantity} & \multicolumn{4}{c}{Prescribed $\bar{N}_{1}$} & \multirow{2}{*}{Baseline}\\
    \cmidrule(lr){2-5}
    & $0$ & $1$ & $2$ & $3$ & \\
    \midrule
    Weights $(\mu_{0},\mu_{1},\mu_{2})$ & \multicolumn{4}{c}{$(100,100,100)$} & $(0,0,0)$ \\
    Compliance & $62.88$ & $62.77$ & $62.24$ & $60.86$ & $57.14$\\
    $(N_{1},N_{2})$ & $(0, 0)$ & $(1, 0)$ & $(2, 0)$ & $(3, 0)$ & $(4, 0)$\\
    \bottomrule
  \end{tabular}
\end{table}

\subsubsection{Michell structure}

\begin{figure}[!htb]
  \centering
  \begin{subfigure}[t]{0.40\textwidth}
    \centering
    \includegraphics[width=\linewidth]{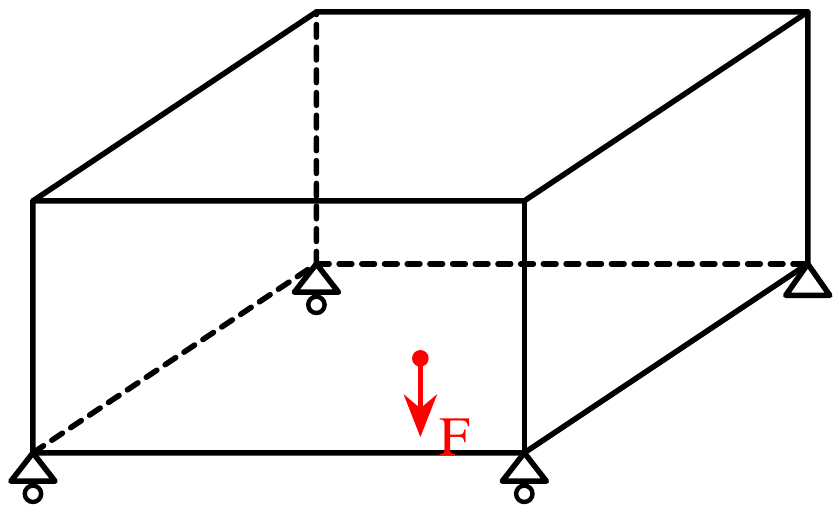}
    \caption{}
    \label{fig:michell3d-model}
  \end{subfigure}
  \hfill
  \begin{subfigure}[t]{0.42\textwidth}
    \centering
    \includegraphics[width=\linewidth]{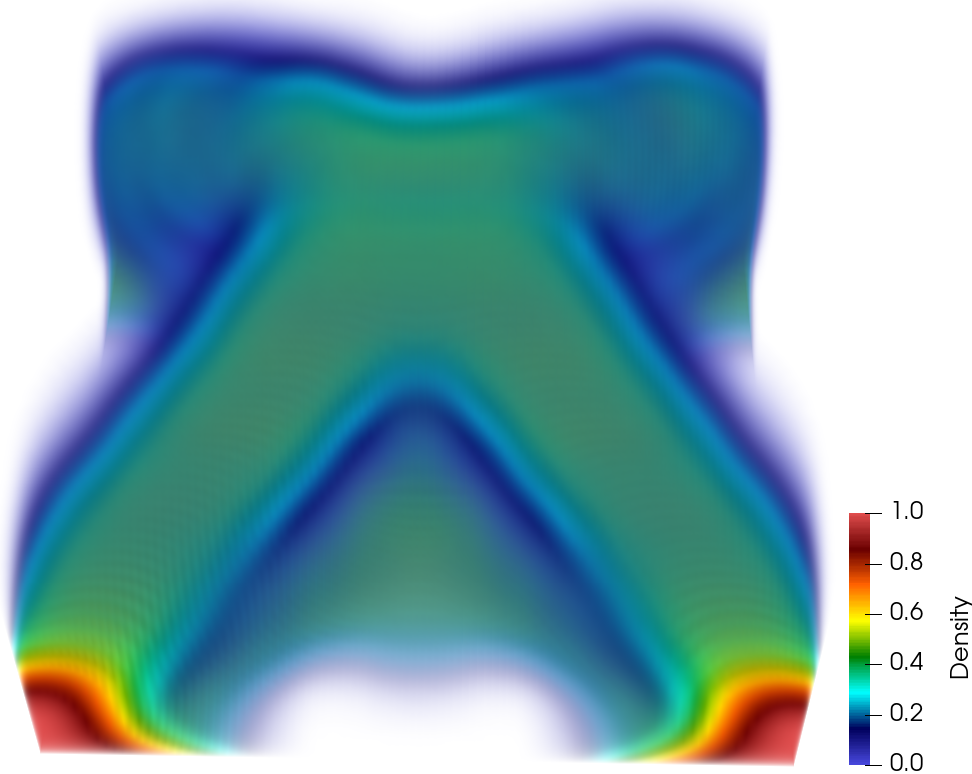}
    \caption{}
    \label{fig:michell3d-density}
  \end{subfigure}\\[8pt]
  \begin{subfigure}[t]{0.40\textwidth}
    \centering
    \includegraphics[width=\linewidth]{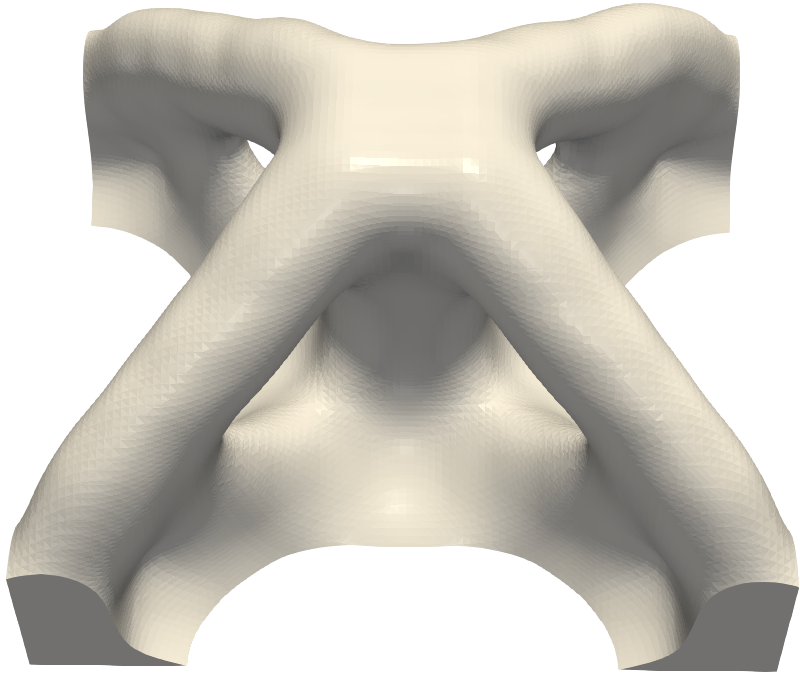}
    \caption{}
    \label{fig:michell3d-topology}
  \end{subfigure}
  \hfill
  \begin{subfigure}[t]{0.40\textwidth}
    \centering
    \includegraphics[width=\linewidth]{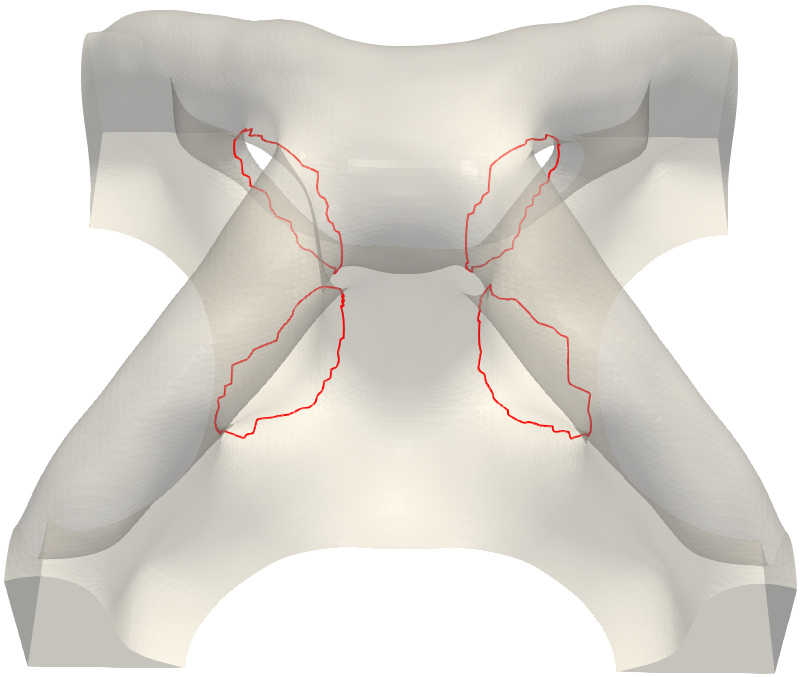}
    \caption{}
    \label{fig:michell3d-topology-transparent}
  \end{subfigure}
  \caption{(a) The 3D Michell structure model.
   (b) Optimized density field without topology control.
   (c) The corresponding binary structure with $N_1=4, N_2=0$.
   (d) Semi-transparent view of (c), with the tunnels marked by red loops.}
\end{figure}

In this example, the 3D Michell benchmark is defined on the cuboid
$[0,3]\times[0,3]\times[0,2]$, represented by a triquadratic NURBS volume with
$18\times 18\times 12$ control points. As illustrated in
Fig.~\ref{fig:michell3d-model}, three corners of the bottom face are supported by
vertical rollers, while the fourth corner is fully fixed, and a unit downward
point load is applied at the centre of the bottom face.
The prescribed maximum volume fraction is $V_{0}=0.2$.
Without topological constraint imposed, 
the baseline structure is obtained and shown in
Figs.~\ref{fig:michell3d-density}--\ref{fig:michell3d-topology-transparent}, which
contains four tunnels and no cavities. \par

\begin{figure}[!htb]
  \centering
  \begin{subfigure}[b]{0.40\textwidth}
    \centering
    \includegraphics[width=\textwidth]{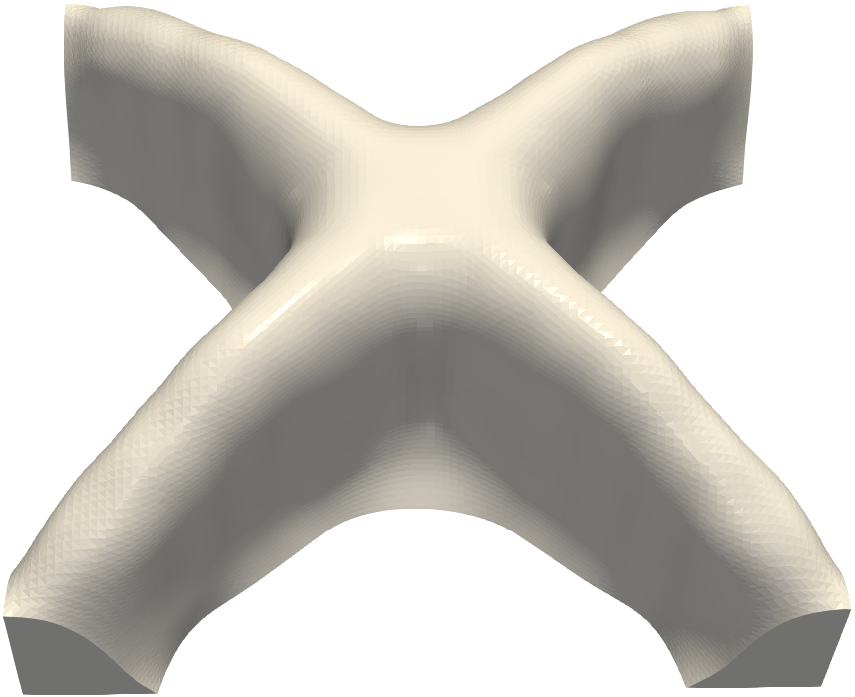}
    \caption{Opaque view}
    \label{fig:michell3d-h1-0-topology}
  \end{subfigure}
  \hfill
  \begin{subfigure}[b]{0.40\textwidth}
    \centering
    \includegraphics[width=\textwidth]{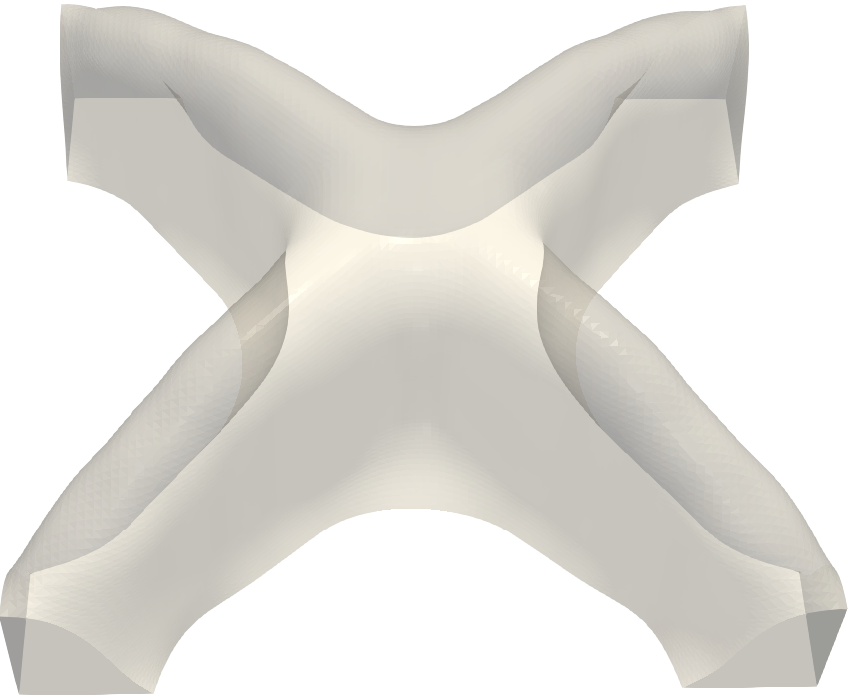}
    \caption{Semi-transparent view}
    \label{fig:michell3d-h1-0-transparent}
  \end{subfigure}
  \caption{Optimized 3D Michell structure for $\bar{N}_{1}=0, \bar{N}_{2}=0$.}
  \label{fig:michell3d-h1-0}
\end{figure}

\begin{figure}[!htb]
  \centering
  \begin{subfigure}[b]{0.40\textwidth}
    \centering
    \includegraphics[width=\textwidth]{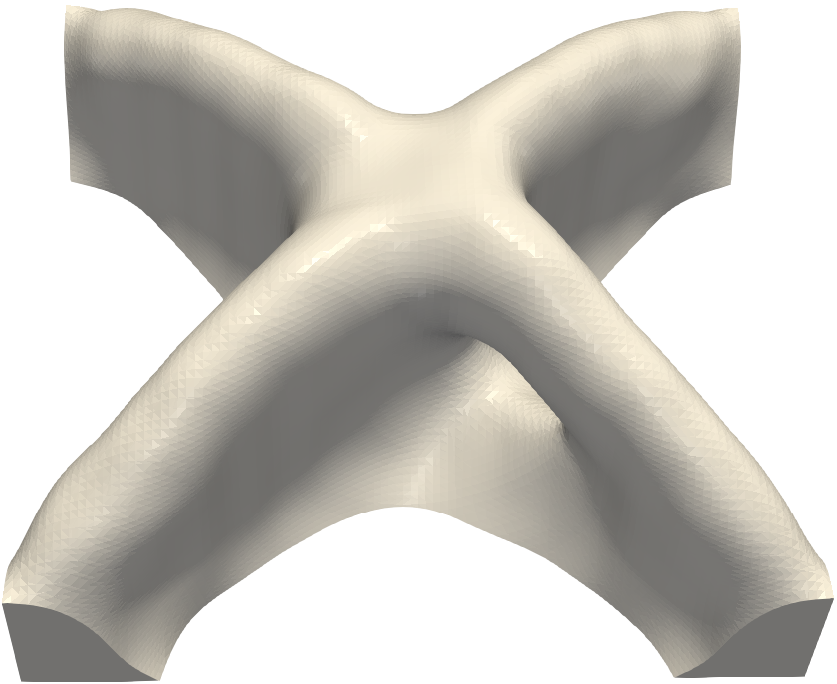}
    \caption{Opaque view}
    \label{fig:michell3d-h1-1-topology}
  \end{subfigure}
  \hfill
  \begin{subfigure}[b]{0.40\textwidth}
    \centering
    \includegraphics[width=\textwidth]{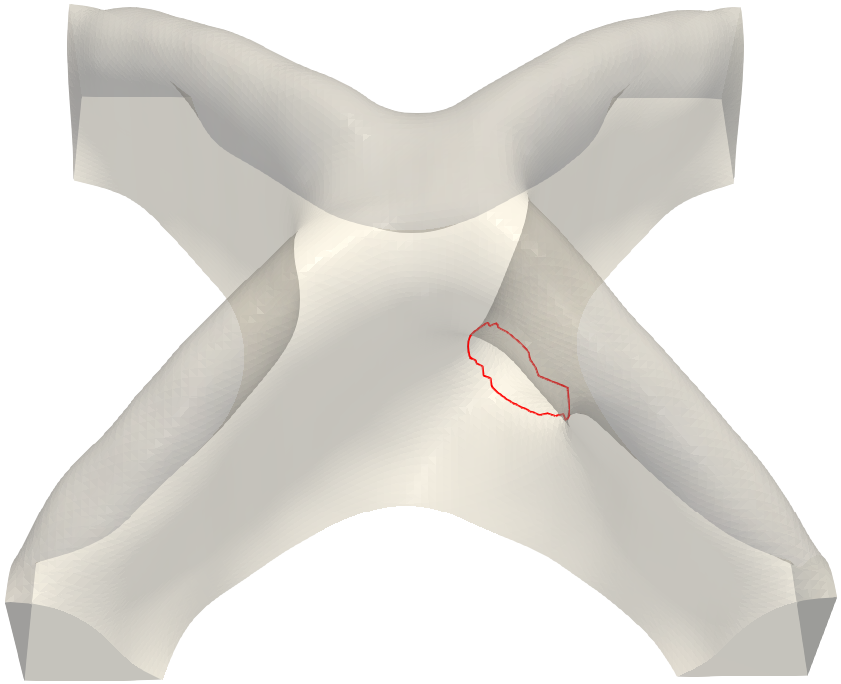}
    \caption{Semi-transparent view}
    \label{fig:michell3d-h1-transparent}
  \end{subfigure}
  \caption{Optimized 3D Michell structure for $\bar{N}_{1}=1, \bar{N}_{2}=0$.}
  \label{fig:michell3d-h1-1}
\end{figure}

\begin{figure}[!htb]
  \centering
  \begin{subfigure}[b]{0.40\textwidth}
    \centering
    \includegraphics[width=\textwidth]{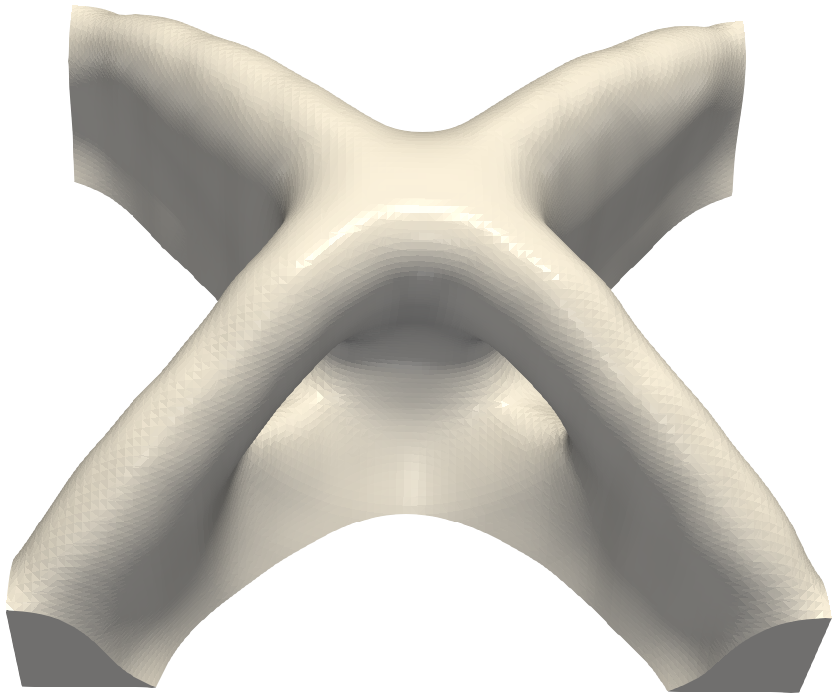}
    \caption{Opaque view}
    \label{fig:michell3d-h1-2-topology}
  \end{subfigure}
  \hfill
  \begin{subfigure}[b]{0.40\textwidth}
    \centering
    \includegraphics[width=\textwidth]{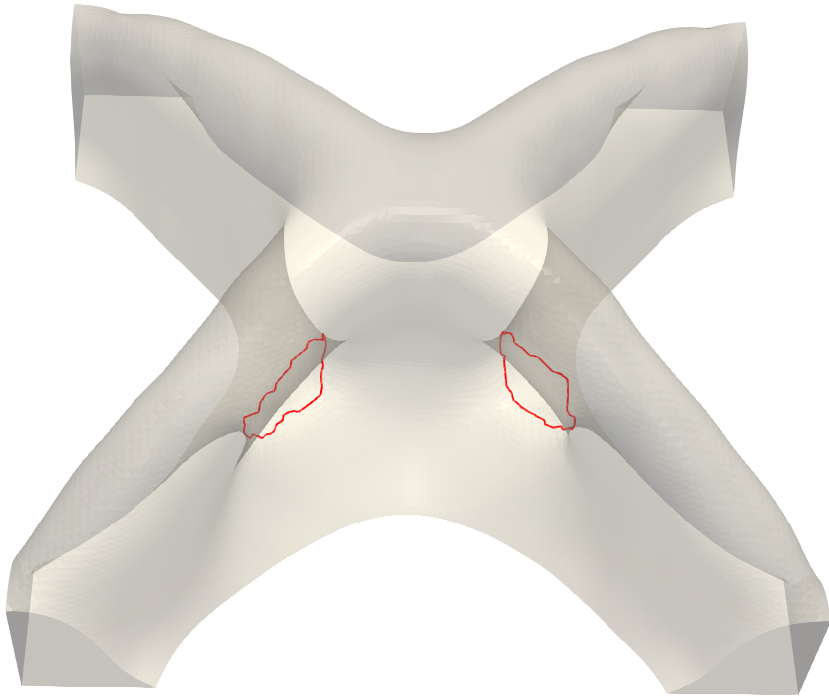}
    \caption{Semi-transparent view}
    \label{fig:michell3d-h1-2-transparent}
  \end{subfigure}
  \caption{Optimized 3D Michell structure for $\bar{N}_{1}=2, \bar{N}_{2}=0$.}
  \label{fig:michell3d-h1-2}
\end{figure}

\begin{figure}[!htb]
  \centering
  \begin{subfigure}[b]{0.40\textwidth}
    \centering
    \includegraphics[width=\textwidth]{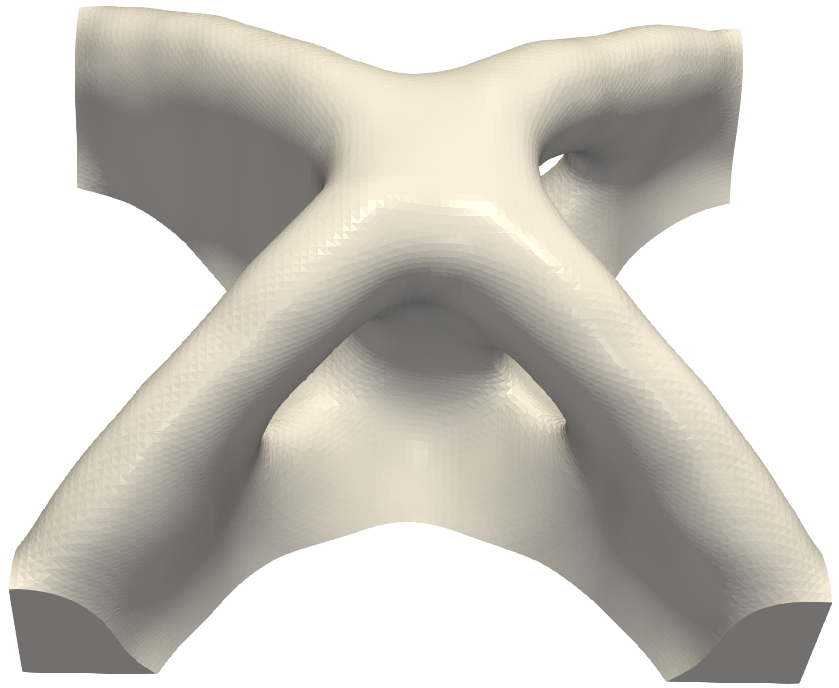}
    \caption{Opaque view}
    \label{fig:michell3d-h1-3-topology}
  \end{subfigure}
  \hfill
  \begin{subfigure}[b]{0.40\textwidth}
    \centering
    \includegraphics[width=\textwidth]{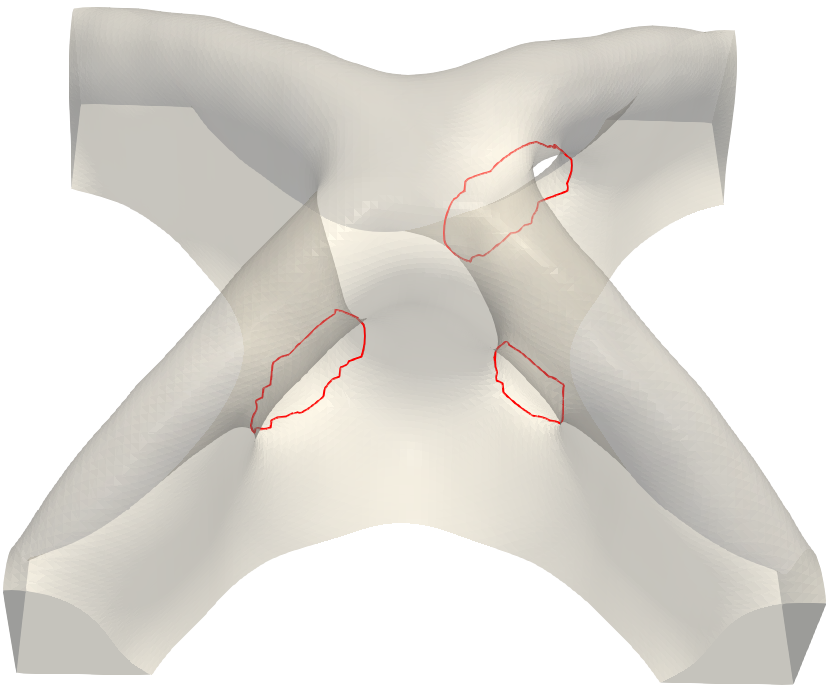}
    \caption{Semi-transparent view}
    \label{fig:michell3d-h1-3-transparent}
  \end{subfigure}
  \caption{Optimized 3D Michell structure for $\bar{N}_{1}=3, \bar{N}_{2}=0$.}
  \label{fig:michell3d-h1-3}
\end{figure}

Topology control is subsequently applied under the same settings, keeping
$\bar{N}_{2}=0$ and setting $\bar{N}_{1}$ to $0,1,2,3$ in turn. The optimized
structures are presented in
Figs.~\ref{fig:michell3d-h1-0}--\ref{fig:michell3d-h1-3}, and the quantitative comparison is reported in
Table~\ref{tab:michell3d}.
Filling removes the redundant tunnels and leaves the final structures smooth. As
the prescribed bound $\bar{N}_{1}$ is tightened, the compliance increases gradually,
which is consistent with the overall trend observed in the 2D examples.\par

\begin{table}[!htb]
  \centering
  \caption{Comparison of the compliance values and of the resulting numbers of tunnels
($N_{1}$) and cavities ($N_{2}$) for the 3D Michell structure under different
prescribed bounds $\bar{N}_{1}$, with $\bar{N}_{2} = 0$.}
  \label{tab:michell3d}
  \small
  \setlength{\tabcolsep}{5pt}
  \begin{tabular}{@{}lccccc@{}}
    \toprule
    \multirow{2}{*}{Quantity} & \multicolumn{4}{c}{Prescribed $\bar{N}_{1}$} & \multirow{2}{*}{Baseline}\\
    \cmidrule(lr){2-5}
    & $0$ & $1$ & $2$ & $3$ & \\
    \midrule
    Weights $(\mu_{0},\mu_{1},\mu_{2})$ & \multicolumn{4}{c}{$(50,50,50)$} & $(0,0,0)$ \\
    Compliance & $45.87$ & $44.62$ & $44.37$ & $43.12$ & $40.43$\\
    $(N_{1},N_{2})$ & $(0, 0)$ & $(1, 0)$ & $(2, 0)$ & $(3, 0)$ & $(4, 0)$\\
    \bottomrule
  \end{tabular}
\end{table}

\subsubsection{Bridge}

\begin{figure}[!htb]
  \centering
  \begin{subfigure}[t]{0.40\textwidth}
    \centering
    \includegraphics[width=\linewidth]{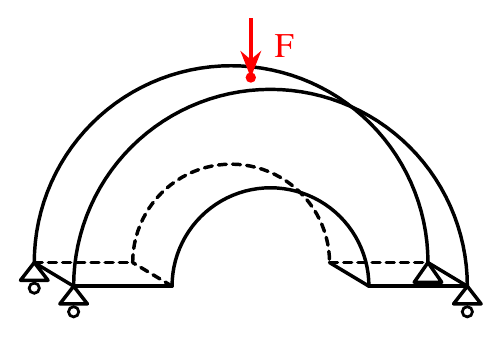}
    \caption{}
    \label{fig:bridge3d-model}
  \end{subfigure}
  \hfill
  \begin{subfigure}[t]{0.42\textwidth}
    \centering
    \includegraphics[width=\linewidth]{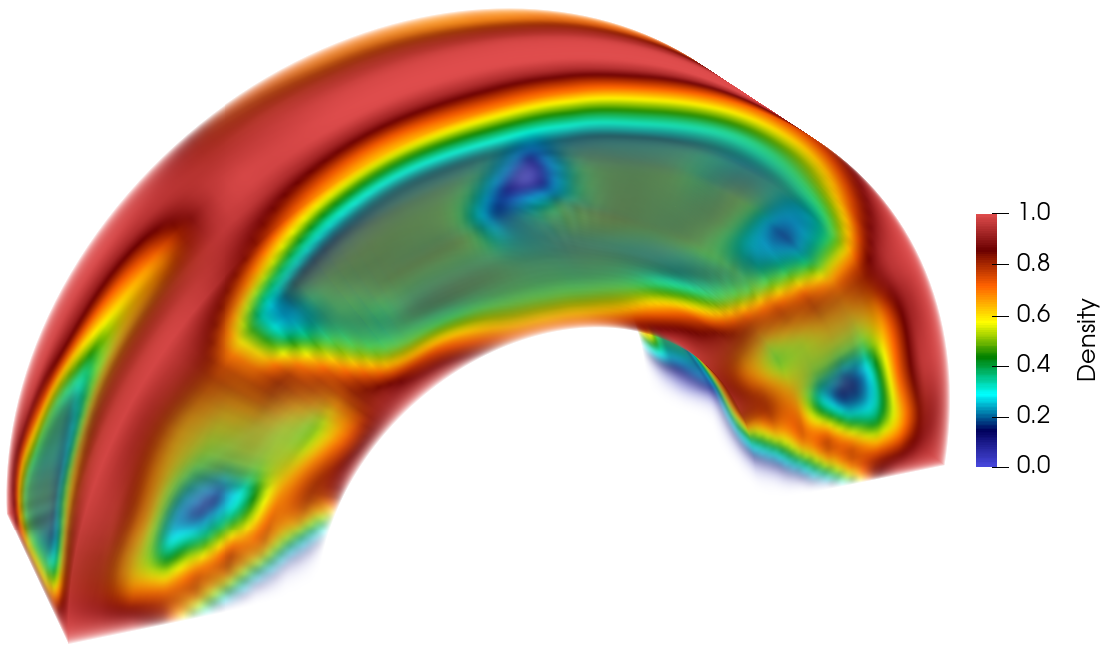}
    \caption{}
    \label{fig:bridge3d-density}
  \end{subfigure}\\[8pt]
  \begin{subfigure}[t]{0.40\textwidth}
    \centering
    \includegraphics[width=\linewidth]{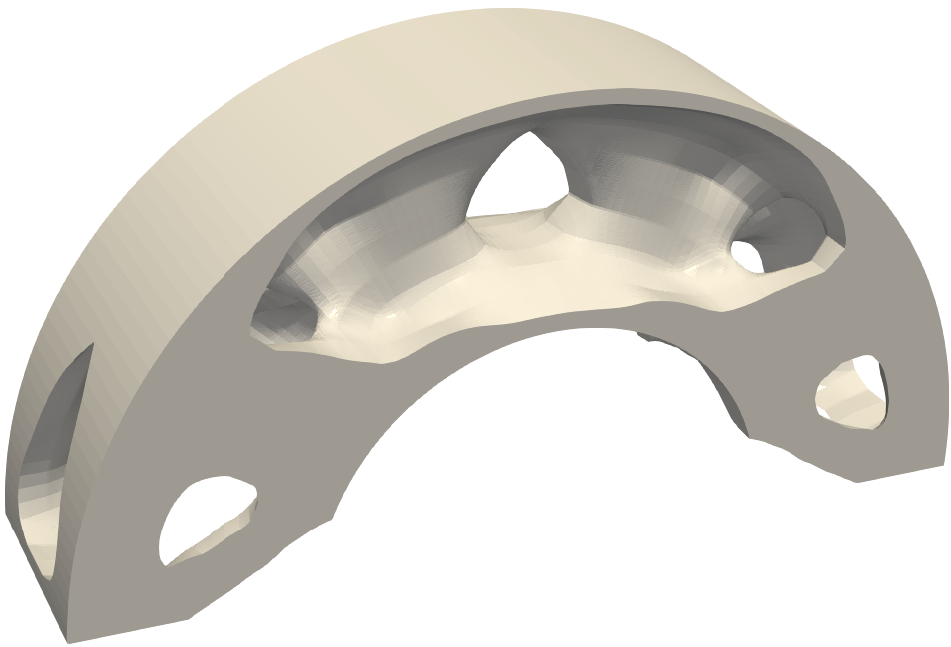}
    \caption{}
    \label{fig:bridge3d-topology}
  \end{subfigure}
  \hfill
  \begin{subfigure}[t]{0.40\textwidth}
    \centering
    \includegraphics[width=\linewidth]{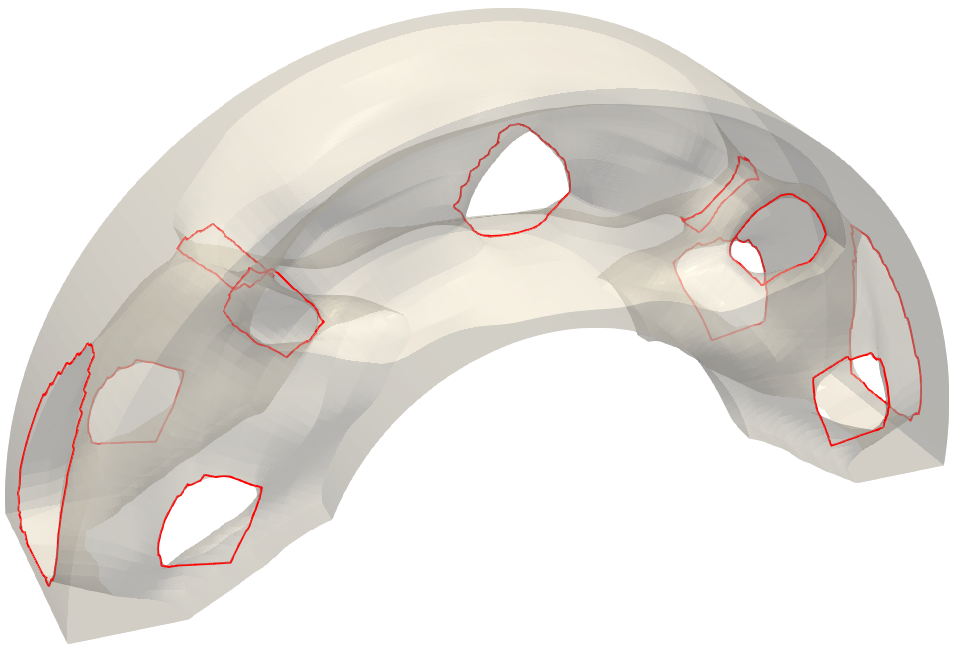}
    \caption{}
    \label{fig:bridge3d-topology-transparent}
  \end{subfigure}
  \caption{(a) The 3D bridge model.
   (b) Optimized density field without topology control.
   (c) The corresponding binary structure with $N_1=11, N_2=0$.
   (d) Semi-transparent view of (c), with the tunnels marked by red loops.}
\end{figure}

\begin{figure}[!htb]
  \centering
  \begin{subfigure}[b]{0.40\textwidth}
    \centering
    \includegraphics[width=\textwidth]{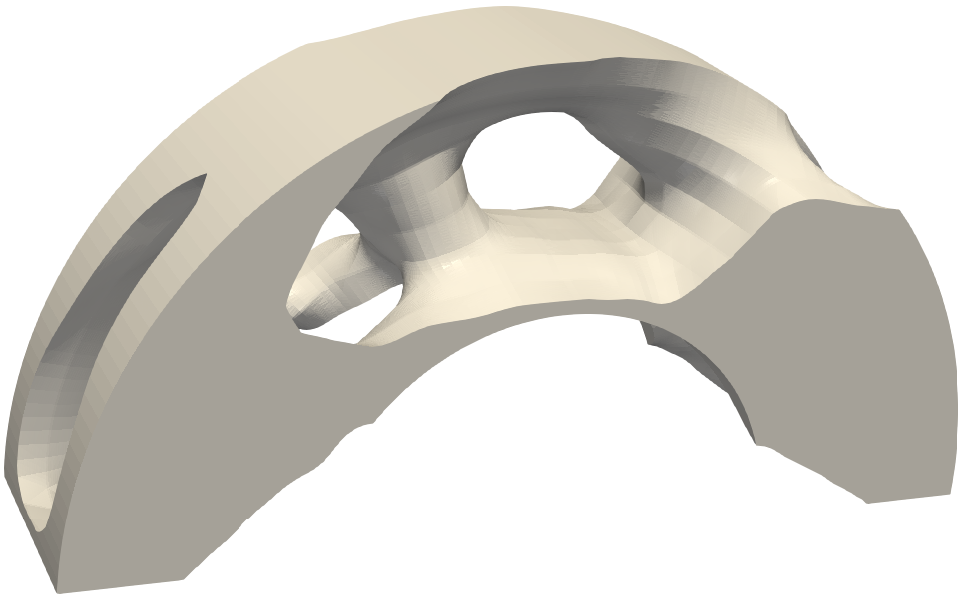}
    \caption{Opaque view}
    \label{fig:bridge3d-h1-5-topology}
  \end{subfigure}
  \hfill
  \begin{subfigure}[b]{0.40\textwidth}
    \centering
    \includegraphics[width=\textwidth]{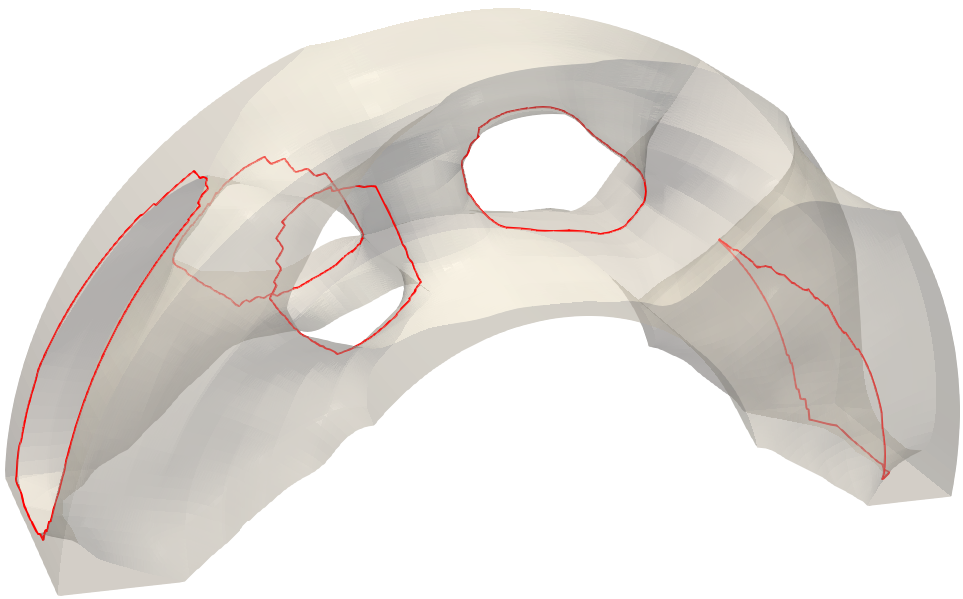}
    \caption{Semi-transparent view}
    \label{fig:bridge3d-h1-5-transparent}
  \end{subfigure}
  \caption{Optimized 3D bridge structure for $\bar{N}_{1}=5, \bar{N}_{2}=0$.}
  \label{fig:bridge3d-h1-5}
\end{figure}

\begin{figure}[!htb]
  \centering
  \begin{subfigure}[b]{0.40\textwidth}
    \centering
    \includegraphics[width=\textwidth]{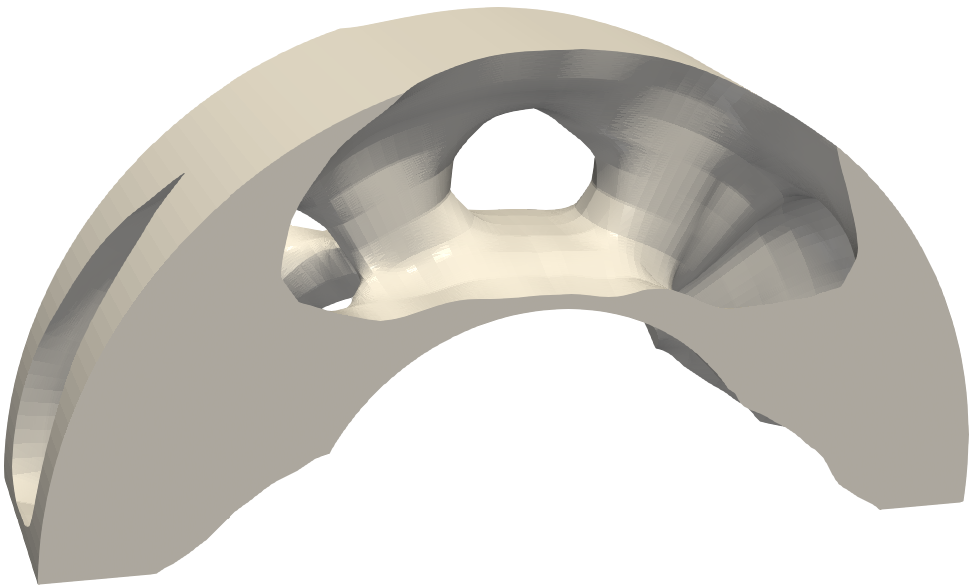}
    \caption{Opaque view}
    \label{fig:bridge3d-h1-6-topology}
  \end{subfigure}
  \hfill
  \begin{subfigure}[b]{0.40\textwidth}
    \centering
    \includegraphics[width=\textwidth]{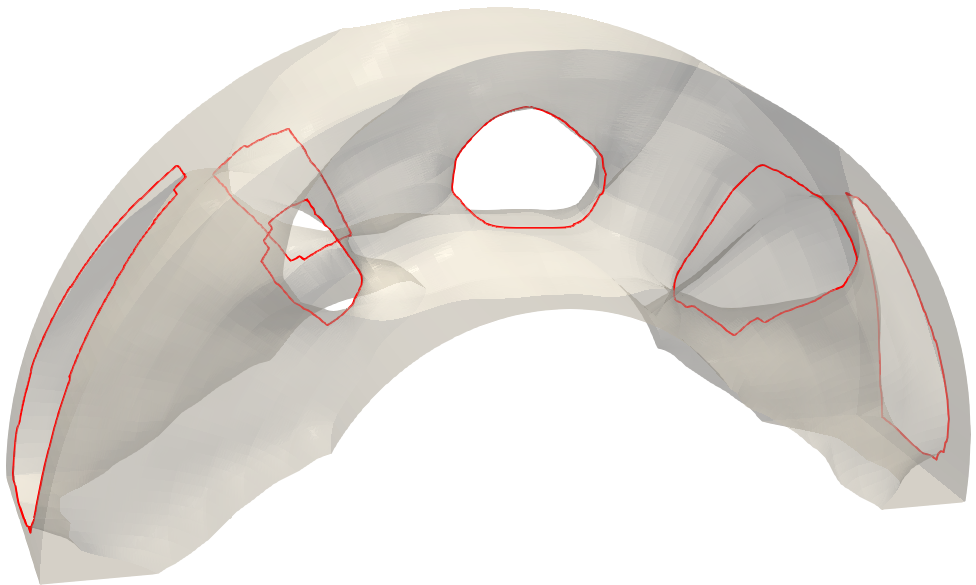}
    \caption{Semi-transparent view}
    \label{fig:bridge3d-h1-6-transparent}
  \end{subfigure}
  \caption{Optimized 3D bridge structure for $\bar{N}_{1}=6, \bar{N}_{2}=0$.}
  \label{fig:bridge3d-h1-6}
\end{figure}

\begin{figure}[!htb]
  \centering
  \begin{subfigure}[b]{0.40\textwidth}
    \centering
    \includegraphics[width=\textwidth]{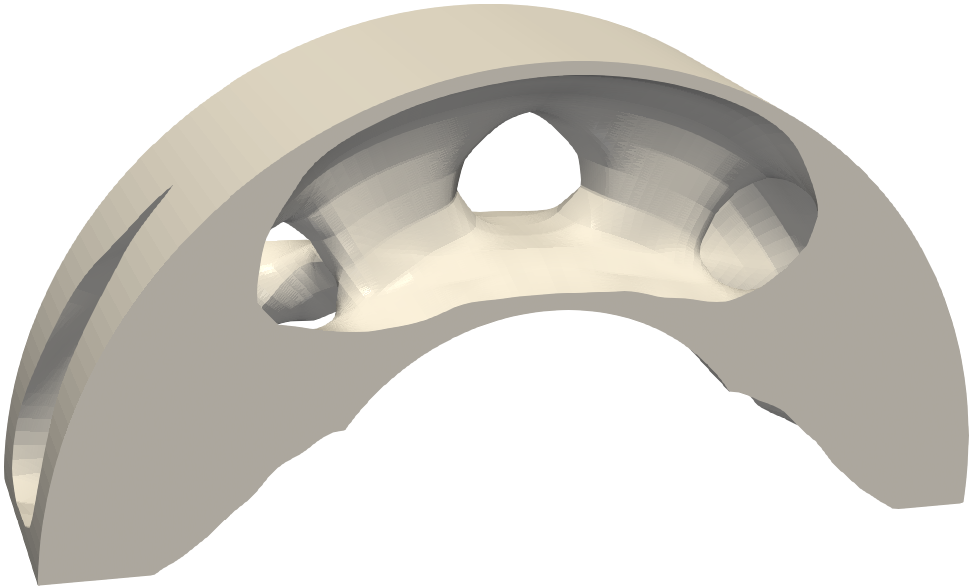}
    \caption{Opaque view}
    \label{fig:bridge3d-h1-7-topology}
  \end{subfigure}
  \hfill
  \begin{subfigure}[b]{0.40\textwidth}
    \centering
    \includegraphics[width=\textwidth]{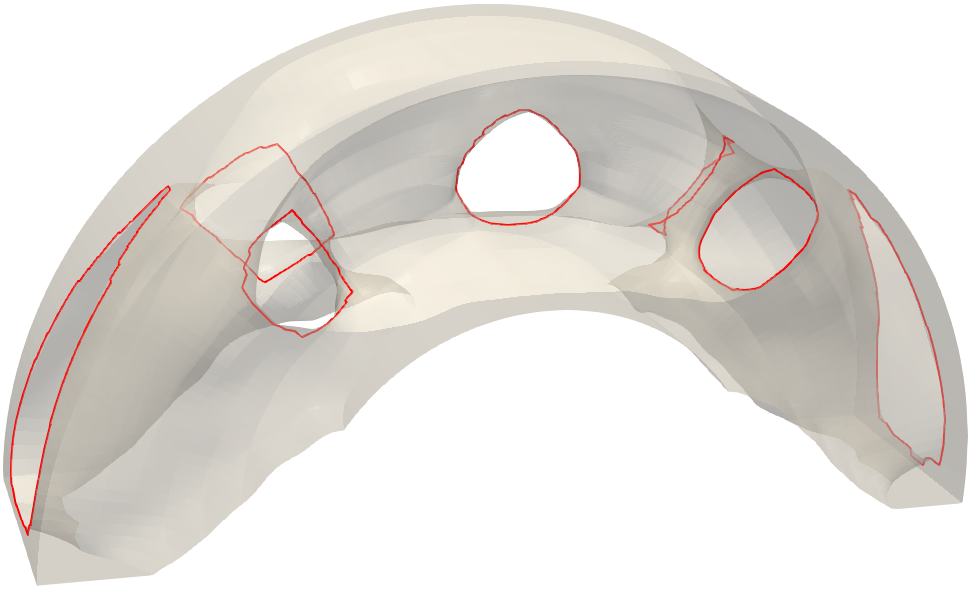}
    \caption{Semi-transparent view}
    \label{fig:bridge3d-h1-7-transparent}
  \end{subfigure}
  \caption{Optimized 3D bridge structure for $\bar{N}_{1}=7, \bar{N}_{2}=0$.}
  \label{fig:bridge3d-h1-7}
\end{figure}

\begin{figure}[!htb]
  \centering
  \begin{subfigure}[b]{0.40\textwidth}
    \centering
    \includegraphics[width=\textwidth]{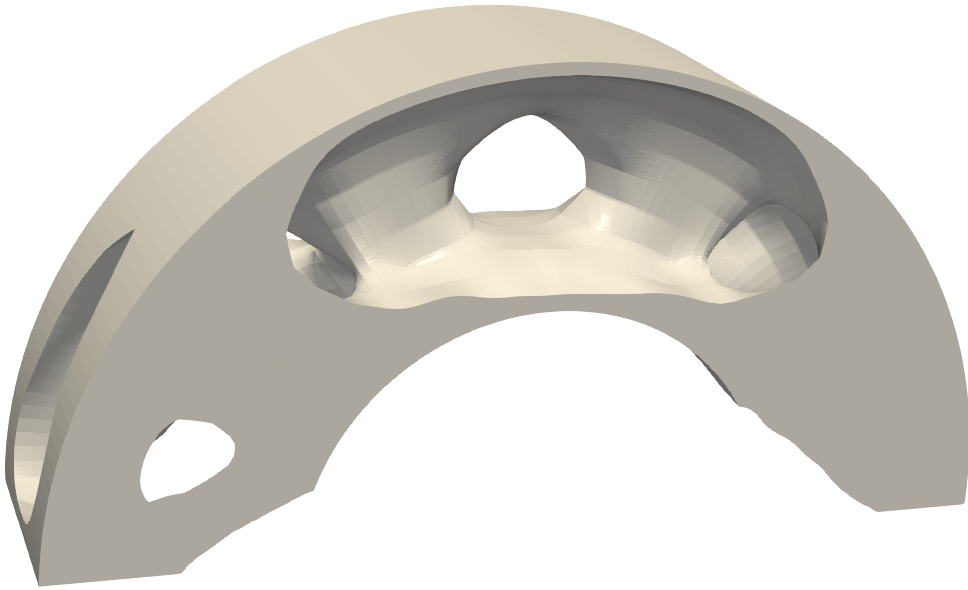}
    \caption{Opaque view}
    \label{fig:bridge3d-h1-8-topology}
  \end{subfigure}
  \hfill
  \begin{subfigure}[b]{0.40\textwidth}
    \centering
    \includegraphics[width=\textwidth]{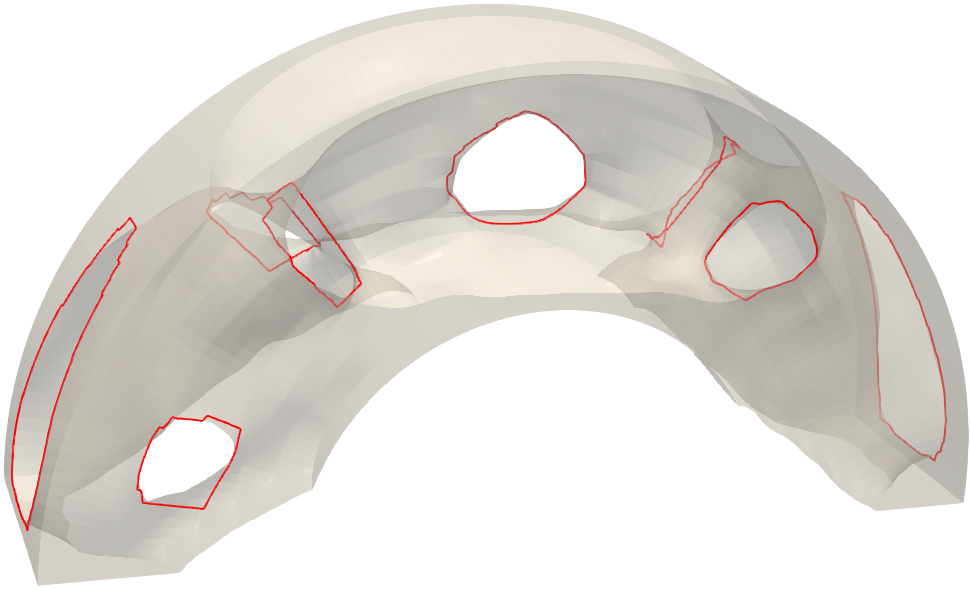}
    \caption{Semi-transparent view}
    \label{fig:bridge3d-h1-8-transparent}
  \end{subfigure}
  \caption{Optimized 3D bridge structure for $\bar{N}_{1}=8, \bar{N}_{2}=0$.}
  \label{fig:bridge3d-h1-8}
\end{figure}

\begin{figure}[!htb]
  \centering
  \begin{subfigure}[b]{0.40\textwidth}
    \centering
    \includegraphics[width=\textwidth]{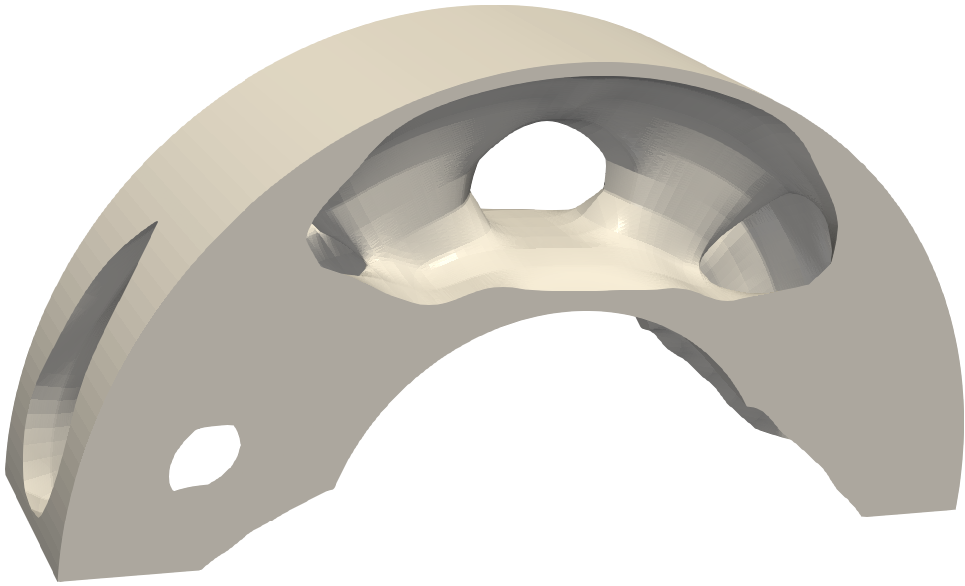}
    \caption{Opaque view}
    \label{fig:bridge3d-h1-9-topology}
  \end{subfigure}
  \hfill
  \begin{subfigure}[b]{0.40\textwidth}
    \centering
    \includegraphics[width=\textwidth]{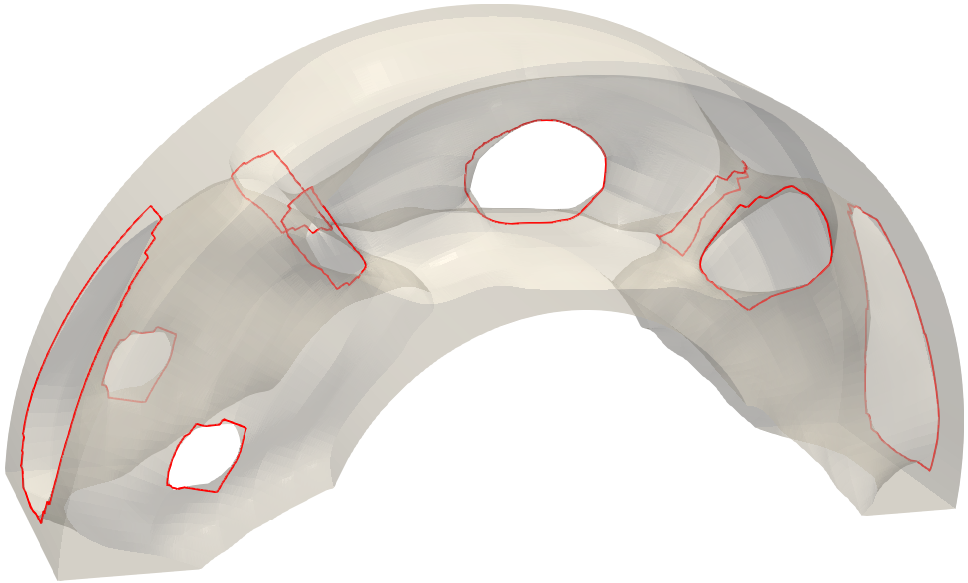}
    \caption{Semi-transparent view}
    \label{fig:bridge3d-h1-9-transparent}
  \end{subfigure}
  \caption{Optimized 3D bridge structure for $\bar{N}_{1}=9, \bar{N}_{2}=0$.}
  \label{fig:bridge3d-h1-9}
\end{figure}

\begin{figure}[!htb]
  \centering
  \begin{subfigure}[b]{0.40\textwidth}
    \centering
    \includegraphics[width=\textwidth]{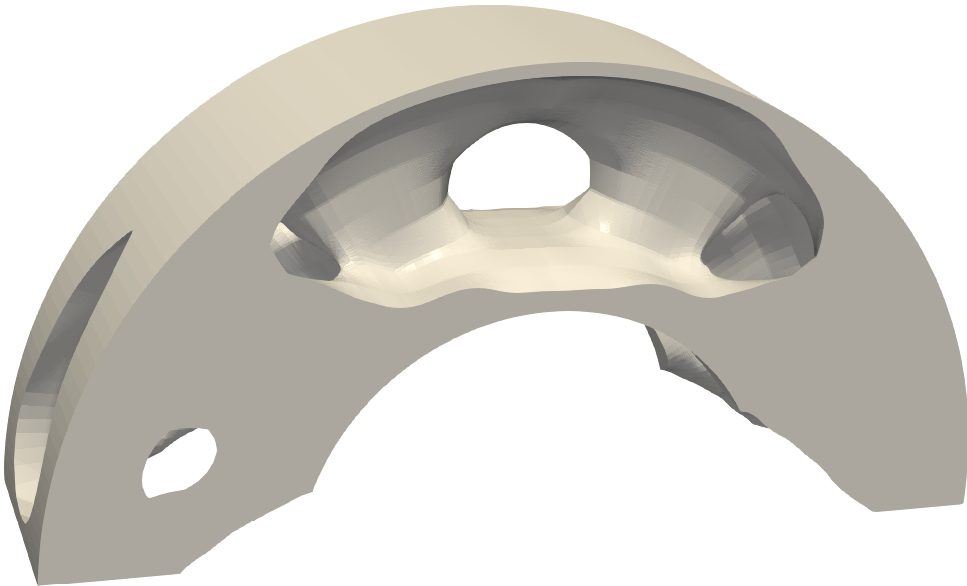}
    \caption{Opaque view}
    \label{fig:bridge3d-h1-10-topology}
  \end{subfigure}
  \hfill
  \begin{subfigure}[b]{0.40\textwidth}
    \centering
    \includegraphics[width=\textwidth]{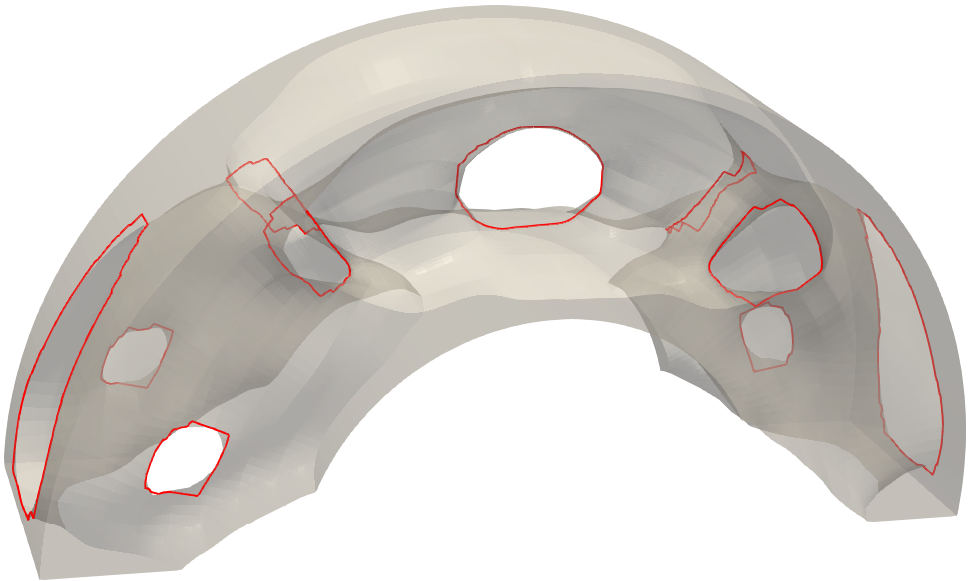}
    \caption{Semi-transparent view}
    \label{fig:bridge3d-h1-10-transparent}
  \end{subfigure}
  \caption{Optimized 3D bridge structure for $\bar{N}_{1}=10, \bar{N}_{2}=0$.}
  \label{fig:bridge3d-h1-10}
\end{figure}

In this subsection, the 3D bridge benchmark defined on a curved domain is considered,
which is represented by a NURBS volume of degree $(2,1,1)$ with $49\times 9\times 9$
control points. As depicted in Fig.~\ref{fig:bridge3d-model}, three corners of the
bottom face are restrained in the vertical direction, while the remaining corner is
fully fixed, and a unit downward point load is applied at the midpoint of the
outer crown. The maximum volume fraction is set to $V_{0}=0.4$.
The unconstrained baseline structure shown in
Figs.~\ref{fig:bridge3d-density}--\ref{fig:bridge3d-topology-transparent}
contains eleven tunnels and no cavities.
Next, the same problem is solved under topology control, with $\bar{N}_{2}$ fixed to $0$ and $\bar{N}_{1}$ ranging from $5$ to $10$. The optimized structures
are shown in Figs.~\ref{fig:bridge3d-h1-5}--\ref{fig:bridge3d-h1-10}, and the quantitative comparison is reported in Table~\ref{tab:bridge3d}.\par

\begin{table}[!htb]
  \centering
  \caption{Comparison of the compliance values and of the resulting numbers of tunnels
($N_{1}$) and cavities ($N_{2}$) for the 3D bridge under different
prescribed bounds $\bar{N}_{1}$, with $\bar{N}_{2} = 0$.}
  \label{tab:bridge3d}
  \small
  \setlength{\tabcolsep}{4pt}
  \renewcommand{\arraystretch}{1.12}
  \begin{tabular}{@{}l *{7}{c} @{}}
    \toprule
    \multirow{2}{*}{Quantity} & \multicolumn{6}{c}{Prescribed $\bar{N}_{1}$} & \multirow{2}{*}{Baseline}\\
    \cmidrule(lr){2-7}
     & $5$ & $6$ & $7$ & $8$ & $9$ & $10$ & \\
    \midrule
    Weights $(\mu_{0},\mu_{1},\mu_{2})$ & $(50,50,50)$ & \multicolumn{5}{c}{$(30,30,30)$} & $(0,0,0)$ \\
    Compliance & $94.61$ & $89.43$ & $86.37$ & $85.26$ & $84.41$ & $84.11$ & $81.91$ \\
    $(N_{1},N_{2})$
      & $(5,0)$ & $(6,0)$ & $(7,0)$ & $(8,0)$ & $(9,0)$ & $(10,0)$ & $(11,0)$ \\
    \bottomrule
  \end{tabular}
\end{table}

\subsubsection{Platform}

\begin{figure}[!htb]
  \centering
  \begin{subfigure}[t]{0.40\textwidth}
    \centering
    \includegraphics[width=\linewidth]{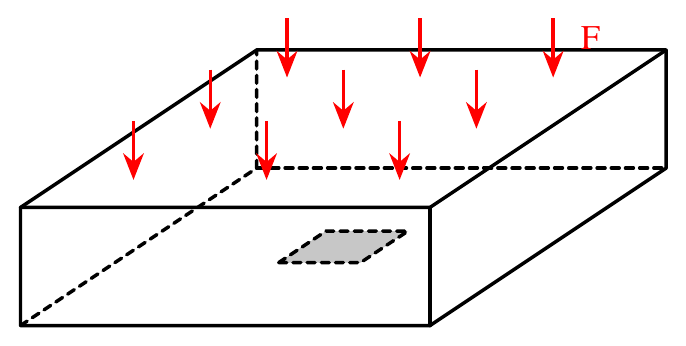}
    \caption{}
    \label{fig:platform3d-model}
  \end{subfigure}
  \hfill
  \begin{subfigure}[t]{0.42\textwidth}
    \centering
    \includegraphics[width=\linewidth]{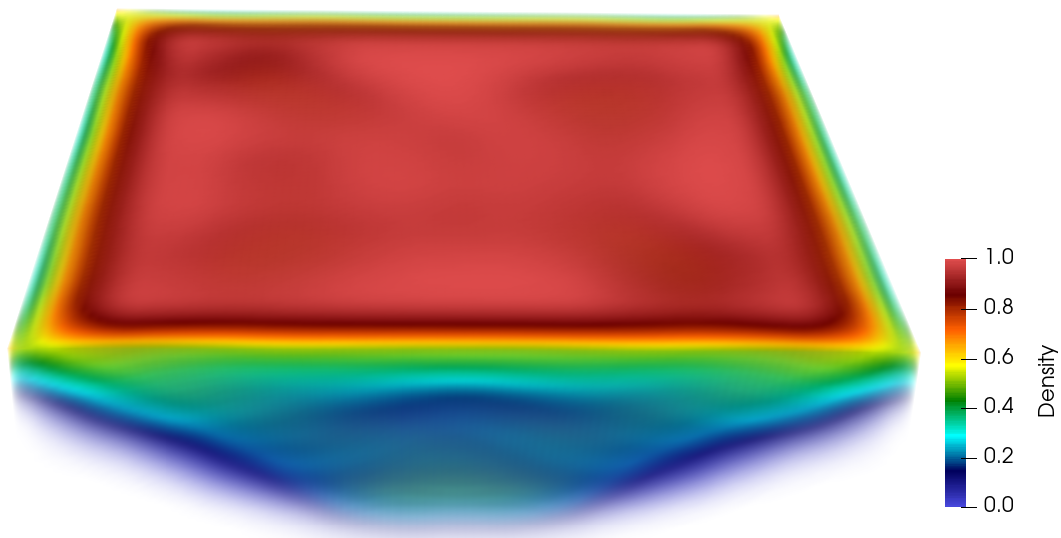}
    \caption{}
    \label{fig:platform3d-density}
  \end{subfigure}\\[8pt]
  \begin{subfigure}[t]{0.40\textwidth}
    \centering
    \includegraphics[width=\linewidth]{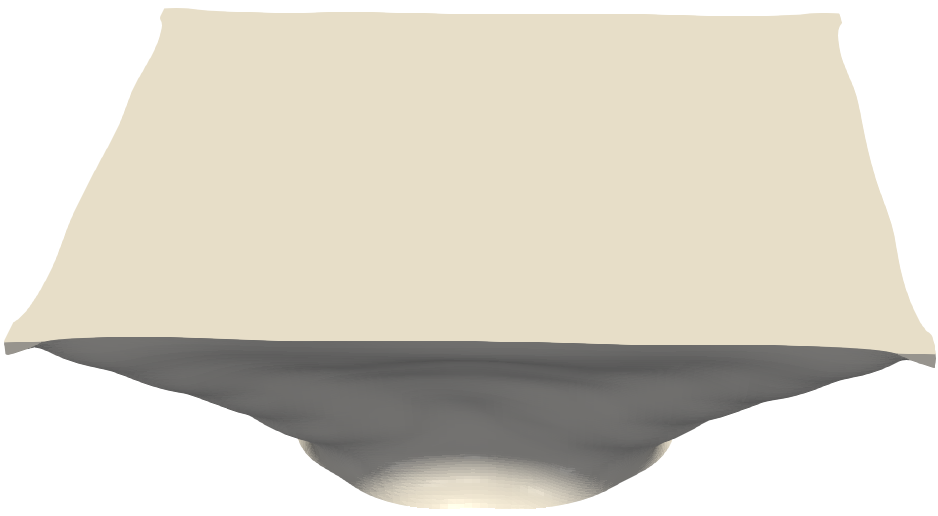}
    \caption{}
    \label{fig:platform3d-topology}
  \end{subfigure}
  \hfill
  \begin{subfigure}[t]{0.40\textwidth}
    \centering
    \includegraphics[width=\linewidth]{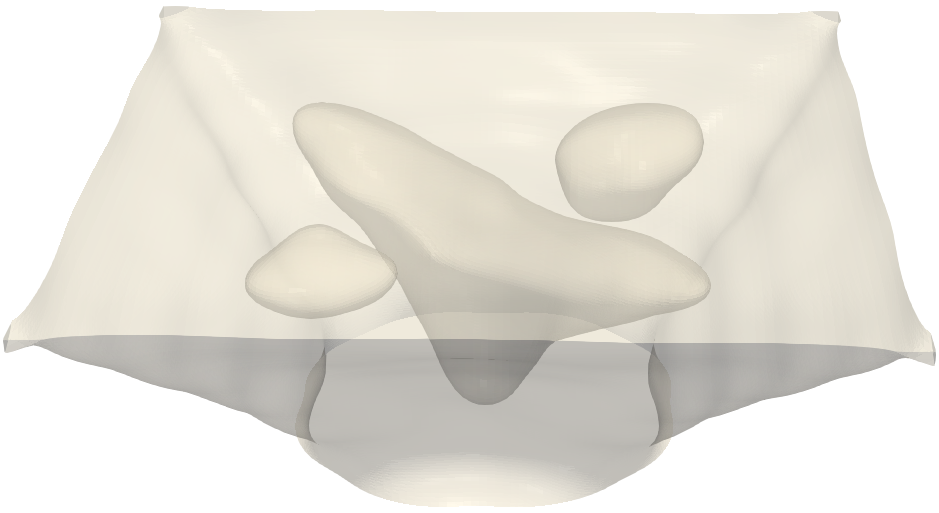}
    \caption{}
    \label{fig:platform3d-topology-transparent}
  \end{subfigure}
  \caption{(a) The 3D platform model.
   (b) Optimized density field without topology control.
   (c) The corresponding binary structure with $N_1=0, N_2=3$.
   (d) Semi-transparent view of (c) for visualization.}
\end{figure}

As the final 3D example, the platform benchmark occupies the cuboid
$[0,4]\times[0,4]\times[0,1.5]$, represented by a triquadratic NURBS volume
with $17\times 17\times 13$ control points. The boundary conditions are
shown in Fig.~\ref{fig:platform3d-model}, with a central square region of the
bottom face fully fixed and a uniform unit downward traction applied over the top
face.
The maximum volume fraction is set to $V_0=0.4$.
Solving the minimum compliance problem without topology
control yields a baseline that, in contrast to the preceding examples, contains
three enclosed cavities but no tunnels, as illustrated in Figs.~\ref{fig:platform3d-density}--\ref{fig:platform3d-topology-transparent}.\par

\begin{figure}[!htb]
  \centering
  \begin{subfigure}[b]{0.40\textwidth}
    \centering
    \includegraphics[width=\textwidth]{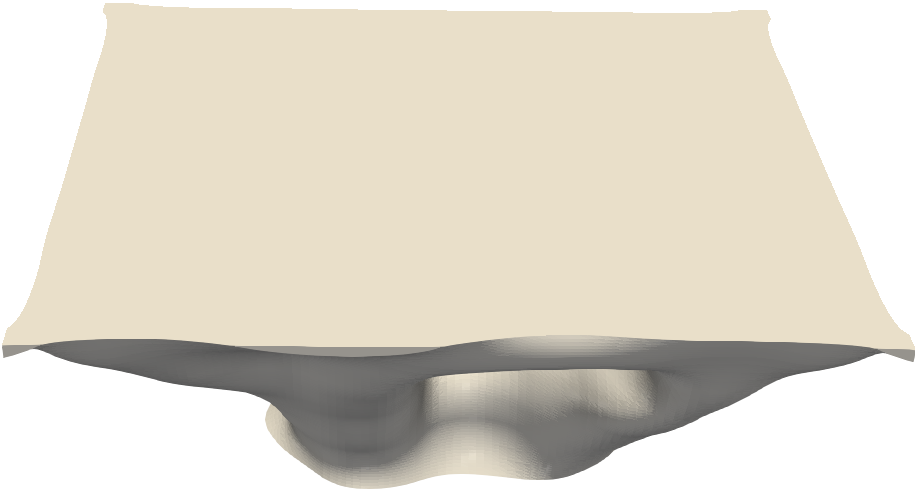}
    \caption{Opaque view}
    \label{fig:platform3d-h1-0}
  \end{subfigure}\hfill
  \begin{subfigure}[b]{0.40\textwidth}
    \centering
    \includegraphics[width=\textwidth]{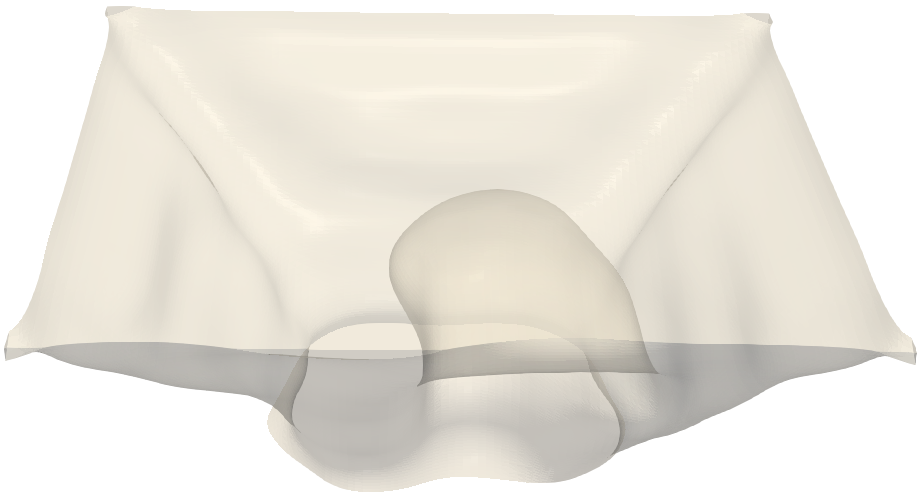}
    \caption{Semi-transparent view}
    \label{fig:platform3d-h1-0-transparent}
  \end{subfigure}
  \caption{Optimized 3D platform structure for $\bar{N}_{1}=0, \bar{N}_2=0$.}
  \label{fig:platform3d-control}
\end{figure}

Topology control is then imposed by setting $\bar{N}_{2}=0$ so as to eliminate
these cavities, while $\bar{N}_{1}=0$ is retained since the baseline already
contains no tunnels. The optimized structure is shown in
Fig.~\ref{fig:platform3d-control}, and the quantitative comparison is reported in Table~\ref{tab:platform3d}.\par

\begin{table}[!htb]
  \centering
  \caption{Comparison of the compliance values and of the resulting numbers of tunnels
($N_{1}$) and cavities ($N_{2}$) for the 3D platform under different
prescribed bounds $\bar{N}_{2}$, with $\bar{N}_{1}=0$.}
  \label{tab:platform3d}
  \small
  \begin{tabular}{@{}lcc@{}}
    \toprule
    Quantity & $\bar{N}_{2}=0$ & Baseline\\
    \midrule
    Weights $(\mu_{0},\mu_{1},\mu_{2})$ & $(100,100,200)$ & $(0, 0, 0)$\\
    Compliance & 419.75 & $374.11$\\
    $(N_1, N_2)$ & $(0,0)$& $(0,3)$\\
    \bottomrule
  \end{tabular}
\end{table}

\section{Conclusion}
\label{sec:conclusion}

In this study, we present a persistent homology-based isogeometric topology
optimization framework for explicit control of structural topology.  
A single sublevel filtration of the density field encodes the topological features as persistence
diagrams, from which the connected components, holes in 2D, and tunnels and cavities in 3D are
rigorously quantified. The resulting topological objectives are differentiable
and can be readily incorporated into gradient-based solvers such as MMA.
In addition, the proposed method removes each topological
feature either by cutting or by filling, according to the region of its
corresponding persistence pair. Constraints on different homology
dimensions are further enforced in a staged manner, resolving the conflict
between connectivity and hole control and, in 3D, the ambiguity introduced into
the tunnel count by cavities.\par

Numerical examples in both 2D and 3D demonstrate that the
framework provides explicit and reliable topology control
while maintaining competitive mechanical performance and stable convergence.
The present study, however, is confined to the minimum compliance problem under
a volume constraint. Extending the framework to other topology optimization
problems, such as heat conduction and compliant mechanisms,
therefore constitutes an important direction for future work. 
In addition, the topological weights are currently set empirically, so multi-objective
formulations that determine them automatically are worth further study. 
From a computational standpoint, it is also worthwhile to improve the efficiency of
persistent homology computations for high-resolution discretizations.\par

\section*{Acknowledgements}
This work is supported in part by the National Natural Science Foundation of China under Grant Nos. 12671457, 62272406.\par

\bibliographystyle{plainnat}
\bibliography{ref}

\end{document}